\documentclass[a4paper,12pt]{amsbook} 
\usepackage{amsmath,amsthm,geometry,graphicx} 
\usepackage{amscd}
\usepackage{psfrag,xspace,graphicx,geometry} 
\usepackage[latin1]{inputenc}
\usepackage[all,dvips,cmtip]{xy}
\usepackage[francais]{babel}

\normalfont  
  
\def\refhautdepage{\ref{conj:invariant,dur}} 

\newcommand{\RedefinitSymbole}[1]{%
  \expandafter\let\csname old\string#1\endcsname=#1 \let#1=\relax 
\newcommand{#1}{\csname old\string#1\endcsname\,}%
} 
\RedefinitSymbole{\forall} \RedefinitSymbole{\exists} 

\def\cf{\textit{cf.}\kern.3em}  
 \def\resp{\textit{resp.}\kern.3em} 
\renewcommand{\k}{\kern2pt} 
 
\numberwithin{equation}{section} 
 
\DeclareMathOperator{\rg}{rang} 
\DeclareMathOperator{\cst}{cst} 
\DeclareMathOperator{\g}{\psi} 
\DeclareMathOperator{\p}{p} 
 
\DeclareMathOperator{\antidiag}{antidiag}

\DeclareMathOperator{\End}{End} 
\DeclareMathOperator{\age}{age} 
 \DeclareMathOperator{\red}{red} 
\DeclareMathOperator{\ev}{ev} 
 \DeclareMathOperator{\orb}{orb} 
\DeclareMathOperator{\pr}{pr}  
  
\DeclareMathOperator{\diag}{diag} 
 \DeclareMathOperator{\pgcd}{pgcd} 
\DeclareMathOperator{\ppcm}{ppcm} 
  
\DeclareMathOperator{\id}{id}  
\DeclareMathOperator{\im}{Im} 
\DeclareMathOperator{\GL}{GL}

\renewcommand{\ker}{\mathop{\rm 
    Ker}\nolimits}
\newcommand{\reg}{\mathrm{reg}} 
\newcommand{\sing}{\mathrm{sing}} 

\newtheorem{conj}[equation]{{Conjecture}} 
\newtheorem{defi}[equation]{{Définition}} 
\newtheorem{prop}[equation]{{Proposition}} 
\newtheorem{thm}[equation]{{Théorème}}
\newtheorem{cor}[equation]{{Corollaire}} 
\newtheorem{lem}[equation]{{Lemme}} 
\theoremstyle{definition} 
 
\newtheorem{rem}[equation]{{Remarque}} 
\newtheorem{notation}[equation]{{Notation}} 
\theoremstyle{definition} 
\newtheorem{expl}[equation]{{Exemple}}

\newcommand{\RR}{\mathbb{R}} \newcommand{\NN}{\mathbb{N}} 
\newcommand{\ZZ}{\mathbb{Z}} \newcommand{\QQ}{\mathbb{Q}} 
\newcommand{\CC}{\mathbb{C}} \newcommand{\PP}{\mathbb{P}} 
\newcommand{\ds}{\displaystyle} \newcommand{\Inv}[1]{\{1-#1\}} 
\newcommand{\bs}{\boldsymbol} \newcommand{\gr}{\mathrm{gr}} 
\newcommand{\V}{S_{w}} 
 
\begin{document} 


\begin{titlepage}
\vspace*{1cm}
  \begin{center}
    \hspace{10cm} \large \textbf{\'Etienne Mann}
  \end{center}
 
\vspace*{7cm}
\hspace{-5mm}
\rule{\textwidth}{0.5mm}
\vspace*{1cm}

  \begin{center}
\Huge \textbf{COHOMOLOGIE QUANTIQUE ORBIFOLDE DES ESPACES PROJECTIFS \`A
  POIDS} 
\end{center}

\vspace*{1cm}

\hspace{-5mm}\rule{\textwidth}{0.5mm}

\end{titlepage}

\begin{titlepage} 
\vspace*{3cm}
E. Mann

 SISSA, Via Beirut 2-4, I--34014 Trieste, Italie

\emph{E-mail :} { mann@sissa.it}

\emph{Url :}{ http://www-irma.u-strasbg.fr/\verb$~$mann/}
  
\vspace{10cm}
\hspace{-0.7cm}\rule{\textwidth}{0.1mm}
\vspace{0.3cm}

\hspace{0.5mm}\emph{\textbf{Key words and phrases.}}--- {Cohomologie
  orbifolde, invariant de Gromov-Witten orbifold, cohomologie
  quantique orbifolde,
  vari\'et\'e de Frobenius, r\'eseau de Brieskorn, syst\`eme de
  Gauss-Manin}

\vspace{-3mm}
\hspace{-7mm}\rule{\textwidth}{0.1mm}
\end{titlepage}

\begin{titlepage}\vspace{3cm}
  \begin{center}
    \Large \textbf{ ORBIFOLD QUANTUM COHOMOLOGY OF WEIGHTED PROJECTIVE
      SPACES} \\ \  \\ \large \textbf{Etienne Mann}
  \end{center}
\vspace{1cm}
  \begin{center}
    \textbf{Abstract:}
  \end{center}
In $2001$, S.\k Barannikov showed that the Frobenius manifold coming
from the quantum cohomology of the complex projective space of
dimension $n$ is isomorphic to the Frobenius manifold associated to
the Laurent polynomial $x_{1}+ \ldots+x_{n}+1/x_{1}\ldots x_{n}$.
 
 The purpose of this thesis is to generalize this result. More
 precisely, given some positive integers $w_{0}, \ldots ,w_{n}$, we
 show, up to a conjecture on the value of some orbifold Gromov-Witten
 invariants, that the Frobenius structure obtained on the orbifold
 quantum cohomology of the weighted projective spaces with weights
 $w_{0}, \ldots ,w_{n}$ is isomorphic to the one obtained from the
 Laurent polynomial $f(u_{0}, \ldots ,u_{n}):=u_{0}+\ldots + u_{n}$
 restricted to $U:=\{(u_{0}, \ldots ,u_{n})\in \CC^{n+1}\mid \prod_{i}
 u_{i}^{w_{i}}=1\}$.
\vspace{2cm}
 \begin{center}
   \textbf{R\'esum\'e :}
 \end{center}
En $2001$, S.\k Barannikov a montré que la variété de Frobenius
provenant de la cohomologie quantique de l'espace projectif complexe
de dimension $n$ est isomorphe à la variété de Frobenius associée au
polyn\^{o}me de Laurent $x_{1}+ \ldots+x_{n}+1/x_{1}\ldots x_{n}$.
 
 L'objectif de cette thèse est de généraliser ce résultat. Plus
 précisément, étant des entiers strictement positifs $w_{0}, \ldots
 ,w_{n}$, nous montrons, modulo une conjecture sur la valeur de
 certains invariants de Gromov-Witten orbifold, que la structure de
 Frobenius obtenue sur la cohomologie quantique orbifolde de l'espace
 projectif de poids $w_{0}, \ldots ,w_{n}$ est isomorphe à celle
 obtenue à partir du polyn\^ome $f(u_{0}, \ldots ,u_{n}):=u_{0}+\ldots
 + u_{n}$ restreint à $U:=\{(u_{0}, \ldots ,u_{n})\in \CC^{n+1}\mid
 \prod_{i} u_{i}^{w_{i}}=1\}$.
\end{titlepage}

\begin{titlepage}
  \ 
\end{titlepage}

\mainmatter
\tableofcontents 

\chapter{Introduction}
\label{cha:introduction}  

\section{Historique du problème}  
 \label{sec:Hist-du-problme}  

En $1991$, B.\kern2pt Dubrovin a défini la structure de Frobenius sur 
une variété complexe pour axiomatiser une partie de la riche structure 
mathématique de la théorie topologique des champs (cf. introduction de  
{\cite{Dtft}}). Il s'est 
inspiré des travaux des quatre physiciens E.\kern2pt Witten \cite{Wtp2dg}, R.\k Dijk\-graaf, 
E.\kern2pt Verlinde et H.\kern2pt Verlinde \cite{DVVtst2dqg}.  Les variétés de Frobenius 
sont des variétés complexes munies d'une forme bilinéaire non 
dégénérée plate et d'un 
produit sur le fibré tangent complexe qui satisfont certaines 
conditions de compatibilité. Elles possèdent aussi un potentiel qui 
satisfait le système d'équations aux dérivées partielles WDVV du nom 
des quatre physiciens cités plus haut.  B.\k Dubrovin a montré qu'une 
variété de Frobenius semi-simple est déterminée par des \og conditions 
initiales \fg en un point.  Ces variétés ont une structure assez riche 
et elles apparaissent naturellement dans différents domaines des 
mathématiques : notamment en théorie des singularités et en 
cohomologie quantique.   
 
La cohomologie quantique, découverte par les physiciens E.\k Witten, 
R.\k Dijk\-graaf et C.\k Vafa dans la théorie des modèles sigma, a été 
axiomatisée en géométrie algébrique par M.\k Kontsevich et Y.\k Manin 
en $1994$. En $1995$, Y.\k Ruan et G.\k Tian ont donné une version de 
la cohomologie quantique en géométrie symplectique. 
 
Dans ces deux approches toute la difficulté réside dans la définition 
des invariants de Gromov-Witten. Ces invariants \og comptent \fg le 
nombre de courbes de genre et de degré fixé qui vérifient certaines 
conditions d'incidence dans une variété compacte (symplectique ou 
projective).  Si l'on ne considère que les invariants de Gromov-Witten 
qui comptent les courbes de genre $0$ alors on obtient la cohomologie 
quantique et elle est naturellement munie d'une structure de variété 
de Frobenius. Son potentiel, appelé potentiel de Gromov-Witten, est la 
série génératrice formée par ces invariants de Gromov-Witten.  La 
forme bilinéaire non dégénérée plate est donnée par la dualité de Poincaré et le produit, 
qui est une déformation du cup produit, est défini à l'aide du 
potentiel de Gromov-Witten et la dualité de Poincaré.

Une autre source d'exemples de variété de Frobenius vient de la 
théorie des singularités. Dans les années soixante-dix et 
quatre-vingt, K.\k Saito a étudié les déploiements universels des 
singularités isolées d'hypersurfaces.  En $1983$, K. Saito a introduit 
la notion de structure plate, et a conjecturé qu'une telle structure 
existe de manière naturelle sur la base du déploiement universel d'une 
singularité isolée d'hypersurface. Cette conjecture a été montrée par 
M. Saito en $1983$. Par ailleurs, la notion de structure plate s'est 
révélée être identique à la notion de structure de Frobenius, 
considérée plus tard par B. Dubrovin.

La symétrie miroir peut se formuler en termes d'isomorphisme 
entre les variétés de Frobenius provenant de la cohomologie quantique 
(c\^{o}té A) et celle provenant de la théorie des singularités (c\^{o}té B). 
 
Inspir\'e par un article de A.\k Givental \cite{Ghgms}, S.\k
Barannikov a montré en $2000$ dans l'article \cite{Bms} que la variété de Frobenius 
provenant de la cohomologie quantique des espaces projectifs complexes 
est isomorphe à une variété de Frobenius naturelle sur la base du 
déploiement universel du polyn\^{o}me de Laurent $x_{1}+ \cdots+x_{n}+1/x_{1}\dots x_{n}$.

L'objectif de cette thèse est de généraliser ce résultat aux espaces
projectifs à poids. Pour cela, 
nous utilisons la théorie des orbifolds et les constructions qui s'y 
rattachent.  Dans les articles \cite{CRnco} et \cite{CRogw}, W.\k Chen 
et Y.\k Ruan définissent l'anneau de cohomologie orbifolde via les 
invariants de Gromov-Witten orbifolds.  Le cup produit orbifold est 
défini comme la partie de degré zéro du produit quantique orbifold et il se 
calcule via la classe d'Euler d'un fibré obstruction.  Le produit 
 quantique orbifold est défini par le potentiel de Gromov-Witten. 
Ainsi, comme dans le cas des variétés, la cohomologie quantique 
orbifolde est naturellement munie d'une structure de Frobenius.

D'un autre c\^{o}té, A.\k Douai et C.\k Sabbah (cf. \cite{DSgm1}) ont 
expliqué comment construire une variété de Frobenius canonique sur la 
base d'un déploiement universel de tout polyn\^{o}me de Laurent commode 
et non dégénéré par rapport à son polyèdre de Newton.  En particulier, 
dans l'article \cite{DSgm2}, les auteurs ont mis cette construction en 
pratique sur le polyn\^ome $w_{0}u_{0}+ \cdots+ w_{n}u_{n}$ restreint 
à $U:=\{(u_{0}, \ldots ,u_{n})\in \CC^{n+1}\mid \prod_{i} 
u_{i}^{w_{i}}=1\}$ où $w_{0}, \ldots ,w_{n}$ sont des entiers 
strictement positifs et  premiers entre eux. Les travaux de A.\k Douai et C.\k Sabbah ou ceux 
de S.\k Barannikov portent sur le polyn\^{o}me de Laurent lui-m\^{e}me 
et non sur une étude locale de ses singularités. 
  
Dans cette thèse, nous comparons les structures de Frobenius, dont
l'existence est assurée par les résultats généraux rappelés ci-dessus, 
obtenues sur la cohomologie quantique orbifolde des espaces projectifs 
à poids $\PP(w_{0}, \ldots ,w_{n})$ (c\^{o}té A) et celles obtenues à partir du 
polyn\^ome $f(u_{0}, \ldots ,u_{n}):=u_{0}+\cdots + u_{n}$ restreint à 
$U$ (c\^{o}té B).

Nous démontrons d'abord une correspondance entre \og les limites 
classiques \fg. Pour expliquer cela, introduisons quelques notations. Pour 
le c\^{o}té A, nous notons $H^{2\star}_{\orb}(\PP(w_{0}, \ldots 
,w_{n}),\CC)$ la cohomologie orbifolde de $\PP(w_{0}, \ldots ,w_{n})$, 
$\cup$ le cup produit orbifold et $\langle\cdot,\cdot\rangle$ la 
dualité de Poincaré orbifolde. Pour le c\^{o}té B, nous considérons 
l'espace vectoriel $\Omega^{n}(U)/df\wedge\Omega^{n-1}(U)$. Il est 
naturellement muni d'une filtration croissante, appelée filtration de 
Newton et notée $\mathcal{N}_{\bullet}$, et d'une 
forme bilinéaire non dégénérée. Le choix d'une forme volume sur $U$ 
nous permet de définir un produit sur cet espace vectoriel. Comme le 
produit et la forme bilinéaire non dégénérée respectent la filtration 
de Newton, nous avons un produit, noté $\cup$, et une forme 
bilinéaire non dégénérée, notée $[\![g]\!]$, sur le gradué de 
$\Omega^{n}(U)/df\wedge\Omega^{n-1}(U)$ par rapport à la filtration de 
Newton.  Le théorème suivant est démontré au paragraphe 
\ref{section:correspondance}. 
   
\begin{thm}[Correspondance classique] 
\label{thm:correspondance,classique} On a un isomorphisme   
d'algèbres de Frobenius graduées entre
   \begin{align*} 
\left(H^{2\star}_{\orb}(\PP(w),\CC),\cup,\langle\cdot,\cdot\rangle\right) 
\end{align*} et 
   \begin{align*} 
\left(\gr^{\mathcal{N}}_{\star}\left(\Omega^{n}(U)/df\wedge\Omega^{n-1}(U)\right),\cup,[\![g]\!](\cdot,\cdot)\right). 
\end{align*} 
\end{thm} 
Signalons que dans un contexte algébrique et plus général, A.\k 
Borisov, L.\k Chen et G.\k Smith \cite{BCSocdms} ont calculé l'anneau 
de cohomologie orbifolde pour un champ de Deligne-Mumford associé à 
une variété torique simpliciale et, dans le cas des espaces projectifs 
à poids, trouvent le m\^{e}me résultat, sans la forme bilinéaire 
non dégénérée cependant. 
  
 Puis, nous énonçons une conjecture (cf. \ref{conj:invariant,dur}) sur 
 la valeur de certains invariants de Gromov-Witten orbifolds et nous 
 montrons au paragraphe \ref{sec:corres-quantique} que cette 
 conjecture implique un isomorphisme entre les variétés de Frobenius 
 provenant du c\^{o}té A et du c\^{o}té B. 
 
Les deux énoncés de correspondance sont démontrés, modulo la 
conjecture \ref{conj:invariant,dur}, au chapitre \ref{cha:correspondance}. 
Ils utilisent les résultats des chapitres \ref{sec:struct-orbif-des} et 
\ref{cha:la-cohom-quant} pour le c\^{o}té A et les résultats du chapitre 
\ref{cha:les-singularites-du} pour le c\^{o}té B.  
 
Au paragraphe \ref{sec:cote,A} (resp. \ref{sec:cote,B}) de l'introduction, nous 
reviendrons plus en détails sur les contenus des chapitres  \ref{sec:struct-orbif-des} et 
\ref{cha:la-cohom-quant} (resp. \ref{cha:les-singularites-du}). 

\bigskip

\section{Les variétés de Frobenius} 
\label{sec:Les-varits-de-Frob} 
 
Rappelons d'abord la définition d'une variété de Frobenius, introduite 
par Dubrovin dans \cite{Dtft}. 
 
Soit $M$ une variété complexe.  Notons  $\mathcal{O}_M$ le faisceau 
des fonctions holomorphes sur $M$ et $\Theta_M$ le faisceau des 
champs de vecteurs holomorphes sur $M$.
 
\begin{defi}[\cite{Dtft} Lecture $1$, voir aussi \cite{Mfm} p.$19$, \cite{Hfm} p.$146$, \cite{Sdivf} 
  p.$240$] Étant donné  un champ de vecteurs $\mathfrak{E}$, 
  appel\'e champ d'Euler, une forme bilin\'eaire non d\'eg\'en\'er\'ee $g$, un produit 
  $\star$ associatif et commutatif d'\'el\'ement unit\'e $e$ sur le
  fibré tangent complexe $TM$.
 On dit 
  que $(M,g,\star,e,\mathfrak{E})$ est une \emph{structure de Frobenius} si les 
  conditions suivantes sont satisfaites  : 
\begin{enumerate} 
\item \label{item:1} la forme bilinéaire non dégénérée $g$ est plate et $\nabla (e)=0$, 
  o\`u $\nabla$ est la connexion sans torsion associée \`a $g$ ; 
\item \label{item:2} pour tous champs de vecteurs $\xi,\eta,\zeta$, on  
  a 
  \begin{align*} 
    \nabla_{\xi}(\eta\star\zeta)-\eta\star\nabla_{\xi}\zeta-\nabla_{\eta}(\xi\star\zeta)  
    + \xi\star \nabla_{\eta}\zeta-[\xi,\eta]\star \zeta& =0 \,; 
  \end{align*} 
   
\item\label{item:6} pour tous champs de vecteurs $\xi,\eta,\zeta$, on 
  a $g(\xi\star\eta,\zeta)=g(\xi,\eta\star\zeta)$ ;

 \item \label{item:5} pour tous champs de vecteurs 
    $\xi,\eta$, on a ${\mathcal 
      L_\mathfrak{E}}(\xi\star\eta)-{\mathcal 
      L_\mathfrak{E}}\xi\star\eta - \xi\star{\mathcal 
      L_\mathfrak{E}}\eta=\xi\star\eta$.
   
\item \label{item:3} il existe un nombre complexe $D$ tels qu'on a 
      ${\mathcal 
      L}_{\mathfrak{E}}(g(\xi,\eta))-g({\mathcal 
      L_\mathfrak{E}}\xi,\eta) - g(\xi,{\mathcal 
      L_\mathfrak{E}}\eta)=D \cdot g (\xi,\eta)$ pour tous champs de 
    vecteurs $\xi,\eta$ où $\mathcal{L}$ est la dérivée de Lie; 
\end{enumerate} 
\end{defi}

\begin{rem} 
   \begin{enumerate} 
  \item D'après les conditions  \eqref{item:1} et \eqref{item:3}, l'endomorphisme 
    $\nabla\mathfrak{E}$ de $\Theta_M$ est une section 
    $\nabla$-horizontale du faisceau End$_{{\mathcal O}_M}(\Theta_M)$. 
   
  \item Supposons que $M$ soit simplement connexe. Soit 
    $(t_{1}, \ldots ,t_{n})$ un système de coordonnées plates sur $M$. 
    D'après le lemme $1.2$ du cours $1$ de \cite{Dtft} (voir aussi le paragraphe $VII.2.b$ de \cite{Sdivf}), il existe une 
    fonction holomorphe, appelée potentiel, $F:M\rightarrow \CC$ 
    telle que pour tous $i,j,k$ dans $\{1, \ldots ,n\}$, nous ayons 
    \begin{align*}  
    F_{ijk}&:=\frac{\partial^{3}F}{\partial t_{i}\partial t_{j}\partial t_{k}}=g(\partial_{t_{i}}\star\partial_{t_{j}},\partial_{t_{k}}).   
    \end{align*}  
    Le potentiel n'est déterminé qu'à l'addition près d'un 
    polyn\^{o}me de degré~$2$. 
     
    Comme le produit $\star$ est associatif, le potentiel est solution des 
    équations WDVV suivantes, qu'on appelle aussi équations 
    d'associativité, Pour tout $i,j,k,\ell$ dans $\{1, \ldots ,n\}$, 
    nous avons 
\begin{align*} 
 \sum_{a,b}F_{ija} 
g^{ab}F_{bk\ell}&=\sum_{a,b}F_{jka}g^{ab}F_{bi\ell}  \ \ \ \  \forall i,j,k,\ell\in \{1, \ldots ,n\}
\end{align*}  
où $(g^{ab})$ est l'inverse de la matrice de la forme bilinéaire non
dégénérée $g$ dans les 
coordonnées $(t_{1}, \ldots ,t_{n})$.  De plus, ce potentiel vérifie 
la condition d'homogénéité suivante par rapport au champ d'Euler : la fonction $ 
\mathfrak{E} (F) -(D+1)F$ est un polyn\^{o}me de degré inférieur ou
égal à deux. 
 
\item La condition (\ref{item:5}) implique que $[e,\mathfrak{E}]=e$. 
  Si nous prenons $\xi=\zeta=e$ et $\eta=\mathfrak{E}$ dans 
  l'égalité (\ref{item:2}) nous obtenons $\nabla_{e}\mathfrak{E}=e$.  Nous 
  avons aussi les égalités 
\begin{align*} 
  \nabla\mathfrak{E}+(\nabla\mathfrak{E})^{\ast}&=D. 
  \id  \mbox{ et } (\mathfrak{E}  \star)^{\ast}=\mathfrak{E}\star 
\end{align*} où $(\cdot)^{\ast}$ est l'adjoint par rapport à la métrique $g$. 
  \end{enumerate} 
\end{rem} 
 
Pour montrer un isomorphisme entre deux variétés de Frobenius, nous 
utilisons le théorème suivant. 
 
\begin{thm}[\cite{Dtft}, lecture $3$ 
  ; voir aussi \cite{Sdivf}, p.$250$] \label{thm:iso,frobenius,dubrovin}
  Soit $g^{\circ}:\CC^{\mu}\times\CC^{ \mu}\rightarrow \CC$ une forme 
   bilinéaire non dégénérée.  Soit $A_{0}^{\circ}$ une matrice complexe de 
  taille $\mu\times\mu$ semi-simple régulière telle que 
  $(A_{0}^{\circ})^{\ast}=A_{0}^{\circ}$. Soit $A_{\infty}$ une 
  matrice complexe de taille $\mu\times\mu$ telle que $A_{\infty}+A_{\infty}^{ 
    \ast}=w\cdot\id$ avec $w\in \ZZ$.  Soit $e^{\circ}$ un vecteur propre 
  de $A_{\infty}$ pour la valeur propre $q$ tel que $(e^{\circ},
  A_{0}^{\circ}e^{\circ}, \ldots , {A_{0}^{\circ}}^{\mu-1}e^{\circ})$
  soit une base de $\CC^{\mu}$.  Le quadruplet 
  $(A^{\circ}_{0},A_{\infty},g^{\circ},e^{\circ})$ détermine  un unique 
germe de variété de Frobenius 
  $((M,0),\star,e,\mathfrak{E},g)$ (avec la constante 
  $D:=2q+2-w$) tels que via l'isomorphisme entre $T_{0}M$ et 
  $\CC^{\mu}$ on ait $g^{\circ}={g}(0)$, 
  $A^{\circ}_{0}=\mathfrak{E}\star$, 
  $A_{\infty}=(q+1)id-\nabla{\mathfrak{E}}$ et $e^{\circ}=e(0)$. 
\end{thm} 
Signalons que ce théorème a été généralisé par C.\k Hertling et Y.\k 
 Manin dans \cite{HMumc}(cf. théorème $4.5$).    
 
Pour montrer un isomorphisme entre la variété de Frobenius provenant 
de $\PP(w_{0}, \ldots ,w_{n})$ et celle provenant du polyn\^{o}me de 
Laurent $f$, nous allons montrer que leurs conditions initiales 
vérifient les hypothèses du  théorème ci-dessus et qu'elles sont égales. 
 
\section{Le c\^{o}té A} 
\label{sec:cote,A} 
 
Nous construisons la structure de Frobenius sur l'espace vectoriel 
complexe $H^\star_{\orb}(\PP(w_{0}, \ldots ,w_{n}),\CC)$ de dimension 
$\mu:=w_{0}+\cdots+w_{n}$. La forme bilinéaire non dégénérée est 
donnée par la dualité de Poincaré pour les orbifolds, nous la notons 
$\langle \cdot,\cdot \rangle$.  Au paragraphe 
\ref{sec:coho-espace}, nous définirons une base 
$(\eta_{0},\ldots,\eta_{\mu-1})$ de l'espace vectoriel 
$H^\star_{\orb}(\PP(w_{0}, \ldots ,w_{n}),\CC)$.  Notons 
$(t_{0},\ldots,t_{\mu-1})$ les coordonnées sur 
$H^\star_{\orb}(\PP(w_{0}, \ldots ,w_{n}),\CC)$ dans cette base.  Le 
champ d'Euler est donné par la formule suivante 
\begin{align*} 
  \mathfrak{E}&:=\mu\partial_{t_{1}}+\sum_{i=0}^{\mu-1}(1-\sigma(i))t_{i}\partial_{t_{i}} 
\end{align*} 
où $\sigma(i)$ est la moitié du degré orbifold de $\eta_{i}$.  Le 
produit quantique, noté $\star$, est défini à l'aide du potentiel de 
Gromow-Witten orbifold, noté $F^{GW}$, par la formule suivante 
\begin{align*} 
 \frac{\partial^3 F^{GW}(t_{0},\ldots,t_{\mu-1})}{\partial 
 t_{i}\partial t_{j} \partial t_{k}}&=\langle \partial t_{i}\star 
 \partial t_{j}, \partial t_{k}\rangle  
\end{align*} 
Les conditions initiales de la variété de Frobenius sont les données 
$(A_{0}^{\circ},A_{\infty},\langle\cdot,\cdot\rangle,\eta_{0})$ où 
$A^\circ_{0}:=\mathfrak{E}\star\mid_{\mathbf{t}=0}$ et 
$A_{\infty}:=\id-\nabla \mathfrak{E}$ ($q=0$ et $w=n$). On peut calculer facilement la 
matrice $A_{\infty}$, mais pour calculer la matrice $A^\circ_{0}$, il 
faut calculer le cup produit orbifold et certains invariants de 
Gromow-Witten orbifolds. 
 
Le chapitre \ref{chap:cohomologie} est composé de rappels sur les 
orbifolds complexes et commutatives que nous utiliserons aux 
chapitres   \ref{sec:struct-orbif-des} et  \ref{cha:la-cohom-quant}. 
 
Au chapitre  \ref{sec:struct-orbif-des}, nous étudions la 
cohomologie orbifolde de $\PP(w_{0},\ldots ,w_{n})$ et nous en donnons 
une base naturelle. Puis, nous exprimons la dualité de Poincaré, notée 
$\langle \cdot,\cdot\rangle$, et le cup produit orbifold dans cette 
base. En particulier, le fibré obstruction est calculé dans le 
théorème \ref{thm:fibre,obstruction}. 
 
Au chapitre \ref{cha:la-cohom-quant}, nous étudions la 
cohomologie quantique orbifolde des espaces projectifs à poids. Nous 
définissons le champ d'Euler et nous calculons la matrice $A_{\infty}$ 
à la proposition \ref{prop:Potentiel-de-Gromov}.  Puis, nous 
décrivons précisément les invariants de Gromov-Witten orbifolds qui 
nous permettent de calculer la matrice $A_{0}^{\circ}$.  A l'aide du 
cup produit orbifold, du théorème \ref{thm:invariant,0} et de la 
conjecture \ref{conj:invariant,dur}, nous pouvons calculer la matrice $A^\circ_{0}$.

\section{Le c\^{o}té B} 
\label{sec:cote,B} 
 
Soit $U:=\{(u_{0}, \ldots ,u_{n})\in \CC^{n+1}\mid \prod_{i} 
u_{i}^{w_{i}}=1\}$. 
Dans l'article \cite{DSgm2}, le polyn\^{o}me considéré est $w_{0}u_{0}+\cdots+w_{n}u_{n}$ 
restreint à $U$ et les poids sont premiers entre eux dans leur 
ensemble. 
Dans notre cas, nous ne faisons pas  cette hypothèse sur les poids et 
le polyn\^{o}me $f$ est $u_{0}+\cdots+u_{n}$ restreint à $U$. 
Néanmoins, nous pouvons utiliser les m\^{e}mes techniques et nous 
démontrons le théorème suivant. 
 
\begin{thm} Il existe une structure de Frobenius canonique sur tout 
  germe de déploiement universel du polyn\^{o}me de Laurent 
  restriction de $f(u_{0}, \ldots ,u_{n})=u_{0}+\cdots+u_{n}$ à $U$. 
\end{thm} 
 
Le chapitre \ref{cha:les-singularites-du} est consacré à l'étude de 
cette structure de Frobenius associée au polyn\^{o}me de Laurent $f$. 
Au paragraphe \ref{sec:le-systeme-de}, nous calculons les 
conditions initiales de cette structure de Frobenius.  
 
Au deuxième 
paragraphe, nous montrerons que le potentiel de la structure est 
définie par certaines conditions initiales, en fait les m\^{e}mes que 
celles pour le c\^{o}té $A$, si la correspondance est exacte. 
 
Au 
troisième paragraphe, nous considérons l'espace vectoriel 
$\Omega^{n}(U)/df\wedge\Omega^{n-1}(U)$ muni d'une filtration 
croissante, appelée filtration de Newton,  et d'une forme 
bilinéaire non dégénérée. Le choix d'une forme volume sur $U$ nous 
permet de définir un produit sur cet espace vectoriel. 
Nous montrerons que le gradué de cet espace 
vectoriel par rapport à la filtration $V_{\bullet}$ est muni d'une 
structure d'algèbre de Frobenius. 
 
\numberwithin{equation}{section} 
\chapter{Préliminaires combinatoires}\label{cha:prel-comb} 

\section{Notations}\label{sec:notations} 
 
Dans la suite de ce travail, nous utiliserons les notations ci-dessous. 

Soient $n$ et $w_{0},\ldots,w_{n}$ des entiers strictement positifs. 
Posons $\mu=w_{0}+\cdots+w_{n}$. 
Pour tout sous-ensemble $I=\{i_{1}, 
\ldots ,i_{\delta}\}$ de $\{0, \ldots ,n\}$, notons $w_{I}=(w_{i_{1}}, 
\ldots ,w_{i_{\delta}})$.

Pour tout $\gamma\in [0,1[$, posons 
$I(\gamma):=\{i\in\{0,\ldots,n\}\mid \gamma w_{i}\in \NN\}$ et notons 
$\delta(\gamma)$ son cardinal.  Posons $a(\gamma):=\{\gamma w_{0}\}+ 
\cdots + \{\gamma w_{n}\}$ où $\{\cdot\}$ 
désigne la partie fractionnaire.  Nous avons les relations suivantes : 
\begin{align} 
I(\gamma)&=I(\{1-\gamma\}) \,; \nonumber\\ 
\label{eq:fondamentale} 
a(\gamma)+a(\{1-\gamma\})&=n+1-\delta(\gamma)\,;\\ 
\label{eq:triplet} 
\gamma_{0}+\gamma_{1}+\gamma_{\infty}\in \NN & \Leftrightarrow 
e^{2i\pi\gamma_{0}}e^{2i\pi\gamma_{1}}e^{2i\pi\gamma_{\infty}}=1 
\Leftrightarrow \gamma_{\infty}=\{1-\{\gamma_{0}+\gamma_{1}\}\}. 
\end{align} 
 
\section{La combinatoire des nombres $\sigma$}\label{sec:sigma} 
  
Considérons l'ensemble $\bigsqcup_{i=0}^{n}\{\ell/w_{i} \mid 
\ell\in\{0,\ldots,w_{i}-1\}\}$ où $\bigsqcup$ désigne la réunion disjointe. 
Soit $\mathcal{V} $ l'application naturelle  
\begin{align*} 
\bigsqcup_{i=0}^{n}\{\ell/w_{i} \mid \ell\in\{0,\ldots,w_{i}-1\}\}\rightarrow 
\QQ\cap[0,1[  
\end{align*} 
 d'image notée $S_{w}$.  Choisissons une bijection 
\begin{align*} 
s:\{0, \ldots ,\mu-1\} & \rightarrow \bigsqcup_{i=0}^{n}\{\ell/w_{i} \mid 
\ell\in\{0,\ldots,w_{i}-1\}\} 
\end{align*} 
 telle que $\mathcal{V}\circ s$ soit croissante. 
 Pour tout $\gamma$ 
dans $\V$, posons $k^{\max}(\gamma):=\max\{i\in\{0, \ldots ,\mu-1\}\mid s(i)=\gamma\}$.  Nous avons 
\begin{align}\label{eq:kw} 
1+k^{\max}(\gamma)&=\#\{i\in\{0, \ldots ,\mu-1\} \mid s(i)\leq \gamma\}. 
\end{align} 
 
Dans la suite, nous allons utiliser la notation suivante 
\begin{align} 
  \label{eq:min,gamma} 
  k^{\min}(\gamma)&:=\min\{i\in \{0, \ldots ,\mu-1\}\mid s(i)=\gamma\}. 
\end{align} 
 Pour tout $\gamma\in S_{w}$, le cardinal de 
$\mathcal{V}^{-1}(\gamma)$ est $\delta(\gamma)$. 
Il est clair que nous avons l'égalité 
\begin{align} 
  \label{eq:min,egalite} 
  k^{\min}(\gamma)&=k^{\max}(\gamma)-\delta(\gamma)+1. 
\end{align} 
 
\begin{prop}\label{prop:position} 
  Soit $\gamma$ dans $\V$. Alors, on a : 
  \begin{align*} 
  k^{\max}(\gamma)&=n+\left[ \gamma w_{0}\right]+\cdots+ \left[\gamma w_{n} 
  \right] 
  \end{align*} 
  où $[x]$ désigne la partie entière de $x$. 
\end{prop}

\begin{proof} Soit $\gamma\in\V$. 
  Les éléments de $S_{w}$ inférieurs ou égaux à $\gamma$ sont 
  \begin{align*} 
  0,\frac{1}{w_{0}},\ldots,\frac{[\gamma 
    w_{0}]}{w_{0}},0,\frac{1}{w_{1}},\ldots,\frac{[\gamma 
    w_{1}]}{w_{1}},\ldots,0,\frac{1}{w_{n}},\ldots,\frac{[\gamma 
    w_{n}]}{w_{n}}. 
  \end{align*} 
   
  Ainsi, $\#\{i\in\{0, \ldots ,n\}\mid s(i)\leq 
  \gamma\}=n+1+\sum_{i=0}^{n}[\gamma w_{i}]$. Puis, l'égalité 
  (\ref{eq:kw}) démontre la proposition. 
\end{proof} 
 
\begin{cor}\label{cor:position,inverse} 
  Soit $\gamma>0$ dans $\V$. On a l'égalité 
\begin{align*} 
k^{\max}(\gamma)+k^{\max}(\{1-\gamma\})&=n+\mu+\delta(\gamma)-1. 
\end{align*} 
\end{cor}

\begin{proof}[Démonstration du corollaire \ref{cor:position,inverse}] 
  Si $\gamma>0$ alors nous avons $\Inv{\gamma}=1-\gamma$.  La 
  proposition \ref{prop:position} implique que $k^{\max}(\{1-\gamma\})=k 
  (1-\gamma)=n+\mu+\sum_{i=0}^{n}[-\gamma w_{i}]$.  Comme on a 
  \begin{align*} 
  [\gamma w_{i}]+[-\gamma w_{i}]=\begin{cases} 0 & 
      \mbox{si } i\in 
      I(\gamma)\,;\\ 
      -1 & \mbox{sinon,} \end{cases} 
  \end{align*} 
  nous en déduisons que $k^{\max}(\gamma)+k^{\max}(\{1-\gamma\})=2n+\mu-\# 
  I(\gamma)^{c}$ où $I(\gamma)^{c}$ est le complémentaire de 
  $I(\gamma)$ dans $\{0, \ldots ,n\}$. 
 \end{proof} 
  
 Considérons, comme dans l'article \cite{DSgm2}, les nombres rationnels suivants 
  :  
\begin{align*} 
\sigma(i)&=i-\mu s(i) \mbox{ pour } i\in\{0,\ldots,\mu-1\}. 
\end{align*} 
 
\begin{prop}\label{prop:spectre} 
  Soient $\gamma$ dans $\V$ et $d\in\{0, \ldots ,\delta(\gamma)-1\}$. 
  On a les égalités \begin{displaymath}\sigma(k^{\max}(\gamma)-d)=n-(d+a(\gamma)) \mbox{ et } 
  \sigma(k^{\max}(\{1-\gamma\})-d)=\delta(\gamma)-1-d+a(\gamma).\end{displaymath} 
\end{prop} 
 
\begin{rem} \label{rem:symetrie} 
  La proposition précédente et la formule (\ref{eq:fondamentale}) 
  impliquent l'é\-qui\-va\-lence suivante : 
  \begin{align*} 
  \sigma(k^{\max}(\gamma)-d) + \sigma(k^{\max}(\{1-\gamma\})-d')=n 
  &\Leftrightarrow d+d'=\delta (\gamma)-1. 
  \end{align*} 
\end{rem} 
 
\begin{proof}[Démonstration de la proposition \ref{prop:spectre} ] 
  Par définition, nous avons 
  \begin{align*} 
\sigma(k^{\max}(\gamma)-d)&=k^{\max}(\gamma)-d-\sum\gamma w_{i}. 
\end{align*} 
 Puis, la 
  proposition \ref{prop:position} nous donne la première égalité. 
   
  Pour la seconde égalité, il suffit d'appliquer la première partie de 
  la proposition et d'utiliser la formule (\ref{eq:fondamentale}). 
\end{proof}

\chapter{Orbifolds complexes et \emph{commutatives}}\label{chap:cohomologie} 
 
Dans ce chapitre, nous rappelons des définitions 
et des propriétés générales sur les orbifolds complexes et 
\emph{commutatives} (c'est-à-dire que tous les groupes considérés seront \emph{commutatifs}).

\section{Les cartes, atlas et applications orbifolds}\label{sec:les-orbif-comm} 
 
La notion d'orbifold (ou de $V$-variété) a été introduite par Satake 
dans l'article \cite{Sgm}. Nous allons utiliser les notations de Chen 
et Ruan dans leurs articles \cite{CRogw} et \cite{CRnco}.  Dans ce 
paragraphe, nous allons définir les objets que nous utiliserons par la 
suite et nous ne donnerons pas les énoncés les plus généraux. 
Pour plus de détails, les lecteurs pourront consulter l'article 
original de Satake \cite{Sgb} ou les articles plus récents de Chen et 
Ruan \cite{CRnco} et \cite{CRogw} voire celui de Fukaya et Ono 
\cite{FOgwi}.

Dans cette section, nous allons d'abord définir les cartes orbifoldes puis les atlas 
orbifolds et enfin les applications entre orbifolds. 
 
\subsection{Les cartes orbifoldes}\label{subsec:Les-cart-orbif} 

Soit $U$ un espace topologique connexe. Une \emph{carte} de $U$ est 
un triplet $(\widetilde{U},G,\pi)$ où $\widetilde{U}$ est un ouvert 
connexe de $\CC^{n}$, $G$ est un groupe \emph{fini} qui agit de manière 
holomorphe sur $\widetilde{U}$ et $\pi$ une application de 
$\widetilde{U}$ sur $U$ telle que $\pi$ induise un homéomorphisme 
entre $\widetilde{U}/G$ et $U$.  Dans les exemples que nous allons 
considérer par la suite, les groupes seront {commutatifs}. Ainsi, pour 
simplifier, nous supposons dorénavant que \emph{tous les groupes sont 
commutatifs}. Quand nous n'aurons pas besoin de préciser le groupe ou 
la projection, nous noterons simplement $\widetilde{U}$ pour une 
carte de $U$.

Deux cartes $(\widetilde{U}_{1},G_{1},\pi_{1})$ et 
$(\widetilde{U}_{2},G_{2},\pi_{2})$ d'un même ouvert $U$ sont 
\emph{isomorphes} s'il existe un biholomorphisme 
$\varphi:\widetilde{U}_{1}\rightarrow \widetilde{U}_{2}$ et un 
isomorphisme de groupes $\kappa: G_{1}\rightarrow G_{2}$ tels que 
$\varphi$ soit $\kappa$-équivariant et $\pi_{2}\circ\varphi=\pi_{1}$. 
Si $(\varphi,\kappa)$ est un automorphisme d'une carte 
$(\widetilde{U},G,\pi)$ alors il existe $g\in G$ tel que 
$\varphi(x)=g\cdot x$ et $\kappa=\id$. Un tel automorphisme est noté 
$(\varphi_{g},\id)$. L'élément $g$ est unique si le groupe $G$ agit de  
manière effective sur $\widetilde{U}$. Notons $\ker (G)$ le sous-groupe de 
$G$ qui agit trivialement sur $\widetilde{U}$. Remarquons que 
$\varphi_{g}=\varphi_{g'}$ si et seulement si $gg'^{-1}$ est dans $\ker(G)$. 
 
Soit $U$ un ouvert connexe de $U'$. Soit $(\widetilde{U}',G',\pi')$ 
une carte de $U'$. Une carte $(\widetilde{U },G,\pi)$ de $U$ est 
\emph{induite} par $(\widetilde{U}',G',\pi')$ s'il existe un 
monomorphisme de groupes $\kappa: G \rightarrow G'$ et  un 
plongement ouvert $\kappa$-équivariant  $\alpha$ de $\widetilde{U}$ dans 
$\widetilde{U}'$ tels que $\kappa$ induise un isomorphisme entre 
$\ker (G)$ et $\ker (G')$ et $\pi'=\alpha\circ\pi$.  Satake appelle un tel 
couple 
$(\alpha,\kappa):(\widetilde{U},G,\pi)\hookrightarrow(\widetilde{U}',G',\pi')$ 
une \emph{injection} de cartes.  Quand nous n'aurons pas besoin 
d'expliciter $\kappa$, nous noterons simplement une injection par 
$\alpha:\widetilde{U}\hookrightarrow \widetilde{U}'$. 
 
\begin{lem}\label{lem:injections} 
  Soient $(\alpha_{1},\kappa_{1})$ et $(\alpha_{2},\kappa_{2})$ deux 
  injections de $(\widetilde{U},G,\pi)$ dans 
  $(\widetilde{U}',G',\pi')$.  Alors il existe $g'\in G'$ tel que 
  $\alpha_{1}=\varphi_{g'}\circ \alpha_{2}$. 
\end{lem} 
 
\begin{rem}\label{rem:injections} 
  \begin{enumerate} 
  \item\label{item:9} Si $G'$ agit effectivement c'est-à-dire $\ker (G')=\{\id\}$ alors nous 
avons l'unicité de $g'$ dans le lemme ci-dessus. 
 
\item \label{item:10} Le lemme ci-dessus implique que $\kappa_{1}=\kappa_{2}$. 
  \end{enumerate} 
   
\end{rem} 
 
\begin{proof}[Démonstration du lemme \ref{lem:injections}] 
  Pour tout $\widetilde{x}$ dans $\widetilde{U}$, nous avons 
  $\pi'\circ\alpha_{1}(\widetilde{x})=\pi'\circ\alpha_{2}(\widetilde{x})$. 
  Il existe donc $g'\in G'$ tel que 
  $\alpha_{1}(\widetilde{x})=g'\cdot\alpha_{2}(\widetilde{x})$.  Nous 
  en déduisons que 
  \begin{align*} 
\widetilde{U}&=\bigcup_{g' \in G'}\{\widetilde{x}\in\widetilde{U}\mid 
  \alpha_{1}(\widetilde{x})={g'}\cdot\alpha_{2}(\widetilde{x})\}. 
\end{align*} 
  Ainsi, $\widetilde{U}$ est une réunion finie d'ensembles  fermés.  Il existe 
  $g'\in G'$ tel que l'ensemble $E_{g'}:=\{\widetilde{x}\in\widetilde{U}\mid 
  \alpha_{1}(\widetilde{x})=g'\cdot\alpha_{2}(\widetilde{x})\}$ ne 
  soit pas d'intérieur vide. En particulier,   les applications  
  holomorphes $\alpha_{1}$ et 
  $\varphi_{g'}\circ\alpha_{2}$ co\"incident sur un ouvert inclus dans 
  $E_{g'} \subset \widetilde{U}$. 
D'après le  théorème du prolongement analytique, nous en déduisons que  
$\alpha_{1}=\varphi_{g'}\circ\alpha_{2}$ sur $\widetilde{U}$ car 
$\widetilde{U}$ est connexe. 
\end{proof} 
 
Le lemme suivant est d\^{u} à Chen et Ruan (cf. lemme $4.1.1$ dans 
l'article \cite{CRogw}). 
 
\begin{lem}\label{lem:ruan} Soit $(\widetilde{U}',G',\pi')$ une 
  carte de $U'$. Si $U$ est un ouvert connexe de $U'$, alors il existe 
  une  carte de $U$, unique à isomorphisme près, qui s'injecte dans 
  $(\widetilde{U}',G',\pi')$. 
\end{lem}

\begin{proof} 
  \emph{Existence :} Soit $\widetilde{U}$ une composante connexe de 
  $\pi'^{-1}(U')$.  Soit $G$ le sous-groupe de $G'$ formé des éléments 
  qui laissent stable $\widetilde{U}$.  Soit $\pi$ la restriction de 
  $\pi'$ à $\widetilde{U}$. Nous avons ainsi construit une carte 
  $(\widetilde{U},G,\pi)$ de $U$. 
   
  \emph{Unicité :} 
\begin{itemize} 
   
\item Nous allons d'abord montrer que deux composantes connexes 
  différentes de $\pi'^{-1}(U')$ donnent deux cartes 
  isomorphes.  Soient $\widetilde{U}_{1}$ et $\widetilde{U}_{2}$ deux 
  composantes connexes de $\pi'^{-1}(U')$.  Nous en déduisons deux 
  cartes $(\widetilde{U}_{1},G_{1},\pi_{1})$ et 
  $(\widetilde{U}_{2},G_{2},\pi_{2})$.  Pour tout $g\in G$, 
  l'application de 
  ${\varphi_{g}}{\mid_{\pi'^{-1}(U')}}:\pi'^{-1}(U')\rightarrow 
  \pi'^{-1}(U')$ qui à $\widetilde{x}$ associe $g\widetilde{x}$ est un 
  homéomorphisme.  Ainsi, il existe $g\in G$ tel que $ 
  \varphi_{g}(\widetilde{U}_{1})=\widetilde{U}_{2}$.  Finalement 
  $(\varphi_{g},\id)$ est un isomorphisme entre 
  $(\widetilde{U}_{1},G_{1},\pi_{1})$ et 
  $(\widetilde{U}_{2},G_{2},\pi_{2})$. 
   
\item Soit $(\widetilde{U}_{1},G_{1},\pi_{1})$ une carte de $U$ qui 
  s'injecte dans $(\widetilde{U}',G',\pi')$.  Par définition, il 
  existe un couple $(\alpha,\kappa)$ tel que : 
\begin{itemize} 
\item[(i)] $\kappa$ soit un monomorphisme de groupes de $G_{1}$ dans 
  $G'$ tel que $\ker (G_{1})=\ker (G')$ ;
\item[(ii)] $\alpha$ soit un plongement $\kappa$-équivariant de 
  $\widetilde{U}_{1}$ dans $\widetilde{U}'$. 
 \end{itemize} 
 Comme $\widetilde{U}_{1}$ est connexe, $\alpha(\widetilde{U}_{1})$ 
 est contenu dans une composante connexe de $\pi'^{-1}(U)$. Notons 
 $\widetilde{U}_{2}$ cette composante connexe. Soit 
 $(\widetilde{U}_{2},G_{2},\pi'{\mid_{\widetilde{U}_{2}}})$ la carte 
 de $U$ construite dans la partie existence.  Nous avons le diagramme 
 commutatif suivant : 
 \begin{align*} 
 \xymatrix{ 
   G_{1} \ar [r]^-{\kappa} & G' \supset G_{2} \\ 
   \widetilde{U}_{1}\ar@{^{(}->}[r]^-{\alpha} \ar [rd]_-{\pi_{1}} & 
   \widetilde{U}_{2} \subset \widetilde{U}' \ar [d]^-{\pi'} \\ 
   & U\subset U' } 
 \end{align*}  
\end{itemize} 
Nous allons montrer que $\alpha(\widetilde{U}_{1})=\widetilde{U}_{2}$ 
par un argument de connexité. Par construction, nous avons 
$\alpha(\widetilde{U}_{1})$ ouvert dans $\widetilde{U}_{2}$. Il reste 
à montrer que $\alpha(\widetilde{U}_{1})$ est fermé dans 
$\widetilde{U}_{2}$. Soit $(\widetilde{x}_{n})_{n\in\NN}$ une suite 
dans $\widetilde{U}_{1}$ telle que $\alpha(\widetilde{x}_{n})$ 
converge vers $\widetilde{y}\in \widetilde{U}_{2}$.  Il suffit de montrer que 
$\widetilde{y}\in \alpha(\widetilde{U}_{1})$. La suite 
$\pi_{1}(\widetilde{x}_{n})=\pi'(\alpha(\widetilde{x}_{n}))$ converge 
vers $\pi'(\widetilde{y})\in U'$. Ainsi, il existe $\widetilde{x}\in 
\pi_{1}^{-1}(\pi'(\widetilde{y}))$ et il existe une suite 
$(\widetilde{x}'_{n})_{n\in\NN}$ dans $\widetilde{U}_{1}$ telle que 
 
\begin{itemize} 
\item[(i)] $\pi_{1}(\widetilde{x}'_{n})=\pi_{1}(\widetilde{x}_{n})$ ; 
\item[(ii)] $(\widetilde{x}'_{n})_{n\in\NN}$ converge vers 
  $\widetilde{x}$. 
\end{itemize} 
Comme $\widetilde{x}_{n}$ et $\widetilde{x}'_{n}$ sont dans la m\^eme 
orbite, il existe $g_{n}\in G_{1}$ tel que 
$g_{n}\widetilde{x}'_{n}=\widetilde{x}_{n}$ pour tout $n$. Comme 
$G_{1}$ est un groupe fini, quitte à extraire une sous-suite, nous 
pouvons supposer que $g_{n}=g$. Puis, nous avons 
$\alpha(\widetilde{x}_{n})=\alpha(g\widetilde{x}'_{n}) 
=\kappa(g)\alpha(\widetilde{x}'_{n})$. En passant à la limite nous 
obtenons l'égalité $\widetilde{y}=\kappa(g)\alpha(\widetilde{x})\in 
\alpha (\widetilde{U}_{1})\subset \widetilde{U}_{2}$, ce qui prouve que 
$\widetilde{U}_{2}=\alpha_{1}(\widetilde{U}_{1})$. 
  
Il reste à montrer que $\kappa(G_{1})=G_{2}$.  Comme $\kappa(G_{1})$ 
stabilise $\alpha(\widetilde{U}_{1})=\widetilde{U}_{2}$, nous en 
déduisons que $\kappa(G_{1})\subset G_{2}$.  Inversement, soit 
$g_{2}\in G_{2}$. Pour tout $\widetilde{x}\in \widetilde{U}_{1}$, il 
existe $g_{1}\in G_{1}$ tel que 
$g_{2}\alpha(\widetilde{x})=\alpha(g_{1}\widetilde{x})=\kappa(g_{1})\alpha(\widetilde{x})$. 
Nous en déduisons que $\widetilde{U}_{1}$ est la réunion sur 
$G_{1}\times G_{2}$ des ensembles fermés 
\begin{align*} 
\{\widetilde{x}\in \widetilde{U}_{1}\mid 
g_{2}\alpha(\widetilde{x})&=\kappa(g_{1})\alpha(\widetilde{x}) \}. 
\end{align*} 
Le m\^{e}me raisonnement que dans la démonstration du lemme 
\ref{lem:injections} implique que   
$g_{2}\alpha=\kappa(g_{1})\alpha$.  En d'autres termes, 
$g_{2}^{-1}\kappa(g_{1})$ appartient à $\ker (G_{2})=\ker (G')\simeq 
\ker (G_{1})$.  Nous obtenons que $g_{2}$ est dans $\kappa(G_{1})$. 
\end{proof} 
 
\begin{cor}\label{cor:injections} 
  Soient $(\widetilde{U},G_{U},\pi_{U}), 
  (\widetilde{V},G_{V},\pi_{V}),(\widetilde{W},G_{W},\pi_{W})$ 
  trois car\-tes de respectivement $U,V,W$ telles qu'il existe deux 
  injections $\widetilde{U}\hookrightarrow \widetilde{W}$ et 
  $\widetilde{V}\hookrightarrow \widetilde{W}$.  Si $U$ est inclus dans $V$ alors 
  il existe une injection $\widetilde{U}\hookrightarrow\widetilde{V}$.  
\end{cor} 
 
\begin{proof} 
Comme $U\subset V$, d'après le lemme \ref{lem:ruan} il existe une 
carte $\widetilde{U}_{V}$ de $U$ qui s'injecte dans $\widetilde{V}$. 
Par hypothèse, nous avons deux cartes $\widetilde{U}$ et 
$\widetilde{U}_{V}$ de $U$ qui s'injectent dans $\widetilde{W}$. Ainsi, 
le lemme \ref{lem:ruan} montre que ces deux cartes $\widetilde{U}$ et 
$\widetilde{U}_{V}$ sont isomorphes. Nous en déduisons le corollaire. 
\end{proof} 
\subsection{Les atlas orbifolds}\label{sec:Les-atlas-orbifolds} 
 
Soit $|X|$ un espace topologique. D'après l'article \cite{MPos}, un 
\emph{atlas orbifold} $\mathcal{A}(|X|)$ de $|X|$ est la donnée d'un 
recouvrement de $|X|$ par des ouverts connexes $(U_{i})_{i\in I}$ tels 
que

 \begin{enumerate}\makeatletter 
\renewcommand\theenumi{\theequation} 
\makeatother 
 \addtocounter{equation}{1}  \item\label{enu:atlas,1} chaque ouvert $U_{i}$ de ce recouvrement 
    ait une carte $(\widetilde{U}_{i},G_{i},\pi_{i})$ ; 
 \addtocounter{equation}{1} \item\label{enu:atlas,2} pour tout $x\in U_{i}\cap U_{j}$, il existe 
    un ouvert $U_{k}\subset U_{i}\cap U_{j}$ contenant $x$ et deux injections 
    $\widetilde{U}_{k}\hookrightarrow \widetilde{U}_{i}$ et 
     $\widetilde{U}_{k}\hookrightarrow \widetilde{U}_{j}$. 
   \end{enumerate} 
   
  D'après le lemme \ref{lem:ruan}, nous pouvons toujours affiner 
  l'atlas orbifold c'est-à-dire rajouter toutes les cartes induites.  Deux 
  atlas orbifolds sont dits équivalents s'il existe un troisième atlas 
  orbifold plus fin que chacun d'eux.

\begin{defi}\label{defi:orbifold} 
  Une \emph{orbifold} est un espace topologique séparé,   muni d'une 
  classe d'équivalence d'atlas orbifold.  Pour alléger les notations, nous 
  notons simplement $X$ pour l'orbifold $(|X|,[\mathcal{A}(|X|)])$. 
\end{defi}

Nous dirons qu'une orbifold est compacte, connexe, ... si son espace 
topologique l'est. 
 
 \begin{lem}\label{lem:isotropie} 
   Soit $X$ une orbifold. Soient $x$ dans $|X|$ et 
   $(\widetilde{U},G,\pi)$ une carte d'un voisinage $U$ de $x$. Soit 
   $\widetilde{x}$ un relevé de $x$ dans $\widetilde{U}$. Le groupe 
   $\{g\in G| g\cdot \widetilde{x}=\widetilde{x}\}$ ne dépend que de 
   $x$. 
\end{lem} 
 
Ce groupe est appelé \emph{groupe d'isotropie} au point $x$ et nous 
le notons $G_{x}$. 
 
 \begin{proof}[Démonstration du lemme \ref{lem:isotropie}] 
   Notons $G^{\widetilde{U}}_{\widetilde{x}}$ le sous-groupe de $G$ 
   défini par  
\begin{align*} 
\{g\in G| g\cdot \widetilde{x}&=\widetilde{x}\}. 
\end{align*} 
Montrons que
$G^{\widetilde{U}}_{\widetilde{x}}=G^{\widetilde{U}}_{g\widetilde{x}}$
pour tout $g$ dans $G$.  Soit $h$ dans
$G^{\widetilde{U}}_{\widetilde{x}}$. Comme $G$ est commutatif, nous
avons $h(g\widetilde{x})=g\widetilde{x}$ c'est-à-dire $h\in
G^{\widetilde{U}}_{g\widetilde{x}}$.  L'autre inclusion est directe.
Notons $G_{x}^{\widetilde{U}}:=G_{\widetilde{x}}^{\widetilde{U}}$.
     
   Montrons que s'il existe une injection 
   $(\alpha,\kappa):(\widetilde{V},G_{\widetilde{V}},\pi_{\widetilde{V}})\hookrightarrow 
   (\widetilde{U},G_{\widetilde{U}},\pi_{\widetilde{U}})$ alors nous avons
   $G_{x}^{\widetilde{U}}=G_{x}^{\widetilde{V}}$.  D'après le lemme 
   \ref{lem:ruan}, nous pouvons supposer que $\widetilde{V}$ est une 
   composante connexe de $\pi_{\widetilde{U}}^{-1}(V)$, 
   $G_{\widetilde{V}}$ est le sous-groupe de $G_{\widetilde{U}}$ qui 
   agit sur $\widetilde{V}$ et 
   $\pi_{\widetilde{V}}={\pi_{\widetilde{U}}}{\mid_{\widetilde{V}}}.$ 
   Il est clair que $G_{x}^{\widetilde{V}}\subset 
   G_{x}^{\widetilde{U}}$.  Inversement, soit $g\in 
   G_{x}^{\widetilde{U}}$. Soit $\widetilde{x}$ un relevé de $x$ dans 
   $\widetilde{V}$. Il suffit de montrer que $g\in G_{\widetilde{V}}$. 
   L'application 
   ${\varphi_{g}}{\mid_{\pi_{\widetilde{V}}^{-1}(U)}}:\pi_{\widetilde{V}}^{-1}(U)\rightarrow 
   \pi_{\widetilde{V}}^{-1}(U)$ qui à $\widetilde{y}$ associe 
   $g\widetilde{y}$ est un homéomorphisme donc elle envoie les 
   composantes connexes sur les composantes connexes.  Or nous avons
   $g\widetilde{x}=\widetilde{x}$, donc 
   $\varphi_{g}(\widetilde{V})=\widetilde{V}$, c'est-à-dire $g\in 
   G_{\widetilde{V}}$. 
    
   Soit $x$ dans $U_{1}\cap U_{2}$. Comme $X$ est une orbifold, il 
   existe une carte 
   $(\widetilde{V},G_{\widetilde{V}},\pi_{\widetilde{V}})$ d'un ouvert 
   $V\subset U_{1}\cap U_{2}$ contenant $x$ et il existe deux 
   injections $(\widetilde{V},G_{\widetilde{V}},\pi_{\widetilde{V}}) 
   \hookrightarrow (\widetilde{U}_{1},G_{1},\pi_{1}) $ et 
   $(\widetilde{V},G_{\widetilde{V}},\pi_{\widetilde{V}}) 
   \hookrightarrow (\widetilde{U}_{2},G_{2},\pi_{2})$.  Nous en 
   déduisons que 
   $G_{x}^{\widetilde{V}}=G_{x}^{\widetilde{U}_{1}}=G_{x}^{\widetilde{U}_{2}}$. 
\end{proof}  
 
\begin{rem}\label{rem:germe} 
  Soit $X$ une orbifold. Soit $x$ un point de $X$. D'après le lemme 
  \ref{lem:ruan}, quitte à prendre un ouvert $U_{x}$, contenant $x$, assez
  petit, il existe une carte $(\widetilde{U}_{x},G_{x},\pi_{x})$ de $U_{x}$. 
\end{rem}

\begin{lem}\label{lem:noyau,action} 
  Soit $X$ une orbifold connexe. Le groupe $\ker (G_{x})$ ne dépend pas 
   du point $x$ dans $|X|$. 
\end{lem} 
 
Notons $\ker (X)$ ce groupe qui agit globalement trivialement. 
 
\begin{proof}[Démonstration du lemme \ref{lem:noyau,action}] 
  Soit $(\widetilde{U},G,\pi)$ une carte de $U$.  D'après le lemme 
  \ref{lem:ruan} et la définition d'une injection, pour tous $x,y$ 
  dans $U$, nous avons $\ker (G_{x})=\ker (G_{y})=\ker (G)$.  Ceci montre 
  que l'ensemble $\{p\in |X|\mid \ker (G_{p}) \simeq\ker (G)\}$ est ouvert et 
  fermé dans $|X|$. 
\end{proof}

Une orbifold est dite \emph{réduite} si $\ker (X)$ est réduit à 
l'identité. 
A toute carte $(\widetilde{U},G,\pi)$ de $U$ d'un atlas 
$\mathcal{A}(|X|)$, nous lui associons la carte réduite 
\begin{align*} 
  (\widetilde{U},G/\ker(X),\pi_{\red}) 
\end{align*} 
de $U$ où 
$\pi_{\red}:\widetilde{U}\rightarrow U$ induit un homéomorphisme entre  
$\widetilde{U}/(G/\ker(X))\simeq \widetilde{U}/G$ et $U$. 
Ainsi, à l'atlas $\mathcal{A}(|X|)$,  nous lui associons un unique atlas 
réduit. Nous obtenons alors une orbifold réduite noté $X_{\red}$.

La partie \emph{régulière} de $X$, notée $X_{\reg}$, est $\{x\in |X| 
\mid G_{x}=\ker (X)\}$. Remarquons que $X_{\reg}$ est une variété 
complexe munie d'une action triviale du groupe $\ker (X)$.  Cette 
définition est différente de celle de Chen et Ruan (cf. définition 
$4.1.2$ de \cite{CRogw}) où ils définissent  
$\widehat{X}_{\reg}:=\{x\in |X| \mid G_{x}=\{\id\}\}$. 
 
\begin{expl}\label{expl:orbifold}  
 \begin{enumerate} 
  \item\label{item:15} Une variété complexe $X$ est une orbifold où 
    nous avons $\ker(X)=\{\id\}$, 
    $G_{x}=\{\id\}$ pour tout $x\in X$ et 
    $X_{\reg}=\widehat{X}_{\reg}=X$. 
  \item  
Soit $G$ un groupe commutatif fini qui agit trivialement sur une variété 
    complexe $Y$. Le quotient $X:=Y/G$ est naturellement munie d'une 
    structure orbifolde. 
 Nous avons $\ker(X)=G$, $G_{x}=G$ pour tout $x\in X$, 
    $X_{\reg}=Y$ mais $\widehat{X}=\{\emptyset \}$. 
 
  \item\label{item:16} Notons $D$ le disque unité ouvert de $\CC$. 
Soit $\pi:\widetilde{D}\rightarrow D$ l'application qui à $z$ associe 
$z^{n}$ où $\widetilde{D}=D$. Le triplet 
$(\widetilde{D},\bs{\mu}_{n},\pi)$, où $\zeta\cdot z:=\zeta z$ pour tout 
$(\zeta,z)\in \bs{\mu}_{n}\times \widetilde{D}$, est une carte de $D$. L'ensemble 
formé de la seule carte $(\widetilde{D},\bs{\mu}_{n},\pi)$ est un atlas 
orbifold. Nous avons 
\begin{itemize} 
\item $G_{z}=\{\id\}$ sauf pour $z=0$ où 
$G_{0}=\bs{\mu}_{n}$; 
\item  $D_{\reg}=\widehat{D}_{\reg}=D-\{0\}$. 
\end{itemize} 
  
\item\label{item:17} Soit $\PP^{1}$ la droite projective complexe. Soient 
  $U_{0}:=\{[x,y]\mid x\neq 0\}$ et 
  $U_{1}:=\{[x,y]\mid y\neq 0\}$. 
Soit $(\widetilde{U}_{0},\bs{\mu}_{w_{0}},\pi_{0})$ 
(resp. $(\widetilde{U}_{1},\bs{\mu}_{w_{1}},\pi_{1})$) la carte de 
$U_{0}$ (resp. $U_{1}$) définie par $\widetilde{U}_{0}=\CC$ 
(resp. $\widetilde{U}_{1}=\CC$), $\zeta\cdot z=\zeta z$ pour tout 
$(\zeta,z)\in \bs{\mu}_{w_{0}}\times \widetilde{U}_{0}$ 
(resp.  $\zeta\cdot t=\zeta t$ pour tout $(\zeta,t)\in \bs{\mu}_{w_{1}}\times \widetilde{U}_{}$) et 
$\pi_{0}(z)=z^{w_{0}}$ (resp. $\pi_{1}(t)=t^{w_{1}}$). 
Soit $U$ un ouvert connexe de $\PP^{1}$. Une carte 
$(\widetilde{U},G_{\widetilde{U}},\pi_{\widetilde{U}})$ de $U$ 
est dite admissible s'il existe $i\in\{0,1\}$ tel que  
\begin{itemize} 
\item $\widetilde{U}$ est une composante connexe de $\pi_{i}^{-1}(U)$ ; 
\item $G_{\widetilde{U}}$ est le sous-groupe de $\bs{\mu}_{w_{i}}$ qui 
  agit sur $\widetilde{U}$, c'est-à-dire que
  $G_{\widetilde{U}}:=\{g\in \bs{\mu}_{w_{i}}\mid g\widetilde{U}\subset \widetilde{U}\}$ ; 
\item $\pi_{\widetilde{U}}:=\pi_{i}\mid_{\widetilde{U}}$. 
\end{itemize} 
 
\begin{prop}\label{prop:P1,twist} 
  L'ensemble des cartes admissibles est un atlas orbifold. Notons 
  $\PP^{1}_{w_{0},w_{1}}$ l'orbifold ainsi construite. 
\end{prop} 
 
\begin{proof}  Le point (\ref{enu:atlas,1}) est évident.  Nous allons montrer le 
  point (\ref{enu:atlas,2}).  Soit $x\in V\cap W$.  
   
  Supposons que $\widetilde{V}$ et $\widetilde{W}$ soient dans 
  $\widetilde{U}_{0}$.  Soit $\widetilde{V\cap W}$ une 
  composante connexe de $\pi_{0}^{-1}(V\cap W)$ qui contient 
  un relevé de $x$. Soit $G_{\widetilde{V\cap W}}$ le 
  sous-groupe de $\bs{\mu}_{w_{0}}$ qui agit sur $\widetilde{V\cap 
    W}$. Ainsi, le triplet $(\widetilde{V\cap 
    W},G_{\widetilde{V\cap W}},\pi_{\widetilde{V\cap 
      W}})$ est une carte de $V\cap W$ où 
  $\pi_{\widetilde{V\cap W}}=\pi_{0}\mid_{\widetilde{V\cap 
      W}}$. Il est clair qu'il existe $g_{0},g_{1} \in 
  \bs{\mu}_{w_{0}}$ tels que $\varphi_{g_{0}}:\widetilde{V\cap 
    W}\hookrightarrow \widetilde{V}$ et 
  $\varphi_{g_{1}}:\widetilde{V\cap W}\hookrightarrow 
  \widetilde{W}$. 
   
  Supposons que $\widetilde{V}\subset\widetilde{U}_{0}$ et 
  $\widetilde{W}\subset\widetilde{U}_{1}$. Comme les applications 
  $\pi_{i}\mid_{\widetilde{U}_{i}-\{0\}}$ sont des rev\^{e}tements, il existe un disque 
  ouvert $D(x,\varepsilon)$ tel que pour $i\in\{0,1\}$,  
$\pi_{i}^{-1}(D(x,\varepsilon))$ soit la réunion disjointe de $w_{i}$ 
disques. Tous ces disques sont biholomorphes à $D(x,\varepsilon)$. 
Soient $\widetilde{D}_{V}$ (resp. $\widetilde{D}_{W}$) une composante connexe de 
  $\pi_{0}^{-1}(D(x,\varepsilon))\cap\widetilde{V}$ (resp. $\pi_{1}^{-1}(D(x,\varepsilon))\cap\widetilde{W}$). 
Ainsi, les triplets  $(\widetilde{D}_{V},\id,\pi_{0}\mid_{\widetilde{D}_{V}})$ et 
  $(\widetilde{D}_{W},\id,\pi_{1}\mid_{\widetilde{D}_{W}})$ sont des 
  cartes de $D(x,\varepsilon)$ et elles sont isomorphes. De plus, 
  elles s'injectent dans respectivement $\widetilde{V}$ et
  $\widetilde{W}$. 
\end{proof} 

Cet exemple est important car les cartes induites sur $U_{0}\cap 
U_{1}$ par celles de $U_{0}$ et $U_{1}$ ne sont pas isomorphes : 
elles sont respectivement 
$(\widetilde{U}_{0}-\{0\},\bs{\mu}_{w_{0}},\pi_{0}\mid_{\widetilde{U}_{0}-\{0\}})$ 
et 
$(\widetilde{U}_{1}-\{0\},\bs{\mu}_{w_{1}},\pi_{1}\mid_{\widetilde{U}_{1}-\{0\}})$. 
Ainsi, pour vérifier la condition (\ref{enu:atlas,2}) de la définition 
d'un atlas orbifold, il faut prendre des ouverts assez petits, ce qui 
complique  les démonstrations. 
  \end{enumerate} 
\end{expl}

\subsection{Les applications entre orbifolds} 
\label{subsec:Applications,orbifolds} 
 
Soient $(\widetilde{U},G,\pi)$ et $(\widetilde{U}',G',\pi')$ deux 
cartes de respectivement $U$ et $U'$.  Soit $f$ une application 
continue de $U$ dans $U'$.  Un \emph{relèvement} holomorphe (resp. 
$C^{\infty}$) de $f$ est une application 
$\widetilde{f}:\widetilde{U}\rightarrow \widetilde{U}'$ holomorphe 
(resp. $C^{\infty}$) telle que 
\begin{enumerate}   
\item\label{enu:relevement,1} $\pi'\circ \widetilde{f}=f\circ \pi$ ; 
\item\label{enu:relevement,2} pour tout $g$ dans $G$, il existe $g'$ 
  dans $G'$ tel que pour tout $\widetilde{x}$ dans $\widetilde{U}$ on 
  ait $g'\cdot 
  \widetilde{f}(\widetilde{x})=\widetilde{f}(g\cdot\widetilde{x})$. 
  \end{enumerate}

\begin{lem}Soient  $(\widetilde{U},G,\pi)$ et  $(\widetilde{U}',G',\pi')$ deux 
  cartes de respectivement $U$ et $U'$. Soit $\widetilde{f}$ un relèvement 
  d'une application continue $f:U\rightarrow U'$. Supposons que 
  l'action de $G'$ soit effective. Alors il existe un morphisme de 
  groupe $\kappa:G\rightarrow G'$ tel que $\widetilde{f}$ soit 
  $\kappa$-équivariante. 
\end{lem} 
 
 \begin{proof} 
   Comme l'action de $G'$ est effective, pour tout $g\in G$ il 
   existe un unique $g'\in G'$ tel que 
   $g'\widetilde{f}(\widetilde{x})=\widetilde{f}(g\widetilde{x})$. 
   Posons $\kappa(g):=g'$. Ceci définit un morphisme de groupes 
   $\kappa:G\rightarrow G'$. 
\end{proof} 
 
Deux relèvements $\widetilde{f}_{1}$ et $\widetilde{f}_{2}$ sont 
\emph{isomorphes} s'il existe deux isomorphismes de cartes 
$(\varphi,\kappa)$ et $(\varphi',\kappa')$ tels que 
$\varphi'\circ\widetilde{f}_{1}=\widetilde{f}_{2}\circ \varphi$.

Soit $\widetilde{f}:\widetilde{U}\rightarrow \widetilde{U}'$ un 
relèvement de $f:U\rightarrow U'$.  Soient $V\subset U$ et $V'\subset 
U'$ deux ouverts connexes tels que $f{\mid_{V}}:V \rightarrow V'$. 
Un relèvement $h:\widetilde{V}\rightarrow \widetilde{V}'$ de 
$f{\mid_{V}}$ est \emph{induit} par le relèvement $\widetilde{f}$ 
si pour toute injection $(\alpha,\kappa):(\widetilde{V},G_{V},\pi_{V}) 
\rightarrow(\widetilde{U},G,\pi)$ il existe une injection 
$(\alpha',\kappa'):(\widetilde{V}',G_{V'},\pi_{V'}) 
\rightarrow(\widetilde{U}',G',\pi')$ telle que 
$\widetilde{h}=(\alpha')^{-1}\circ\widetilde{f} \circ \alpha$. 
 
\begin{lem}\label{lem:relevement,induit} 
  Soit $\widetilde{f}:\widetilde{U}\rightarrow \widetilde{U}'$ un 
  relèvement de $f:U\rightarrow U'$. Soient $V\subset U$ et $V'\subset 
  U'$ tel que $f{\mid_{V}}:V\rightarrow V'$. Alors il existe un 
  unique, à isomorphisme près, relèvement de $f{\mid_{V}}$ induit par 
  $\widetilde{f}$. 
\end{lem} 
 
\begin{proof}\emph{Existence:} Soient $(\widetilde{U},G,\pi)$ une 
  carte de $U$ et $(\widetilde{U}',G',\pi')$ une 
  carte de $U'$. Soient $V\subset U$ et $V'\subset 
  U'$ tels que $f{\mid_{V}}:V\rightarrow V'$.  Soit $\widetilde{V}$ 
  une composante connexe de $\pi^{-1}(V)$. Ainsi, l'ensemble  
  $\widetilde{f}(\widetilde{V})$ est contenu dans une unique 
  composante connexe, notée $\widetilde{V}'$, de $\pi'^{-1}(V')$. 
  Montrons que l'application 
  $\widetilde{f}{\mid_{\widetilde{V}}}:\widetilde{V}\rightarrow 
  \widetilde{V}'$ est un relèvement induit de $\widetilde{f}$. Soit 
  $(\alpha,\kappa):(\widetilde{V},G_{V},\pi_{V}) 
  \hookrightarrow(\widetilde{U},G,\pi)$ une injection. D'après le
  lemme~\ref{lem:injections},
 il existe $g\in G$ tel que 
  $\alpha(\widetilde{x})=g\cdot\widetilde{x}$ et d'après la définition 
  d'un relèvement, il existe $g'\in G'$ tel que 
  $\widetilde{f}(g\cdot\widetilde{x})=g'\cdot 
  \widetilde{f}{\mid_{\widetilde{V}}}(\widetilde{x})$. 
   
  \emph{Unicité:} D'après le lemme \ref{lem:ruan}, nous pouvons 
  supposer que $\widetilde{f}_{1}$ et $\widetilde{f}_{2}$ sont deux 
  relèvements entre les m\^{e}mes cartes.  Ainsi, nous avons 
  $f{\mid_{V}}\circ\pi=\pi'\circ 
  \widetilde{f}_{1}=\pi'\circ\widetilde{f}_{2}$. Pour tout 
  $\widetilde{x}$ dans $\widetilde{V}$ il existe $g'\in G'$ tel que 
  $\widetilde{f}_{1}(\widetilde{x})=g'\cdot\widetilde{f}_{2}(\widetilde{x})$. 
  Nous en déduisons que $\widetilde{V}$ est la réunion sur les 
  éléments  $g'$ de $G'$ des ensembles fermés 
\begin{align*} 
\{\widetilde{x}\in \widetilde{V}\mid \widetilde{f}_{1}(\widetilde{x})=g'\cdot\widetilde{f}_{2}(\widetilde{x}) \}. 
\end{align*} 
Le m\^{e}me raisonnement que dans la démonstration du lemme \ref{lem:injections} montre que  
$\widetilde{f}_{1}=g'\cdot \widetilde{f}_{2}$ sur $\widetilde{V}$. 
\end{proof} 
  
 \begin{defi} Soient $X$ et $Y$ 
   deux orbifolds. Une application holomorphe orbifolde entre $X$ et 
   $Y$ est une application continue $|f|:|X|\rightarrow |Y|$ telle que, pour 
   tout point $x$ de $|X|$ il existe une carte 
   $(\widetilde{U}_{x},G_{x},\pi_{x})$ d'un ouvert $U_{x}$ qui 
   contient $x$ et une carte 
   $(\widetilde{U}_{|f|(x)},G_{|f|(x)},\pi_{|f|(x)})$ d'un ouvert 
   $U_{|f|(x)}$ qui contient $|f|(x)$, satisfaisant aux propriétés 
   suivantes : 
\begin{enumerate}  
\item $|f|(U_{x})$ est inclus dans $U_{|f|(x)}$; 
\item il existe un relèvement 
  $\widetilde{f}_{x}:\widetilde{U}_{x}\rightarrow 
  \widetilde{U}_{|f|(x)}$ de $f{\mid_{U}}$; 
\item si $y\in\pi_{x}(\widetilde{U}_{x})$ alors $\widetilde{f}_{x}$ et 
  $\widetilde{f}_{y}$ induisent des relèvements isomorphes sur un 
  voisinage de $y$. 
   \end{enumerate} 
\end{defi}

Deux applications orbifoldes $f_{1},f_{2}:X\rightarrow Y$ sont 
isomorphes si $|f_{1}|=|f_{2}|$ et si pour tout $x$ dans $|X|$, 
$\widetilde{f}_{1,x}$ et $\widetilde{f}_{2,x}$ induisent des 
relèvements isomorphes sur un voisinage de $x$.

Une application orbifolde $f:X\rightarrow Y$ est \emph{un plongement orbifold} 
si les relève\-ments $\widetilde{f}_{x}$ sont des immersions et 
l'application  continue sous-jacente $|f|:|X|\rightarrow |Y|$ est un 
homéomorphisme sur son image. 
Une \emph{sous-orbifold} de $Y$ est l'image d'un plongement orbifold. 
 
\begin{expl} 
  Soit $\bs{\mu}_{2}$ agissant sur $\CC\times\CC$ de la façon suivante  
  $g(x,y)=(gx,y)$. Le quotient $Y:=\CC\times\CC/\bs{\mu}_{2}$ est une 
  orbifolde. Nous avons $G_{(0,y)}=\bs{\mu}_{2}$ et 
  $G_{(x,y)}=\{\id\}$ si $x\neq 0$. Le lieu singulier du quotient est 
  $(\{0\}\times \CC)/\bs{\mu}_{2}$. 
 
Soit l'orbifold $X_{1}:= \CC/\bs{\mu}_{2}$ où $\bs{\mu}_{2}$ agit 
trivialement sur $\CC$. L'application $\CC \rightarrow \CC\times\CC$ 
qui à $y$ associe $(0,y)$ induit un plongement orbifold 
$X_{1}\rightarrow Y$. Ainsi, $(\{0\}\times \CC) /\bs{\mu}_{2}$ est une 
sous-orbifold de  $Y$. Nous pouvons voir $X_{1}$ comme le lieu singulier de  
$Y$ que nous avons \og sorti \fg de $|Y|$.  
 
Soit l'orbifold $X_{2}:= \CC/\bs{\mu}_{2}$ où $\bs{\mu}_{2}$ agit 
sur $\CC$ par multiplication. L'application $\CC \rightarrow \CC\times\CC$ 
qui à $x$ associe $(x,0)$ induit un plongement orbifold 
$X_{2}\rightarrow Y$. Ainsi, $(\CC\times\{0\}) /\bs{\mu}_{2}$ est une 
sous-orbifold de  $Y$. 
\end{expl} 
   
Nous définissons le faisceau des fonctions holomorphes sur 
 $|X|$.  Notons $\mathcal{O}_{X_{\reg}}$ le faisceau des 
fonctions holomorphes sur la variété complexe $X_{\reg}$. Notons 
$j:X_{\reg}\hookrightarrow X$ l'inclusion naturelle. Nous allons 
définir le faisceau d'anneaux commutatifs et unitaires 
$\mathcal{O}_{|X|}$ comme le sous-faisceau de 
$j_{\ast}\mathcal{O}_{X_{\reg}}$ formé des fonctions localement 
bornées. Ainsi, $(|X|,\mathcal{O}_{|X|})$ est un espace annelé. 
Remarquons que la donnée de l'espace annelé $(|X|,\mathcal{O}_{|X|})$ 
est plus faible que la donnée d'une structure orbifolde sur $|X|$. Sur 
un espace annelé $(|X|,\mathcal{O}_{|X|})$ nous \og oublions \fg les 
actions de groupes.  La proposition suivante explique cette perte.

\begin{expl}\label{expl:espace,annele}  
  Il est clair que l'espace annelé 
  $(|\PP^{1}_{w_{0},w_{1}}|,\mathcal{O}_{|\PP^{1}_{w_{0},w_{1}}|})$ 
  est isomorphe à l'espace annelé $(|\PP^{1}|,\mathcal{O}_{|\PP^{1}|})$. 
\end{expl} 
 
\begin{prop}\label{prop:faisceau,fonction}Soit $X$ une orbifold. 
  Soit $U$ un ouvert de $X$. Soit $f$ dans $\mathcal{O}_{|X|}(U)$. Pour  
  tout $x$ dans $U -U\cap X_{\reg}$, il existe une carte 
  $(\widetilde{U}_{x},G_{x},\pi_{x})$ d'un ouvert $U_{x}$ contenant $x$ et $f_{1}\in 
  \left({\pi_{x}}_{\ast}\mathcal{O}_{\widetilde{U}_{x}}\right)^{G_{x}}(U_{x})$ tel que $f_{1}=f\circ \pi_{x}$. 
\end{prop} 
 
\begin{proof} 
  Posons $f_{1}:=f\circ\pi_{x}\mid_{\pi_{x}^{-1}(U_{x}\cap 
    X_{\reg})}:\pi_{x}^{-1}(U_{x}\cap X_{\reg})\to \CC$. 
La fonction $f_{1}$ est holomorphe sur l'ouvert dense $\pi_{x}^{-1}(U_{x}\cap 
    X_{\reg})\subset \widetilde{U}_{x}$ et elle est 
    $G_{x}$-invariante. 
Comme $f$ est localement bornée, $f_{1}$ est localement bornée  
sur $\widetilde{U}_{x}$. Ainsi, la fonction $f_{1}$ se prolonge en une 
fonction holomorphe  $G_{x}$-invariante sur $\widetilde{U}_{x}$. 
\end{proof} 
 
Pour finir ce paragraphe, nous allons définir le faisceau, noté 
$\mathcal{C}^{\infty}_{|X|}$, des fonctions $C^{\infty}$ sur 
une orbifold complexe $X$.\footnote{Ce faisceau est l'analogue pour les 
orbifolds du faisceau des fonctions $C^{\infty}$ sur une 
variété complexe.} 
Le faisceau $\mathcal{C}^{\infty}_{|X|}$ est le 
sous-faisceau de $j_{\ast}\mathcal{C}^{\infty}_{X_{\reg}}$ défini par 
\begin{align*} 
  \mathcal{C}^{\infty}_{|X|}(U):= 
\left\{\begin{array}{l} 
f\in j_{\ast}\mathcal{C}^{\infty}_{X_{\reg}}(U)\mid \forall x\in 
U_{\sing}, \exists (\widetilde{U}_{x},G_{x},p_{x}) \mbox{ une carte} \\ 
\mbox{de } U_{x}  \mbox{ et } f_{1}\in \left({\pi_{x}}_{\ast} \mathcal{C}^{\infty}_{\widetilde{U}_{x}}\right)^{G_{x}}(U_{x}) 
\mbox{ telles que } f_{1}= f\circ\pi_{x} \mbox{ sur } U_{x,\reg} 
\end{array}\right\} 
\end{align*}

Soit $\mathcal{U}:=(U_{\alpha})_{\alpha\in A}$ un recouvrement ouvert de $|X|$. 
Une \emph{partition de l'unité  sur l'orbifold $X$ subordonnée au recouvrement 
$\mathcal{U}$} est une famille 
$(\rho_{\alpha})_{\alpha\in A}$ de fonctions $C^{\infty}$ telle que 
\begin{enumerate} 
\item pour tout $\alpha\in A$ et pour tout $x\in |X|$, 
  $\rho_{\alpha}(x)\in [0,1]$ ; 
\item pour tout $\alpha\in A$, les supports de $\rho_{\alpha}$ sont 
  inclus dans $U_{\alpha}$ ; 
\item la famille des supports des fonctions $\rho_{\alpha}$ est un recouvrement 
localement fini ; 
\item pour tout $x\in |X|$, nous avons 
$\sum_{\alpha\in A}\rho_{\alpha}(x)=1$. 
\end{enumerate}

\begin{prop}[cf. lemme $4.2.1$ dans \cite{CRogw}]\label{prop:partition,unite} Soit $|X|$ un espace topologique paracompact.  Soit $\mathcal{A}(|X|)$ un 
  atlas orbifold de $|X|$. Soit $(U_{\alpha})_{\alpha\in A}$ le 
  recouvrement de $|X|$ associé à cet atlas.  Il existe une partition 
  de l'unité subordonnée au recouvrement $(U_{\alpha})_{\alpha\in A}$. 
\end{prop} 
 
\begin{proof} 
  Comme $|X|$ est paracompact, il existe un raffinement $(V_{i})_{i\in 
    I}$ localement fini du recouvrement $(U_{\alpha})_{\alpha\in A}$. 
  Chaque $V_{i}$ admet une carte $(\widetilde{V}_{i},G_{i},\pi_{i})$. 
   
Montrons qu'il existe une famille de fonctions $C^{\infty}$ 
  $\widetilde{f}_{i}:\widetilde{V}_{i}\to \CC$ $G_{i}$-invariantes et 
  positives telle que les fonctions induites $f_{i}$ sur $V_{i}$ soient 
  à support dans $V_{i}$ et que la réunion des supports des fonctions 
  $f_{i}$ recouvre $|X|$.  
 Il existe un recouvrement $B_{i}$ de $|X|$ 
  tel que $\overline{B}_{i}\subset V_{i}$. Il existe une famille 
  $(h_{i})_{i\in I}$ de fonctions continues telle que  
  \begin{itemize} 
    \item  $h_{i}\geq 0$; 
\item $h_{i}\equiv 1$ 
  sur $\overline{B}_{i}$; 
\item le support de $h_{i}$ est inclus dans 
  $V_{i}$. 
  \end{itemize} 
  Nous relevons les fonctions $h_{i}$ en des fonctions continues 
  $\widetilde{h}_{i}$ $G_{i}$-invariantes sur $\widetilde{V}_{i}$.
  Remarquons que $\widetilde{h_{i}}\equiv 1 $ sur 
  $\pi_{i}^{-1}(\overline{B}_{i})$.  Nous lissons ses fonctions et 
  nous obtenons des fonctions $C^{\infty}$, notées  $\widetilde{h}_{i}^{\infty}$, telles que 
  $\widetilde{h}_{i}^{\infty}\equiv 1$ sur 
  $\pi_{i}^{-1}(\overline{B}_{i})$.  Puis, nous posons 
  \begin{align*} 
    \widetilde{f}_{i}&:=\frac{1}{\#G_{i}}\sum_{g\in G} 
  \widetilde{h}^{\infty}_{i}\circ\varphi_{g}. 
  \end{align*} 
 Pour tout $i\in I$, la fonction 
  $\widetilde{f}_{i}$ est positive, $G_{i}$-invariante, $C^{\infty}$ et son 
  support contient $\pi_{i}^{-1}(\overline{B}_{i})$.  Notons $f_{i}$ 
  la fonction induite  par $\widetilde{f}_{i}$ sur $V_{i}$. Le support de 
  $f_{i}$ contient $\overline{B}_{i}$ et donc l'ensemble des supports 
  des fonctions $f_{i}$ recouvre $|X|$. 
 
Puis, nous prolongeons les fonctions $f_{i}$ par $0$ en dehors de $U_{i}$. 
Nous en déduisons que $f:=\sum_{i}f_{i}$ est  
  une fonction strictement positive sur $|X|$. 
 La famille des fonctions $\rho_{i}:=f_{i}/f$ est une partition de
 l'unité subordonnée au recouvrement $(V_{i})_{i\in I}$. Comme le
 recouvrement $(V_{i})_{i\in I}$ est un raffinement du recouvrement
 $(U_{\alpha})_{\alpha\in A}$, nous en d\'eduisons une partition de
 l'unit\'e subordonn\'ee au recouvrement $(U_{\alpha})_{\alpha\in A}$.  
\end{proof} 
 
 \section{Les fibrés vectoriels complexes orbifolds}\label{sec:les-fibr-vect} 
  
Dans ce paragraphe, nous allons d'abord définir les fibrés vectoriels 
complexes orbifolds triviaux et définir leurs faisceaux des sections. 
Puis, nous donnerons la définition générale des  fibrés vectoriels 
complexes orbifolds et de leur faisceau des sections.

\subsection{Les fibrés vectoriels complexes orbifolds triviaux} 
\label{subsec:Fibrs-vect-orbif-triviaux} 
 
 Soit $(\widetilde{U},G,\pi)$ une carte de $U$. Soient $E$ une orbifold 
 et $\pr:E\rightarrow U$ une application orbifolde surjective. 
 L'application $\pr:E\rightarrow U$ est un \emph{fibré vectoriel complexe 
   orbifold trivial de rang $r$} sur $U$ si 
\begin{enumerate} \makeatletter 
  \renewcommand\theenumi{\theequation} 
  \makeatother  
\addtocounter{equation}{1} 
\item \label{item:ker} $\ker (E)=\ker (U)$; 
\addtocounter{equation}{1} 
\item\label{item:11} $(\widetilde{U}\times \CC^{r},G,\pi_{E})$ est une carte de $E$; 
\addtocounter{equation}{1} 
\item\label{item:13} l'application $\pr$ se relève en la projection 
  $\widetilde{\pr}:\widetilde{U}\times \CC^{r} \rightarrow 
  \widetilde{U}$ c'est-à-dire $\pi\circ\widetilde{\pr}=\pr\circ\pi_{E}$; 
\addtocounter{equation}{1} 
\item\label{item:14} l'action de $G$ sur $\widetilde{U}\times\CC^{r}$ est donnée par 
  $g(\widetilde{x},w)=(g\widetilde{x},\rho(\widetilde{x},g)w)$ où 
  l'application $\rho:\widetilde{U}\times G\rightarrow GL_{r}(\CC)$ 
  est holomorphe telle que 
  \begin{align*} 
    \rho(\widetilde{x},gh)&=\rho(h\widetilde{x},g)\circ\rho(\widetilde{x},h). 
  \end{align*} 
\end{enumerate}

La condition (\ref{item:ker}) implique que 
\begin{itemize} 
\item $\rho(\widetilde{x},g)=\id$ pour tout $\widetilde{x}\in \widetilde{U}$ 
et pour tout $g\in \ker (U)$; 
\item l'application $\pr{\mid_{E_{\reg}}}:E_{\reg}\rightarrow U_{\reg}$ est un fibré 
vectoriel complexe trivial sur la variété complexe $X_{\reg}$; 
\item l'application $\pr\mid_{E_{\red}}:E_{\red}\rightarrow U_{\red}$ est un 
  fibré vectoriel orbifold trivial. 
\end{itemize} 
 
La fibre $E_{x}:=\pr^{-1}(x)$ est isomorphe à $\CC^{r}/G_{x}$ où
l'action de $G_{x}$ sur $\CC^{r}$ est donnée par l'application
$\rho(x,\cdot)$.  La fibre $E_{x}$ contient l'espace vectoriel
$E_{x}^{G_{x}}:=\{\pi_{E}(\widetilde{x},\widetilde{v})\mid \forall
g\in G_{x},
g(\widetilde{x},\widetilde{v})=(\widetilde{x},\widetilde{v})\}$.

Une section holomorphe d'un fibré orbifold trivial $\pr:E\rightarrow U$ est une 
application orbifolde $s:U\rightarrow E$ qui se relève en une 
application holomorphe $(\id,\widetilde{s}):\widetilde{U}\rightarrow 
\widetilde{U}\times \CC^{r}$ $G$-équivariante.  En d'autres termes, une 
section est la donnée d'une application holomorphe 
$\widetilde{s}:\widetilde{U}\rightarrow \CC^{r}$ telle que 
\begin{align*} 
\widetilde{s}(g\widetilde{x})&=\rho(\widetilde{x},g)\widetilde{s}(\widetilde{x}). 
\end{align*} 
Ainsi, l'ensemble des sections holomorphes est un faisceau de 
$\mathcal{O}_{|U|}$-modules. Nous le notons $\mathcal{E}_{|U|}$. 
Nous définissons l'action de $G$ sur l'ensemble des applications 
holomorphes $\widetilde{s}:\widetilde{U}\rightarrow \CC^{r}$ par la 
formule 
\begin{align} 
  \label{eq:5} 
g\cdot\widetilde{s}(\widetilde{x})&=\rho(g^{-1}\widetilde{x},g)\widetilde{s}(g^{-1}\widetilde{x}). 
\end{align} 
Une application holomorphe $\widetilde{s}:\widetilde{U}\rightarrow 
\CC^{r}$ $G$-invariante induit par passage au quotient une section 
holomorphe $s:U\rightarrow E$.  Le faisceau $\mathcal{E}_{|U|}$ des 
sections holomorphes du fibré $\pr:E\rightarrow U$ est isomorphe au 
faisceau 
$\left({\pi}_{\ast}{\mathcal{O}_{\widetilde{U}}}^{r}\right)^{G}$. 
 
En remplaçant le mot holomorphe par $C^{\infty}$ au 
paragraphe précédent, nous définissons le faisceau des sections 
${C}^{\infty}$ d'un fibré orbifold $\pr:E\rightarrow 
U$.\footnote{Ce faisceau est l'analogue du faisceau des sections 
  $C^{\infty}$ pour un fibré holomorphe sur une variété 
  complexe.} Ce faisceau, noté $\mathcal{E}^{\infty}_{|U|}$, est un 
faisceau de $\mathcal{C}^{\infty}_{|U|}$-module. Il est isomorphe au 
faisceau 
$\left({\pi}_{\ast}{\mathcal{C}^{\infty}_{\widetilde{U}}}^{r}\right)^{G}$.

 Soit $s:U\rightarrow E$ une section holomorphe ou $C^{\infty}$ du fibré 
$\pr:E\rightarrow U$.  Nous en déduisons que 
  \begin{itemize} 
   \item si $\widetilde{s}(\widetilde{x})$ est non nul alors 
$\widetilde{s}(\widetilde{x})$ est un vecteur propre de la matrice 
$\rho(\widetilde{x},g)$ pour la valeur propre $1$ pour tout $g\in G_{x}$ ;  
\item  pour tout $x$ dans $U$,  $s(x)$ appartient à 
l'espace vectoriel $E_{x}^{G_{x}}$.  
\end{itemize}  
Contrairement au cas des variétés, le faisceau des sections d'un fibré 
vectoriel orbifold ne permet de reconstruire le fibré vectoriel 
orbifold car nous perdons l'action du groupe $G$ sur le fibré.

 \begin{expl}\label{expl:fibre} Nous reprenons les notations de l'exemple \ref{expl:orbifold}.(\ref{item:16}). 
   Consi\-dé\-rons le fibré orbifold trivial défini par le diagramme 
   commutatif 
   \begin{align*} 
     \xymatrix{ \widetilde{D}\times \CC \ar[d]_-{\widetilde{\pr}} 
       \ar[r]^-{ \pi_{E}} & (\widetilde{D}\times \CC) /\bs{\mu}_{n} \ar[d]^{\pr}\\ 
\widetilde{D} \ar[r]^-{\pi} & D} 
   \end{align*} 
   où $\zeta\cdot(z,v):=(\zeta z,\zeta v)$. Soit $s$ une section de ce 
   fibré. Ainsi, nous avons 
   \begin{align*} 
\widetilde{s}(\zeta\widetilde{x})&=\zeta\widetilde{s}(\widetilde{x}). 
\end{align*} 
   Ceci implique que $\widetilde{s}(0)=0$. 
 \end{expl} 
 
\begin{rem}\label{rem:fibre,chen,ruan} 
  Au paragraphe $4.1$ de l'article \cite{CRogw}, un fibré 
  vectoriel orbifold trivial vérifie les conditions (\ref{item:11}), 
  (\ref{item:13}) et (\ref{item:14}). Un fibré vectoriel qui vérifie 
  aussi (\ref{item:ker}) est appelé un \emph{bon fibré} (cf. 
  paragraphe $4.3$ de \cite{CRogw}). Nous avons l'équivalence suivante 
  : 
\begin{align*} 
  E\rightarrow U \mbox{ est un bon fibré orbifold trivial} \\ \Leftrightarrow 
  E_{\red}\rightarrow U_{\red} \mbox { est un fibré orbifold trivial}. 
\end{align*} 
Modifions légèrement l'exemple du fibré vectoriel \ref{expl:fibre} de 
la façon suivante : l'action de $\bs{\mu}_{n}$ sur $\widetilde{D}$ est 
trivial et l'action de $\bs{\mu}_{n}$ sur $\widetilde{D}\times \CC$ 
est donnée par $\zeta\cdot(z,v):=(z,\zeta v)$.  Pour ce fibré, nous avons 
$\ker(E)=\{\id\}$ et $\ker(\widetilde{D})=\bs{\mu}_{n}$. Soit $s$ une 
section de ce fibré alors nous avons 
\begin{align*} 
  \widetilde{s}(\widetilde{x})&=\widetilde{s}(\zeta\widetilde{x})=\zeta\widetilde{s}(\widetilde{x}). 
\end{align*} 
Ceci implique que la seule section possible d'un tel fibré est la 
section nulle. C'est pour cette raison que nous ne considérons que les 
bons fibrés orbifolds. 
 \end{rem}

 Soient $(\widetilde{U}_{1},G_{1},\pi_{1})$ et 
 $(\widetilde{U}_{2},G_{2},\pi_{2})$ deux cartes de $U$. Deux fibrés 
 vectoriels orbifolds triviaux $\pr_{1}:E_{1}\rightarrow U$ et 
 $\pr_{2}:E_{2}\rightarrow U$ sont \emph{isomorphes} s'il existe un 
 isomorphisme de cartes 
 $(\varphi,\kappa):(\widetilde{U}_{1},G_{1},\pi_{1})\hookrightarrow 
 (\widetilde{U}_{2},G_{2},\pi_{2})$ et s'il existe une application 
 $\delta:\widetilde{U}_{1}\rightarrow GL_{r}(\CC)$ telle que 
 l'application 
  \begin{align} 
    \widetilde{U}_{1}\times \CC^{r} & \longrightarrow \widetilde{U}_{2}\times \CC^{r} \label{eq:iso,fibre} \\ 
    (\widetilde{x},w) & \longmapsto 
    (\varphi(\widetilde{x}),\delta(\widetilde{x})w) \nonumber 
  \end{align} 
  soit  $\kappa$-équivariante. 
   
  Soient $\pr_{1}:E_{1}\rightarrow U$ et $\pr_{2}:E_{2}\rightarrow U$ 
  deux fibrés orbifolds triviaux. Soient $(\widetilde{U},G,\pi)$ une 
  carte de $U$ et $(\widetilde{U}\times\CC^{r_{1}},G,\pi_{E_{1}})$ 
  (resp. $(\widetilde{U}\times\CC^{r_{2}},G,\pi_{E_{1}})$) une carte de 
  $E_{1}$ (resp. $E_{2}$). 
Un \emph{morphisme entre deux fibrés orbifolds triviaux}  est une 
application orbifolde $\varphi:E_{1} \rightarrow E_{2}$ qui commute 
avec les projections $\pr_{1}$ et $\pr_{2}$ telle qu'il 
existe un relèvement linéaire $(\id,\widetilde{\varphi}): \widetilde{U}\times \CC^{r_{1}}\rightarrow 
\widetilde{U}\times\CC^{r_{2}}$ $G$-équivariant c'est-à-dire que nous 
avons  
\begin{align}\label{eq:morphisme,fibre,trivial} 
  \widetilde{\varphi}(g\widetilde{x})\rho^{E_{1}}(\widetilde{x},g)&=\rho^{E_{2}}(\widetilde{x},g)\widetilde{\varphi}(\widetilde{x})  \hspace{1cm} \forall (\widetilde{x},g)\in\widetilde{U}\times G. 
\end{align}

  Soit $(\widetilde{U},G,\pi)$ une carte de $U$.  Soit 
  $\pr:E\rightarrow U$ un fibré vectoriel orbifold trivial. Soit $V$ 
  un ouvert connexe de $U$.  Soit $(\alpha,\kappa): 
  (\widetilde{V},G_{\widetilde{V}},\pi_{\widetilde{V}}) 
  \hookrightarrow (\widetilde{U},G,\pi)$ une injection.  Soit 
  $\pr_{V}:F\rightarrow V$ un fibré vectoriel orbifold trivial. Le 
  fibré $\pr_{V}:F\rightarrow V$ \emph{s'injecte} dans le fibré 
  $\pr:E\rightarrow U$ s'il existe une application holomorphe 
  $\g_{\alpha}:\widetilde{V}\rightarrow GL_{r}(\CC)$ telle que 
  l'application 
\begin{align*} 
  \widetilde{V} \times \CC^{r} &  \longrightarrow \widetilde{U}\times \CC^{r}  \\ 
  (\widetilde{x},w) & \longmapsto  
  (\alpha(\widetilde{x}),\g_{\alpha}(\widetilde{x})w) 
\end{align*} 
soit une injection de cartes. En particulier, nous avons 
\begin{align*} 
\g_{\alpha}(g\widetilde{x})\circ\rho_{\widetilde{V}}(\widetilde{x},g)&=\rho_{\widetilde{U}}(\alpha(\widetilde{x}),\kappa(g))\circ 
\g_{\alpha}(\widetilde{x}). 
\end{align*}

 \begin{lem} \label{lem:injections,fibre,triviaux} Soit $(\widetilde{U},G,\pi)$ une carte de $U$. 
  Soit $\pr:E\rightarrow U$ un fibré  vectoriel orbifold trivial. Soit $V$ un ouvert 
  connexe de $U$. Il existe un unique, à isomorphisme de fibrés près, 
  fibré  vectoriel orbifold $F\rightarrow V$ qui s'injecte dans $\pr:E\rightarrow U$. 
 \end{lem} 
 
 \begin{proof} 
   \emph{Existence:} Soit $\widetilde{V}$ une composante connexe de 
   $\pi^{-1}(V)$. Notons $G_{\widetilde{V}}$ le sous-groupe de $G$ qui 
   stabilise $\widetilde{V}$.  Ainsi, $\widetilde{V}\times 
   \CC^{r}\subset \widetilde{U}\times \CC^{r}$ est une composante 
   connexe de $\pi_{E}^{-1}\pr^{-1}(V)$. Le groupe $G_{\widetilde{V}}$ 
   agit sur $\widetilde{V}\times \CC^{r}$. Finalement, 
   $(\widetilde{V}\times 
   \CC^{r},G_{\widetilde{V}},{\pi_{E}}{\mid_{\widetilde{V}\times 
       \CC^{r}}})$ est une carte de $\pr^{-1}(V)$ et 
   $\pr{\mid_{\pr^{-1}(V)}}:\pr^{-1}(V)\rightarrow V$ est un fibré 
   trivial qui s'injecte dans $\pr:E\rightarrow U$. 
  
   \emph{Unicité:} Soit $F\rightarrow V$ un fibré trivial qui 
   s'injecte dans $\pr:E\rightarrow U$. Soit 
   $(\widetilde{V}_{F},G_{\widetilde{V}_{F}},\pi_{\widetilde{V}_{F}})$ 
   une carte de $V$. Soit $(\alpha_{F},\kappa): 
   (\widetilde{V}_{F},G_{\widetilde{V}_{F}},\pi_{\widetilde{V}_{F}}) 
   \hookrightarrow (\widetilde{U},G,\pi)$ une injection.  D'après le 
   lemme \ref{lem:ruan}, il existe un isomorphisme de cartes 
   $(\varphi,\kappa'):(\widetilde{V}_{F},G_{\widetilde{V}_{F}},\pi_{\widetilde{V}_{F}})\rightarrow 
   (\widetilde{V},G_{\widetilde{V}},\pi_{\widetilde{V}})\subset 
   (\widetilde{U},G,\pi)$.  Alors, l'application 
   $\widetilde{V}_{F}\times \CC^{r} \rightarrow \widetilde{V} \times 
   \CC^{r}$ qui à $(\widetilde{x},w)$ associe 
   $(\varphi(\widetilde{x}),\g_{\alpha_{F}}(\widetilde{x})w)$ est un 
   isomorphisme de fibrés triviaux. 
  \end{proof}  
 
 \bigskip
 
  \subsection{Le cas général : les fibrés vectoriels complexes orbifolds} 
  \label{subsec:Le-cas-gnral}  
 
\begin{defi}\label{defi:fibre,orbifold}  
  Une application $\pr:E\rightarrow X$ surjective entre deux orbifolds 
  est un \emph{fibré vectoriel complexe orbifold de rang $r$} si pour 
  tout $x$ dans $X$, il existe une carte 
  $(\widetilde{U}_{x},G_{x},\pi_{x})$ d'un ouvert $U_{x}$ contenant 
  $x$ telle que 
\begin{enumerate}  
\item l'application 
  $\pr{\mid_{\pr^{-1}(U_{x})}}:\pr^{-1}(U_{x})\rightarrow U_{x}$ est un 
  fibré vectoriel complexe orbifold trivial de rang $r$ ; 
\item pour toute injection 
  $(\alpha,\kappa):(\widetilde{U}_{y},G_{y},\pi_{y})\hookrightarrow(\widetilde{U}_{x},G_{x},\pi_{x})$, 
  le fibré orbifold trivial $\pr^{-1}(U_{y})\rightarrow U_{y}$ s'injecte dans 
  $\pr^{-1}(U_{x})\rightarrow U_{x}$. 
    \end{enumerate}  
\end{defi}  
 
Soit $\pr:E\rightarrow X$ un fibré vectoriel orbifold de rang $r$. Une 
carte $(\widetilde{U},G,\pi)$ de $U\subset X$ est dite 
\emph{trivialisante} si $(\widetilde{U}\times \CC^{r},G,\pi_{E})$ est 
une carte de $\pr^{-1}(U)$.  Remarquons que 
\begin{itemize} 
\item $\pr{\mid_{E_{\reg}}}:E_{\reg}\rightarrow X_{\reg}$ est un fibré 
  vectoriel complexe de rang $r$ sur la variété complexe $X_{\reg}$ ; 
\item $E_{\red}\rightarrow X_{\red}$ est un fibré vectoriel orbifold. 
\end{itemize}

Une section holomorphe (resp. $C^{\infty}$) d'un fibré 
orbifold $\pr:E\rightarrow X$ est une application orbifolde 
$s:X\rightarrow E$ qui localement se relève en une application 
holomorphe (resp. $C^{\infty}$) de 
$(\id,\widetilde{s}_{x}):\widetilde{U}_{x}\rightarrow 
\widetilde{U}_{x}\times \CC^{r}$ $G_{x}$-invariante c'est-à-dire que 
nous avons pour tout $g\in G_{x}$ 
\begin{align*} 
  g\cdot\widetilde{s}_{x}&=\widetilde{s}_{x}. 
\end{align*} 
Notons $\mathcal{E}_{X_{\reg}}$ (resp. 
$\mathcal{E}_{X_{\reg}}^{\infty}$) le faisceau des sections 
holomorphes (resp. $C^{\infty}$) du fibré vectoriel 
$\pr{\mid_{E_{\reg}}}:E_{\reg}\rightarrow X_{\reg}$.  Nous allons 
définir le faisceau, noté $\mathcal{E}_{|X|}$ (resp. 
$\mathcal{E}_{|X|}^{\infty}$), des sections holomorphes (resp. 
$C^{\infty}$) d'un fibré orbifold $\pr:E \rightarrow X$ comme 
un sous-faisceau de $j_{\ast}\mathcal{E}_{X_{\reg}}$ (resp. 
$j_{\ast}\mathcal{E}^{\infty}_{X_{\reg}}$), où $j$ est l'inclusion 
naturelle de $X_{\reg}$ dans $X$.

Pour tout ouvert $U$ dans $|X|$, nous posons 
\begin{align*} 
\mathcal{E}_{|X|}(U):= 
\left\{\begin{array}{l} 
s\in j_{\ast}\mathcal{E}_{X_{\reg}}(U)\mid \forall x\in U-U_{\reg}, 
\exists (\widetilde{U}_{x},G_{x},\pi_{x}) \mbox{ une carte 
} \\  \mbox{trivialisante d'un voisinage de } U_{x}  
\mbox{ et } \exists 
(\id,\widetilde{s}):\widetilde{U}_{x}\rightarrow \widetilde{U}_{x}\times 
\CC^{r}   \\ G_{x}\mbox{-équivariante tels que } 
\pi_{x}^{E}\circ(\id,\widetilde{s})=s\circ{\pi}_{x}    
\end{array}\right\} 
\end{align*} 
 
\begin{align*} 
\mathcal{E}^{\infty}_{|X|}(U):= 
\left\{\begin{array}{l} 
s\in j_{\ast}\mathcal{E}^{\infty}_{X_{\reg}}(U)\mid \forall x\in U-U_{\reg}, 
\exists (\widetilde{U}_{x},G_{x},\pi_{x}) \mbox{ une carte 
} \\  \mbox{trivialisante d'un voisinage de } U_{x}  
\mbox{ et } \exists 
(\id,\widetilde{s}):\widetilde{U}_{x}\rightarrow \widetilde{U}_{x}\times 
\CC^{r}  \\  G_{x}\mbox{-équivariante tels que } 
\pi_{x}^{E}\circ(\id,\widetilde{s})=s\circ{\pi}_{x}    
\end{array}\right\} 
\end{align*}

Nous pouvons aussi définir le faisceau $\mathcal{E}_{|X|}$ (resp. 
$\mathcal{E}^{\infty}_{|X|}$) en recollant les faisceaux 
$({\pi_{x}}_{\ast}{\mathcal{O}_{\widetilde{U}_{x}}}^{r})^{G_{x}}$ 
(resp. 
$({\pi_{x}}_{\ast}{\mathcal{C}^{\infty}_{\widetilde{U}_{x}}}^{r})^{G_{x}}$). 
 
En général, on ne peut pas reconstruire le fibré à partir du faisceau 
de ses sections c'est-à-dire que la donnée d'un fibré orbifold est plus 
riche que la donnée de son faisceau des sections.

\begin{lem} 
  Soit $E\rightarrow X$ un fibré orbifold. Soient 
  $(\widetilde{U},G_{U},\pi_{U}),(\widetilde{V},G_{V},\pi_{V})$ et 
  $(\widetilde{W},G_{W},\pi_{W})$ trois cartes trivialisantes telles 
  qu'il existe des injections 
  $(\alpha,\kappa_{\alpha}):(\widetilde{U},G_{U},\pi_{U})\hookrightarrow(\widetilde{V},G_{V},\pi_{V})$ 
  et 
  $(\beta,\kappa_{\beta}):(\widetilde{V},G_{V},\pi_{V})\hookrightarrow(\widetilde{W},G_{W},\pi_{W})$. 
  Pour tout $\widetilde{x}$ dans $\widetilde{U}$, nous avons 
  \begin{align*} 
  \g_{\beta\circ\alpha}(\widetilde{x})&=\g_{\beta}(\alpha(\widetilde{x}))\circ 
  \g_{\alpha}(\widetilde{x}). 
  \end{align*} 
\end{lem}

\begin{proof} 
  Par hypothèse, nous en déduisons deux injections de 
  $(\widetilde{U}\times \CC^{r},G_{U},\pi_{E})\hookrightarrow 
  (\widetilde{W}\times \CC^{r},G_{W},\pi_{E})$. La première injection 
  est donnée par la formule suivante $(\widetilde{x},w)\mapsto 
  (\beta\circ\alpha(\widetilde{x}),\g_{\beta\circ\alpha}(x)w)$ et la 
  deuxième injection par $(\widetilde{x},w)\mapsto 
  (\beta\circ\alpha(\widetilde{x}),\g_{\beta}(\alpha(\widetilde{x}))\circ 
  \g_{\alpha}(x)w)$.  D'après le lemme \ref{lem:ruan}, il existe 
  $g\in G_{W}$ tel que pour tout $\widetilde{x}\in \widetilde{U}$ et 
  pour tout $w\in \CC^{r}$ nous ayons 
  \begin{align*} 
  g\cdot(\beta\circ\alpha(\widetilde{x}),\g_{\beta\circ\alpha}(x)w)&=(\beta\circ\alpha(\widetilde{x}),\g_{\beta}(\alpha(\widetilde{x}))\circ 
  \g_{\alpha}(x)w). 
  \end{align*} 
  Nous en déduisons que pour tout $\widetilde{x}$ dans 
  $\widetilde{U}$ 
  \begin{align*} 
  \begin{cases} 
      g
      \cdot\beta\circ\alpha(\widetilde{x})=\beta\circ\alpha(\widetilde{x}) ; \\ 
      \rho_{\widetilde{U}}(\beta\circ\alpha(\widetilde{x}))(g)\g_{\beta\circ\alpha}(x)= 
      \g_{\beta}(\alpha(\widetilde{x}))\circ \g_{\alpha}(x)w. 
\end{cases} 
\end{align*} 
La première condition implique que $g\in \ker (U)$ et donc que 
$\rho_{\widetilde{U}}(\beta\circ \alpha(\widetilde{x}))(g)=\id$ car $\ker(E)=\ker(X)$. 
 \end{proof} 
  
 Deux fibrés orbifolds $\pr_{1}:E_{1}\rightarrow X$ et 
 $\pr_{2}:E_{2}\rightarrow X$ sont \emph{isomorphes} s'il existe une 
 application $\psi:E_{1}\rightarrow E_{2}$ telle que pour tout $x$ dans 
 $X$, il existe un isomorphisme $(\widetilde{U}_{x}^{1}\times 
 \CC^{r},G_{x}^{1},\pi_{x}^{E_{1}})\rightarrow 
 (\widetilde{U}_{x}^{2}\times \CC^{r},G_{x}^{2},\pi_{x}^{E_{2}})$ qui 
 soit linéaire entre les fibres de $\widetilde{\pr}_{1}$ et celles de 
 $\widetilde{\pr}_{2}$ et qui induise un isomorphisme entre les cartes 
 $(\widetilde{U}_{x}^{1},G^{1}_{x},\pi_{x}^{1})$ et 
 $(\widetilde{U}^{2}_{x},G_{x}^{2},\pi_{x}^{2})$.  En d'autres termes, 
 pour chaque $x\in X$, il existe une carte 
 $(\widetilde{U}_{x},G_{x},\pi_{x})$ d'un ouvert contenant $x$ telle 
 que nous ayons un isomorphisme entre les fibrés orbifolds triviaux 
 ${E_{1}}{\mid_{U_{x}}}\rightarrow U_{x}$ et 
 ${E_{2}}{\mid_{U_{x}}}\rightarrow U_{x}$. 
  
 Soient fibrés orbifolds $\pr_{1}:E_{1}\rightarrow X$ et 
 $\pr_{2}:E_{2}\rightarrow X$. Un \emph{morphisme} entre les fibrés 
 $\pr_{1}:E_{1}\rightarrow X$ et $\pr_{2}:E_{2}\rightarrow X$ est une 
 application orbifolde $\varphi:E_{1}\rightarrow E_{2}$ qui commute 
 avec les projections $\pr_{1}$ et $\pr_{2}$ telle que pour tout $x 
 \in X$, il existe une carte $(\widetilde{U}_{x},G_{x},\pi_{x})$ d'un 
 ouvert $U_{x}$ contenant $x$ qui vérifie 
 \begin{enumerate} 
 \item $\varphi\mid_{U_{x}}:\pr_{1}^{-1}(U_{x}) \rightarrow 
   \pr_{2}^{-1}(U_{x})$ est un morphisme de fibrés orbifold trivial; 
 \item pour toute injection $\alpha:\widetilde{U}_{y}\hookrightarrow 
   \widetilde{U}_{x}$, nous avons  
   \begin{align}\label{eq:morphisme,fibre} 
     \psi^{E_{2}}_{\alpha}\circ\widetilde{\varphi}\mid_{U_{y}}&=\widetilde{\varphi}\mid_{U_{x}}\circ\psi^{E_{1}}_{\alpha} 
   \end{align} où 
   $\widetilde{\varphi}\mid_{U_{y}}:\widetilde{E}_{1}\mid_{U_{y}}:=\widetilde{U}_{y}\times
   \CC^{r_{1}}\rightarrow
   \widetilde{E}_{2}\mid_{U_{y}}:=\widetilde{U}_{y}\times \CC^{r_{2}}$
   (resp.  $\widetilde{\varphi}\mid_{U_{x}}$) est un relevé linéaire
   $G_{y}$-équivariant (resp. $G_{x}$-équivariant) de
   ${\varphi}\mid_{U_{x}}$ (resp. $\varphi\mid_{U_{y}}$).
 \end{enumerate} 
  
 Soit $\varphi:E_{1}\rightarrow E_{2}$ un morphisme de fibrés 
 orbifolds. Pour tout $x\in |X|$, nous avons une application 
 $\varphi_{x}:E_{1,x}:=\pr_{1}^{-1}(x)\rightarrow 
 E_{2,x}:=\pr_{2}^{-1}(x)$ qui se relève en une application 
 $\widetilde{\varphi}\mid_{U_{x}}(\widetilde{x}):\{\widetilde{x}\}\times 
 \CC^{r_{1}}\rightarrow\{\widetilde{x}\}\times \CC^{r_{2}}$ 
 $G_{x}$-équivariante où $\widetilde{x}$ est un relevé de $x$. 
 Nous posons 
\begin{align*} 
  \ker \varphi &:=\bigcup_{x\in 
    |X|}\pi_{x}^{E_{1}}(\ker(\widetilde{\varphi}\mid_{U_{x}}(\widetilde{x}))) \,; \\ 
 \im \varphi &:=\bigcup_{x\in 
    |X|}\pi_{x}^{E_{2}}(\im(\widetilde{\varphi}\mid_{U_{x}}(\widetilde{x}))). 
\end{align*} 
Remarquons que $\ker\varphi\subset E_{1}$ et $\im\varphi\subset E_{1}$ ne dépendent pas du choix  
des cartes et des relevés. 
 
 \begin{prop}\label{prop:noyau,image,morphisme,fibre} 
   Soient  $\pr_{1}:E_{1}\rightarrow X$ et 
 $\pr_{2}:E_{2}\rightarrow X$ deux fibrés orbifolds. Soit $\varphi:E_{1}\rightarrow E_{2}$ 
 un morphisme de fibrés orbifolds tel que pour tout $x\in U_{x}$ le 
 rang de $\widetilde{\varphi}\mid_{U_{x}}$ est constant. 
 Les espaces topologiques $\ker\varphi$ et $\im\varphi$ sont des fibrés 
 vectoriels orbifolds. 
 \end{prop}

\begin{proof}
Nous allons décrire les fonctions de transition de ses  
   fibrés orbifolds puis nous appliquerons le théorème \ref{thm:Les-fibrs-vectoriels}.  
    Pour toute injection $\alpha:\widetilde{U}_{y}\hookrightarrow 
   \widetilde{U}_{x}$, nous avons le diagramme commutatif suivant 
   (cf. l'égalité (\ref{eq:morphisme,fibre})) 
   \begin{align*} 
     \xymatrix{\widetilde{E}_{1}\mid_{U_{x}}\ar[rr]^{(\id,\widetilde{\varphi}\mid_{U_{x}})} 
       \ar@{^{(}->}[d]_{(\alpha,\psi_{\alpha}^{E_{1}})} 
       &&\widetilde{E}_{2}\mid_{U_{x}} \ar@{^{(}->}[d]^{(\alpha,\psi_{\alpha}^{E_{2}})}\\ 
       \widetilde{E}_{1}\mid_{U_{y}} \ar[rr]^{(\id,\widetilde{\varphi}\mid_{U_{y}})} 
       &&\widetilde{E}_{2}\mid_{U_{y}}}  
   \end{align*} 
Nous en déduisons que   
\begin{align*} 
  \psi_{\alpha}^{E_{1}}\mid_{\ker\widetilde{\varphi}\mid_{U_{x}}}: \ker\widetilde{\varphi}\mid_{U_{x}}&\longrightarrow \ker\widetilde{\varphi}\mid_{U_{y}}\,;\\ 
  \psi_{\alpha}^{E_{2}}\mid_{\im\widetilde{\varphi}\mid_{U_{x}}}:\im\widetilde{\varphi}\mid_{U_{x}}&\longrightarrow 
  \im\widetilde{\varphi}\mid_{U_{y}}. 
\end{align*} 
D'après les résultats sur les fibrés vectoriels sur une variété, 
$\ker\widetilde{\varphi}\mid_{U_{x}}$ et 
$\im\widetilde{\varphi}\mid_{U_{x}}$ sont des fibrés sur 
$\widetilde{U}_{x}$ car le rang de $\widetilde{\varphi}\mid_{U_{x}}$ 
est constant. 
Nous posons  
\begin{align*} 
   \psi_{\alpha}^{\ker\varphi}&:= 
   \psi_{\alpha}^{E_{1}}\mid_{\ker\widetilde{\varphi}\mid_{U_{x}}}\,;\\ 
 \psi_{\alpha}^{\im\varphi}&:= 
   \psi_{\alpha}^{E_{2}}\mid_{\ker\widetilde{\varphi}\mid_{U_{y}}}. 
\end{align*} 
Ces fonctions de transition satisfont bien les hypothèses du théorème  
\ref{thm:Les-fibrs-vectoriels}. Nous en déduisons que $\ker\varphi$ et
$\im\varphi$ sont des  fibrés vectoriels orbifolds sur $X$. 
 \end{proof}

\begin{thm}\label{thm:Les-fibrs-vectoriels} Soit $X$ une orbifold. Les  
  données suivantes définissent un unique, à isomorphisme près, fibré 
  orbifold  de rang $r$ : 
  \begin{enumerate} 
  \item un atlas orbifold $\mathcal{A}(|X|)$ ; 
  \item pour toute 
  injection $\alpha:\widetilde{U_{i}}\hookrightarrow 
  \widetilde{U}_{j}$ entre deux cartes de l'atlas, on se donne  une 
  application holomorphe $\g_{\alpha}:\widetilde{U}_{i}\rightarrow 
  GL_{r}(\CC)$ telle que 
   \begin{enumerate} \item l'application  
    \begin{align*} 
      \widetilde{U}_{i} \times \CC^{r} & \longrightarrow 
      \widetilde{U}_{j} \times \CC^{r}\\  
     (\widetilde{x},w) & \longmapsto (\alpha(\widetilde{x}),\psi_{\alpha}(\widetilde{x})w) 
     \end{align*} soit un plongement ouvert ; 
   \item pour deux injections successives 
  $\alpha,\beta$ nous ayons 
  \begin{align*} 
  \g_{\beta\circ\alpha}(\widetilde{x})&=\g_{\beta}(\alpha(\widetilde{x}))\circ 
  \g_{\alpha}(\widetilde{x}). 
  \end{align*} 
  \end{enumerate}\end{enumerate} 
 \end{thm}

\begin{rem}\label{rem:section,global} 
  Soit $E \rightarrow X$ un fibré vectoriel orbifold donnée par ses 
  fonctions de transition c'est-\`a-dire donn\'ee par les $\psi_{\alpha}$. Soit $(U_{i})_{i\in I}$ le recouvrement de 
  $|X|$ induit par $\mathcal{A}(|X|)$. Supposons que pour tout $i\in 
  I$, il existe $\widetilde{s}_{i}:\widetilde{U}_{i}\rightarrow 
  \CC^{r}$ $G_{i}$-équivariante tel que pour toute injection 
  $\alpha:\widetilde{U}_{i}\hookrightarrow \widetilde{U}_{j}$ nous ayons 
  \begin{align*} 
    \widetilde{s}_{j}(\alpha(\widetilde{x}))=\g_{\alpha}(\widetilde{x})\widetilde{s}_{i}(\widetilde{x}). 
  \end{align*} 
Ces données se recollent en une section globale du fibré $E \rightarrow X$. 
\end{rem} 
 
\begin{proof} 
  Pour tout $i\in I$, notons $(\widetilde{U}_{i},G_{i},\pi_{i})$ la 
  carte de $U_{i}$ dans l'atlas orbifold $\mathcal{A}(|X|)$.  Nous définissons 
  l'action de $G_{i}$ sur $\widetilde{U}_{i}\times\CC^{r}$ de la façon 
  suivante. Pour tout $g\in G_{i}$, l'application 
  $\varphi_{g}:\widetilde{U}_{i}\rightarrow \widetilde{U}_{i}$ qui à 
  $\widetilde{x}$ associe $g\widetilde{x}$ est une injection. Nous 
  posons 
  $g\cdot(\widetilde{x},w):=(\varphi_{g}(\widetilde{x}),\g_{\varphi_{g}}(\widetilde{x})w)$ 
  pour tout $(\widetilde{x},w)\in \widetilde{U}\times
  \CC^{r}$. Remarquons que nous avons $\rho(\widetilde{x},g)=\psi_{\varphi_{g}}(\widetilde{x})$. 
   
  Considérons l'espace topologique $\bigsqcup_{i\in 
    I}\widetilde{U}_{i}\times \CC^{r}/G_{i}$. Notons 
  $[\widetilde{x}_{i},w_{i}]$ la classe de 
  $(\widetilde{x}_{i},w_{i})$. 
   
  Nous dirons que $[\widetilde{x}_{i},w_{i}]\sim 
  [\widetilde{x}_{j},w_{j}]$ s'il existe deux injections 
  $\alpha_{i}:\widetilde{U}_{ij}\hookrightarrow \widetilde{U}_{i}$ et 
  $\alpha_{j}:\widetilde{U}_{ij}\hookrightarrow \widetilde{U}_{j}$ et 
  $(\widetilde{x}_{ij},w_{ij})\in \widetilde{U}_{ij}\times \CC^{r}$ 
  tels que 
  \begin{align*} 
 [\alpha_{i}(\widetilde{x}_{ij}),\g_{\alpha_{i}}(\widetilde{x}_{ij})w_{ij}]& =[\widetilde{x}_{i},w_{i}]\,; \\ 
[\alpha_{j}(\widetilde{x}_{ij}),\g_{\alpha_{j}}(\widetilde{x}_{ij})w_{ij}]&=[\widetilde{x}_{j},w_{j}]. 
 \end{align*} 
Quitte à modifier les injections, nous pouvons supposer que les 
égalités ci-dessus sont sur les couples et non sur les classes. 
 
Montrons que cette relation est une relation d'équivalence.  Le seul 
point délicat est la transitivité.  Soient $[\widetilde{x}_{i},w_{i}], 
[\widetilde{x}_{j},w_{j}],[\widetilde{x}_{k},w_{k}]$ tels que 
$[\widetilde{x}_{i},w_{i}]\sim [\widetilde{x}_{j},w_{j}]$ et 
$[\widetilde{x}_{j},w_{j}]\sim[\widetilde{x}_{k},w_{k}]$.  Ainsi, il 
existe $\beta_{j}:\widetilde{U}_{jk}\hookrightarrow \widetilde{U}_{j}$ 
et $\beta_{k}:\widetilde{U}_{jk}\hookrightarrow \widetilde{U}_{k}$ et 
$(\widetilde{x}_{jk},w_{jk})\in \widetilde{U}_{jk}\times \CC^{r}$ tels 
que 
\begin{align*} 
(\alpha_{i}(\widetilde{x}_{ij}),\g_{\alpha_{i}}(\widetilde{x}_{ij})w_{ij})&=(\widetilde{x}_{i},w_{i})\,;\\ 
(\alpha_{j}(\widetilde{x}_{ij}),\g_{\alpha_{j}}(\widetilde{x}_{ij})w_{ij})&=(\widetilde{x}_{j},w_{j}). 
\end{align*} 
Comme 
$x:=\varphi_{ij}(\widetilde{x}_{ij})=\varphi_{jk}(\widetilde{x}_{jk})$, 
il existe deux injections 
$\gamma_{ij}:\widetilde{U}_{ijk}\hookrightarrow \widetilde{U}_{ij}$ et 
$\gamma_{jk}:\widetilde{U}_{ijk}\hookrightarrow \widetilde{U}_{jk}$ et 
$\widetilde{x}_{ijk}\in \widetilde{U}_{ijk},w_{ijk},w'_{ijk}\in 
\CC^{r}$ tels que 
$\gamma_{ij}(\widetilde{x}_{ijk})=\widetilde{x}_{ij},\gamma_{jk}(\widetilde{x}_{ijk})=\widetilde{x}_{jk}$ 
et 
$\g_{\gamma_{ij}}(\widetilde{x}_{ijk})w_{ijk}=w_{ij},\g_{\gamma_{jk}}(\widetilde{x}_{ijk})w'_{ijk}=w_{jk}$. 
D'après le lemme \ref{lem:injections}, il existe $g_{j}\in G_{j}$ tel 
que $\alpha_{j}\circ \gamma_{ij}=\varphi_{g_{j}}\circ \beta_{j}\circ 
\gamma_{jk}$.  Pour résumer, nous avons le diagramme commutatif 
suivant :
\begin{align*} 
\xymatrix{ & & \widetilde{U}_{ijk} \ar@{_{(}->}[dl]_-{\gamma_{ij}} \ar@{^{(}->}[dr]^-{\gamma_{jk}} & & \\ 
&  \widetilde{U}_{ij} \ar@{_{(}->}[d]_-{\alpha_{j}} \ar@{_{(}->}[dl]_-{\alpha_{i}}& & \widetilde{U}_{jk} \ar@{^{(}->}[d]^-{\beta_{j}} \ar@{^{(}->}[rd]^-{\beta_{k}}&   \\ 
\widetilde{U}_{i} & \widetilde{U}_{j} & & \widetilde{U}_{j} \ar[ll]_{\simeq}^-{\varphi_{g_{j}}}& \widetilde{U}_{k}\\} 
\end{align*} 
Nous en déduisons que $g_{j}\in G_{x}\subset G_{j}$. Ainsi, il 
existe $g_{ijk}\in G_{x} \subset G_{ijk}$ tel que $\varphi_{g_{j}}\circ 
\beta_{j}\circ \gamma_{jk}= \beta_{j}\circ \gamma_{jk} 
\circ\varphi_{g_{ijk}}$.  Comme nous avons 
\begin{align*} 
w_{j}&=\g_{\alpha_{j}\circ\gamma_{ij}}(\widetilde{x}_{ijk})w_{ijk}=\g_{\varphi_{g_{j}}\circ 
  \beta_{j}\circ \gamma_{jk}}(\widetilde{x}_{ijk})w_{ijk} 
\end{align*} 
et 
\begin{align*}w_{j}=\g_{\beta_{j\circ 
    \gamma_{jk}}\circ\varphi_{g_{ijk}}}(\widetilde{x}_{ijk})w'_{ijk}. 
\end{align*} 
nous obtenons que 
$w'_{ijk}=\g_{\varphi_{g_{ijk}}}(\widetilde{x}_{ijk})w_{ijk}$. 
Finalement, les injections $\beta_{j}\circ 
\gamma_{jk}\circ\varphi_{g_{ijk}}:\widetilde{U}_{ijk}\hookrightarrow 
\widetilde{U}_{k}$ et $\alpha_{j}\circ 
\gamma_{ij}:\widetilde{U}_{ijk}\hookrightarrow \widetilde{U}_{j}$ et 
$(\widetilde{x}_{ijk},w_{ijk})\in \widetilde{U}_{ijk}\times \CC^{r}$ 
montre que $[\widetilde{x}_{i},w_{i}]=[\widetilde{x}_{k},w_{k}]$.

Posons $|E|:=\left( \bigsqcup_{i\in I} \widetilde{U}_{i}\times 
  \CC^{r}/G_{i}\right) /\sim$.  Notons 
$[\![\widetilde{x}_{i},w_{i}]\!]$ la classe de 
$[\widetilde{x}_{i},w_{i}]$ pour $(\widetilde{x}_{i},w_{i})\in 
\widetilde{U}_{i}\times \CC^{r}$.  L'application $|\pr|:|E|\rightarrow 
|X|$ qui à $[\![\widetilde{x}_{i},w_{i}]\!]$ associe 
$\pi_{i}(\widetilde{x}_{i})$ est bien définie et elle est continue. 
De plus, $(\widetilde{U}_{i}\times \CC^{r},G_{i},\pi_{E,i})$, où 
$\pi_{E,i}$ est la projection de $\widetilde{U}_{i}\times \CC^{r}$ 
dans $\widetilde{U}_{i}\times \CC^{r}/G_{i}$ , est une carte de 
$\pr^{-1}(U_{i})$.  Nous en déduisons un atlas orbifold sur $|E|$. 
Finalement, $\pr:E\rightarrow X$ est un fibré orbifold de rang $r$ sur 
$X$. 
\end{proof} 
 
\begin{rem} 
  Soit $\pr:E\rightarrow X$ un fibré orbifold sur $X$. Soit $Y$ une 
  sous-orbifold de $X$.  L'application $\pr_{Y}: \pr^{-1}(Y)\rightarrow 
  Y$ est naturellement munie d'une structure de fibré orbifold. De 
  plus, la restriction $\mathcal{E}_{|X|}\mid_{|Y|}$,  comme faisceau de  
  $\mathcal{O}_{|X|}$-modules, du faisceau des sections 
  $\mathcal{E}_{|X|}$s est le faisceau des sections de  $\pr_{Y}: \pr^{-1}(Y)\rightarrow 
  Y$. 
\end{rem} 
 
 \begin{expl}[Le fibré tangent complexe orbifold] 
  Soit $X$ une orbifold complexe de dimension $n$. Nous allons définir le fibré  
  tangent  complexe orbifold.  Soit $\alpha:(\widetilde{U},G,\pi)\hookrightarrow 
  (\widetilde{U}',G',\pi')$ une injection. Soient $(u_{i})$ des 
  coordonnées complexes  sur $\widetilde{U}$ et $(u'_{i})$ des 
  coordonnées complexes sur 
  $\widetilde{U}'$. Nous notons $\g_{\alpha}(\widetilde{x})$ la 
  matrice suivante  
  \begin{align} 
    \label{eq:tangent} 
  \left(\frac{\partial 
        u'_{i}\circ\alpha}{\partial u_{j}}\right)_{i,j\in\{1, \ldots ,n\}}. 
  \end{align} 
  Les conditions du théorème \ref{thm:Les-fibrs-vectoriels} sont 
  vérifiés et nous obtenons un fibré orbifold qu'on appelle \emph{fibré 
  tangent complexe orbifold} de $X$. Le faisceau des sections du fibré  
tangent est noté $\Theta_{|X|}$.

  Une section du fibré tangent complexe orbifold est appelée un \emph{champ de vecteurs 
    complexe}.  Soit l'injection 
  $\varphi_{g}:\widetilde{U}\rightarrow\widetilde{U}$ qui à 
  $\widetilde{x}$ associe $g\widetilde{x}$. Nous avons les égalités suivantes 
  \begin{align*}
\rho(\widetilde{x},g)=\g_{\varphi_{g}}(\widetilde{x})=d\varphi_{g}(\widetilde{x}).
\end{align*} 
  Ainsi, localement un champ de vecteurs est une application 
  $\mathcal{X}:U\rightarrow TU$ qui se relève en un champ de vecteurs 
  $\widetilde{\mathcal{X}}:\widetilde{U}\rightarrow T\widetilde{U}$ où 
  $T\widetilde{U}$ est le fibré vectoriel tangent complexe de $\widetilde{U}$ tel que 
  $\widetilde{\mathcal{X}}(g\widetilde{x})=d\varphi_{g}(\widetilde{x})(\widetilde{\mathcal{X}}(\widetilde{x}))$. 
   
  Remarquons, que si l'on se restreint à $|X_{\reg}|$, alors le fibré 
  tangent complexe orbifold n'est rien d'autre que le fibré tangent 
  complexe de la variété $|X_{\reg}|$. 
\end{expl} 
 
\begin{expl}[Le fibré cotangent complexe orbifold]\label{expl:Les-fibrs-vectoriels} 
  Nous allons définir le fibré cotangent complexe. 
  Soit $\alpha:(\widetilde{U},G,\pi)\hookrightarrow 
  (\widetilde{U}',G',\pi')$ une injection. Nous posons 
  \begin{align} 
    \label{eq:cotangent} 
    \g_{\alpha}^{T^{\ast}X}(\widetilde{x})&:= \ ^{t} \g_{\alpha}^{TX}(\widetilde{x})^{-1}.
  \end{align} 
Les conditions du théorème \ref{thm:Les-fibrs-vectoriels} sont 
  vérifiés et nous obtenons un fibré orbifold qu'on appelle \emph{fibré 
  cotangent complexe orbifold} de $X$. Nous notons ce fibré   
$T^{\ast}X$. 
Ainsi, nous avons $\rho_{T^{\ast}X}(\widetilde{x},g)=\ 
^{t}\rho_{TX}(\widetilde{x},g)^{-1}= \ ^{t}d\varphi_{g}(\widetilde{x})^{-1}$. 
Le faisceau des sections du fibré orbifold $T^{\ast}X$ est noté 
$\Omega^{1}_{|X|}$. 
Une section du fibré  $T^{\ast}X$ est appelée une $1$-forme holomorphe.  
Localement une $1$-forme holomorphe est une application 
$\omega:U\rightarrow T^{\ast}U$ qui se relève en une application 
$\widetilde{\omega}:\widetilde{U}\rightarrow T^{\ast}\widetilde{U}$ 
telle que 
\begin{align}\label{eq:G,equivariant} 
\widetilde{\omega}(g\widetilde{x})(X)&=\widetilde{\omega}(\widetilde{x})(d 
\varphi_{g}(\widetilde{x})^{-1}(X)) 
\end{align}  
où $X\in T_{g\widetilde{x}}\widetilde{U}$. 
La condition (\ref{eq:G,equivariant}) est équivalente à 
$g\cdot\widetilde{\omega}=\widetilde{\omega}$ où l'action de $G$ sur 
les $1$-formes holomorphes de $\widetilde{U}$ est donnée par la formule 
\begin{align}\label{eq:6} 
  {g}\cdot\widetilde{\omega}&= 
  \left(\varphi_{g}^{-1}\right)^{\ast}\widetilde{\omega}=\left(\varphi_{g^{-1}}\right)^{\ast}\widetilde{\omega} . 
\end{align}

 Soit $\alpha:(\widetilde{U},G,\pi)\hookrightarrow 
  (\widetilde{U}',G',\pi')$ une injection. Nous posons 
  \begin{align} 
    \label{eq:cotangent,k} 
    \g_{\alpha}^{\wedge^{k}T^{\ast}X}&:= \wedge^{k} 
    \g_{\alpha}^{T^{\ast}X}. 
  \end{align} 
Les conditions du théorème \ref{thm:Les-fibrs-vectoriels} sont 
  vérifiées et nous obtenons un fibré orbifold qu'on appelle \emph{fibré 
  des $k$-formes holomorphes} de $X$. Nous notons ce fibré   
$\wedge^{k}T^{\ast}X$. 
\end{expl} 
 
\section[Formes différentielles et intégrale orbifolde]{Faisceaux des formes différentielles sur une 
  orbifold et intégrale orbifolde}\label{appendice:A} 
Dans ce paragraphe, nous allons définir le faisceau des $k$-formes 
différentielles $C^{\infty}$ sur une orbifold. Puis, nous 
montrerons que ces faisceaux permettent d'interpréter une classe de cohomologie de 
l'espace topologique sous-jacent à une orbifold comme la classe d'une 
forme différentielle fermée. Ceci nous permettra de définir 
l'intégrale orbifolde.  
  
Dans ce paragraphe $X$ est une orbifold complexe connexe de dimension $n$. 
 
Notons $\mathcal{E}^{k}_{|X_{\reg}|}$ le faisceau des $k$-formes 
différentielles $C^{\infty}$ \'a valeur complexe sur la variété complexe $X_{\reg}$. Soit 
$j$ l'inclusion de $|X_{\reg}|$ dans $|X|$. 
Soit $(\widetilde{U},G,\pi)$ une carte de $U$. Nous notons 
$\mathcal{E}^{k}_{\widetilde{U}}$ le faisceau des $k$-formes 
différentielles \`a valeur complexe sur $\widetilde{U}$. Notons $(\pi_{\ast} 
\mathcal{E}^{k}_{\widetilde{U}})^{G}$ le sous-faisceau $G$-invariant 
de $\pi_{\ast}\mathcal{E}^{k}_{\widetilde{U}}$.  Nous allons définir 
le faisceau des $k$-formes différentielles sur l'orbifold $X$ comme 
un sous-faisceau de $j_{\ast}\mathcal{E}_{|X_{\reg}|}$.  Pour tout 
ouvert $U$ dans $|X|$, posons 
\begin{align*} 
\mathcal{E}^{k}_{|X|}(U)&:=\left\{ 
\begin{array}{l} 
\omega\in j_{\ast}\mathcal{E}^{k}_{X_{\reg}}(U)\mid \forall x\in 
U_{sing},\exists (\widetilde{U}_{x},G_{x},\pi_{x}) \mbox{ carte de  
  } U_{x}\\ 
 \mbox{ et } \omega'\in ({\pi_{x}}_{\ast}\mathcal{E}^{k}_{\widetilde{U}_{x}})^{G_{x}}(U_{x}) 
\mbox{ telles que } \omega'={\omega} \mbox{ sur }U_{x,reg} 
\end{array}\right\}. 
\end{align*} 
 
De la m\^{e}me mani\`ere que dans l'exemple \ref{expl:Les-fibrs-vectoriels},
nous pouvons d\'efinir le complexifi\'e du fibr\'e cotangent r\'eel
orbifold de $X$. Nous le notons $T^{\star}_{\CC}X$. Son faisceau des
sections est $\mathcal{E}^{1}_{|X|}$. De mani\`ere g\'en\'erale, le
faisceau des sections du fibr\'e orbifold $\wedge^{k} T^{\star}_{\CC}X$ est $\mathcal{E}^{k}_{|X|}$.

\begin{prop}\label{prop:cohomologie} Supposons $|X|$ paracompact. 
\begin{enumerate} \item \label{prop:cohomologie,acyclique} Les faisceaux $\mathcal{E}^{k}_{|X|}$ sont 
  acycliques pour le foncteur $\Gamma(|X|,\cdot)$. 
\item Il existe une différentielle $d$ telle que 
  \begin{align*} 
  \xymatrix{\mathcal{E}_{|X|}^{\bullet}: & \mathcal{E}^{0}_{|X|} 
    \ar[r]^-{d} & \mathcal{E}^{1}_{|X|} \ar[r]^-{d} & \cdots \ar[r]^-{d} 
    & \mathcal{E}^{n}_{|X|}\ar[r]^-{d} & 0} 
  \end{align*} 
  soit une résolution du faisceau constant $\underline{\CC}_{|X|}$. 
  \end{enumerate} 
\end{prop}

\begin{proof} 
\begin{enumerate} \item   
  D'après la fin du paragraphe \ref{sec:les-orbif-comm}, le faisceau 
  $\mathcal{E}^{0}_{|X|}$ est le faisceau des fonctions $C^{\infty}$ 
  sur $X$.  Comme $|X|$ est paracompact, il existe des 
  partitions de l'unité dans $\mathcal{E}^{0}_{|X|}$ (cf. proposition 
  \ref{prop:partition,unite}).  Ainsi, le faisceau 
  $\mathcal{E}^{0}_{|X|}$ est fin. Comme $\mathcal{E}^{k}_{|X|}$ 
  est un faisceau de $\mathcal{E}^{0}_{|X|}$-modules, nous en déduisons que le 
  faisceau $\mathcal{E}^{k}_{|X|}$ est fin. 
   
\item Comme la condition est locale, il suffit de la vérifier pour 
  tout ouvert assez petit. Soit 
  $(\widetilde{U},G,\pi)$ une carte de $U$.  Comme $\widetilde{U}$ est 
  un ouvert de $\CC^{n}$, le complexe de De Rham 
  \begin{align*} 
  \xymatrix{\mathcal{E}^{\bullet}_{\widetilde{U}}  :& 
    \mathcal{E}^{0}_{\widetilde{U}} \ar[r]^-{d}& 
    \mathcal{E}^{1}_{\widetilde{U}} \ar[r]^-{d} & \cdots \ar[r]^-{d} & 
    \mathcal{E}^{n}_{\widetilde{U}} \ar[r]^-{d}& 0} 
  \end{align*} 
  est une résolution de $\underline{\CC}_{\widetilde{U}}$.   Puis, 
  nous appliquons le foncteur exact à gauche $\pi_{\ast}$ à 
  $\mathcal{E}^{\bullet}_{\widetilde{U}}$. Pour tout 
  ${x}\in{U}$, nous avons 
  \begin{align*}(R^{i}\pi_{\ast} 
  \mathcal{E}_{\widetilde{U}}^{k})_{{x}}&=\lim_{V\mid x\in 
    V}H^{i}(\pi^{-1}(V), \mathcal{E}_{\widetilde{U}}^{k}\mid_{V}). 
  \end{align*} 
  Comme $ \mathcal{E}_{\widetilde{U}}^{k}$ est un faisceau fin, nous 
  avons  
  $H^{i}(\pi^{-1}(V),\mathcal{E}_{\widetilde{U}}^{k}\mid_{V})=0$ pour 
  $i>0$. Ainsi, $(R^{i}\pi_{\ast} 
  \mathcal{E}_{\widetilde{U}}^{k})_{{x}}$ est nul c'est-à-dire que le complexe 
  \begin{align*} 
  \xymatrix{\pi_{\ast}\mathcal{E}^{\bullet}_{\widetilde{U}}  :& 
    \pi_{\ast}\mathcal{E}^{0}_{\widetilde{U}} \ar[r]^-{\pi_{\ast}d}& 
    \pi_{\ast}\mathcal{E}^{1}_{\widetilde{U}} \ar[r]^-{\pi_{\ast}d} & 
    \cdots \ar[r]^-{\pi_{\ast}d} & 
    \pi_{\ast}\mathcal{E}^{n}_{\widetilde{U}} \ar[r]^-{\pi_{\ast}d}& 
    0} 
  \end{align*} 
  est une résolution de 
  $\pi_{\ast}\underline{\CC}_{\widetilde{U}}$. Les différentielles 
  $\pi_{\ast}d$ sont les restrictions de la différentielle $d$. 
   
  Puis, nous appliquons le foncteur covariant exact à gauche qui à un 
  faisceau $\mathcal{F}$ sur lequel $G$ agit, associe le sous-faisceau 
  $G$-invariant, noté $\mathcal{F}^{G}$.  Nous en déduisons le 
  complexe 
  \begin{align*} 
  \xymatrix{(\pi_{\ast}\mathcal{E}^{\bullet}_{\widetilde{U}})^{G} 
     :& (\pi_{\ast}\mathcal{E}^{0}_{\widetilde{U}})^{G} 
    \ar[r]^-{\pi_{\ast}d}&(\pi_{\ast}\mathcal{E}^{1}_{\widetilde{U}})^{G} 
    \ar[r]^-{\pi_{\ast}d} & \cdots \ar[r]^-{\pi_{\ast}d} & 
    (\pi_{\ast}\mathcal{E}^{n}_{\widetilde{U}})^{G}}. 
  \end{align*} 
   
  Il reste à montrer que ce complexe est encore exact. 
   
  Soit ${x}$ un point de ${U}$. Soit 
  $\widetilde{\omega}\in 
  (\pi_{\ast}\mathcal{E}^{k}_{\widetilde{U}})^{G}_{{x}}$ 
  tel que $(\pi_{\ast}d)\widetilde{\omega}=0$ où $k$ est un entier 
  strictement positif. Le 
  complexe $\pi_{\ast}\mathcal{E}^{\bullet}_{\widetilde{U}}$ est exact, 
  il existe $\widetilde{\eta}$ dans 
  $\pi_{\ast}\mathcal{E}^{k-1}_{\widetilde{U},{x}}$ tel que 
  $(\pi_{\ast}d)\widetilde{\eta}=\widetilde{\omega}$. Puis on pose 
  $\widetilde{\alpha}:=\frac{1}{\# G}\sum_{g\in 
    G}g^{\ast}\widetilde{\eta}$.  Il est clair que 
  $\widetilde{\alpha}$ est $G$-invariante et 
  \begin{align*} (\pi_{\ast}d) \widetilde{\alpha} &= 
    \frac{1}{ \# G}\sum_{g\in G} (\pi_{\ast}d)(g^{\ast} \widetilde{\eta}) \\ 
    &= \frac{1}{\# G}\sum_{g\in G} g^{\ast}(\pi_{\ast}d) \widetilde{\eta} \\ 
    & = \widetilde{\omega}. \end{align*} 
   
  Nous en déduisons aisément que le complexe 
  $(\pi_{\ast}\mathcal{E}^{\bullet}_{\widetilde{U}})^{G}$ est une 
  résolution de $(\pi_{\ast}\CC_{\widetilde{U}})^{G}$.  Or nous avons 
  les égalités suivantes 
  $(\pi_{\ast}\CC_{\widetilde{U}})^{G}=\underline{\CC}_{|U|}=\underline{\CC}_{|X|}\mid_{|U|}$. 
 \end{enumerate} 
 \end{proof} 
  
 La proposition précédente implique directement le corollaire suivant. 
 
\begin{cor}\label{cor:egalite} 
  Pour $k\in\{0,\ldots,n\}$, on a un isomorphisme d'espace vectoriel entre 
  $H^{k}(|X|,\CC)$ et 
  $H^{k}\left(\Gamma(|X|,\mathcal{E}^{\bullet}_{|X|})\right)$. 
\end{cor} 
 
Supposons que $|X|$ soit compacte.  Soit $\omega$ dans 
$H^{n}(|X|,\CC)$. Grâce au corollaire \ref{cor:egalite}, nous voyons 
$\omega$ comme une classe de forme différentielle et nous définissons 
\begin{align}\label{eq:defi,integrale} 
\int^{\orb}_{X}\omega & :=\frac{1}{\# \ker (X)}\int_{|X_{\reg}|}\omega. 
\end{align}

 \section{Classes de Chern  orbifoldes par la théorie de Chern-Weil} 
 \label{sec:classe,Chern}  

 Dans cette section, nous allons définir les classes de Chern pour les 
 fibrés vectoriels complexes orbifolds. Nous adaptons l'appendice C du  
 livre de Milnor et Stasheff \cite{MScc} aux orbifolds. 
 
Nous commençons par une étude locale.  Soit $\pr : 
F\rightarrow U$ un fibré vectoriel complexe orbifold trivial sur une 
orbifold complexe $U$. Soient $(\widetilde{U},G,\pi)$ une carte de $U$ 
et $(\widetilde{U}\times\CC^{r},G,\pi_{F})$ une carte de $F$. 
L'action de $G$ sur $\widetilde{F}:=\widetilde{U}\times\CC^{r}$ est donnée par 
\begin{align*} 
  g\cdot(x,w)&=(gx,\rho(x,g)w),\  \forall (g,x,w)\in G\times\widetilde{U}\times\CC^{r}  
\end{align*} 
où $\rho: \widetilde{U}\times G \rightarrow \GL (r,\CC)$ est holomorphe. 
Rappelons que l'action de $G$ sur une section ${s}$ est 
donnée par 
\begin{align*} 
  g\cdot{s}(\widetilde{x}):=\rho(g^{-1}\widetilde{x},g){s}(g^{-1}\widetilde{x}). 
\end{align*}   
 
Soit $s$ une section $C^{\infty}$ 
du fibré trivial $\widetilde{F}\rightarrow \widetilde{U}$. Soit 
$\omega$ une $1$-forme $C^{\infty}$ sur $\widetilde{U}$.  
Nous définissons l'action de $G$ sur $\omega\otimes s$ par la formule suivante 
\begin{align*} 
  g\cdot (\omega\otimes s):= (g\cdot\omega)\otimes (g\cdot s). 
\end{align*} 
Une \emph{connexion $\widetilde{\nabla}$ $G$-équivariante} sur le 
fibré trivial $\widetilde{F}\rightarrow \widetilde{U}$ est une 
connexion qui vérifie $\widetilde{\nabla} (g\cdot{s})=g\cdot 
(\widetilde{\nabla}{s})$ où ${s}$ est une section $C^{\infty}$ 
du fibré trivial $\widetilde{F}\rightarrow \widetilde{U}$\footnote{Les connexions 
  $G$-\'equivariantes existent.  En effet, si nous avons une connexion 
  $\nabla$ sur le fibré trivial $\widetilde{F}\rightarrow 
  \widetilde{U}$, la connexion 
\begin{align*} 
      \widetilde{\nabla}&:= \sum_{g\in G} g\cdot \nabla 
      (g^{-1}\cdot ) 
\end{align*}  
 est  $G$-équivariante.}. 
 
\begin{lem}\label{lem:connexion,changement,base} 
  Soit $\widetilde{\nabla}$ une connexion $G$-équivariante sur le 
  fibré trivial $\widetilde{U}\times\CC^{r}\rightarrow \widetilde{U}$. 
  Soit $\widetilde{\omega}$ la matrice de cette connexion dans une 
  base $s_{1}, \ldots ,s_{r}$ de sections du fibré 
  $\widetilde{U}\times\CC^{r}\rightarrow \widetilde{U}$.  Pour tout 
  $g\in G$, nous avons 
\begin{align*} 
  \varphi_{g^{-1}}^{\ast}\widetilde{\omega}&=\rho(\cdot,g)^{-1}\widetilde{\omega}\rho(\cdot,g)+\rho(\cdot,g)^{-1}d\rho(\cdot,g). 
\end{align*} 
\end{lem} 
 
\begin{proof} 
  Nous allons calculer $\widetilde{\nabla}(gs_{i})$ de deux manières. 
  D'abord, nous utilisons la $G$-équivariance de la connexion 
  $\widetilde{\nabla}$. Pour tout $g\in G$, nous avons 
 \begin{align*} 
   \widetilde{\nabla}(gs_{i})&=g\widetilde{\nabla}s_{i}\\ 
  &=g \sum_{i=1}^{r}\widetilde{\omega}_{ij}\otimes s_{j}\\ 
&= \sum_{i=1}^{r}\varphi_{g^{-1}}^{\ast}\widetilde{\omega}_{ij}\otimes gs_{j}. 
 \end{align*} 
 D'un autre c\^{o}té, nous savons comment la matrice d'une connexion se 
 comporte par un changement de base. Ainsi dans la base $g s_{1}, 
 \ldots ,g s_{r}$, la matrice de la connexion est 
 \begin{align*} 
   \rho(\cdot,g)^{-1}\widetilde{\omega}\rho(\cdot,g)+\rho(\cdot,g)^{-1}d\rho(\cdot,g). 
\end{align*} 
Nous en déduisons le lemme. 
\end{proof}  
 
Soit $(\widetilde{U},G,\pi)$ une carte de $U$. Le fibré orbifold trivial $\pr: F\rightarrow U$ est muni d'une 
connexion $\nabla$ s'il existe 
$(\widetilde{F}:=\widetilde{U}\times\CC^{r},G,\pi_{F})$ une carte de $F$ et une 
connexion $G$-équivariante telles que  
l'application 
$\widetilde{pr}:\widetilde{F}\rightarrow \widetilde{U}$ 
relève $\pr$.

Pour $i\in\{1,2\}$, soient $(F_{i},\nabla_{i})$ deux fibrés orbifolds 
triviaux  sur $U$ munis d'une connexion. 
Nous dirons que $(F_{1},\nabla_{1})$ et $(F_{2},\nabla_{2})$ 
sont \emph{isomorphes} s'il existe un isomorphisme  $\varphi$ entre les 
cartes $(\widetilde{F}_{1},G_{1},\pi_{F_{1}})$ et 
$(\widetilde{F}_{2},G_{2},\pi_{F_{2}})$ de respectivement $F_{1}$ et $F_{2}$ tel que 
$\varphi^{\ast}\widetilde{\nabla}_{2}=\widetilde{\nabla}_{1}$.

Le fibré trivial $F_{1}\rightarrow U_{1}$ 
muni d'une connexion ${\nabla}_{1}$ 
\emph{s'injecte} dans le fibré ${F}_{2}\rightarrow {U}_{2}$ 
muni d'une connexion ${\nabla}_{2}$ s'il existe une 
injection de fibrés orbifolds $(\alpha,\psi_{\alpha})$ c'est-à-dire un diagramme commutatif 
\begin{align*} 
  \xymatrix{\widetilde{F}_{1} \ar@{^{(}->}[rr]^{(\alpha,\psi_{\alpha})} 
    \ar[d]& &  
    \widetilde{F}_{2} \ar[d] \\ 
    \widetilde{U}_{1} \ar@{^{(}->}[rr]^{\alpha}&& \widetilde{U}_{2}} 
\end{align*} 
telle que  
$(\alpha^{\ast}\widetilde{F}_{2},\alpha^{\ast}\widetilde{\nabla}_{2})$ 
soit isomorphe à $(\widetilde{F}_{1},\widetilde{\nabla}_{1})$.

\begin{lem}\label{lem:injection,fibre,connexion} 
   Soit $(F_{2},\nabla_{2})$ un fibré orbifold trivial sur $U_{2}$ muni 
  d'une con\-nexion. Soit $F_{1}$ un fibré orbifold trivial sur $U_{1}$ 
  qui s'injecte dans $F_{2}$.  Il existe une connexion 
  $\nabla_{1}$ unique, à isomorphisme près, sur $F_{1}$ telle que 
  $(F_{1},\nabla_{1})$ s'injecte dans $(F_{2},\nabla_{2})$. 
\end{lem} 
 
\begin{proof}Pour $i\in\{1,2\}$, soient $(\widetilde{U}_{i},G_{i},\pi_{i})$ 
  deux cartes de $U_{i}$.  Soient 
  $(\widetilde{F}_{i},G_{i},\pi_{F_{i}})$ deux cartes de $F_{i}$. 
  Soit $\widetilde{\nabla}_{2}$ la connexion sur le fibré trivial 
  $\widetilde{F}_{2}\rightarrow \widetilde{U}_{2}$.  D'après le lemme 
  \ref{lem:injections,fibre,triviaux}, nous pouvons supposer que 
  $\widetilde{U}_{1}\subset \widetilde{U}_{2}$ et 
  $\widetilde{F}_{1}\subset \widetilde{F}_{2}$.   
 
\textbf{Existence :} 
  Nous posons 
  $\widetilde{\nabla}_{1}:=\widetilde{\nabla}_{2}\mid_{\widetilde{F}_{1}}$. 
  Ceci définit bien une connexion $G_{1}$-équivariante sur 
  $\widetilde{F}_{1}\rightarrow \widetilde{U}_{1}$. 
 
\textbf{Unicité:} Soient $(\widetilde{F}_{3},G_{3},\pi_{F_{3}})$ une 
carte de $F_{1}$ et $\widetilde{\nabla}_{3}$ une con\-nexion 
$G_{3}$-équiva\-ri\-ante sur le fibré trivial $\widetilde{F}_{3}\rightarrow 
\widetilde{U}_{1}$ telles qu'on ait le diagramme commutatif  
 \begin{align*} 
  \xymatrix{\widetilde{F}_{3} \ar@{^{(}->}[rr]^{(\alpha,\psi_{\alpha})} 
    \ar[d]& &  
    \widetilde{F}_{2} \ar[d] \\ 
    \widetilde{U}_{1} \ar@{^{(}->}[rr]^{\alpha}&& \widetilde{U}_{2}} 
\end{align*} 
Supposons que $(\alpha^{\ast}\widetilde{F}_{2},\alpha^{\ast}\widetilde{\nabla}_{2})$ 
soit isomorphe à $(\widetilde{F}_{3},\widetilde{\nabla}_{3})$. 
D'après le lemme \ref{lem:injections,fibre,triviaux}, le fibré 
$\widetilde{F}_{3}\rightarrow \widetilde{U}_{1}$ est isomorphe au 
fibré 
$\widetilde{F}_{1}:=(\alpha,\psi_{\alpha})(\widetilde{F}_{3})\rightarrow 
\alpha(\widetilde{U}_{1})$. 
Pour finir, il suffit de remarquer que 
$\alpha^{\ast}\widetilde{F}_{2}=\alpha^{\ast}\widetilde{F}_{1}$ et $\alpha^{\ast}(\widetilde{\nabla}_{2}\mid_{\widetilde{F}_{1}})=\alpha^{\ast}\widetilde{\nabla}_{2}$. 
\end{proof} 
 
\begin{defi} 
  Soit $\pr:F\rightarrow X$ un fibré orbifold. Une \emph{connexion orbifolde}  
  sur le fibré 
  orbifold  $\pr:F\rightarrow X$ est la donnée pour chaque $x\in X$ d'une carte 
  $(\widetilde{U}_{x},G_{x},\pi_{x})$ d'un ouvert $U_{x}$ contenant 
  $x$ telle que  
  \begin{enumerate} 
  \item $\pr^{-1}(U_{x})\rightarrow U_{x}$ est un fibré orbifold 
    trivial muni d'une connexion $\nabla_{x}$; 
  \item pour toute injection $\alpha:\widetilde{U}_{x}\hookrightarrow 
    \widetilde{U}_{y}$, il existe une injection 
    $(\pr^{-1}(U_{x}),\nabla_{x})\hookrightarrow(\pr^{-1}(U_{y}),\nabla_{y})$. 
  \end{enumerate} 
\end{defi} 
 
Rappelons que le fibr\'e orbifold $T^{\ast}_{\CC}X$ est le
complexifi\'e du fibr\'e cotangent r\'eel de $X$ et que son faisceau
des sections est $\mathcal{E}^{1}_{|X|}$ (cf. paragraphe \ref{appendice:A}).
 
\begin{prop}\label{prop:existence,connexion} 
  Soit $X$ une orbifold paracompacte. Soit $F\rightarrow X$ un fibré 
  orbifold. 
  \begin{enumerate} 
  \item Il existe une connexion sur $F$. 
  \item La différence entre deux connexions orbifoldes est une section  
    du fibré $T^{\ast}_{\CC}X\otimes_{\mathcal{C}^{\infty}_{|X|}}\End(F)$. 
  \end{enumerate} 
\end{prop}

\begin{proof} La première partie de la proposition est démontrée dans 
  le lemme $4.3.2$ de \cite{CRogw}. 
  Soit $\widetilde{F}\rightarrow \widetilde{U}$ un 
  fibré trivial orbifold muni de deux connexions $G$-équivariantes 
  notées $\widetilde{\nabla}_{1}$ et $\widetilde{\nabla}_{2}$. Il  
  suffit de voir que $\widetilde{\nabla}_{1}-\widetilde{\nabla}_{2}$ 
  est une section orbifold du fibré 
  $T^{\ast}_{\CC}\widetilde{U}\otimes_{\mathcal{C}^{\infty}_{\widetilde{U}}}\End(\widetilde{F})$. 
Soit $s_{1}, \ldots ,s_{r}$ une base de sections du fibré 
$\widetilde{F}\rightarrow \widetilde{U}$.  
Notons $\widetilde{\eta}$ la matrice de 
$\widetilde{\nabla}_{1}-\widetilde{\nabla}_{2}$ dans cette base. 
Montrons que $g\cdot \widetilde{\eta}=\widetilde{\eta}$, où l'action de $G$ sur la matrice de 
$1$-forme différentielle est donnée par  
\begin{align*} 
  g\cdot \widetilde{\eta} & = \rho(\cdot,g)\left(\varphi_{g^{-1}}^{\ast}\widetilde{\eta}\right)\rho(\cdot,g)^{-1}  
   \mbox{ où } (\varphi_{g^{-1}}^{\ast}\widetilde{\eta})_{ij}=\varphi_{g^{-1}}^{\ast}\widetilde{\eta}_{ij}. 
\end{align*} 
L'égalité ci-dessus est une conséquence directe du lemme \ref{lem:connexion,changement,base}. 
\end{proof}

\begin{lem}\label{lem:courbure,fibre,trivial}  Soit 
  $(\widetilde{U},G,\pi)$ une carte de $U$. Soit $\pr:F\rightarrow U$ 
  un fibré orbifold trivial muni d'une connexion $\nabla$.  
  Soit 
  $\widetilde{\pr}:\widetilde{F}:=\widetilde{U}\times\CC^{r}\rightarrow 
  \widetilde{U}$ le relevé du fibré orbifold trivial $\pr:F\rightarrow 
  U$  muni 
  d'une connexion $\widetilde{\nabla}$ $G$-équivariante. La courbure, notée $K(\widetilde{\nabla})$, de la connexion 
  $\widetilde{\nabla}$ est une section $G$-invariante du fibré 
  $\bigwedge^{2}T^{\ast}_{\CC}\widetilde{U}\otimes_{\mathcal{C}^{\infty}_{\widetilde{U}}} 
  \End(\widetilde{F})\rightarrow \widetilde{U}$. 
\end{lem} 
 
\begin{proof} 
Soit $s_{1}, \ldots ,s_{r}$ une base de sections du fibré trivial 
$\widetilde{U}\times \CC^{r}\rightarrow \widetilde{U}$. 
Notons $\widetilde{\Omega}$ la matrice de la courbure dans cette base. 
 Montrons que $g\cdot \widetilde{\Omega}=\widetilde{\Omega}$, où l'action du groupe $G$  
sur la matrice de $2$-formes est donnée par 
\begin{align}\label{eq:29} 
  g\cdot\widetilde{\Omega}&=\rho(\cdot,g)(\varphi_{g^{-1}}^{\ast}\widetilde{\Omega})\rho(\cdot,g)^{-1} \mbox{ où }(\varphi_{g^{-1}}^{\ast}\widetilde{\Omega})_{ij}=\varphi_{g^{-1}}^{\ast}\widetilde{\Omega}_{ij}. 
\end{align} 
 
Avant de démontrer cette égalité, nous allons montrer que  
\begin{align*} 
  K(\widetilde{\nabla})(g s)&=gK(\widetilde{\nabla})(s), 
\end{align*} 
pour tout $(g,s)\in G\times 
\widetilde{\mathcal{F}}_{\widetilde{U}}(\widetilde{U})$, où 
$\widetilde{\mathcal{F}}_{\widetilde{U}}$ est le faisceau des 
sections $C^{\infty}$ du fibré $\widetilde{F}\rightarrow \widetilde{U}$. 
La courbure $K(\widetilde{\nabla})$ est la composée de 
$\widehat{\widetilde{\nabla}}$ avec $\widetilde{\nabla}$ où  
\begin{align*} 
  \widehat{\widetilde{\nabla}}: 
  \widetilde{\mathcal{F}}_{\widetilde{U}}\otimes_{C^{\infty}_{\widetilde{U}}}\mathcal{E}^{1}_{\widetilde{U}} \longrightarrow \widetilde{\mathcal{F}}_{\widetilde{U}}\otimes_{C^{\infty}_{\widetilde{U}}}\mathcal{E}^{2}_{\widetilde{U}} 
\end{align*} est définie par  
\begin{align*} 
   \widehat{\widetilde{\nabla}}(\theta\otimes 
   s)&:=d\theta\otimes s-\theta\wedge\widetilde{\nabla}s. 
\end{align*} 
Un calcul direct montre que 
$\widehat{\widetilde{\nabla}}(g\cdot(\theta\otimes 
s))=g\widehat{\widetilde{\nabla}}(\theta\otimes s)$ pour tout $g\in G$ 
et pour tout $(\theta,s)\in 
\widetilde{\mathcal{F}}_{\widetilde{U}}(\widetilde{U})\times 
\mathcal{E}^{1}_{\widetilde{U}}(\widetilde{U})$. 
Nous en déduisons que $K(\widetilde{\nabla})(g 
s)=gK(\widetilde{\nabla})(s)$ pour tout $(g,s)\in G\times 
\widetilde{\mathcal{F}}_{\widetilde{U}}(\widetilde{U})$. 
 
Pour démontrer l'égalité (\ref{eq:29}), nous allons calculer 
$K(\widetilde{\nabla})(gs_{i})$ de deux façons. 
D'un c\^{o}té nous avons 
\begin{align} 
  K(\widetilde{\nabla})(gs_{i})&=gK(\widetilde{\nabla})(s_{i})\nonumber \\ 
&= g\left(\sum_{j=1}^{r}\widetilde{\Omega}_{ij}\otimes s_{j}\right) \nonumber\\\label{eq:30} 
&= \sum_{j=1}^{r} \varphi_{g^{-1}}^{\ast}\widetilde{\Omega}_{ij}\otimes gs_{j}. 
\end{align} 
D'un autre c\^{o}té, le changement de base entre $(s_{1}, \ldots ,s_{r})$ 
et  $(gs_{1}, \ldots ,gs_{r})$ montre que la matrice de 
$K(\widetilde{\nabla})$ dans la base $(gs_{1}, \ldots ,gs_{r})$ est  
\begin{align}\label{eq:31} 
  \rho(\cdot,g)\widetilde{\Omega}\rho(\cdot,g)^{-1}. 
\end{align} 
En comparant les égalités (\ref{eq:30}) et (\ref{eq:31}), nous en 
déduisons l'égalité (\ref{eq:29}) c'est-à-dire que 
$K(\widetilde{\nabla})$ est une section $G$-équivariante du fibré  
$\bigwedge^{2}T^{\ast}_{\CC}\widetilde{U}\otimes_{\mathcal{C}^{\infty}_{\widetilde{U}}} 
\End{\widetilde{F}}\rightarrow \widetilde{U}$. 
 \end{proof}

 Soit $\pr: F\rightarrow X$ un fibré orbifold de rang $r$ muni d'une 
 connexion $\nabla$.  Pour tout $x\in|X|$, il existe une carte 
 $(\widetilde{U}_{x},G_{x},\pi_{x})$ d'un ouvert $U_{x}$ contenant $x$ 
 telle que le fibré trivial 
 $\widetilde{F}_{x}:=\widetilde{U}_{x}\times \CC^{r}\rightarrow 
 \widetilde{U}_{x}$ soit muni d'une connexion, noté 
 $\widetilde{\nabla}_{x}$, $G_{x}$-équivariante où 
 $\widetilde{U}_{x}\times \CC^{r}$ est une carte de $\pr^{-1}(U_{x})$. 
 D'après le lemme \ref{lem:courbure,fibre,trivial}, la courbure 
 $K(\widetilde{\nabla}_{x})$ est une section $G_{x}$-invariante du 
 fibré $\bigwedge^{2}T^{\ast}_{\CC}\widetilde{U}_{x}\otimes_{\mathcal{C}^{\infty}_{\widetilde{U}_{x}}} 
 \End(\widetilde{F}_{x})$.  
  
 Pour $i\in\{1, \ldots ,n\}$, notons $\sigma_{i}$ le $i$-ième polyn\^{o}me 
 symétrique élémentaire en $r$ variables.  Soit $A$ une matrice carrée 
 de taille $r$. Notons $\sigma_{i}(A)$ le polyn\^{o}me $\sigma_{i}$ 
 évalué en les valeurs propres de $A$. Nous avons 
\begin{align*} 
  \det(\id+tA)=1+t\sigma_{1}(A)+\cdots+t^{n}\sigma_{n}(A). 
\end{align*}  
Pour tout $x\in |X|$, soit $s_{1,x}, \ldots ,s_{r,x}$ une base de 
sections du fibré $\widetilde{F}_{x}\rightarrow \widetilde{U}_{x}$. 
Notons $\widetilde{\Omega}_{x}$ la matrice de la courbure 
$K(\widetilde{\nabla}_{x})$ dans cette base.  Comme les $2$-formes 
différentielles commutent entre elles, nous pouvons évaluer les 
polyn\^{o}mes $\sigma_{i}$ en $\widetilde{\Omega}_{x}$.  La forme différentielle 
$\sigma_{i}(\widetilde{\Omega}_{x})$ est de degré $2i$. Dans une nouvelle base de 
sections, la matrice de la courbure est $P^{-1}\widetilde{\Omega}_{x} P$ où $P$ 
est la matrice de changement de base.  Comme 
$\sigma_{i}(P^{-1}AP)=\sigma_{i}(A)$, la forme différentielle 
$\sigma_{i}(\widetilde{\Omega}_{x})$ ne dépend pas de la base choisie.  Dorénavant, 
nous notons $\sigma_{i}(K(\widetilde{\nabla}_{x}))$ cette forme 
différentielle. Pour 
tout $g\in G_{x}$, nous avons 
\begin{align*} 
  \varphi_{g}^{\ast}\sigma_{i}(\widetilde{\Omega}_{x}):=&\sigma_{i}(\varphi_{g}^{\ast}\widetilde{\Omega}_{x}) \\ 
=&\sigma_{i}(\widetilde{\Omega}_{x}) \hspace{1cm} \mbox{ d'après l'égalité (\ref{eq:29})}. 
\end{align*} 
Ainsi, $\sigma_{i}(K(\widetilde{\nabla}_{x})$ définie une forme différentielle de 
degré $2i$ $G_{x}$-invariante sur $\widetilde{U}_{x}$ qui, par passage 
au quotient, définit 
une forme différentielle, notée  $\sigma_{i}(K({\nabla}_{x}))$, dans $\mathcal{E}^{2i}_{|U_{x}|}(U_{x})$.  
 
\begin{lem}\label{lem:forme,diff,globale} 
  Soit $F\rightarrow X$ un fibré orbifold muni d'une connexion 
  $\nabla$. 
Nous gardons les notations ci-dessus. Pour tout $i\in\{1, \ldots 
,r\}$, les formes différentielles  $\sigma_{i}(K({\nabla}_{x})$ dans 
$\mathcal{E}^{2i}_{|U_{x}|}(U_{x})$ se recollent  en une forme 
différentielle, notée $\sigma_{i}(K({\nabla}))$, dans $\mathcal{E}^{2i}_{|X|}(|X|)$. 
\end{lem}

\begin{proof}
Soit $\alpha:\widetilde{U}_{x}\hookrightarrow \widetilde{U}_{y}$. Il 
suffit de montrer que 
\begin{align*} 
  \sigma_{i}(K({\nabla}_{y}))\mid_{U_{x}}&=\sigma_{i}(K({\nabla}_{x})). 
\end{align*} 
Par 
définition d'un fibré orbifold muni d'une connexion, il existe deux 
applications $\psi_{\alpha}$ et $\varphi:\widetilde{F}_{x}\rightarrow 
\alpha^{\ast}\widetilde{F}_{y}$ telles que  
\begin{itemize} 
\item  le diagramme suivant soit commutatif 
\begin{align*} 
  \xymatrix{\widetilde{F}_{x} \ar@{^{(}->}[rr]^{(\alpha,\psi_{\alpha})} 
    \ar[d]& &  
    \widetilde{F}_{y} \ar[d] \\ 
    \widetilde{U}_{x} \ar@{^{(}->}[rr]^{\alpha}&& \widetilde{U}_{y}} 
\end{align*} 
\item  $\varphi$ soit un isomorphisme de cartes ; 
\item  $\varphi^{\ast}\alpha^{\ast}\widetilde{\nabla}_{y}=\widetilde{\nabla}_{x}$. 
\end{itemize} 
Nous en déduisons l'égalité 
$\varphi^{\ast}\alpha^{\ast}K(\widetilde{\nabla}_{y})=K(\widetilde{\nabla}_{x})$. 
Puis, nous obtenons l'égalité 
$\sigma_{i}(K({\nabla}_{y}))\mid_{U_{x}}=\sigma_{i}(K({\nabla}_{x}))$. 
\end{proof} 
  
 \begin{prop}\label{prop:poly,invariant,forme,fermee} Soit 
   $F\rightarrow X$ un fibré orbifold muni d'une connexion $\nabla$. 
  Pour tout $i\in\{1, \ldots ,n\}$, la  forme différentielle 
    $\sigma_{i}(K(\nabla))$ est fermée. 
\end{prop}  
 
\begin{proof} Soit 
   $F\rightarrow X$ un fibré orbifold muni d'une connexion $\nabla$.  
Nous en déduisons une connexion, notée $\nabla_{\reg}$, sur le fibré vectoriel $F_{\reg}\rightarrow X_{\reg}$. 
Nous en obtenons l'égalité   
\begin{align*} 
\sigma_{i}(K(\nabla_{\reg}))=\sigma_{i}(K(\nabla))\mid_{X_{\reg}}. 
\end{align*} 
Dans le cas des variétés, on sait que 
$\sigma_{i}(K(\nabla_{\reg}))$ est fermée. Ainsi $d\sigma_{i}(K(\nabla))$ est 
nulle sur un ouvert dense. Nous en déduisons que 
$\sigma_{i}(K(\nabla))$ est fermée. 
\end{proof} 
 
Cette proposition montre que les formes différentielles 
$\sigma_{0}(K(\nabla)), \ldots ,\sigma_{r}(K(\nabla))$ sont fermées de degré  
respectivement $0, \ldots ,2r$. Nous noterons 
$[\sigma_{i}(K(\nabla))]$ la classe de cette forme différentielle dans 
$H^{2i}(|X|,\CC)$ (cf. le corollaire \ref{cor:egalite}). 
 
\begin{cor}\label{cor:Chern,independance,connexion}Soit $F\rightarrow X$ un fibré orbifold muni d'une connexion. 
   La classe de la forme différentielle $\sigma_{i}(K(\nabla))$ ne dépend pas 
    du choix de la connexion sur $F$. 
  \end{cor} 
 
\begin{proof} 
    Soit $\pr:F\rightarrow X$ un fibré vectoriel orbifold muni de deux 
    connexions $\nabla_{0}$ et $\nabla_{1}$.  
Pour tout $x\in X$, nous avons deux connexions 
$\widetilde{\nabla}_{0,x}$ et $\widetilde{\nabla}_{1,x}$ sur 
$\widetilde{F}_{x}\rightarrow \widetilde{U}_{x}$ où 
$\widetilde{U}_{x}$ est une carte d'un ouver $U_{x}$ contenant $x$ et 
$\widetilde{F}_{x}$ une carte de $\pr^{-1}(U_{x})$. 
Posons 
$\widetilde{\eta}_{x}:=\widetilde{\nabla}_{1,x}-\widetilde{\nabla}_{0,x}$. 
D'après la proposition \ref{prop:existence,connexion}, 
$\widetilde{\eta}_{x}$ est une section $G_{x}$-invariante de 
$T^{\ast}_{\CC}\widetilde{U}_{x}\otimes_{\mathcal{C}^{\infty}_{\widetilde{U}_{x}}} \End(\widetilde{F}_{x})$. 
Pour $t\in[0,1]$, nous définissons une connexion sur le fibré 
$\widetilde{F}_{x}\rightarrow \widetilde{U}_{x}$ par la formule 
 
    \begin{align*} 
      \widetilde{\nabla}_{t,x}:= \widetilde{\nabla}_{0,x}+ t\widetilde{\eta}_{x}. 
    \end{align*} 
  
D'après les calculs dans \cite{GHag} p.$405$, nous avons  
\begin{align}\label{eq:32} 
  d\left(\int_{0}^{1}\widetilde{\sigma}_{i}\left(\widetilde{\eta}_{x},K(\widetilde{\nabla}_{t,x}),
      \ldots
      ,K(\widetilde{\nabla}_{t,x})\right)dt\right)&=\sigma_{i}(K(\widetilde{\nabla}_{1,x}))-\sigma_{i}(K(\widetilde{\nabla}_{0,x}))
\end{align} 
où $\widetilde{\sigma}_{i}$ est l'application $i$-linéaire de 
$M_{r}(\CC)\times \cdots \times M_{r}(\CC)$ dans $\CC$ telle que  
\begin{align*} 
  \widetilde{\sigma}_{i}(A, \ldots ,A)=\sigma_{i}(A) \mbox{ pour toute matrice} A\in M_{r}(\CC). 
\end{align*} 
  L'égalité (\ref{eq:32}) montre 
que 
$\sigma_{i}(K(\widetilde{\nabla}_{1,x}))-\sigma_{i}(K(\widetilde{\nabla}_{0,x}))$ 
est une forme exacte.   
 
Montrons que pour toute injection $\alpha: 
\widetilde{U}_{y}\hookrightarrow \widetilde{U}_{x}$, nous avons 
$$\alpha^{\ast}\widetilde{\sigma}_{i}\left(\widetilde{\eta}_{x},K(\widetilde{\nabla}_{t,x}), 
  \ldots 
  ,K(\widetilde{\nabla}_{t,x})\right)=\widetilde{\sigma}_{i}\left(\widetilde{\eta}_{y},K(\widetilde{\nabla}_{t,y}), 
  \ldots ,K(\widetilde{\nabla}_{t,y})\right).$$ 
D'après le lemme 
\ref{lem:injection,fibre,connexion}, nous pouvons supposer que  
$\alpha^{\ast}(\widetilde{\nabla}_{i,x})= \widetilde{\nabla}_{i,y}$ 
pour $i\in\{1,2\}$. 
Ceci implique que $\alpha^{\ast}\eta_{x}=\eta_{y}$ et que 
$\alpha^{\ast}K(\widetilde{\nabla}_{t,x})=K(\widetilde{\nabla}_{t,y})$. 
Finalement, nous avons 
\begin{align*} 
  \alpha^{\ast}\widetilde{\sigma}_{i}\left(\widetilde{\eta}_{x},K(\widetilde{\nabla}_{t,x}), 
  \ldots 
  ,K(\widetilde{\nabla}_{t,x})\right)=&\widetilde{\sigma}_{i}\left(a^{\ast}\widetilde{\eta}_{x},\alpha^{\ast}K(\widetilde{\nabla}_{t,x}), 
  \ldots 
  ,\alpha^{\ast}K(\widetilde{\nabla}_{t,x})\right)\\ 
=&\widetilde{\sigma}_{i}\left(\widetilde{\eta}_{y},K(\widetilde{\nabla}_{t,y}), 
  \ldots ,K(\widetilde{\nabla}_{t,y})\right). 
\end{align*} 
 
Les égalités ci-dessus nous montrent que les formes différentielles 
\begin{align*} 
  \widetilde{\sigma}_{i}\left(\widetilde{\eta}_{x},K(\widetilde{\nabla}_{t,x}), 
  \ldots ,K(\widetilde{\nabla}_{t,x})\right) 
\end{align*} 
de degré $2i-1$ sur 
$\widetilde{U}_{x}$ sont $G_{x}$-invariantes. Elles définissent par 
passage au quotient une forme différentielle de degré $2i-1$, notée 
$\widetilde{\sigma}_{i}\left({\eta}_{x},K({\nabla}_{t,x}), \ldots 
  ,K({\nabla}_{t,x})\right)$, sur $|X|$. Nous en déduisons que 
$\sigma_{i}(K(\nabla_{1}))-\sigma_{i}(K(\nabla_{0}))$ est une forme 
exacte. 
\end{proof}  
 
Soit $F\rightarrow X$ un fibré orbifold. Nous définissons la $i$-ième 
classe de Chern 
 du fibré $F\rightarrow X$ par 
\begin{align*} 
  c _{i}(F) :=\left[\sigma_{i}(K(\nabla))\right]/(2\sqrt{-1}\pi)^{i} 
\end{align*}  
où $\nabla$ est une connexion sur $F$. D'après le corollaire \ref{cor:Chern,independance,connexion}, les classes de Chern de $F$ ne 
dépendent pas de la connexion choisie et la $i$-ième classe de Chern 
définit une classe de cohomologie dans $H^{2i}(|X|,\CC)$.  La classe 
de Chern totale est  
\begin{align*} 
  c ({F})=c _{0}({F})+tc _{1}({F})+\cdots 
  +t^{n}c _{n}({F}). 
\end{align*}

\begin{prop}\label{prop:suite,exacte,fibre,chern} 
   Soit 
    \begin{align*}\xymatrix{ 
      0 \ar[r]& {E}\ar[r]^-{\alpha} & 
      {F} \ar[r]^-{\beta}& 
      {G} \ar[r]& 0} 
    \end{align*} 
    une suite exacte de fibrés orbifolds sur l'orbifold $X$. Nous avons l'égalité 
    \begin{align*} 
      c ({F})=c ({E})c ({G}). 
    \end{align*} 
\end{prop} 
 
 \begin{proof}  
   Cette propriété est locale.  Soit $p$ (resp.  $r$) le rang du 
   fibré $E$ (resp. $F$).  Pour tout $x\in |X|$, il existe une carte 
   $(\widetilde{U}_{x},G_{x},\pi_{x})$ d'un ouvert $U_{x}$ contenant $x$ qui 
   trivialise les trois fibrés. Soient $\widetilde{E}\mid_{U_{x}}$, 
   $\widetilde{F}\mid_{U_{x}}$ et $\widetilde{G}\mid_{U_{x}}$ trois 
   cartes de respectivement 
   $\pr_{E}^{-1}(U_{x})$,$\pr_{F}^{-1}(U_{x})$ et 
   $\pr_{G}^{-1}(U_{x})$. Nous avons une suite exacte de fibrés 
   triviaux sur $\widetilde{U}_{x}$ 
    \begin{align*}\xymatrix{ 
      0 \ar[r]& \widetilde{E}\mid_{U_{x}}\ar[r]^-{\widetilde{\alpha}\mid_{U_{x}}} & 
      \widetilde{F}\mid_{U_{x}} \ar[r]^-{\widetilde{\beta}\mid_{U_{x}}}& 
      \widetilde{G}\mid_{U_{x}} \ar[r]& 0.} 
    \end{align*} 
     Soit 
    $(s^{E}_{1}, \ldots ,s^{E}_{p})$ une base  des 
    sections du fibré $\widetilde{E}\mid_{U_{x}}\rightarrow\widetilde{U}_{x}$. Notons 
    $s^{F}_{i}:=\widetilde{\alpha}\mid_{U_{x}}(s^{E}_{i})$.  On complète $(s^{F}_{1}, \ldots 
    ,s^{F}_{p})$ en une base $(s^{F}_{1}, \ldots ,s^{F}_{r})$  
    des sections du fibré trivial $\widetilde{F}\mid_{U_{x}}\rightarrow\widetilde{U}_{x}$. Pour $i\in\{p+1, \ldots ,r\}$, les sections 
    $s^{G}_{i}:=\widetilde{\beta}\mid_{U_{x}}(s^{F}_{i})$ forment une base  
    des sections du fibré trivial $\widetilde{G}\mid_{U_{x}}\rightarrow\widetilde{U}_{x}$. 
  
 Nous choisissons une connexion $\widetilde{\nabla}_{F,x}$ $G_{x}$-équivariante de matrice $\widetilde{\omega}_{x}^{F}$ sur $\widetilde{F}\mid_{U_{x}}$ telle que sa 
 matrice soit triangulaire supérieure sur $\widetilde{U}_{x}$ dans la base $(s^{F}_{1}, \ldots ,s^{F}_{r})$. Nous en déduisons deux 
 connexions $\widetilde{\nabla}_{E,x}$ de matrice $\widetilde{\omega}_{x}^{E}$ et $\widetilde{\nabla}_{G,x}$ de matrice $\widetilde{\omega}_{x}^{G}$ sur $\widetilde{E}\mid_{U_{x}}$ et sur 
 $\widetilde{G}\mid_{U_{x}}$ qui sont définis par les matrices  
 \begin{align*} 
   \left(\widetilde{\omega}_{x}^{E}\right)_{ij}&:=\left(\widetilde{\omega}_{x}^{F}\right)_{ij} & & \mbox{ pour } i,j\in\{1, \ldots 
   ,p\}\,; \\ \left(\widetilde{\omega}_{x}^{G}\right)_{ij}&:=\left(\widetilde{\omega}_{x}^{F}\right)_{i+p,j+p}&& \mbox{ pour } i,j\in\{1, \ldots 
   ,r-p\}. 
 \end{align*} 
 La courbure $K(\widetilde{\nabla}_{F,x})$ sur $\widetilde{F}\mid_{U_{x}}$ est triangulaire supérieure. Nous en 
 déduisons que  
 \begin{align*} 
   \det(\id+tK(\widetilde{\nabla}_{F,x}))=\det(\id+tK(\widetilde{\nabla}_{E,x}))\det(\id+tK(\widetilde{\nabla}_{G,x})). 
 \end{align*} 
 \end{proof}

\section[Bonne application et image inverse de fibrés orbifolds]{Bonne application orbifolde et image inverse de fibrés 
  vectoriels orbifolds} 
\label{sec:bonne-appl-orbif} 
Dans la théorie générale des orbifolds, l'image inverse d'un fibré 
vectoriel orbifold n'existe pas toujours. C'est pour cette raison que 
Chen et Ruan ont défini la notion de bonne application orbifolde (cf. 
paragraphe $4.4$ de \cite{CRogw}).  Soit $|X|$ un espace topologique 
séparé. D'après Satake \cite{Sgb}, un recouvrement 
$\mathcal{U}=(U_{i})_{i\in I}$ est dit \emph{compatible} si 
 
\begin{enumerate}\makeatletter 
\renewcommand\theenumi{\theequation} 
\makeatother 
 \addtocounter{equation}{1} 
\item\label{item:12} chaque ouvert $U_{i}$ de ce recouvrement a une 
  carte $(\widetilde{U}_{i},G_{i},\pi_{i})$; 
  \addtocounter{equation}{1}\item \label{item:18} pour tout $x\in 
  U_{i}\cap U_{j}$, il existe un ouvert $U_{k}\subset U_{i}\cap U_{j}$ 
  tel que $x\in U_{k}$; \addtocounter{equation}{1}\item 
  \label{item:19} si $U_{i}\subset U_{j}$ alors il existe une 
  injection de $(\widetilde{U}_{i},G_{i},\pi_{i})$ dans 
  $(\widetilde{U}_{j},G_{j},\pi_{j})$. 
\end{enumerate} 
 
Un recouvrement compatible sur $|X|$ est un atlas orbifold (cf. le début 
du paragraphe \ref{sec:Les-atlas-orbifolds}).  
  
Nous avons la proposition suivante. 
\begin{prop}[cf. paragraphe $4.1$ p.$67$ de \cite{CRogw}]\label{prop:atlas,recouvrement} 
  Soit $|X|$ un espace topologique paracompact. Étant donné un atlas 
  orbifold sur $|X|$, il existe un recouvrement compatible plus fin sur $|X|$. 
\end{prop} 
 
\begin{rem}\label{rem:variete,recouvrement} 
  Si $X$ est une variété complexe, un atlas orbifold est aussi un atlas de 
  $X$ en tant que variété.  Or un atlas de variété vérifie bien les 
  conditions (\ref{item:12}), (\ref{item:18}) et (\ref{item:19}). 
  C'est-à-dire que pour une variété, la proposition ci-dessus est 
  vraie. 
 \end{rem}  
  
 \begin{proof}[Démonstration de la proposition \ref{prop:atlas,recouvrement}]

Soit $\mathcal{A}(|X|)$ un atlas orbifold sur $|X|$. Comme $|X|$ est paracompact, on peut supposer 
que le recouvrement associé à $\mathcal{A}(|X|)$ est localement 
fini. Notons-le $(U_{i})_{i\in I}$. 
Soit $x\in |X|$. Le point $x$ appartient à un nombre fini d'ouverts de 
$(U_{i})_{i\in I}$. Notons-les $U_{i_{1}}, \ldots ,U_{i_{n(x)}}$. 
Ainsi, il existe une carte $(\widetilde{V}_{x},G_{x},\pi_{x})$ d'un 
ouvert $V_{x}$ contenant $x$ tel que  
\begin{enumerate} 
\item $V_{x}$ est inclus dans $U_{i_{1}}\cap\ldots\cap U_{i_{n(x)}}$; 
\item la carte $\widetilde{V}_{x}$ s'injecte dans 
  $\widetilde{U}_{i_{\alpha}}$ pour tout $\alpha\in\{1, \ldots ,n(x)\}$. 
\end{enumerate} 
Notons $\mathcal{U}$ le recouvrement formé par de tels ouverts $V_{x}$  
et les ouverts inclus dans $V_{x}$ et contenant $x$. 
Montrons que $\mathcal{V}$ est un recouvrement compatible. 
Les  conditions (\ref{item:12}) et (\ref{item:18}) sont clairement 
vérifiées. 
\begin{figure}[bthp] 
\begin{center} 
  \psfrag{x}{$x$}  
\psfrag{ux}{$V_{k}$} 
\psfrag{vx}{$V_{x}$} 
\psfrag{y}{$y$} 
\psfrag{uy}{$V_{\ell}$} 
\psfrag{vy}{$V_{y}$} 
 \includegraphics[width=0.6\linewidth]{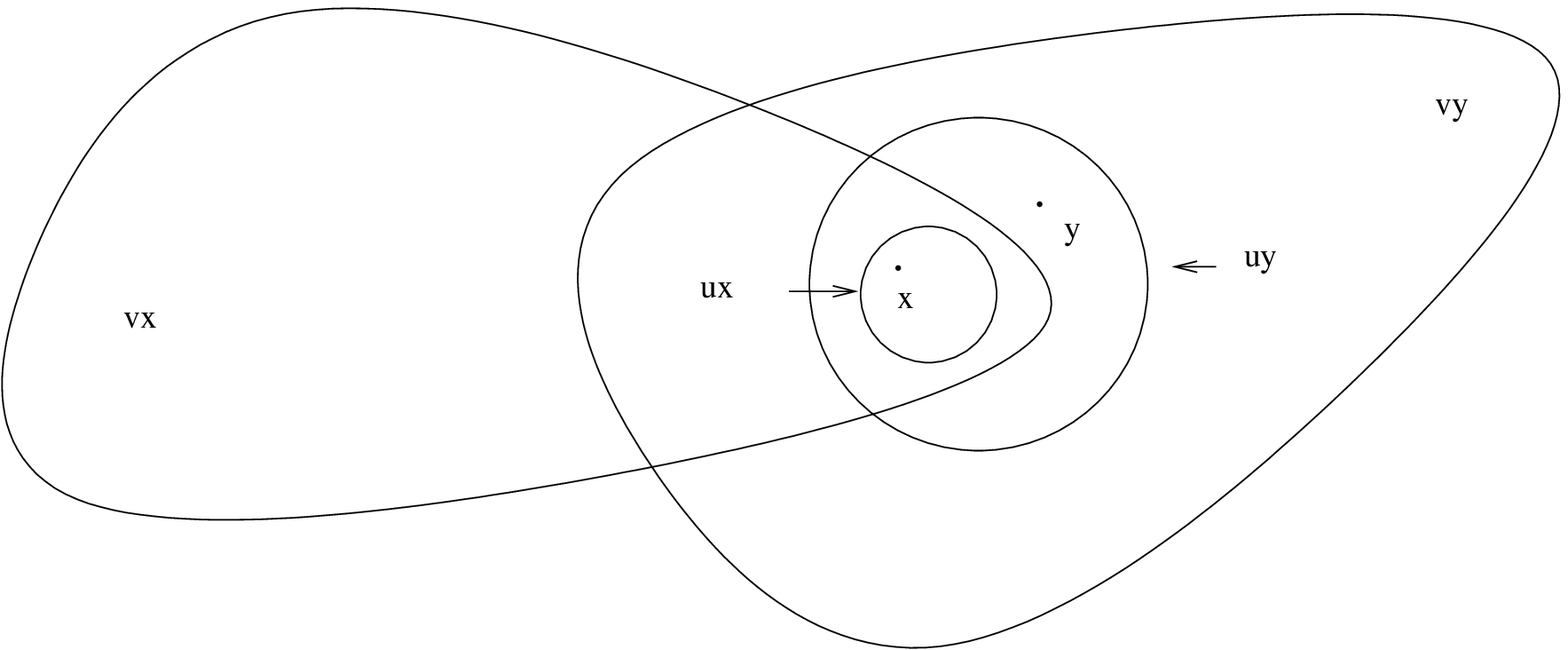} 
\end{center}\caption{}\label{fig:3} 
\end{figure} 
  
Il nous reste à montrer la condition (\ref{item:19}).  Soient
$V_{k},V_{\ell}$ deux ouverts de $\mathcal{V}$ tels que $V_{k}\subset
V_{\ell}$.  Montrons qu'il existe une injection entre une carte de
${V}_{k}$ et une carte de $V_{\ell}$. Il existe $x,y\in |X|$ tels que
$V_{k}\subset V_{x}$ et $V_{\ell}\subset V_{y}$ (cf. la figure
\ref{fig:3}).  Par d\'efinition, s'il existe $i\in I$ tel que nous
ayons $V_{k}\subset U_{i}$ (resp. $V_{k}\subset U_{i}$) alors $i\in
\{i_{1}, \ldots ,i_{n(x)}\}$ (resp. $i\in \{j_{1}, \ldots
,j_{n(y)}\}$).  L'inclusion de $V_{k}$ dans $V_{\ell}$ implique que
l'ensemble $\{j_{1}, \ldots ,j_{n(y)}\}$ est un sous-ensemble de
$\{i_{1}, \ldots ,i_{n(x)}\}$. Soit $\widetilde{V}_{k,x}$ (resp.
$\widetilde{V}_{\ell,y}$) une carte de $V_{k}$ (resp. $V_{\ell}$)
induite par $\widetilde{V}_{x}$ (resp. $\widetilde{V}_{y}$).  Nous en
d\'eduisons qu'il existe $i\in I$ tel que $\widetilde{V}_{x}$ et
$\widetilde{V}_{y}$ s'injectent dans $\widetilde{U}_{i}$.  Finalement,
le corollaire \ref{cor:injections} appliqué avec $U:=V_{k}$,
$V:=V_{\ell}$ et $W:=U_{i}$ montre qu'il existe une injection de
$\widetilde{V}_{k,x}$ dans $\widetilde{V}_{\ell,y}$.
 \end{proof}

\begin{defi}[cf. paragraphe $4.4$ de \cite{CRogw}]\label{defi bonne application} 
Soit $f:X\rightarrow Y$ une application orbifolde. On dit que $f$ est 
\emph{une bonne application} s'il existe deux recouvrements 
compatibles $\mathcal{U}_{X}$ et 
$\mathcal{V}_{Y}$  de respectivement $|X|$  
et $|Y|$ tels que  
\begin{enumerate} 
\item \label{item:7} il existe une correspondance bijective, notée $\mathfrak{F}$, 
  entre les ouverts de $\mathcal{U}_{X}$ et ceux de $\mathcal{V}_{Y}$ 
  telle que pour tout $U$ ouvert de $\mathcal{U}_{X}$, nous ayons $f(U)\subset \mathfrak{F}(U)$ et que si 
  $U_{1}\subset U_{2}$ alors $\mathfrak{F}(U_{1})\subset 
  \mathfrak{F}(U_{2})$; 
\item \label{item:8} l'ensemble des relèvements 
  $\widetilde{f}_{UV}:\widetilde{U}\rightarrow \widetilde{V}$ vérifie 
  la condition suivante : à chaque injection $\alpha: 
  (\widetilde{U}_{1},G_{1},\pi_{1})\hookrightarrow 
  (\widetilde{U}_{2},G_{2},\pi_{2})$,  correspond une injection, 
  notée $\mathfrak{F}(\alpha): (\widetilde{V_{1}},H_{1},\p_{1}) 
  \hookrightarrow (\widetilde{V_{2}},H_{2},\p_{2})$ où 
  $\p_{i}(\widetilde{V}_{i}):=\mathfrak{F}(U_{i})$ telle que 
\begin{itemize} \item le diagramme suivant soit commutatif 
  \begin{align*} 
  \xymatrix{ ( \widetilde{U}_{1},G_{1},\pi_{1}) \ar 
    @{^{(}->}[rr]^-{\alpha} 
    \ar[d]_-{\widetilde{f}_{U_{1}V_{1}}} & &( \widetilde{U}_{2},G_{2},\pi_{2}) \ar[d]^-{\widetilde{f}_{U_{2}V_{2}}} \\ 
    ( \widetilde{V}_{1},H_{1},\p_{1}) \ar@{^{(}->}[rr]^-{\mathfrak{F}(\alpha)}&& ( 
    \widetilde{V}_{2},H_{2},\p_{2})} 
  \end{align*} 
\item  pour 
  toutes composées d'injections $\alpha\circ \beta$, on ait 
  $\mathfrak{F}(\alpha\circ\beta)=\mathfrak{F}(\alpha)\circ\mathfrak{F}(\beta)$. 
\end{itemize} 
\end{enumerate} 
\end{defi} 
 
Soit $f:X\to Y$ une application orbifolde.
Nous appelons \emph{système compatible} l'ensemble des données supplémentaires 
$\{\widetilde{f}_{{U}{V}},\mathfrak{F}\}$ tel que l'application $f:X\rightarrow Y$ 
soit une bonne application.  Remarquons que la notation $\mathfrak{F}$ 
est utilisée pour deux correspondances : la première au niveau des 
ouverts des recouvrements compatibles et la deuxième au niveau des 
injections.

Chen et Ruan ont démontré les deux propositions suivantes dans l'article 
\cite{CRogw}. 
 
 \begin{prop}[cf. lemme $4.4.3$ dans \cite{CRogw}]\label{prop:pullback,fibre} 
  Soit $f:X\rightarrow Y$ une bonne application entre deux orbifolds. 
  Supposons qu'on ait un fibré vectoriel orbifold $E$ sur $Y$ 
  alors on peut définir un fibré vectoriel, noté 
  $f^{\ast}E$, sur $X$ par les fonctions de transition 
  suivantes : pour toute injection
  $\alpha:\widetilde{U}_{1}\hookrightarrow \widetilde{U}_{2}$ entre
  deux cartes du recouvrement compatible de $|X|$, on pose  
\begin{align*}
  \g_{\alpha}^{f^{\ast}E}(\widetilde{x}):=\g_{\mathfrak{F}(\alpha)}^{E}(\widetilde{f}(\widetilde{x})).
\end{align*}
\end{prop}

\begin{prop}[cf. lemme $4.4.3$ de \cite{CRogw}]\label{prop:chern,fonctoriel} 
  Les classes de Chern orbifoldes sont fonctorielles c'est-à-dire que 
   pour tout fibré vectoriel orbifold complexe $E\rightarrow Y$ et
   pour toute bonne application orbifolde  $f:X\rightarrow Y$, on a 
  $f^{\ast}c(\mathcal{E}^{\infty}_{|Y|})=c(f^{\ast}\mathcal{E}^{\infty}_{|Y|})$.
\end{prop} 
 
\begin{proof}[Démonstration de la proposition 
  \ref{prop:chern,fonctoriel}]  
Soient $\mathcal{U}_{X}$ et $\mathcal{V}_{Y}$ deux recouvrements 
compatibles de respectivement $|X|$ et $|Y|$. 
Notons $\mathfrak{F}$ la bijection  
\begin{align*} 
  \left\{\mbox{ ouverts de }\mathcal{U}_{X}\right\} & \rightarrow  \left\{\mbox{ ouverts de }\mathcal{V}_{Y}\right\}. 
\end{align*} 
Soit $\nabla$ une connexion sur le fibré $\pr:E\rightarrow Y$.  Soit $V$ 
un ouvert du recouvrement $\mathcal{V}_{Y}$. Comme $\mathfrak{F}$ est bijective, il existe 
un unique $U$ tel que $\mathfrak{F}(U)=V$.  Soit 
$\widetilde{f}_{U,V}:\widetilde{U}\rightarrow \widetilde{V}$ un 
relèvement de $f\mid_{U}: U \rightarrow V$. 
Soit $\widetilde{\nabla}_{\widetilde{V}}$ la connexion sur 
le fibré $\widetilde{E}_{V}:=\widetilde{V}\times\CC^{r}\rightarrow \widetilde{V}$ où 
$\widetilde{E}_{V}$ est une carte de $\pr^{-1}(V)$. 
Notons  $\widetilde{\omega}_{\widetilde{V}}$ la matrice de la connexion 
$\widetilde{\nabla}_{\widetilde{V}}$ dans une base de sections du fibré $\widetilde{E}_{V}\rightarrow \widetilde{V}$. 
Considérons la matrice $\widetilde{f}_{U,V}^{\ast}\widetilde{\omega}_{\widetilde{V}}$ de $1$-formes  
sur $\widetilde{U}$. Le diagramme commutatif appliqué avec 
$\alpha:=\varphi_{g}$ de la définition~\ref{defi bonne application} 
montre que la matrice de $1$-formes  
$\widetilde{f}_{U,V}^{\ast}\widetilde{\omega}_{\widetilde{V}}$ est 
$G_{\widetilde{U}}$-invariante. Cette matrice définit une connexion, 
notée $\widetilde{f}_{U,V}^{\ast}\widetilde{\nabla}_{\widetilde{V}}$, 
sur le fibré $\widetilde{U}\times\CC^{r}\rightarrow \widetilde{U}$. 
Le diagramme commutatif de la  définition \ref{defi bonne application}  
montre que nous avons bien défini  une connexion, notée $f^{\ast}\nabla$, sur le 
fibré $f^{\ast}E\rightarrow X$. 
Nous en déduisons l'égalité suivante, qui implique la proposition, 
\begin{align*} 
\widetilde{f}_{U,V}^{\ast}K(\widetilde{\nabla}_{\widetilde{V}})=K(\widetilde{f}_{U,V}^{\ast}\widetilde{\nabla}_{\widetilde{V}}). 
\end{align*} 
\end{proof} 
 
\begin{defi}[cf. définition $4.4$ de \cite{CRogw}] 
  Soit $f:X\rightarrow Y$ une application orbifolde. Deux systèmes 
  compatibles $\xi$ et $\xi'$ sont \emph{isomorphes} si pour tout fibré vectoriel orbifold 
  $E\rightarrow Y$, les deux fibrés vectoriel orbifold 
  $f^{\ast}_{\xi}E$ et $f^{\ast}_{\xi'}E$ sont isomorphes.  
\end{defi} 
 
Nous rappelons que la partie régulière d'une orbifold $X$ d'après Chen  
et Ruan est  
$\widehat{X}_{\reg}=\{x\in |X'|\mid G_{x}=\{\id\}\}$. 
 
\begin{prop}[cf. lemme $4.4.11$ dans \cite{CRogw}]\label{prop:application,regulier} 
Soit $f:X\rightarrow Y$ une application entre deux orbifolds. Si 
  $f^{-1}(\widehat{X}_{\reg})$ est un ouvert dense et connexe de $X$ alors il 
  existe un unique système compatible, à isomorphisme près. 
\end{prop}

\chapter{Cohomologie orbifolde des espaces projectifs à poids}\label{sec:struct-orbif-des} 
Dans ce chapitre, nous calculons l'anneau de cohomologie orbifolde des  
espaces projectifs à poids muni de la dualité de Poincaré orbifolde.

Dans le premier paragraphe, nous définissons la structure orbifolde 
sur les espaces projectifs à poids. Dans le deuxi\`eme, nous
d\'efinissons les fibr\'es $\mathcal{O}_{\PP(w)}(\cdot)$ sur $\PP(w)$. 
 
Les trois derniers suivants se présentent de la façon suivante. Dans 
chaque début de paragraphe, nous rappelons rapidement les définitions 
et les propriétés générales puis nous les appliquerons à notre exemple  
préféré $\PP(w)$.  
 
 Le 
paragraphe \ref{sec:coho-espace} traite de la cohomologie orbifolde en  
temps qu'espace vectoriel et la proposition \ref{prop:base} donne une 
base, notée $\bs{\eta}$, de la cohomologie orbifolde des espaces projectifs à poids. 
 
Le troisième paragraphe \ref{sec:dualite-de-poincare} est consacré à la dualité 
de Poincaré orbifolde et la proposition \ref{prop:dualite} calcule la 
dualité de Poincaré de $\PP(w)$ dans la base $\bs{\eta}$.  
 
Le dernier paragraphe concerne le cup produit orbifold. Le 
théorème  \ref{thm:fibre,obstruction} explicite le fibré obstruction 
et le corollaire \ref{cor:cup} donne une formule explicite pour le cup produit 
orbifold de $\PP(w)$ dans la base $\bs{\eta}$.

\section{La structure orbifolde sur $\PP(w)$}\label{sec:La-struct-orbif} 
Dans ce paragraphe, nous définissons l'atlas orbifold des espaces 
projectifs à poids.  
 
Nous définissons l'action du groupe multiplicatif $\CC^{\star}$ sur 
$\CC^{n+1}-\{0\}$ par 
\begin{align}\label{eq:action}\lambda\cdot(y_{0}, \ldots ,y_{n})&:=(\lambda^{w_{0}}y_{0}, \ldots 
,\lambda^{w_{n}}y_{n}). 
\end{align} L'espace projectif à poids est le 
quotient de $\CC^{n+1}-\{0\}$ par cette action. Notons $|\PP(w)|$ cet 
espace topologique quotient et $\pi_{w}$ l'application de passage au 
quotient. Notons $[y_{0}:\ldots:y_{n}]_{w}$ la classe de 
$\pi_{w}(y_{0},\ldots ,y_{n})$ dans $|\PP(w)|$.  Pour tout 
$y:=[y_{0}:\ldots:y_{n}]_{w}$ dans $|\PP(w)|$, posons 
\begin{align}\label{eq:groupe,local} 
I_{y}&:=\{k\mid y_{k}\neq 0\}.  
\end{align} 
On a le diagramme commutatif suivant : 
\begin{align*} 
\xymatrix{(z_{0}, \ldots ,z_{n}) \ar @{|->}[d] & \CC^{n+1}-\{0\} 
  \ar[d]_-{\widetilde{f}_{w}} \ar[r]^-{\pi} & \PP^{n}\ar[d]^-{f_{w}} & 
  [z_{0}: \ldots :z_{n}]\ar @{|->}[d] \\ (z_{0}^{w_{0}}, \ldots 
  ,z_{n}^{w_{n}})& \CC^{n+1}-\{0\} \ar[r]^-{\pi_{w}} & |\PP(w)|& 
  [z_{0}^{w_{0}}: \ldots :z_{n}^{w_{n}}]_{w}} 
\end{align*} 
où $\pi$ est l'application standard de passage au quotient pour les 
espaces projectifs complexes.  Notons $\bs{\mu}_{k}$ le groupe des 
racines $k$-ièmes de l'unité.  L'espace projectif à poids peut \^{e}tre 
muni de deux structures orbifoldes différentes. 
\begin{enumerate}  
\item[(i)]  
   Le groupe $\bs{\mu}_{w_{0}}\times\cdots\times\bs{\mu}_{w_{n}}$ agit 
  sur $\PP^n$ de la façon suivante : 
  \begin{align*} 
\bs{\mu}_{w_{0}}\times\cdots\times\bs{\mu}_{w_{n}}\times \PP^n & \longrightarrow \PP^n\\ 
\left((\lambda_{0},\ldots,\lambda_{n}),[z_{0}:\ldots:z_{n}]\right)&\longmapsto 
[\lambda_{0}z_{0}:\ldots:\lambda_{n}z_{n}]\end{align*} 
 
L'application $f_{w}$ induit un homéomorphisme entre 
$\PP^{n}/\bs{\mu}_{w_{0}}\times\cdots\times\bs{\mu}_{w_{n}}$ et $|\PP(w)|$. Ainsi, l'espace topologique 
$|\PP(w)|$ est muni d'une structure orbifolde, dite globale. 
\item[(ii)] L'espace topologique $|\PP(w)|$ peut aussi \^etre muni 
  d'une structure orbifolde via l'application $\pi_{w}$. L'atlas orbifold qui 
  définit cette structure est décrit ci-dessous. 
 \end{enumerate} 
 Dans la suite de la thèse, nous nous intéressons uniquement à la 
 structure orbifolde provenant de (ii).  Pour $i \in \{0,\ldots,n\}$, 
 notons $U_{i}:=\{[y_{0}:\ldots:y_{n}]_{w} \mid y_{i}\neq 0\}\subset 
 |\PP(w)|$.  Soit $\widetilde{U}_{i}$ l'ensemble des points de 
 $\CC^{n+1}-\{0\}$ tels que $y_{i}=1$.  Le sous-groupe de $\CC^\star$ 
 qui stabilise $\widetilde{U}_{i}$ est $\bs{\mu}_{w_{i}}$. 
 L'application $\pi_{i}:=\pi_{w}\mid_{\widetilde{U}_{i}}: 
 \widetilde{U}_{i} \longrightarrow U_{i}$ induit un homéomorphisme 
 entre $\widetilde{U}_{i}/\bs{\mu}_{w_{i}}$ et $U_{i}$. 
  
 Soit $U$ un ouvert connexe de $|\PP(w)|$. Une carte 
 $(\widetilde{U},G_{\widetilde{U}},\pi_{\widetilde{U}})$ de $U$ est 
 dite admissible s'il existe $i\in\{0, \ldots ,n\}$ tel que 
 
\begin{enumerate} \makeatletter 
  \renewcommand\theenumi{\theequation} 
  \makeatother \addtocounter{equation}{1} 
\item \label{cond:a} 
  $\widetilde{U}$ soit une composante connexe de $\pi_{i}^{-1}(U)$\,;
  \addtocounter{equation}{1}
\item \label{cond:b} $G_{\widetilde{U}}$ 
  soit le sous-groupe de $\bs{\mu}_{w_{i}}$ qui stabilise 
  $\widetilde{U}$\,; \addtocounter{equation}{1}
\item \label{cond:c} 
  $\pi_{\widetilde{U}}= \pi_{i}\mid_{\widetilde{U}_{i}}$. 
\end{enumerate}
En particulier, les cartes 
$(\widetilde{U_{i}},\bs{\mu}_{w_{i}},\pi_{i})$ de $U_{i}$ sont des 
cartes admissibles.  Notons $\mathcal{A}(|\PP(w)|)$ l'ensemble des 
cartes admissibles. L'ensemble des cartes de $\mathcal{A}(|\PP(w)|)$ 
induit un recouvrement, noté $\mathcal{U}_{w}$, de $|\PP(w)|$. 
 
\bigskip
   
 \begin{lem}\label{lem:injection}  
   Soit $\alpha : 
   (\widetilde{U},G_{\widetilde{U}},\pi_{\widetilde{U}})\hookrightarrow 
   (\widetilde{V},G_{\widetilde{V}},\pi_{\widetilde{V}})$ une 
   injection  entre deux cartes de 
   $\mathcal{A}(|\PP(w)|)$. Il existe une unique fonction holomorphe 
   $\lambda_{\alpha}:\widetilde{U}\rightarrow \CC$ telle que 
   $\alpha(y_{0},\ldots,y_{n})=(y_{0}\lambda_{\alpha}^{w_{0}/\pgcd(w)}(y), 
   \ldots ,y_{n}\lambda_{\alpha}^{w_{n}/\pgcd(w)}(y))$.  Pour deux 
   injections successives $\alpha,\beta$, on a 
   \begin{align*}\lambda_{\beta\circ\alpha}(y)=\lambda_{\beta}(\alpha(y))\lambda_{\alpha}(y). 
   \end{align*} 
 \end{lem}

 \begin{rem} \label{rem:orbite}Soit $\widetilde{p}=(p_{0}, \ldots 
   ,p_{n})$ dans $\CC^{n+1}-\{0\}$. Notons 
   $D^{n+1}(\widetilde{p},\varepsilon)$ le produit des disques 
   $D^{1}(p_{i},\varepsilon)$ de centre $p_{i}$ et de rayon 
   $\varepsilon$. Pour tout $k\in I_{\pi(\widetilde{p})}$ et pour 
   toute détermination de l'argument $\arg_{k}(\cdot)$ dans 
   $D^{1}(p_{k},\varepsilon)$, considérons l'application 
    \begin{align*}  \psi_{\widetilde{p},\arg_{k}} :   
  D^{n+1}(\widetilde{p},\varepsilon) 
  & \longrightarrow  \widetilde{U}_{k} \\ 
  (y_{0},\ldots,y_{n}) & \longmapsto  (y_{0}/y_{k}^{w_{0}/w_{k}},\ldots,1_{k},\ldots,y_{n}/y_{k}^{w_{n}/w_{k}}) 
 \end{align*}   
  
 où 
 $y_{k}^{1/w_{k}}:=|y_{k}|^{1/w_{k}}\exp(i\arg_{k}(y_{k})/w_{k})$. 
  
 Soit $\arg_{k}'$ une autre détermination de l'argument sur 
 $D^{1}(p_{k},\varepsilon)$. Il existe un unique $\alpha\in\ZZ$ tel 
 que $\arg_{k}-\arg_{k}'=2\pi\alpha$.  Nous en déduisons que 
 \begin{align}\label{eq:changement,argument} 
 \psi_{\widetilde{p},\arg_{k}}&=e^{2i\pi\alpha/w_{k}}\cdot 
 \psi_{\widetilde{p},\arg_{k}'}.  
  \end{align}   
 \end{rem}

 \begin{proof}[Démonstration du lemme \ref{lem:injection}]  

\begin{itemize} \item    
  Il existe une fonction  
  \begin{align*} 
    \widetilde{\lambda}_{\alpha} 
  :\widetilde{U}\cap \cap_{k=0}^{n}\{y\in \CC^{n+1}\mid y_{k}\neq 
  0\}&\longrightarrow \CC^{\star} 
  \end{align*} 
telle que 
  \begin{align*} 
    \widetilde{\lambda}_{\alpha}(y)\cdot y=\alpha(y)=(\alpha_{0}(y), \ldots ,\alpha_{n}(y)). 
  \end{align*} 
Ainsi pour tout 
  $k\in\{0, \ldots ,n\}$, nous avons 
  $\widetilde{\lambda}_{\alpha}(y)^{w_{k}}y_{k}=\alpha_{k}(y)$. 
  Nous en déduisons que la fonction 
  $\widetilde{\lambda}_{\alpha}^{w_{k}}=\alpha_{k}/y_{k}$ est holomorphe sur 
  $\widetilde{U}\cap \{y_{k}\neq 0 \}$. 
   
  Montrons que $\widetilde{\lambda}^{w_{k}}_{\alpha}$ est holomorphe 
  sur $\widetilde{U}$. Il existe un unique $j$ tel que 
  $\widetilde{V}\subset \widetilde{U}_{j}$.  Pour tout 
  $\widetilde{p}\in \widetilde{U} \cap \{ y_{k}=0 \}$, la remarque 
  \ref{rem:orbite} implique qu'il existe $\arg_{j}$  tel que nous 
  ayons une application holomorphe 
  \begin{align*} 
    \psi_{\widetilde{p},\arg_{j}}\mid_{\widetilde{U}\cap 
  D^{n+1}(\widetilde{p},\varepsilon)} : \widetilde{U}\cap 
  D^{n+1}(\widetilde{p},\varepsilon) \rightarrow \widetilde{U}_{j} 
  \end{align*} 
 Nous vérifions sans peine que cette application est une injection.  
 Nous obtenons alors deux injections $\alpha\mid_{\widetilde{U}\cap 
  D^{n+1}(\widetilde{p},\varepsilon)}$ et $\psi_{\widetilde{p},\arg_{j}}\mid_{\widetilde{U}\cap 
  D^{n+1}(\widetilde{p},\varepsilon)}$ de $\widetilde{U}\cap 
  D^{n+1}(\widetilde{p},\varepsilon)$ dans $\widetilde{U}_{j}$. 
 D'après le lemme \ref{lem:injections} et l'égalité 
 (\ref{eq:changement,argument}), nous pouvons choisir une 
 détermination de l'argument $\arg'_{j}$ tel que 
  $\psi_{\widetilde{p},\arg'_{j}} = \alpha$ sur $\widetilde{U}\cap 
  D^{n+1}(\widetilde{p},\varepsilon)$. 
  Nous obtenons que  
  $\widetilde{\lambda}_{\alpha}^{w_{k}}$ est holomorphe sur 
  $\widetilde{U}$.  D'après Bézout, 
  $\lambda_{\alpha}:=\widetilde{\lambda}_{\alpha}^{\pgcd(w)}$ est 
  aussi holomorphe sur $\widetilde{U}$. L'unicité découle de Bézout. 
\item Nous avons l'égalité 
  \begin{align*} 
  \beta\circ\alpha(y)=\widetilde{\lambda}_{\beta\circ\alpha}(y)\cdot y 
  = 
  \left(\widetilde{\lambda}_{\beta}(\alpha(y))\widetilde{\lambda}_{\alpha}(y)\right)\cdot 
  y. 
  \end{align*} 
  Ainsi, il existe $\zeta\in \bs{\mu}_{\pgcd(w)}$ tel 
  que $ 
  \widetilde{\lambda}_{\beta\circ\alpha}(y)\zeta=\widetilde{\lambda}_{\beta}(\alpha(y))\widetilde{\lambda}_{\alpha}(y)$. 
  Nous élevons cette équation à la puissance $\pgcd(w)$ et nous en 
  déduisons l'égalité voulue. 
\end{itemize} 
\end{proof} 
 
\begin{notation}\label{not:injection} 
  Soit $\alpha : 
  (\widetilde{U},G_{\widetilde{U}},\pi_{\widetilde{U}})\hookrightarrow 
  (\widetilde{V},G_{\widetilde{V}},\pi_{\widetilde{V}})$ une injection 
  entre deux cartes de $\mathcal{A}(|\PP(w)|)$. Supposons que 
  $\widetilde{U}\subset \widetilde{U}_{i}$ et que 
  $\widetilde{V}\subset \widetilde{U}_{j}$. 
 Nous avons alors deux possibilités  
 \begin{enumerate} 
 \item soit $i=j$ et $\alpha$ est simplement l'action d'un élément de 
   $\bs{\mu}_{w_{i}}$ c'est-à-dire que $\alpha(y)=\zeta\cdot y$ avec $\zeta\in  
   \bs{\mu}_{w_{i}}$. Dans ce cas, nous avons 
   $\lambda_{\alpha}(y)\cdot y=\zeta^{\pgcd(w)}\cdot y$. 
 \item Soit $i\neq j$ et alors nous avons 
   $\lambda_{\alpha}^{w_{j}/\pgcd(w)}(y)=1/y_{j}$. Ainsi, nous notons 
   $1/y_{j}^{\pgcd(w)/w_{j}}:=\lambda_{\alpha}(y)$. 
 \end{enumerate} 
 \end{notation}

\begin{prop}\label{prop:orbifold,poids}  
  L'ensemble $\mathcal{A}(|\PP(w)|)$ est un atlas orbifold. 
\end{prop} 
 
Nous noterons $\PP(w)$ l'orbifold $(|\PP(w)|,\mathcal{A}(|\PP(w)|))$.

\begin{proof}[Démonstration de la proposition \ref{prop:orbifold,poids}]  
Il faut vérifier les conditions 
(\ref{enu:atlas,1}) et (\ref{enu:atlas,2}) 
p.\pageref{enu:atlas,1}. 
La première condition est trivialement vraie, il reste à 
montrer la seconde condition. 
 
  Considérons des cartes 
  $(\widetilde{U},G_{\widetilde{U}},\pi_{\widetilde{U}})$ de $U$ et 
  $(\widetilde{V},G_{\widetilde{V}},\pi_{\widetilde{V}})$ de $V$ dans 
  $\mathcal{A}(|\PP(w)|)$. Soit $p$ dans $U\cap V$. Il existe $i$ et 
  $j$ tels que $\widetilde{U}\subset \widetilde{U}_{i}$ et 
  $\widetilde{V}\subset \widetilde{U}_{j}$.

  Puisque $\pi_{i}^{-1}(p)\cap \widetilde{U}_{i}$ est fini, il existe 
  un réel $\varepsilon>0$ tel que les polydisques $D^{n+1}(\widetilde{p},\varepsilon)$ 
  avec $\widetilde{p}\in \pi_{i}^{-1} (p)$ soient disjoints (cf. 
  figure \ref{fig:2}). 
 
\begin{figure}[tbhp] 
\begin{center} 
  \psfrag{U'}{$U$} \psfrag{V}{$V$} \psfrag{W}{$W$} 
  \psfrag{a}{$\widetilde{p}^{U}_{1}$} \psfrag{d}{$\widetilde{p}^{U}$} 
  \psfrag{e}{$\widetilde{e}_{j}$} 
  \psfrag{f}{$\widetilde{e}_{j}(\widetilde{p}^{U})$} \psfrag{p}{$ p$} 
  \psfrag{A}{$|\PP(w)|$} 
  \psfrag{U}{$\widetilde{U}_{i}\subset \CC^{n+1}-\{ 0\}$} 
  \psfrag{g}{$\pi_{i}$} 
  \psfrag{D}{$D^{n+1}(\widetilde{p}^{U},\varepsilon)$} 
  \psfrag{D'}{$D^{n+1}(\widetilde{p}^{U}_{1},\varepsilon)$} 
  \psfrag{D''}{$D^{n+1}(\widetilde{p}^{U}_{\ell},\varepsilon)$} 
  \psfrag{E}{$D^{1}(\widetilde{e}_{j}(\widetilde{p}^{U}),\varepsilon)\subset 
    \CC$} \includegraphics[width=0.7\linewidth]{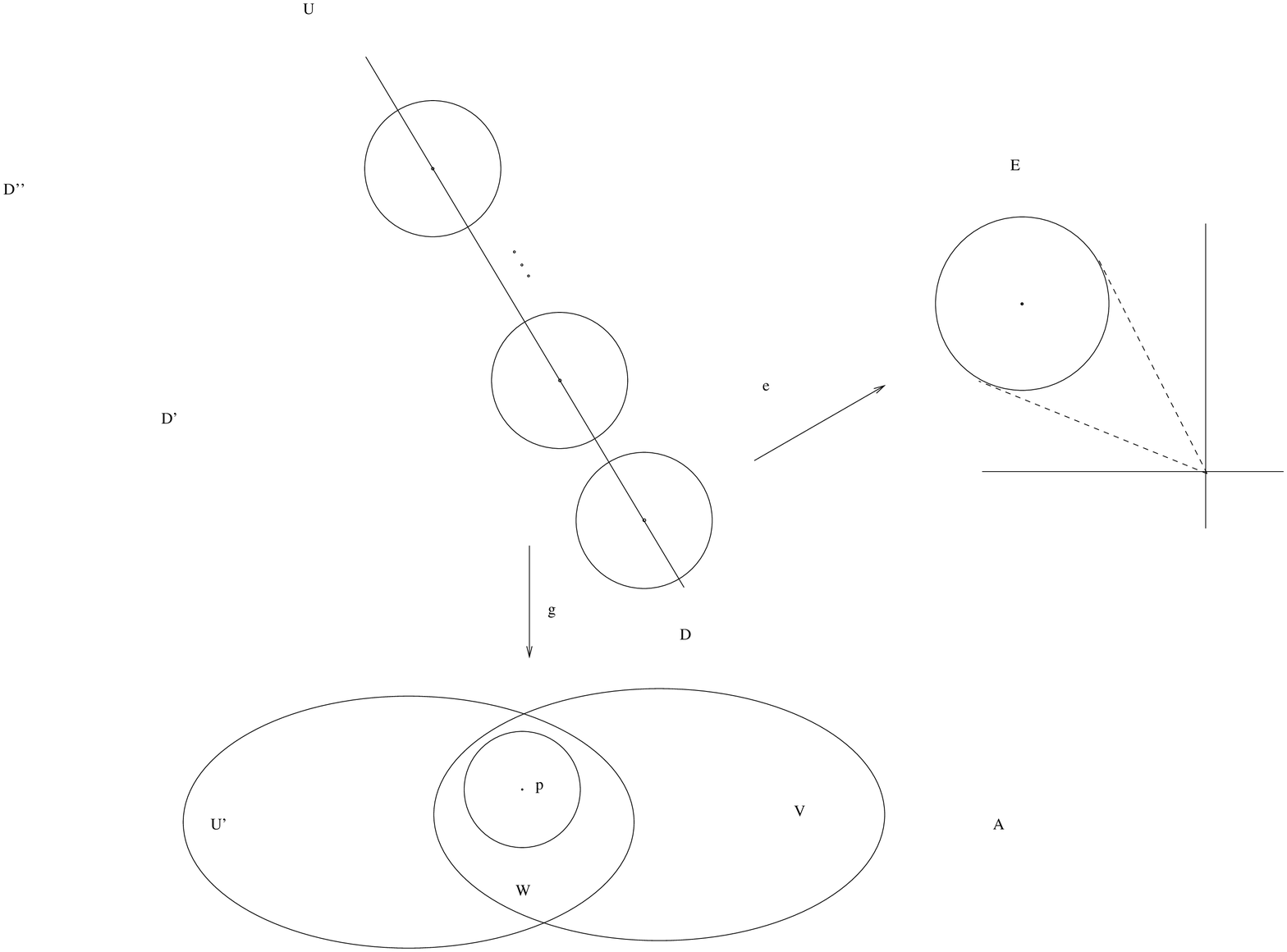} 
\end{center}\caption{$\widetilde{e}_{j} : 
  \CC^{n+1}-\{0\} \rightarrow \CC $ est la projection sur la $j$-ième 
  coordonnée.}\label{fig:2}\end{figure}    
 
Fixons $\widetilde{p}^{U}=(p^{U}_{0}, \ldots ,1_{i}, 
\ldots,p^{U}_{n})$ un relevé de $p$ dans $\widetilde{U}$.  Quitte à 
diminuer $\varepsilon$, on peut supposer que 
$D^{n+1}(\widetilde{p}^{U},\varepsilon)\cap \widetilde{U}_{i}\subset 
\widetilde{U}$ et que 
$\pi_{i}(D^{n+1}(\widetilde{p}^{U},\varepsilon)\cap 
\widetilde{U}_{i})\subset U\cap V.$ Posons 
$\widetilde{W}:=D^{n+1}(\widetilde{p}^{U},\varepsilon)\cap 
\widetilde{U}_{i}$ et $W:=\pi_{i}(\widetilde{W})$. Nous avons 
$\pi_{i}^{-1}(W)=\bigsqcup_{\widetilde{p}\in\pi_{i}^{-1}(p)} 
D^{n+1}(\widetilde{p},\varepsilon)\cap \widetilde{U}_{i}$ car 
$\bs{\mu}_{w_{i}}$ agit transitivement sur 
$\{D^{n+1}(\widetilde{p},\varepsilon)\cap \widetilde{U}_{i}\mid 
\widetilde{p}\in \pi_{i}^{-1}(p)\}$. Ainsi, $\widetilde{W}$ est une 
composante connexe de $\pi_{i}^{-1}(W)$.

Posons $G_{\widetilde{W}}:=\cap_{k\in I_{p}}\bs{\mu}_{w_{k}}$ où 
$I_{p}$ est défini par (\ref{eq:groupe,local}). Le groupe 
$G_{\widetilde{W}}$ fixe $\widetilde{p}^{U}$ et il stabilise 
$\widetilde{W}$. De plus, comme $\pi_{i}^{-1}(p)\cap 
\widetilde{W}=\{\widetilde{p}^{U}\}$, un élément qui stabilise 
$\widetilde{W}$ doit fixer $\widetilde{p}^{U}$. Finalement, 
$G_{\widetilde{W}}$ est le sous-groupe de $\bs{\mu}_{w_{i}}$ qui 
stabilise $\widetilde{W}$. 
 
Posons $\pi_{\widetilde{W}}:=\pi_{i}\mid_{\widetilde{W}}$. Ainsi, 
$(\widetilde{W},G_{\widetilde{W}},\pi_{\widetilde{W}})$ est bien une 
carte de $W$ et elle est dans $\mathcal{A}(|\PP(w)|)$.

Montrons que $(\widetilde{W},G_{\widetilde{W}},\pi_{\widetilde{W}})$ 
s'injecte dans $(\widetilde{U},G_{U},\pi_{\widetilde{U}})$. Pour cela 
il suffit de montrer que $G_{\widetilde{W}}$ est inclus dans 
$G_{\widetilde{U}}$. Supposons qu'il existe $g\in 
G_{\widetilde{W}}-G_{\widetilde{U}}$. Comme la carte 
$(\widetilde{U},G_{\widetilde{U}},\pi_{\widetilde{U}})$ est dans 
$\mathcal{A}(|\PP(w)|)$, $g$ envoie $\widetilde{U}$ sur une 
autre composante connexe de $\pi_{i}^{-1}(U)$ or ceci est impossible 
car $g$ stabilise $\widetilde{W}\subset\widetilde{U}$.

Finalement, nous avons montré que $\widetilde{W}\subset 
\widetilde{U}$, $G_{\widetilde{W}}\subset G_{\widetilde{U}}$ et 
$\pi_{\widetilde{W}}=\pi_{\widetilde{U}}\mid_{\widetilde{W}}$
c'est-à-dire
 que nous avons une injection 
$(\widetilde{W},G_{\widetilde{W}},\pi_{\widetilde{W}})\hookrightarrow 
(\widetilde{U},G_{\widetilde{U}},\pi_{\widetilde{U}})$ qui est donnée 
par l'inclusion. 
 
 Montrons qu'il existe une injection $(\widetilde{W},G_{\widetilde{W}},\pi_{\widetilde{W}})\hookrightarrow(\widetilde{V},G_{\widetilde{V}},\pi_{\widetilde{V}})$. %
  
 Comme $p^{U}_{j}$ est non nul, nous appliquons la remarque 
 \ref{rem:orbite} pour $k=j$.  Notons 
 $\alpha:=\psi_{\widetilde{p}^{U},\arg_{j}}\mid_{\widetilde{W}}$. 
 D'après l'égalité (\ref{eq:changement,argument}), on peut choisir 
 $\arg_{j}$ tel que $\alpha(\widetilde{W})\subset \widetilde{V}$.

 Montrons que $\alpha$ est injective.  Soient $\widetilde{y}$ et 
 $\widetilde{y}'$ dans $\widetilde{W}$ tels que 
 $\alpha(\widetilde{y})=\alpha(\widetilde{y}')$. Pour tout $k$, nous 
 avons 
 \begin{align*} 
 \frac{y_{k}}{y_{j}^{w_{k}/w_{j}}}=\frac{y_{k}'}{{y'}_{j}^{w_{k}/w_{j}}}. 
 \end{align*} 
 En particulier, pour $k=i$, nous obtenons 
 \begin{align*} 
 \left(\frac{y_{j}^{1/w_{j}}}{{y'}_{j}^{1/w_{j}}}\right)^{w_{i}}=1. 
 \end{align*} 
 Ceci implique que 
 $\left(\arg_{j}(y_{j})-\arg_{j}(y_{j}')\right)w_{i}/w_{j}$ est dans 
 $2\pi\ZZ$. Quitte à diminuer $\varepsilon$ (cf. figure \ref{fig:2}), 
 on peut supposer que $\arg_{j}(y_{j})=\arg_{j}(y_{j}')$.  Nous 
 obtenons $y_{j}=y'_{j}$ et nous en déduisons que $\alpha$ est 
 injective. 
  
 Il reste à montrer que $G_{\widetilde{W}}$ est un sous-groupe de 
 $G_{\widetilde{V}}$.  Pour tout $g\in G_{\widetilde{W}}$, nous 
 avons $g\cdot\alpha(\widetilde{y})=\alpha(g\cdot 
 \widetilde{y})$.  Ainsi, $G_{\widetilde{W}}$ stabilise l'ouvert 
 $\alpha(\widetilde{W})$. Nous utilisons le m\^eme raisonnement que 
 précédemment pour montrer que $G_{\widetilde{W}}\subset 
 G_{\widetilde{V}}$. 
\end{proof} 
 
\begin{rem}\label{rem:groupe,isotropie} 
 \begin{enumerate} 
 \item\label{rem:isotropie} D'après la démonstration précédente, le groupe d'isotropie  
   au point $p\in \PP(w)$ est  
   $G_{p}=\cap_{k\in I_{p}}\bs{\mu}_{w_{k}}$. 
    
 \item\label{rem:element,trivial} Nous avons 
   $\ker (\PP(w))=\cap_{i=0}^{n}\bs{\mu}_{w_{i}}$ et $\# 
   \ker (\PP(w))=\pgcd(w)$. 
    
 \item\label{rem:regulier} Nous avons  $|\PP(w)_{\reg}|=\{p\in 
   |\PP(w)|\mid \cap_{k\in 
     I_{p}}\bs{\mu}_{w_{k}}=\cap_{i=0}^{n}\bs{\mu}_{w_{i}}\}$. Par 
   contre $|\widehat{\PP(w)}_{\reg}|=\{\emptyset\}$ si $\pgcd(w)>1$.

 \item \label{rem:item:epp,P,twist} 
 Soient $w_{0},w_{1}$ deux entiers premiers entre eux. Il est 
 facile de voir que les orbifolds $\PP(w_{0},w_{1})$ et 
 $\PP^{1}_{w_{0},w_{1}}$ (cf. proposition \ref{prop:P1,twist}) sont isomorphes. 
 \item\label{rem:P(w)I} 
    
   Nous définissons l'espace topologique $|\PP(w)_{I}|:=\{ p\in 
   |\PP(w)| \mid \forall i\in I^{c}, p_{i}=0 \}$.  Montrons que 
   l'orbifold $\PP(w)$ induit une structure orbifolde sur 
   $|\PP(w)_{I}|$.  Pour tout $i\in I$, posons 
   $U_{i}\mid_{I}:=U_{i}\cap |\PP(w)_{I}|$ et 
   $\widetilde{U}_{i}\mid_{I}:=\{(y_{0}, \ldots ,y_{n})\in 
   \CC^{n+1}-\{0\}\mid y_{i}=1\mbox{ et } \forall k\in I^{c}, 
   y_{k}=0\}$.  Nous avons 
   $\widetilde{U}_{i}\mid_{I}=\pi_{i}^{-1}(U_{i}\mid_{I})$. Posons 
   $\pi_{i}\mid_{I}:=\pi_{i}\mid_{\widetilde{U}_{i}\mid_{I}}$.  Ainsi, 
   $(\widetilde{U}_{i}\mid_{I},\bs{\mu}_{w_{i}},\pi_{i}\mid_{I})$ est 
   une carte de $U_{i}\mid_{I}$.  Soit $U\mid_{I}$ un ouvert connexe 
   de $|\PP(w)_{I}|$.  Une carte 
   $(\widetilde{U}\mid_{I},G_{\widetilde{U}\mid_{I}},\pi_{\widetilde{U}\mid_{I}})$ 
   de $U\mid_{I}\subset |\PP(w)_{I}|$ est dite $I$-admissible s'il existe 
   $i\in I$ tel que :
\begin{enumerate}  \makeatletter 
  \renewcommand\theenumii{\theequation} 
  \makeatother \addtocounter{equation}{1}
\item \label{cond:I,a} 
  $\widetilde{U}\mid_{I}$ soit une composante connexe de 
  $\left(\pi_{i}\mid_{I}\right)^{-1}(U\mid_{I})$\,;
  \addtocounter{equation}{1}
\item \label{cond:I,b} $G_{\widetilde{U}}$ 
  soit le sous-groupe de $\bs{\mu}_{w_{i}}$ qui stabilise 
  $\widetilde{U}\mid_{I}$\,; \addtocounter{equation}{1}
\item 
  \label{cond:I,c} $\pi_{\widetilde{U}\mid_{I}}= 
  \left(\pi_{i}\mid_{I}\right)\mid_{\widetilde{U}_{i}\mid_{I}}$. 
 \end{enumerate} 
  
 Notons $\mathcal{A}(|\PP(w)|)\mid_{I}$ l'ensemble des cartes 
 $I$-admissibles et notons $\mathcal{U}_{w}\mid_{I}$ le recouvrement de 
 $|\PP(w)_{I}|$ induit par les cartes de 
 $\mathcal{A}(|\PP(w)|)\mid_{I}$. Le m\^eme type de raisonnement que 
 dans la démonstration de la proposition \ref{prop:orbifold,poids} 
 montre que $\mathcal{A}(|\PP(w)|)\mid_{I}$ est un atlas orbifold. 
 Notons $\PP(w)_{I}$ l'orbifold 
 $(|\PP(w_{I})|,\mathcal{A}(|\PP(w)|)\mid_{I})$. 
 \end{enumerate}  
\end{rem} 
 
Pour tout sous-ensemble $I:=\{i_{0}, \ldots ,i_{\delta}\}$ de $\{0, 
\ldots ,n\}$, considérons les applications $\CC^{\star}$-équivariantes 
: 
 \begin{align*} 
\widetilde{\iota}_{I} : \CC^{\# I}-\{0\} & \longrightarrow 
 \CC^{n+1}-\{0\}  \\ 
 (y_{i_{0}},\ldots,y_{i_{\delta}}) & \longmapsto 
 (0,\ldots,0,y_{i_{0}},0,\ldots,0,y_{i_{\delta}},0,\ldots,0) \\ 
 \widetilde{\pi}_{I} : \CC^{n+1}-\{0\} & \longrightarrow  
\CC^{\# I} \\ 
 (y_{0},\ldots,y_{n}) & \longmapsto (y_{i_{0}}, \ldots ,y_{i_{\delta}}) 
\end{align*} 
 
Notons $\iota_{I}:|\PP(w_{I})| \rightarrow |\PP(w)|$ l'application 
quotient.  
 
\begin{prop}\label{prop:bonne,application}L'application  
 \begin{align*} 
\iota_{I} : \PP(w_{I})& \longrightarrow  
\PP(w)\\ 
 [z_{0}:\ldots:z_{\delta}]_{w_{I}}& \longmapsto  [0:\ldots:0:z_{i_{0}}:0:\ldots:0:z_{i_{\delta}}:0:\ldots:0]_{w} 
\end{align*} 
est une bonne application orbifolde. 
\end{prop}

\begin{proof}Pour simplifier la démonstration, nous allons démontrer la 
  proposition pour $I=\{0, \ldots ,\delta\}$. Il faut 
  vérifier les conditions (\ref{item:7}) et (\ref{item:8}) de la 
  définition \ref{defi bonne application}. Mais, nous allons d'abord 
  expliciter les recouvrements compatibles que nous allons considérer. 
 
 La structure orbifolde sur 
  $|\PP(w_{I})|$ est donnée par l'atlas orbifold 
  $\mathcal{A}(|\PP(w_{I})|)$.  D'après la proposition 
  \ref{prop:atlas,recouvrement}, il existe un recouvrement compatible, 
  noté $\mathcal{U}_{I}$, associé à cet atlas.  Pour tout ouvert 
  $U_{I}$ dans $\mathcal{U}_{I}$ nous considérons l'ouvert de $\PP(w)$ 
  défini par $\mathfrak{F}(U_{I}):=\{[y_{0}:\ldots:y_{n}]_{w}\mid 
  [y_{0}:\ldots:y_{\delta}]_{w_{I}}\in U_{I} \}$. 
  Soit 
  $(\widetilde{U}_{I},G_{{U}_{I}},\pi_{{U}_{I}})$ 
  une carte de $U_{I}$ dans $\mathcal{U}_{I}$. Il existe $i\in I$ tel 
  que $\widetilde{U}_{I}\subset \widetilde{U}_{I,i}$.  Posons 
  $\widetilde{\mathfrak{F}(U_{I})}:=\{y\in\CC^{n+1}-\{0\}\mid \widetilde{\pi}_{I}(y)\in 
  \widetilde{U}_{I}\}=\widetilde{U}_{I}\times \CC^{n-\delta} 
  \subset\widetilde{U}_{i}$.  L'ensemble $\widetilde{\mathfrak{F}(U_{I})}$ est une 
  composante connexe de $\pi_{i}^{-1}(\mathfrak{F}(U_{I}))$. Notons 
  $G_{\mathfrak{F}(U_{I})}$ le sous-groupe de $\bs{\mu}_{w_{i}}$ qui 
  stabilise $\widetilde{\mathfrak{F}(U_{I})}$ et posons 
  $\pi_{\mathfrak{F}(U_{I})}:=\pi_{i}\mid_{\widetilde{\mathfrak{F}(U_{I})}}$.  Ainsi 
  $(\widetilde{\mathfrak{F}(U_{I})},G_{\mathfrak{F}(U_{I})},\pi_{\mathfrak{F}(U_{I})})$ est une carte de 
  $\mathfrak{F}(U_{I})$.  L'ensemble des ouverts $\mathfrak{F}(U_{I})$ 
  avec $U_{I}$ dans $\mathcal{U}_{I}$ forme un recouvrement 
  compatible de $\PP(w)$. 
   
  Soit $\mathfrak{F}$ la correspondance qui à $U_{I}$ associe 
  $\mathfrak {F}(U_{I})$. Cette correspondance vérifie clairement le 
  point $(\ref{item:7})$ de la définition \ref{defi bonne 
    application}. 
 Nous avons le diagramme commutatif suivant 
 : 
 \begin{align*} 
 \xymatrix{\widetilde{U}_{I} 
   \ar[rr]^-{\widetilde{\iota}_{I}\mid_{\widetilde{U}_{I}}} 
   \ar[d]_-{\pi_{{U}_{I}}}& & \widetilde{\mathfrak{F}(U_{I})} 
   \ar[d]^-{\pi_{\mathfrak{F}(U_{I})}} \\ U_{I} 
   \ar[rr]^-{\iota_{I}\mid_{U_{I}}}& & \mathfrak{F}(U_{I}) } 
 \end{align*} 
 Le m\^{e}me raisonnement que dans la démonstration de la 
 proposition \ref{prop:orbifold,poids} montre que 
 $G_{{U}_{I}}=G_{\mathfrak{F}(U_{I})}$. 
 Soit $\alpha : 
 (\widetilde{U}_{I},G_{{U}_{I}},\pi_{{U}_{I}}) 
 \hookrightarrow 
 (\widetilde{V}_{I},G_{{V}_{I}},\pi_{{V}_{I}})$ 
 une injection. 
 Soient $(\widetilde{\mathfrak{F}(U_{I})},G_{\mathfrak{F}(U_{I})},\pi_{\mathfrak{F}(U_{I})})$ et 
 $(\widetilde{\mathfrak{F}(V_{I})},G_{\mathfrak{F}(V_{I})},\pi_{\mathfrak{F}(V_{I})})$ les cartes de 
 $\mathfrak{F}(U_{I})$ et $\mathfrak{F}(V_{I})$ cons\-trui\-tes ci-dessus. 
 Posons $\mathfrak{F}(\alpha):=(\alpha,\id)$ c'est-à-dire qu'on applique 
 $\alpha$ aux $\delta+1$ premières coordonnées et l'identité aux 
 autres. 
  
 L'application $\mathfrak{F}(\alpha):\widetilde{\mathfrak{F}(U_{I})}\rightarrow\widetilde{\mathfrak{F}(V_{I})}$ 
 est injective. Comme $G_{{U}_{I}}$ est un sous-groupe de 
 $G_{{V}_{I}}$, $G_{\mathfrak{F}(U_{I})}$ est un sous-groupe de 
 $G_{\mathfrak{F}(V_{I})}$. Finalement, $\mathfrak{F}(\alpha)$ est une injection 
 qui satisfait la condition $(\ref{item:8})$ de la définition \ref{defi bonne application}. 
\end{proof}

\begin{cor}\label{cor:identification} 
  Pour tout sous-ensemble $I$ de $\{0,\ldots,n\}$, l'application 
  orbifolde $\iota_{I}:\PP(w_{I})\rightarrow \PP(w)$ induit un 
  isomorphisme entre $\PP(w_{I})$ et $\PP(w)_{I}$. 
\end{cor} 
  
Dans la suite, nous identifions les orbifolds 
$\PP(w_{I})\hookrightarrow \PP(w)$ et $\PP(w)_{I}\subset \PP(w)$. 
 
\begin{proof}[Démonstration du corollaire \ref{cor:identification}] 
  L'application  
  \begin{align*} 
    \iota_{I}:|\PP(w_{I})|&\longrightarrow |\PP(w)_{I}| 
  \end{align*} 
est 
  un homéomorphisme.  Pour toute carte 
  $(\widetilde{U}_{I},G_{{U}_{I}},\pi_{{U}_{I}})$ 
  de $U_{I}$ dans $\mathcal{A}(|\PP(w_{I})|)$, 
  $\widetilde{\iota}_{I}\mid_{\widetilde{U}_{I}}: 
  (\widetilde{U}_{I},G_{{U}_{I}},\pi_{{U}_{I}}) 
  \rightarrow 
  (\widetilde{\iota}_{I}(\widetilde{U}_{I}),G_{{U}_{I}},\pi\mid_{\widetilde{\iota}_{I}(\widetilde{U}_{I})})$ 
  est un isomorphisme de carte et 
  $(\widetilde{\iota}_{I}(\widetilde{U}_{I}),G_{{U}_{I}},\pi\mid_{ 
    \widetilde{\iota}_{I}(\widetilde{U}_{I})})$ est une carte de 
  $\iota_{I}(U_{I})\cap |\PP(w)_{I}|$ dans 
  $\mathcal{A}(|\PP(w)|)\mid_{I}$. 
\end{proof} 
 
Le corollaire suivant est une conséquence de la proposition \ref{prop:pullback,fibre}. 
 
\begin{cor}\label{cor:fibre,inverse} 
  Soit $I$ un sous-ensemble de $\{0,\ldots,n\}$. Pour tout fibré 
  vectoriel orbifold $E$ sur $\PP(w)$, l'image inverse de $E$, noté 
  $\iota_{I}^\ast E$, par l'application $\iota_{I}$ est un fibré 
  vectoriel orbifold. 
\end{cor}

\begin{prop}\label{prop:application,orbifold} 
  L'application 
  \begin{align*} 
 f_{w}  :  \PP^{n} & \longrightarrow  \PP(w) \\ 
  [z_{0}: \ldots:z_{n}]  &\longmapsto  [z_{0}^{w_{0}}: \ldots:z_{n}^{w_{n}}]_{w}\\ 
\end{align*} 
est une bonne application orbifolde. 
\end{prop}    
 
\begin{rem}\label{rem:degre} 
  Le degré de $f_{w}$, vue comme une application entre espaces 
  topologiques, est $\prod w_{i}/ \pgcd(w_{0}, \ldots ,w_{n})$. 
\end{rem} 
 
\begin{proof}[Démonstration de la proposition \ref{prop:application,orbifold}] 
  Rappelons que l'application  
  \begin{align*} 
    \widetilde{f}_{w}: \CC^{n+1}-\{0\} & \longrightarrow   \CC^{n+1}-\{0\} \\ 
(z_{0}, \ldots ,z_{n}) & \longmapsto (z_{0}^{w_{0}}, \ldots ,z_{n}^{w_{n}}) 
  \end{align*} est $\CC^{\star}$-équivariante et qu'elle relève 
  $|f_{w}|:|\PP^{n}|\rightarrow |\PP(w)|$.  L'application $f_{w}$ est 
  surjective et ouverte. 
   
Remarquons que si les poids sont premiers entre eux, nous avons 
   $\PP(w)_{\reg}=\widehat{\PP(w)}_{\reg}$ (cf. remarque \ref{rem:groupe,isotropie}.(\ref{rem:regulier})). 
   Ainsi, l'application $f_{w}$ est une application orbifolde régulière 
    c'est-à-dire 
   que $f_{w}^{-1}(\widehat{|\PP(w)|}_{\reg})$ est un ouvert 
   connexe et dense. Puis, la proposition \ref{prop:application,regulier} 
   montre que l'application $f_{w}$ est bonne. 
 
Dans la suite, nous démontrons la proposition dans le cas général. 
  Soit $U$ un ouvert de $\PP^{n}$. Soit $(\widetilde{U},\id,\pi_{U})$ 
  une carte de ${U}$ où $\widetilde{U}$ est un ouvert de 
  $\widetilde{U}_{i}=\{(y_{0}, \ldots ,y_{n})\in \CC^{n+1}-\{0\}\mid 
  y_{i}=1\}$. Soit $\widetilde{f_{w}(U)}\subset \widetilde{U}_{i}$ la 
  composante connexe de $\pi_{i}^{-1}(f_{w}(U))$ qui contient 
  $\widetilde{f}_{w}(\widetilde{U})$. Nous obtenons une carte 
  $(\widetilde{f_{w}(U)},G_{f_{w}(U)},\pi_{f_{w}(U)})$ de $f_{w}(U)$ 
  et nous avons le diagramme commutatif  
  \begin{align*} 
  \xymatrix{\widetilde{U} \ar[rr]^-{\widetilde{f}_{w}} 
    \ar[d]^-{\pi_{U}}& &\widetilde{f}_{w}(\widetilde{U}) 
    \ar[d]^-{\pi_{f_{w}(U)}}\\ U \ar[rr]^-{f_{w}}&& f_{w}(U)} 
\end{align*}   
  
Pour tout $p\in \PP^{n}$, il existe un ouvert connexe $U_{p}$ 
contenant $p$ tel que 
\begin{enumerate} 
\item  la carte $\widetilde{f_{w}(U_{p})}$ s'injecte 
dans $\widetilde{U}_{i_{1}}, \ldots ,\widetilde{U}_{i_{n(p)}}$ où 
$\{i_{1}, \ldots ,i_{n(p)}\}=I_{f_{w}(p)}$ (cf. définition 
(\ref{eq:groupe,local})); 
\item  la carte $\widetilde{f_{w}(U_{p})}$ soit incluse dans 
  $D^{n+1}(\widetilde{f_{w}(p)},\varepsilon)\cap \widetilde{U}_{i}$ où  
  $\varepsilon$ est choisi pour qu'on ait une détermination de 
  l'argument (cf. remarque \ref{rem:orbite}). 
\end{enumerate} 
La famille $\mathcal{U}_{\PP^{n}}:=(U_{p})_{p\in \PP^{n}}$ de tels 
ouverts est un recouvrement compatible de $\PP^{n}$. Montrons que la 
famille $\mathcal{U}_{\PP(w)}:=(f_{w}(U_{p}))_{p\in \PP^{n}}$ est un 
recouvrement compatible de $\PP(w)$. Le point (\ref{item:12}) est 
trivialement vrai.  Soit $f_{w}(y)\in f_{w}(U_{p})\cap f_{w}(U_{q})$. 
D'après la preuve de la proposition \ref{prop:orbifold,poids}, il 
existe une carte $\widetilde{U}_{f_{w}(y)}$ de $U_{f_{w}(y)}$ qui 
s'injecte dans $\widetilde{f_{w}(U_{p})}$ et dans 
$\widetilde{f_{w}(U_{q})}$. Nous prenons un ouvert $U_{y}\subset 
\PP^{n}$ tel que $f_{w}(U_{y})\subset U_{f_{w}(y)}$. Ainsi, nous avons 
l'injection $\widetilde{f_{w}(U_{y})}\hookrightarrow 
\widetilde{U}_{f_{w}(y)}$, ce qui nous montre le point (\ref{item:18}). 
Supposons que $f_{w}(U_{p})\subset f_{w}(U_{q})$. Ainsi, la carte 
$\widetilde{f_{w}(U_{p})}$ s'injecte dans $\widetilde{U}_{i}$ pour 
$i\in I_{f_{w}(p)}$ et la carte $\widetilde{f_{w}(U_{q})}$ s'injecte 
dans $\widetilde{U}_{j}$ pour $j\in I_{f_{w}(q)}$. Comme 
$f_{w}(U_{p})\subset f_{w}(U_{q})$, nous avons $I_{f_{w}(q)}\subset 
I_{f_{w}(p)}$. Soit $\widetilde{f_{w}(U_{p})}_{q}\subset 
\widetilde{f_{w}(U_{q})} $ une carte de $f_{w}(U_{p})$ induite par la 
carte $\widetilde{f_{w}(U_{q})}$. Ainsi, 
$\widetilde{f_{w}(U_{p})}_{q}$ et $\widetilde{f_{w}(U_{p})}$ 
s'injectent dans les cartes $\widetilde{U}_{j}$ pour $j\in 
I_{f_{w}(q)}$. Le lemme \ref{lem:ruan} montre que les cartes 
$\widetilde{f_{w}(U_{p})}_{q}$ et $\widetilde{f_{w}(U_{p})}$ sont 
isomorphes. Nous en concluons que $\widetilde{f_{w}(U_{p})}$ s'injecte 
dans $\widetilde{f_{w}(U_{q})}$. Ce qui montre le point 
(\ref{item:19}). 
 
Nous allons maintenant vérifier que les conditions de la définition 
\ref{defi bonne application} sont satisfaites.  La correspondance 
$\mathfrak{F}:\mathcal{U}_{\PP^{n}}\rightarrow\mathcal{U}_{\PP(w)}$ 
définie par $U_{p}\mapsto f_{w}(U_{p})$ vérifie la condition 
(\ref{item:7}) de la définition \ref{defi bonne application}. 
 
Soit $\alpha:\widetilde{U}_{p}\subset\widetilde{U}_{i}\hookrightarrow 
\widetilde{U}_{q}\subset\widetilde{U}_{j}$ une injection. Nous avons 
deux cas 
\begin{enumerate} 
\item soit $i=j$ et $\alpha$ est simplement une inclusion; 
\item soit $i\neq j$ et nous avons 
  \begin{align*} 
    \alpha(z_{0}, \ldots ,1_{i}, \ldots ,z_{n})=(z_{0}/z_{j}, \ldots ,1_{j}, \ldots ,z_{n}/z_{j}). 
  \end{align*} 
\end{enumerate} 
 
Dans le premier cas, $\mathfrak{F}(\alpha)$ est simplement l'inclusion 
de $\widetilde{f_{w}(U_{p})}$ dans $\widetilde{f_{w}(U_{q})}$. Nous 
avons le diagramme commutatif 
\begin{align*} 
\xymatrix{ \widetilde{U}_{p} \ar@{^{(}->}[rr]^-{\alpha} 
  \ar[d]^-{\widetilde{f}_{w}}&&\widetilde{U}_{q} \ar[d]^-{\widetilde{f}_{w}}\\ 
  \widetilde{f_{w}(U_{p})}\ar@{^{(}->}[rr]^-{\mathfrak{F}(\alpha)}&&\widetilde{f_{w}(U_{q})}} 
\end{align*} 
Nous avons unicité de $\mathfrak{F}(\alpha)$ car d'après le lemme 
\ref{lem:injections} une autre injection aurait la forme 
$\varphi_{g}\circ\mathfrak{F}(\alpha)$ avec $g\in G_{f_{w}(U_{q})}$ 
mais le diagramme ci-dessus ne serait commutatif que si $g\in 
\ker(\PP(w))$ c'est-à-dire $\varphi_{g}=\id$. 
 
Dans le second cas, $\mathfrak{F}(\alpha)$ est l'injection 
\begin{align*} 
\mathfrak{F}(\alpha):    \widetilde{f_{w}(U_{p})}&\hookrightarrow \widetilde{f_{w}(U_{q})} \\ 
  (y_{0}, \ldots ,1_{i}, \ldots ,y_{n})& \mapsto 
  (y_{0}/y_{j}^{w_{0}/w_{j}}, \ldots ,1_{j}, \ldots 
  ,y_{n}/y_{j}^{w_{n}/w_{j}}) 
 \end{align*} 
 Cette injection est bien définie car 
nous avons choisi $\varepsilon$ assez petit pour avoir une 
détermination de l'argument. Nous avons un diagramme 
commutatif comme le précédent et l'unicité de $\mathfrak{F}(\alpha)$ se 
déduit de la m\^{e}me manière.  
 
L'unicité montre que pour deux 
injections successives, nous avons 
$\mathfrak{F}(\alpha\circ\beta)=\mathfrak{F}(\alpha)\circ\mathfrak{F}(\beta)$. 
Ce qui termine la démonstration. 
 \end{proof}

\section{Les fibrés vectoriels orbifolds $\mathcal{O}_{\PP(w)}(k)$ sur $\PP(w)$}\label{sec:Les-fibrs-vectoriels} 
Avant de définir les fibrés vectoriels orbifolds 
$\mathcal{O}_{\PP(w)}(\pgcd(w))$, nous allons faire une remarque sur 
les notations. Remarquons que la notation $\mathcal{O}_{\PP^{n}}(k)$ 
est ambigu\"{e} car elle désigne à la fois un fibré vectoriel 
holomorphe sur $\PP^{n}$ et son faisceau des sections. Cette 
confusion n'est pas gênante car nous avons une équivalence de 
catégories entre les faisceaux de $\mathcal{O}_{M}$-modules localement 
libres de rang $r$ et les fibrés vectoriels holomorphes de rang $r$ 
sur une variété $M$.  Pour les orbifolds, nous n'avons pas cette 
équivalence de catégories : le fibré vectoriel orbifold contient plus 
d'informations que son faisceau des sections (cf. paragraphe 
\ref{subsec:Le-cas-gnral}). 
 
\begin{prop}\label{prop:fibre,O(1)} Il existe un fibré vectoriel 
  complexe orbifold de rang $1$, noté 
  $\mathcal{O}_{\PP(w)}(\pgcd(w))$, sur $\PP(w)$ tel que 
  $f_{w}^{\ast}\mathcal{O}_{\PP(w)}(\pgcd(w))$ soit isomorphe au fibré 
  $\mathcal{O}_{\PP^{n}}(\pgcd(w))$ sur $\PP^{n}$. 
\end{prop}

\begin{proof}
Pour alléger les notations, notons $d(w):=\pgcd(w)$. 
   Nous utilisons les résultats du paragraphe 
   \ref{sec:bonne-appl-orbif}. 
    
   Nous allons définir le fibré vectoriel orbifold 
   $\mathcal{O}_{\PP(w)}(d(w))$ sur $\PP(w)$ par ses fonctions de 
   transition. D'après le lemme \ref{lem:injection} et la notation 
   \ref{not:injection}, pour toute injection $\alpha 
   :(\widetilde{U},G_{\widetilde{U}},\pi_{\widetilde{U}})\hookrightarrow 
   (\widetilde{V},G_{\widetilde{V}},\pi_{\widetilde{V}})$ entre deux 
   cartes de $\mathcal{A}(|\PP(w)|)$, notons 
   $\psi^{\mathcal{O}_{_{\PP(w)}}(d(w))}_{\alpha}(y)$ l'application de 
   $\CC$ dans $\CC$ qui à $t$ associe $\lambda_{\alpha}(y) t$.  Plus 
   précisément (cf. notation \ref{not:injection}), supposons que $\widetilde{U}\subset \widetilde{U}_{i}$ 
   et $\widetilde{V}\subset \widetilde{U}_{j}$ alors nous avons 
   \begin{enumerate} 
   \item soit $i=j$ et 
     $\psi^{\mathcal{O}_{_{\PP(w)}}(d(w))}_{\alpha}(y)(t)=\zeta^{d(w)}t$  
     où $\zeta\in\bs{\mu}_{w_{i}}$ ; 
   \item soit $i\neq j$ et $\psi^{\mathcal{O}_{_{\PP(w)}}(d(w))}_{\alpha}(y)(t)=t/y_{j}^{d(w)/w_{j}}$. 
   \end{enumerate} 
D'après le lemme 
   \ref{lem:injection}, nous avons la relation de cocycle orbifolde 
   suivante 
    
   \begin{align*}\psi^{\mathcal{O}_{{\PP(w)}}(d(w))}_{\beta\circ\alpha}(y)(t) = 
   \psi^{\mathcal{O}_{_{\PP(w)}}(d(w))}_{\beta}\left(\alpha(y)\right)\left(\psi^{\mathcal{O}_{_{\PP(w)}}(d(w))}_{\alpha}(t)\right).\end{align*} 
    
   Nous avons ainsi défini un fibré vectoriel orbifold de rang $1$, 
   noté $\mathcal{O}_{\PP(w)}(d(w))$, sur $\PP(w)$ et $\ker(\mathcal{O}_{\PP(w)}(d(w)))=\ker(\PP(w))$.

   D'après la proposition \ref{prop:application,orbifold}, 
   l'application $f_{w}$ est une bonne application orbifolde. D'après la proposition \ref{prop:pullback,fibre} 
    le fibré $f_{w}^{\ast}(\mathcal{O}_{\PP(w)}(d(w)))$ existe. 
    
   Montrons que les fonctions de transition de 
   $f_{w}^{\ast}(\mathcal{O}_{\PP(w)}(d(w)))$ sont les m\^emes que les 
   fonctions de transition de $\mathcal{O}_{\PP^{n}}(d(w))$.  Soit 
   $p\in U_{i}^{\PP^{n}}\cap U_{j}^{\PP^{n}}$ avec $i\neq j$. Ainsi, $f_{w}(p)$ est 
   dans $U_{i}^{\PP(w)}\cap U_{j}^{\PP(w)}$. Soit 
   $\left(\widetilde{W}_{f_{w}(p)},G_{\widetilde{W}_{f_{w}(p)}},\pi_{\widetilde{W}_{f_{w}(p)}}\right)$ 
   une carte de $W_{f_{w}(p)}$ de $\mathcal{A}(|\PP(w)|)$ telle que 
\begin{itemize} \item $f_{w}(p)\in W_{f(p)}\subset U_{i}^{\PP(w)}\cap 
  U_{j}^{\PP(w)}$ \,;
\item il existe une injection 
  \begin{align*}\alpha:\left(\widetilde{W}_{f_{w}(p)},G_{\widetilde{W}_{f_{w}(p)}},\pi_{\widetilde{W}_{f_{w}(p)}}\right) 
  \hookrightarrow 
  (\widetilde{U}_{j}^{\PP(w)},\bs{\mu}_{w_{j}},\pi_{j})\end{align*} 
   \end{itemize}  
    
   Soit $U_{p}^{\PP^{n}}$ la composante connexe de 
   $f_{w}^{-1}\left(W_{f_{w}(p)}\right)$ qui contient $p$.  Posons 
   \begin{align*} 
     \widetilde{U}_{p}^{\PP^{n}}&:=\widetilde{f}_{w}^{-1}\left(\widetilde{W}_{f_{w}(p)}\right)\cap  \widetilde{U}_{i}^{\PP^{n}}. 
   \end{align*} 
Ainsi, 
   $(\widetilde{U}_{p}^{\PP^{n}},\id,\pi\mid_{\widetilde{U}_{p}^{\PP^{n}}})$ 
   est une carte de $U_{p}^{\PP^{n}}$. 
    
   Nous avons le diagramme commutatif 
   \begin{align*} 
   \xymatrix{ \widetilde{U}_{p}^{\PP^{n}} 
     \ar[rr]^-{\widetilde{f}_{w}\mid_{\widetilde{U}_{p}^{\PP^{n}}}} 
     \ar[d]_-{\beta} && \widetilde{W}_{f_{w}(p)} 
     \ar[d]^-{\alpha} \\ 
     \widetilde{U}_{j}^{\PP^{n}} 
     \ar[rr]^-{\widetilde{f}_{w}\mid_{\widetilde{U}_{j}^{\PP^{n}}}}&& 
     \widetilde{U}_{j}^{\PP(w)} } 
   \end{align*} 
   où $\beta((z_{0}, \ldots ,1_{i}, \ldots ,z_{n}))=(z_{0}/z_{j}, 
   \ldots ,1_{j}, \ldots ,z_{n}/z_{j})$. 
      Ainsi, nous avons 
   \begin{align*}\psi^{f_{w}^{\ast}\mathcal{O}_{\PP(w)}(d(w))}_{\beta}(z)(t)=\psi^{\mathcal{O}_{\PP(w)}(d(w))}_{\alpha}(f_{w}(z))(t)=t/z^{d(w)}_{j}.\end{align*} 
\end{proof} 
  
\begin{rem} \label{rem:fibre}  
\begin{enumerate} \item \label{rem:fibre,O(m)}  
  Soit $m$ un entier. Nous définissons le fibré orbifold de rang $1$, 
  noté $\mathcal{O}_{\PP(w)}(m.\pgcd(w))$, en prenant $m$ fois le 
  produit tensoriel du fibré $\mathcal{O}_{\PP(w)}(\pgcd(w))$ avec 
  lui-m\^{e}me.  Avec les notations ci-dessus, les fonctions de 
  transition du fibré $\mathcal{O}_{\PP(w)}(m.\pgcd(w))$ sont 
\begin{align*} \g_{\alpha}^{\mathcal{O}_{\PP(w)}(m.\pgcd(w))}(y) 
    :  \CC 
    & \longrightarrow \CC \\ 
     t &\longmapsto t\lambda_{\alpha}^{m}(y) 
\end{align*} 
 
\item A l'orbifold $\PP(w)$, on peut associer l'orbifold, dite 
  réduite, $\PP(w/\pgcd(w))$ en quotientant par le groupe 
  $\ker (\PP(w))$.  Comme nous avons 
  \begin{align*} 
    \ker (\PP(w))=\ker (\mathcal{O}_{\PP(w)}(\pgcd(w))), 
  \end{align*} le fibré 
  orbifold $\mathcal{O}_{\PP(w)}(\pgcd(w))\rightarrow \PP(w)$ peut 
  aussi \^{e}tre réduit en un fibré vectoriel orbifold 
  $\mathcal{O}_{\PP(w/\pgcd(w))}(1) \rightarrow\PP(w/\pgcd(w))$.  
   
  On peut définir un fibré vectoriel orbifold, au sens de Chen et Ruan  
  (cf. remarque \ref{rem:fibre,chen,ruan}), 
  $\mathcal{O}_{\PP(w)}(1)$ sur $\PP(w)$. Mais ce n'est pas un fibré 
  vectoriel orbifold au sens de la définition 
  \ref{defi:fibre,orbifold} car  
  \begin{align*} 
    \ker (\mathcal{O}_{\PP(w)}(1))=\{ \id \} \neq \ker (\PP(w)). 
  \end{align*} 
   
\item D'après le corollaire \ref{cor:fibre,inverse}, pour tout 
  sous-ensemble $I$ de $\{0, \ldots ,n\}$, le fibré orbifold 
  $\iota_{I}^{\ast}\mathcal{O}_{\PP(w)}(\pgcd(w_{I}))$ est isomorphe 
  au fibré vectoriel orbifold 
  $\mathcal{O}_{_{\PP(w_{I})}}(\pgcd(w_{I}))$. 
\end{enumerate} 
\end{rem} 
 
\begin{lem}\label{lem:section,O(k)} 
  Pour toute carte $(\widetilde{U},G_{U},\pi_{U})$ de $U$ où 
  $\widetilde{U}\subset \widetilde{U}_{i}$, nous posons 
  \begin{align*} 
    \widetilde{s}_{k,\widetilde{U}}:  \widetilde{U} & \longrightarrow   
    \CC \\  
(y_{0}, \ldots ,1_{i}, \ldots ,y_{n})  &\longmapsto y_{k} 
  \end{align*} 
Ces sections locales sont compatibles avec les 
  changements de cartes  et elles défi\-nis\-sent une section, notée $s_{k}$, du fibré $\mathcal{O}_{\PP(w)}(w_{k})\rightarrow \PP(w)$. 
\end{lem} 
 
\begin{proof} D'après la remarque \ref{rem:section,global}, il suffit 
  de montrer que ces sections locales sont compatibles avec les 
  changements de cartes c'est-à-dire que 
  \begin{align}\label{eq:2} 
    \widetilde{s}_{k,\widetilde{V}}(\alpha(y))&=\g_{\alpha}^{\mathcal{O}_{\PP(w)}(w_{k})}(y)(\widetilde{s}_{k,\widetilde{U}}(y))=\lambda_{\alpha}^{w_{k}/\pgcd(w)}(y)\widetilde{s}_{k,\widetilde{U}}(y) 
  \end{align} 
  pour toute injection $\alpha:\widetilde{U}\hookrightarrow 
  \widetilde{V}$.  Nous utilisons les notations de 
  ~\ref{not:injection}, supposons que $\widetilde{U}\subset 
  \widetilde{U}_{i}$ et $\widetilde{V}\subset \widetilde{U}_{j}$, nous  
  avons  
  \begin{itemize} 
  \item soit $i=j$ et $\alpha(y)=\zeta\cdot y$ où $\zeta\in 
    \bs{\mu}_{w_{i}}$ et nous avons  
    \begin{align*} 
      \widetilde{s}_{k,\widetilde{V}}(\zeta\cdot y)=\zeta^{w_{k}}y_{k}=\lambda_{\alpha}^{w_{k}/d(w)}(y)\widetilde{s}_{k,\widetilde{U}}(y). 
    \end{align*} 
  \item Soit $i\neq j$ et nous avons 
    \begin{align*} 
       \widetilde{s}_{k,\widetilde{V}}(y_{0}/y_{j}^{w_{0}/w_{j}}, 
       \ldots ,1_{j}, \ldots 
       ,y_{n}/y_{j}^{w_{n}/w_{j}})=y_{k}/y_{j}^{w_{k}/w_{j}}=\lambda_{\alpha}^{w_{k}/d(w)}(y) \widetilde{s}_{k,\widetilde{U}}(y). 
    \end{align*} 
  \end{itemize} 
\end{proof} 
  
\section[La cohomologie orbifolde de 
$\PP(w)$ vue comme  espace vectoriel]{La cohomologie orbifolde de 
$\PP(w)$ vue comme  $\CC$-espace vectoriel gradué} \label{sec:coho-espace} 
 
Dans ce paragraphe, nous allons d'abord donner la définition de la 
structure de $\CC$-espace vectoriel gradué de la 
cohomologie orbifolde pour les orbifolds complexes et 
commutatives. Nous nous appuierons sur l'article de Chen et Ruan \cite{CRnco}. 
Puis, nous appliquerons cette définition à notre exemple favori 
$\PP(w)$. En particulier, nous expliciterons une base (cf. proposition  
~\ref{prop:base}) de la cohomologie orbifolde des espaces projectifs à  
poids. 
 
\subsection{Définition de la cohomologie orbifolde}\label{subsec:Dfinition-de-la} 
 
Nous allons rappeler la définition de la cohomologie orbifolde pour 
une orbifold complexe et commutative. Pour plus de précisions, nous 
renvoyons le lecteur à l'article de Chen et Ruan \cite{CRnco} 
paragraphe $3.2$. 
 
Soit $X$ une orbifold  complexe, compacte et commutative de dimension $n$. 
Notons $\Sigma X:=\bigsqcup_{x\in X} G_{x}$.  
Nous définissons une topologie sur $\Sigma X$ à l'aide d'une base 
d'ouverts. 
Pour tout $(x,g)\in \Sigma X$ et 
pour toute carte $(\widetilde{U}_{x},G_{x},\pi_{x})$ d'un ouvert 
$U_{x}$ contenant $x$, nous posons  
\begin{align}\label{eq:ouvert} 
\mathcal{V}_{(x,g)}(\widetilde{U}_{x})&:=\left\{ 
  \begin{array}{l} 
(y,h)\in \Sigma X\mid \exists 
    (\alpha,\kappa):(\widetilde{U}_{y},G_{y},\pi_{y})\hookrightarrow 
    (\widetilde{U}_{x},G_{x},\pi_{x}) \\\mbox{ qui vérifie } \kappa(h)=g 
  \end{array} 
\right\}. 
\end{align}

\begin{lem}\label{lem:topologie,separe} 
  \begin{enumerate} 
  \item La collection de ces sous-ensembles de $\Sigma X$ définit une 
topologie sur $\Sigma X$.  
\item L'ensemble $\Sigma X$ muni de cette topologie est séparé. 
  \end{enumerate} 
\end{lem} 
 
\begin{proof} 
  \begin{enumerate} 
  \item  
  Il est clair que $(x,g)$ appartient à 
  $\mathcal{V}_{(x,g)}(\widetilde{U}_{x})$. 
 
Il reste à montrer que si $(z,k)\in 
  \mathcal{V}_{(x,g)}(\widetilde{U}_{x})\cap \mathcal{V}_{(y,h)}(\widetilde{U}_{y})$ 
  alors il existe une carte $\widetilde{W}_{z}$ d'un ouvert contenant 
  $z$ telle que  
  \begin{align*} 
    \mathcal{V}_{(z,k)}(\widetilde{W}_{z})\subset 
  \mathcal{V}_{(x,g)}(\widetilde{U}_{x})\cap 
  \mathcal{V}_{(y,h)}(\widetilde{U}_{y}). 
  \end{align*} 
Soit $(z,k)\in 
  \mathcal{V}_{(x,g)}(\widetilde{U}_{x})\cap 
  \mathcal{V}_{(y,h)}(\widetilde{U}_{y})$. Il existe deux injections  
  \begin{align*} 
    (\alpha,\kappa_{\alpha}): \widetilde{U}_{z} & \hookrightarrow 
    \widetilde{U}_{x} \,;\\ 
 (\alpha',\kappa_{\alpha'}): \widetilde{U}'_{z} & \hookrightarrow 
    \widetilde{U}_{y} 
  \end{align*} 
 telles que $\kappa_{\alpha}(k)=g$ et $\kappa_{\alpha'}(k)=h$. 
Il existe une carte $(\widetilde{W}_{z},G_{z},\pi_{z})$ d'un ouvert 
$W_{z}\subset U_{z}\cap U'_{z}$ qui s'injecte dans respectivement 
$ \widetilde{U}_{z}$ et $  \widetilde{U}'_{z}$. Nous en concluons que 
$\mathcal{V}_{(z,k)}(\widetilde{W}_{z})\subset 
  \mathcal{V}_{(x,g)}(\widetilde{U}_{x})\cap 
  \mathcal{V}_{(y,h)}(\widetilde{U}_{y})$. 
 
\item Soient $(x_{1},g_{1})$ et $(x_{2},g_{2})$ deux éléments différents dans $\Sigma 
  X$. 
Si $x_{1}\neq x_{2}$ alors il existe deux cartes $\widetilde{U}_{x_{1}}$ et 
$\widetilde{U}_{x_{2}}$ de  $U_{x_{1}}$  et $U_{x_{2}}$ contenant respectivement 
$x_{1}$ et $x_{2}$ telles que 
$U_{x_{1}}\cap U_{x_{2}}=\{\emptyset\}$. Nous en déduisons que $\mathcal{V}_{(x_{1},g_{1})}(\widetilde{U}_{x_{1}})\cap 
  \mathcal{V}_{(x_{2},g_{2})}(\widetilde{U}_{x_{2}})=\{\emptyset \}$. 
Supposons que $x:=x_{1}=x_{2}$ et $g_{1}\neq g_{2}$. 
Soit $\widetilde{U}_{x}$ une carte d'un ouvert contenant $x$. 
Supposons qu'il existe $(y,h)\in \mathcal{V}_{(x,g_{1})}(\widetilde{U}_{x})\cap 
  \mathcal{V}_{(x,g_{2})}(\widetilde{U}_{x})$. 
Il existe deux injections  
 \begin{align*} 
    (\alpha,\kappa_{\alpha}): (\widetilde{U}_{y},G_{y},\pi_{y}) & \hookrightarrow 
    (\widetilde{U}_{x},G_{x},\pi_{x}) \,;\\ 
 (\alpha',\kappa_{\alpha'}): (\widetilde{U}'_{y},G_{y},\pi_{y}) & \hookrightarrow 
    (\widetilde{U}_{x},G_{x},\pi_{x}) 
  \end{align*} 
telles que $\kappa_{\alpha}(h)=g_{1}$ et $\kappa_{\alpha'}(h)=g_{2}$. 
Quitte à diminuer $U'_{y}$, on peut supposer que $U'_{y}\subset 
U_{y}$. 
Puis le corollaire \ref{cor:injections} implique qu'il existe une 
injection de $\widetilde{U}'_{y}$ dans   $\widetilde{U}_{y}$. Le 
morphisme de groupes $\kappa:G_{y}\rightarrow G_{y}$ associé à cette 
injection est l'identité. Nous en concluons que 
$\kappa_{\alpha}(h)=\kappa_{\alpha'}(h)$. Ce qui contredit l'hypothèse  
$g_{1}\neq g_{2}$. 
  \end{enumerate} 
\end{proof} 
 
L'application $P:\Sigma X \rightarrow |X|$ qui à $(x,g)$ associe $x$ est une application 
continue et $\#P^{-1}(x)=\# G_{x}$.   
 
\begin{lem}\label{lem:composantes,connexes} Soit $X$ une orbifold complexe, compacte et commutative. 
 L'espace topologique $\Sigma X$ n'a qu'un nombre fini de composantes connexes. 
\end{lem} 
 
\begin{proof} 
  Nous allons montrer que $P:\Sigma X \rightarrow |X|$ est une 
  application propre.  Il faut montrer que, pour tout $x$ dans $X$,
  $P^{-1}(x)$ est quasi-compact, c'est-\`a-dire que de tout
  recouvrement ouvert on peut en extraire un sous-recouvrement
  fini\footnote{Rappelons qu'en France un espace topologique est compact s'il
    est quasi-compact et s\'epar\'e.},
  et que l'application $P$ est fermée.  Le premier point est évident 
  car $P^{-1}(x)$ est un ensemble fini. 
   
  Montrons que l'application $P$ est fermée.  Soit $F$ un fermé de 
  $\Sigma X$. Soit $x\in \overline{P(F)}$.  Pour toute carte   
  $(\widetilde{U}_{x},G_{x},\pi_{x})$  de $U_{x}$, il existe $y\in 
  P(F)\cap U_{x}$. Nous en déduisons qu'il existe $h\in G_{y}$ tel que 
  $(y,h)\in F\subset \Sigma X$ et qu'il existe une injection 
\begin{align*} 
  (\alpha,\kappa):(\widetilde{U}_{y},G_{y},\pi_{y})\hookrightarrow (\widetilde{U}_{x},G_{x},\pi_{x}) 
\end{align*} 
Considérons $(U_{x,n})_{n\in \NN}$ une suite d'ouverts emboîtés 
contenant $x$ telle que  
\begin{align*} 
  \cap_{n\in\NN}U_{x,n}=\{x\}. 
\end{align*} 
Nous en 
déduisons une suite $(y_{n},h_{n})_{n\in\NN}$ d'éléments de $F$ telle 
que $y_{n}$ converge vers $x$. Comme $G_{x}$ est fini, quitte à 
prendre une sous-suite, nous pouvons supposer que la suite 
$(\kappa_{n}(h_{n}))_{n\in \NN}$ dans $G_{x}$ est constamment égale à 
un élément $\kappa(h)\in G_{x}$. Nous en déduisons que la suite 
$(y_{n},h_{n})_{n\in\NN}$ converge vers $(x,\kappa(h))$. Comme $F$ est 
fermé, l'élément $(x,\kappa(h))$ appartient à $F$.  Nous en concluons 
que $P$ est une application fermée puis qu'elle est propre. 
 
 Comme $\Sigma X$ est un espace topologique séparé (cf. lemme 
 \ref{lem:topologie,separe}) et que $X$ est 
 compacte, nous en déduisons (cf. Proposition $7$ du chapitre $I$  
 paragraphe $10.3$ de \cite{Btg}) que pour tout  
 compact $K$ de $|X|$, l'ensemble $P^{-1}(K)$ est compact. 
Finalement, l'espace topologique $\Sigma X$ est compact et n'a qu'un nombre fini de 
composantes connexes. 
\end{proof} 
 
Nous dirons que $(x,g)\sim (y,h)$ si $(x,g)$ et $(y,h)$ sont dans la 
m\^{e}me composante connexe de $\Sigma X$. Soit $T$ l'ensemble des 
classes d'équivalence. Pour tout $g\in G_{x}$, nous notons $(g)$ la 
classe d'équivalence de $(x,g)$. Cette notation $(g)$ sous-entend que 
$g\in G_{x}$ pour un certain $x\in|X|$.  Notons $X_{(g)}$ la composante 
connexe de $\Sigma X$ qui contient $(x,g)$.  Nous avons la 
décomposition suivante de $\Sigma X$ 
\begin{align} 
  \label{eq:decomposition} 
  \Sigma X &= \bigsqcup_{(g)\in T} X_{(g)}.  
\end{align}  
 
Les espaces topologiques $X_{(g)}$ sont naturellement munis d'une 
structure orbifolde (cf. lemme $3.1.1$ de \cite{CRnco} ou l'article \cite{Ksto}). 
Dans le cas des espaces projectifs à poids, cette 
structure orbifolde sera claire. 
 
Remarquons que pour $g\in \ker (X)$, nous avons $|X_{(g)}|=|X|$. 
 
Dans l'article \cite{Ksto}, Kawasaki donne une stratification 
canonique d'une orbifold.  
Heuristiquement, nous pouvons voir les composantes connexes de $\Sigma X$ comme des 
adhérences de strates de l'orbifold $X$ qu'on a \og sorties \fg de 
$|X|$. Regardons sur l'exemple suivant ce qu'il se passe. 
 
\begin{expl} 
  Considérons l'orbifold $\PP^{1}_{w_{0},w_{1}}$. L'ensemble $T$ est 
  en bijection avec $\bs{\mu}_{w_{0}}\sqcup\bs{\mu}_{w_{1}}$. 
L'ensemble $\Sigma  \PP^{1}_{w_{0},w_{1}}$ se décompose de la façon 
suivante 
\begin{align*} 
  \Sigma  \PP^{1}_{w_{0},w_{1}}= \PP^{1}\times\{\id\} \bigsqcup 
  \{[1:0]\}\times \left(\bs{\mu}_{w_{0}}-\{\id\}\right)\bigsqcup \{[0:1]\}\times\left(\bs{\mu_{w_{1}}}-\{\id\}\right). 
\end{align*} 
\end{expl} 
 
Soit $x$ dans $|X|$. Soit $(\widetilde{U}_{x},G_{x},\pi_{x})$ une 
carte d'un voisinage $U_{x}$ de $x$. L'action du groupe d'isotropie 
$G_{x}$ sur l'espace vectoriel tangent 
$T_{\widetilde{x}}\widetilde{U}_{x}$ induit une représentation 
$\rho_{x}:G_{x}\rightarrow GL(n,\CC)$ qui ne dépend pas de la carte 
choisie. Comme $g$ est d'ordre fini, la matrice $\rho_{x}(g)$ est 
diagonalisable. Dans une telle base, la matrice $\rho_{x}(g)$ s'écrit 
\begin{align*} 
  \diag (\exp(2i\pi r_{1}), \ldots ,\exp(2i\pi r_{n})) \end{align*} où
  les $r_{i}$ sont dans l'intervalle $[0,1[\cap \QQ$.

 \emph{L'\^{a}ge} de $g\in G_{x}$ est défini $\age(g,x):=r_{1}+\cdots+r_{n}$.   
\bigskip  
 \begin{lem}[cf. lemme $3.2.1$ de \cite{CRnco}] 
   \begin{enumerate} 
   \item L'application $\age:X_{(g)}\rightarrow \QQ$ est constante. 
   \item Nous avons $\age(g,x)+\age(g^{-1},x)=n-\dim X_{(g)}$. 
  \end{enumerate} 
\end{lem} 
D'après le lemme, l'\^{a}ge d'un élément $g\in G_{x}$ ne dépend que de la 
composante connexe $X_{(g)}$.  Dorénavant, nous le notons $\age(g)$.

\begin{defi}[cf. définition $3.2.3$ de \cite{CRnco}]\label{defi:cohomologie,orbifold} Soit $X$ une orbifold complexe et commutative. 
  Nous définissons les espaces vectoriels complexes de cohomologie 
  orbifolde de $X$ par 
\begin{align*} 
  H^{\star}_{\orb}(X,\CC):=\bigoplus_{(g)\in T} H^{\star-2 \age(g)}(X_{(g)},\CC). 
\end{align*} 
\end{defi}  

\begin{rem}
  Remarquons que $|X_{(\id)}|$ est simplement $|X|$. Ainsi, la cohomologie ordinaire de
  $|X|$ est un sous-espace vectoriel de la cohomologie orbifolde de $X$.
\end{rem}

\subsection{Cohomologie orbifolde des espaces projectifs à poids}\label{subsec:Cohom-orbif-des} 
 
Pour l'espace projectif à poids, l'ensemble  $T$ est en bijection avec  
$\bigcup_{i=0}^{n}\bs{\mu}_{w_{i}}$ (ou avec $S_{w}$ défini dans le
paragraphe \ref{sec:sigma}). 
Pour $\gamma\in S_{w}$, nous avons $\PP(w)_{(e^{2i\pi\gamma})}=\{p\in |\PP(w)|\mid e^{2i\pi\gamma}\in 
G_{p}\}\times\{e^{2i\pi\gamma}\}$ (nous utilisons la notation définie juste
 avant l'égalité  (\ref{eq:decomposition}) avec $X=\PP(w)$). 
 
Soit $p$ dans $|\PP(w)|$. Soit $(\widetilde{U}_{p},G_{p},\pi_{p})$ une 
carte d'un voisinage $U_{p}$ de $p$. Notons $\widetilde{p}$ le relevé 
de $p$ dans $\widetilde{U}_{p}$. L'action de $G_{p}$ sur l'espace 
tangent $T_{\widetilde{p}}\widetilde{U}_{p}$ induit la représentation 
de groupe suivante 
 
\begin{align*} G_{p} \longrightarrow & GL(n,\CC) \\ 
\exp(2i\pi\gamma) \longmapsto & \diag (e^{2i\pi\gamma 
  w_{0}}, \ldots ,e^{2i\pi\gamma w_{n}})   
\end{align*}

 \begin{prop}\label{prop:cohomologie,orbifold}  
  La structure de $\CC$-espace vectoriel gradué de la cohomologie 
  orbifolde de $\PP(w)$ est donnée par  
  \begin{align*} H^{2\star}_{\orb}(\PP(w),\CC) & =  
    \bigoplus_{\gamma\in 
      \V}H^{2(\star-a(\gamma))} ( |\PP(w)_{(e^{2i\pi\gamma})}|,\CC) \\ 
    & \simeq \bigoplus_{\gamma\in 
      \V}H^{2(\star-a(\gamma))}(|\PP(w_{I(\gamma)})|,\CC)  
\end{align*} 
où $a(\gamma)=\{\gamma w_{0}\}+\cdots+\{\gamma w_{n}\}$ et 
$I(\gamma)=\{i\mid \gamma w_{i} \in \NN\}$ (cf. le paragraphe \ref{sec:notations}). 
\end{prop} 
 
\begin{rem}\label{rem:amrani}
\begin{enumerate} 
\item Pour $\gamma=0$, nous avons
  $|\PP(w_{I(\gamma)})|=|\PP(w)|$. Ainsi, la cohomologie ordinaire de
  $|\PP(w)|$ est un sous-espace vectoriel de la cohomologie orbifolde
  de $\PP(w)$.
   \item 
    \label{rem:alamrani} Sur la deuxième page de l'article \cite{Kcwps}, 
    T.\ Kawasaki  
    démontre le résultat suivant 
    \begin{align*} 
    H^{2i}(|\PP(w)|,\CC)=\begin{cases} \CC & \mbox{si } 
        i\in\{0, \ldots ,n\} \,; \\ 0 & \mbox{sinon.}  \end{cases} 
    \end{align*} 
    Ainsi, d'après la proposition précédente, nous avons une 
    description explicite du $\CC$-espace vectoriel 
    $H^{2\star}_{\orb}(\PP(w),\CC)$. 
     
  \item\label{rem:iso} Nous avons 
    $|\PP(w)_{(e^{2i\pi\gamma})}|=|\PP(w)_{I(\gamma)}|$. D'après la  
    remarque \ref{rem:groupe,isotropie}.(\ref{rem:P(w)I}), l'espace 
    topologique $|\PP(w)_{I(\gamma)}|$ est muni d'une structure 
    orbifolde. 
 En effet, 
    nous avons les équivalences suivantes : 
\begin{align*}  
[p_{0}:\ldots:p_{n}]_{w} \in |\PP(w)_{(e^{2i\pi\gamma})}| 
& \Leftrightarrow  e^{2i\pi\gamma} \cdot (p_{0},\ldots,p_{n})=(p_{0},\ldots,p_{n})\\ 
&\Leftrightarrow  p_{i}=0 \mbox{ pour tout }i \mbox{ tels que } \gamma 
w_{i}\notin \NN \\ 
&\Leftrightarrow  p_{i}=0 \mbox{ pour tout } i \mbox{ dans } I^{c}(\gamma). 
\end{align*}

\item Prenons l'exemple des poids $w=(1,2,4)$. Nous avons alors  
  \begin{align*} 
    \Sigma \PP(1,2,4) &= |\PP(1,2,4)|\times\{1\} \bigsqcup 
    |\PP(2,4)|\times\{-1\} \\  & \bigsqcup 
    |\PP(4)|\times\{\sqrt{-1}\} \bigsqcup |\PP(4)|\times\{-\sqrt{-1}\}. 
  \end{align*} 
Les \og strates \fg définies par Kawasaki dans \cite{Ksto} sont 
$\{[0:0:1]\}$, $\{[0,y,z]\mid y\neq 0\}$ et $\{[x:y:z]\mid x\neq 0\}$. 
Remarquons que l'adhérence des strates est respectivement $\PP(4)$, 
$\PP(2,4)$ et $\PP(1,2,4)$, et qu'elle est isomorphe comme orbifolde aux composantes  
connexes de $\Sigma \PP(1,2,4)$. 
\end{enumerate} 
\end{rem}  
 
\begin{proof}[Démonstration de la proposition \ref{prop:cohomologie,orbifold}] 
  L'age de l'élément $e^{2i\pi\gamma}$ est $a(\gamma)$. Nous en 
  déduisons l'égalité de la proposition. D'après la remarque 
  \ref{rem:amrani} (\ref{rem:iso}), l'espace topologique 
  $|\PP(w)_{(e^{2i\pi\gamma})}|$ est égal à $|\PP(w)_{I(\gamma)}|$. Puis 
  le corollaire \ref{cor:identification} permet d'identifier 
  $\PP(w_{I(\gamma)})$ et $\PP(w)_{I(\gamma)}$. 
\end{proof}

De la remarque \ref{rem:amrani} (\ref{rem:alamrani}), nous en 
déduisons directement le corollaire suivant. 
 
\begin{cor} Soit $\gamma$ dans $S_{w}$.
  Les degrés de la cohomologie orbifolde qui proviennent de l'espace 
  topologique $|\PP(w)_{(e^{2i\pi\gamma})}|$ sont les nombres rationnels suivants 
  : 
  \begin{align*} 
  2\left(d+a(\gamma) \right) \mbox{ avec } 
  d\in\{0,\ldots,\delta(\gamma)-1\} 
  \end{align*} où $\delta(\gamma)=\# I(\gamma)$ (cf. paragraphe \ref{sec:notations}). 
\end{cor} 
 
 \begin{cor} La dimension du $\CC$-espace vectoriel $H^{2\star}_{\orb}(\PP(w),\CC)$ est $\mu$. 
\end{cor} 
 
\begin{proof} 
  D'après la proposition \ref{prop:cohomologie,orbifold}, la dimension 
  de l'espace vectoriel $H^{2\star}_{\orb}(\PP(w),\CC)$ est 
  \begin{align*} 
  \sum_{\gamma\in\V}\dim_{\CC}H^{\star}(|\PP(w)_{(e^{2i\pi\gamma})}|,\CC)=\sum_{\gamma\in 
    \V} \delta(\gamma)=\mu. 
  \end{align*} 
\end{proof} 
 
Pour tout $\gamma$ dans $\V$ posons 
\begin{align*} 
\eta_{\gamma}^{d} := 
\left(\frac{c_{1}(\mathcal{O}_{\PP(w_{I(\gamma)})}(\pgcd(w_{I(\gamma)})))}{\pgcd(w_{I(\gamma)})}\right)^{d} 
\in H^{2d}(|\PP(w_{I(\gamma)})|,\CC).\end{align*} 
Remarquons que 
$\eta_{\gamma}^{d}$ est nul pour $d\geq\delta(\gamma)$. 
 
Nous renvoyons à la 
formule (\ref{eq:defi,integrale}) pour la définition de l'intégrale orbifolde. 
 
\begin{prop}\label{prop:integrale} 
  On a l'égalité suivante 
  \begin{align*} 
  \int^{\orb}_{\PP(w)}\eta_{0}^{n}=\left({\prod_{i=0}^{n} w_{i}}\right)^{-1}. 
  \end{align*} 
\end{prop} 
  
\begin{proof} 
  D'a\-près la définition d'une intégrale orbifolde 
  (\ref{eq:defi,integrale}) et la remarque \ref{rem:groupe,isotropie}, 
  nous avons 
  \begin{align*} 
  \ds{\int^{\orb}_{\PP(w)}\eta_{0}^{n}=\frac{1}{\pgcd(w)}\int_{|\PP(w)_{\reg}|}\eta_{0}^{n}}. 
  \end{align*} 
  Puis les propositions \ref{prop:fibre,O(1)} et 
  \ref{prop:chern,fonctoriel} impliquent les égalités 
  suivantes : 
  \begin{align*}\ds{ 
    \int_{|\PP(w)_{\reg}|}\eta_{0}^{n}=\frac{1}{\deg(f_{w})}\int_{\PP^{n}}\left(\frac{c_{1}(\mathcal{O}_{\PP^{n}}(\pgcd(w)))}{\pgcd(w)}\right)^{n}=\frac{1}{\deg(f_{w})}}. 
  \end{align*} 
  La remarque \ref{rem:degre} termine la démonstration. 
\end{proof} 
 
 \begin{prop}\label{prop:base} 
  L'ensemble $\bs{\eta} :=\{\eta^{d}_{\gamma}\mid \gamma\in \V, 
  d\in\{0, \ldots ,\delta(\gamma)-1\}\}$ est une base du $\CC$-espace 
  vectoriel $H^{\star}_{\orb}(\PP(w),\CC)$. Le degré orbifold de 
  $\eta^{d}_{\gamma}$ est $2(d + a(\gamma))$. 
\end{prop} 
 
\begin{proof}  
  D'après la remarque \ref{rem:amrani}, la classe $\eta^{d}_{\gamma}$ 
  engendre le $\CC$-espace vectoriel 
  $H^{2d}(|\PP(w)_{(e^{2i\pi\gamma})})|,\CC)$. Ainsi, les éléments 
  $\eta^{0}_{\gamma}, \ldots ,\eta^{\delta(\gamma)-1}_{\gamma}$ 
  forment une base de $H^{\star}(|\PP(w)_{(e^{2i\pi\gamma})}|,\CC) $. 
  Nous en déduisons que $\bs{\eta}$ est une base de 
  l'espace vectoriel $H_{\orb}^{2\star}(\PP(w),\CC)$ et que 
  $\deg^{\orb}(\eta^{d}_{\gamma})=2(d+a(\gamma))$. 
\end{proof} 
 
\begin{expl}\label{expl:base,coho,122333}Considérons les poids $w=(1,2,2,3,3,3)$ (cet exemple est 
  considéré dans l'article \cite{Jocwps}).  
 Nous avons 
 \begin{align*} 
&  n=5, \,  \mu=14\, , \\  
&S_{w}=\left\{0,\frac{1}{3},\frac{1}{2},\frac{2}{3}\right\}, \\ 
&  \delta(0)=6,\, 
  \delta(1/3)=\delta(2/3)=3,\, \delta(1/2)=2 \\ 
& I(0)=\{0,1,2,3,4,5\},\, I(1/3)=I(2/3)=\{3,4,5\},\, I(1/2)=\{1,2\}. 
  \end{align*} 
Les ages sont $a(0)=0$, $a(1/3)=5/3$, $a(1/2)=2$ et $a(2/3)=4/3$. 
La cohomologie orbifolde  $H^{\ast}_{\orb}(\PP(1,2,2,3,3,3))$ est 
\begin{align*} 
   H^{2\ast}(|\PP(1,2,2,3,3,3)|)\oplus H^{2\ast-10/3}(|\PP(3,3,3)|) \\ 
 \oplus  H^{2\ast-4}(|\PP(2,2)|)\oplus H^{2\ast-8/3}(|\PP(3,3,3)|). 
\end{align*} 
Nous visualisons cette cohomologie \og en ligne \fg. Chaque ligne 
correspond à la cohomologie de $|\PP(w_{I(\gamma)})|$ et chaque groupe  
de cohomologie, noté $H^{i}$, qui apparaît dans le tableau ci-dessous est $\CC$.  
\begin{align*} 
&\gamma=0, \ |\PP(w_{I(0)})|=|\PP(1,2,2,3,3,3)| & & H^{0}\oplus 
  H^{2}\oplus H^{4}\oplus H^{6} \oplus H^{8}\oplus H^{10}\\ 
&\gamma = 1/3,\ |\PP(w_{I(1/3)})|=|\PP(3,3,3)| & & H^{10/3}\oplus 
H^{2+10/3}\oplus H^{4+10/3}\\ 
&\gamma = 1/2,\ |\PP(w_{I(1/2)})|=|\PP(2,2)| & & H^{4}\oplus 
H^{6}\\ 
&\gamma = 2/3,\ |\PP(w_{I(2/3)})|=|\PP(3,3,3)|& & H^{8/3}\oplus 
H^{2+8/3}\oplus H^{4+8/3} 
\end{align*} 
Nous visualisons la base $\bs{\eta}$ de la m\^{e}me manière. 
\begin{align*} 
&\gamma=0,  & & \CC\eta_{0}^{0}\oplus\CC\eta_{0}^{1}\oplus\CC\eta_{0}^{2}\oplus\CC\eta_{0}^{3}\oplus\CC\eta_{0}^{4}\oplus\CC\eta_{0}^{5}\\ 
&\gamma = 1/3, & & \CC\eta_{1/3}^{0}\oplus\CC\eta_{1/3}^{1}\oplus\CC\eta_{1/3}^{2}\\ 
&\gamma = 1/2, & & \CC\eta_{1/2}^{0}\oplus\CC\eta_{1/2}^{1} \\ 
&\gamma = 2/3, & & \CC\eta_{2/3}^{0}\oplus\CC\eta_{2/3}^{1}\oplus\CC\eta_{2/3}^{2} 
\end{align*} 
 
\end{expl}

\section{La dualité de Poincaré orbifolde des espaces projectifs à poids} 
\label{sec:dualite-de-poincare}  
 
Dans un premier paragraphe, nous donnerons la définition générale de 
la dualité de Poincaré pour les orbifolds complexes et commutatives 
(cf. l'article de \cite{CRnco}). 
Puis, nous appliquerons cette définition sur notre exemple fétiche $\PP(w)$. 
La proposition \ref{prop:dualite} nous donne une formule explicite de cette 
dualité dans la base $\bs{\eta}$. 
 
\subsection{Définition générale de la dualité de Poincaré orbifolde}\label{subsec:Dfinition-gnrale-de} 
 
Nous suivons le paragraphe $3.3$ de \cite{CRnco}.  Soit $X$ une 
orbifold complexe, commutative et compacte de dimension complexe $n$. 
L'application $I:X_{(g)}\rightarrow X_{(g^{-1})}$ qui à $(x,g)$ 
associe $(x,g^{-1})$ est un isomorphisme entre orbifolds. 
 
Pour tout $g$ et pour tout $0\leq d \leq n$, posons 
  \begin{align*} 
    \langle \cdot,\cdot\rangle_{g} : H^{2(d-\age(g))}(X_{(g)},\CC) 
    \times 
    H^{2(n-d-\age(g^{-1}))}(X_{(g^{-1})},\CC) & \longrightarrow  \CC \\ 
    (\alpha,\beta) & \longmapsto \int^{\orb}_{X_{(g)}}\alpha\wedge 
    I^{\ast}\beta 
  \end{align*}

  La dualité de Poincaré orbifolde est un accouplement 
\begin{align*} 
  \langle \cdot,\cdot\rangle : H^{d}_{\orb}(X,\CC) \times 
  H^{2n-d}_{\orb}(X,\CC) \rightarrow \CC 
\end{align*} 
pour tout $0\leq d \leq n$, défini par la somme directe des 
accouplements $\langle \cdot,\cdot\rangle_{g}$.

\begin{prop}[cf. proposition $3.3.1$ de \cite{CRnco}] 
  La dualité de Poincaré orbifolde est une forme bilinéaire non 
  dégénérée. 
\end{prop} 
 
\begin{rem}  
  D'après l'article de \cite{Sgm}, l'espace topologique sous-jacent à 
  une orbifold a une dualité de Poincaré qui est donnée par 
  l'intégrale orbifolde. Ainsi, la restriction de la dualité de 
  Poincaré orbifolde à $X_{(\id)}=|X|$ est la dualité de Poincaré de 
  Satake sur $H^{\star}(|X|,\CC)$. Remarquons que si l'orbifold n'est 
  pas réduite c'est-à-dire que $\ker(X)$ n'est pas réduit à $\{\id\}$, 
  alors d'après la formule (\ref{eq:defi,integrale}) l'intégrale 
  orbifolde est un multiple de l'intégrale ordinaire. 
\end{rem} 
 
\begin{expl} 
  \begin{enumerate} 
  \item Soit $Y$ une variété complexe vue comme orbifold.  La dualité 
    de Poincaré orbifolde est exactement la dualité de Poincaré 
    ordinaire sur $Y$ car l'intégrale orbifolde est l'intégrale 
    ordinaire. 
  \item Soit $G$ un groupe commutatif qui agit trivialement sur une variété 
complexe $Y$ de dimension $n$. Le quotient $X:=Y/G$ est une orbifold. Nous avons 
$|X|=|Y|$ mais pour tout $\alpha\in H^{2n}(Y,\CC)$, nous avons 
\begin{align*} 
\int^{\orb}_{X}\alpha=\frac{1}{\# G}\int_{Y}\alpha.   
\end{align*} 
Ainsi, la dualité de Poincaré orbifolde restreint à $X_{(id)}=X$ est 
à un multiple près la dualité de Poincaré de la variété $Y$.   
 
  \end{enumerate} 
\end{expl} 
 
\subsection{Dualité de Poincaré orbifolde de $\PP(w)$}\label{subsec:Dualit-de-Poincar} 
 
La proposition suivante  exprime la dualité de Poincaré orbifolde 
de $\PP(w)$ dans la base $\bs{\eta}$. 
 
\begin{prop}\label{prop:dualite} 
  Soient $\eta^{d}_{\gamma}$ et $\eta^{d'}_{\gamma'}$ deux éléments de 
  la base $\bs{\eta}$. 
\begin{enumerate} \item  Si  $\gamma'\neq\{1-\gamma\}$, alors on a $\langle \eta^{d}_{\gamma}, 
  \eta^{d'}_{\gamma'} \rangle = 0$. 
   
\item Si $\gamma'=\{1-\gamma\}$ alors on a $I(\gamma')=I(\gamma)$ et 
  \begin{align*} 
  \langle \eta^{d}_{\gamma},\eta^{d'}_{\{1-\gamma\}}\rangle= 
  \begin{cases}\left({\prod_{i\in I(\gamma)}w_{i}}\right)^{-1} & \mbox{si 
      } \deg^{\orb}( 
    \eta^{d}_{\gamma})+\deg^{\orb}(\eta^{d'}_{\{1-\gamma\}})=2n \,; \\0& 
    \mbox{sinon.}  \end{cases} 
  \end{align*} 
  \end{enumerate} 
\end{prop} 
 
\begin{rem}\label{rem:dual} 
\begin{enumerate} \item \label{item:25}  La remarque \ref{rem:symetrie} et la 
  proposition \ref{prop:base} impliquent que la condition $\deg^{\orb}( 
  \eta^{d}_{\gamma})+\deg^{\orb}(\eta^{d'}_{\{1-\gamma\}})=2n$ est 
  équivalente à la condition $d+d'=\delta(\gamma)-1$. 
\item Soient $\gamma\in \V$ et $d\in\{0,\ldots,\delta(\gamma)\}$. Le 
  dual de $\eta^{d}_{\gamma}$ pour la forme bilinéaire 
  $\langle\cdot,\cdot\rangle$ est $\left(\prod_{i\in 
      I(\gamma)}w_{i}\right)\eta^{\delta(\gamma)-1-d}_{\{1-\gamma\}}$. 
\end{enumerate} 
\end{rem} 
 
\begin{proof}[Démonstration de la proposition \ref{prop:dualite}] 
  D'après la définition de la dualité de Poincaré, si 
  $\gamma'\neq\{1-\gamma\}$ ou $d+d'\neq \delta(\gamma)-1$ alors 
  $\langle \eta^{d}_{\gamma},\eta^{d'}_{\{1-\gamma\}}\rangle =0$.  Si 
  $d+d'= \delta(\gamma)-1$, alors nous avons 
\begin{align*}\langle 
  \eta^{d}_{\gamma},\eta^{d'}_{\{1-\gamma\}}\rangle & = 
  {\int_{\PP(w_{I(\gamma)})}^{\orb} 
    \eta^{d}_{\gamma}\wedge\eta^{d'}_{\{1-\gamma\}}}\\ 
  &= 
  {\int_{\PP(w_{I(\gamma)})}^{\orb}\left(\frac{c_{1}^{\orb}(\mathcal{O}_{\PP(w_{I(\gamma)})}(\pgcd(w_{I(\gamma)})))}{\pgcd(w_{I(\gamma)})}\right)^{d+d'}}.\end{align*} 
 
Puis la proposition \ref{prop:integrale} implique que cette intégrale 
vaut $\left({\prod_{i\in I(\gamma)}w_{i}}\right)^{-1}$. Finalement, la 
remarque \ref{rem:dual}.(\ref{item:25}) nous permet de conclure. 
\end{proof} 
 
\begin{expl} 
  Pour les poids $w=(1,2,2,3,3,3)$, la dualité de Poincaré dans la 
  base $\bs{\eta}$ (cf. l'exemple \ref{expl:base,coho,122333})
  s'exprime par 
  \begin{align*} 
    \left( 
      \begin{array}{c|c|c|c} 
2^{-2}3^{-3} \antidiag_{6} & 0 &0 & 0\\ 
 \hline 0&0 &0 & 3^{-3}\antidiag_{3}\\ 
\hline 0&0 & 2^{-2}\antidiag_{2}& \\ 
\hline 0& 3^{-3}\antidiag_{3} &0 &0  
      \end{array} 
\right) 
  \end{align*} 
où $\antidiag_{n}$ est la matrice antidiagonale de taille $n\times n$ 
avec des $1$ sur l'antidiagonale. 
\end{expl} 
 
\section{Cup produit orbifold des espaces projectifs à poids} 
\label{sec:cup,produit} 
 
Dans un premier paragraphe, nous donnons la définition du cup  
produit orbifold pour les orbifolds complexes, compactes et commutatives 
(cf. l'article \cite{CRnco}).  
 
Dans un second paragraphe, nous expliciterons le cup produit orbifold 
de $\PP(w)$ dans la base $\bs{\eta}$. Ce résultat est donné 
explicitement dans le corollaire \ref{cor:cup}. En particulier, le 
fibré obstruction est calculé dans le théorème \ref{thm:fibre,obstruction}. 
 
\subsection{La définition du cup produit orbifold}\label{subsec:La-dfinition-du} 
Avant de commencer la définition, a priori surprenante, du cup produit 
orbifold, nous allons expliquer l'heuristique du problème. Si on 
veut couper une classe d'homologie de $X_{(g)}$ avec une classe de 
$X_{(h)}$, on ne peut pas forcément choisir des représentants de 
ces classes d'homologie de façon qu'ils soient transverses. En effet ces 
représentants doivent rester dans leur espace tordu respectif, ce qui 
est une contrainte supplémentaire.  Ainsi, il se peut que 
l'intersection des représentants n'ait pas la bonne dimension. Pour 
remédier à ce phénomène, on doit introduire un fibré correcteur. Dans
la littérature, on l'appelle aussi fibré obstruction.

Nous suivons le paragraphe $4.1$ de \cite{CRnco}.  Soit $X$ une 
orbifolde complexe, compacte et commutative. 
 
Considérons l'ensemble 
\begin{align*} 
  \Sigma_{3} X= \{(x,g_{1},g_{2},g_{3})\in \sqcup_{x\in  
  X} G_{x}\times G_{x}\times G_{x} \mid g_{1}g_{2}g_{3}=\id \}. 
\end{align*} 
 
Comme au paragraphe \ref{sec:coho-espace}, nous définissons une 
topologie sur $\Sigma_{3} X$ à l'aide d'une base d'ouverts.  
Pour tout $(x,g_{1},g_{2},g_{3})\in 
\Sigma_{3} X$ et pour toute carte $(\widetilde{U}_{x},G_{x},\pi_{x})$ d'un ouvert 
$U_{x}$ contenant $x$, nous posons 
\begin{align*}\left\{ 
  \begin{array}{l} 
 (y,h_{1},h_{2},h_{3})\in\Sigma_{3} X\mid \exists 
  (\alpha,\kappa):(\widetilde{U}_{y},G_{y},\pi_{y})\hookrightarrow(\widetilde{U}_{x},G_{x},\pi_{x}) \\ \mbox{ qui vérifie }\kappa(h_{i})=g_{i}. 
  \end{array}\right\} 
\end{align*} 
 Le lemme suivant se démontre de la m\^{e}me façon que le lemme \ref{lem:topologie,separe}. 
 
\begin{lem} 
  \begin{enumerate} 
  \item La collection de ces sous-ensembles de $\Sigma_{3} X$ définit une 
topologie sur $\Sigma_{3} X$.  
\item L'ensemble $\Sigma_{3} X$ muni de cette topologie est séparé. 
  \end{enumerate} 
\end{lem} 
 
L'application $P_{3}:\Sigma_{3}X\rightarrow |X|$ qui à 
$(x,g_{1},g_{2},g_{3})$ associe $x$ est continue. 
Le lemme suivante se démontre de la m\^{e}me façon que le lemme \ref{lem:composantes,connexes}. 
 
\begin{lem} 
  Soit $X$ une orbifolde complexe, compacte et commutative. 
  L'espace topologique $\Sigma_{3}X$ a un nombre fini de composantes connexes. 
\end{lem} 
 
Nous dirons que $(x,g_{1},g_{2},g_{3})\sim (y,h_{1},h_{2},h_{3})$ si 
$(x,g_{1},g_{2},g_{3})$ et $(y,h_{1},h_{2},h_{3})$ sont dans la 
m\^{e}me composante connexe de $\Sigma_{3}X$. Soit $T_{3}$ l'ensemble 
des classes d'équivalence. Nous 
notons $(g_{1},g_{2},g_{3})$ la classe d'équivalence de 
$(x,g_{1},g_{2},g_{3})$. La notation $(g_{1},g_{2},g_{3})$ sous-entend  
que $g_{1},g_{2},g_{3}$ sont dans un groupe $G_{x}$ pour un certain point $x\in |X|$. 
 Notons $X_{(g_{1},g_{2},g_{3})}$ la 
composante connexe de $\Sigma X$ qui contient 
$(x,g_{1},g_{2},g_{3})$. Nous avons la décomposition suivante de 
$\Sigma X$ 
\begin{align*} 
  \Sigma_{3} X=\bigsqcup_{(g_{1},g_{2},g_{3})\in T_{3}}X_{(g_{1},g_{2},g_{3})}. 
\end{align*} 
 
D'après le lemme $4.1.1$ de \cite{CRnco}, les espaces topologiques 
$X_{(g_{1},g_{2},g_{3})}$ sont naturellement munis d'une structure 
d'orbifolde. Nous ne donnerons pas de précision supplémentaire 
concernant le cas général mais dans le cas des espaces projectifs à poids, cette 
structure orbifolde est donnée par la remarque \ref{rem:amrani}.(\ref{rem:iso}). 
 
Soit $\PP^{1}$ la sphère de Riemann avec trois points marqués 
$(0,1,\infty)$.   
Soit $U_{0}$ la carte affine  de $\PP^{1}$ contenant $0$. Notons $z$ une coordonnée 
sur $U_{0}$. Soit $z=1/2$ le point base noté $\ast$.  
Nous considérons les lacets suivants :
\begin{align*} 
  \lambda_{0}:[0,2\pi]&\rightarrow U_{0} \\ 
t & \mapsto \exp(it)/2\,; \\ 
 \lambda_{1}:[0,2\pi]&\rightarrow U_{0} \\ 
t & \mapsto 1-\frac{1}{2}\exp(it)\,;\\ 
\lambda_{\infty}&:=(\lambda_{0}\lambda_{1})^{-1}. 
\end{align*} 
Le groupe fondamental de $\PP^{1}-\{0,1,\infty\}$ est engendré par les 
lacets $\lambda_{0},\lambda_{1},\lambda_{\infty}$ (nous notons de la 
m\^{e}me façon le lacet et sa classe dans le groupe fondamental). 
Soit $(x,g_{0},g_{1},g_{\infty})$ un élément de $\Sigma_{3}X$. Nous 
avons un morphisme de groupes surjectif 
$\rho:\pi_{1}(\PP^{1}-\{0,1,\infty \},\ast)\rightarrow H$ qui à 
$\lambda_{i}$ associe $g_{i}$ où $H$ est le sous-groupe de $G_{x}$ 
engendré par $g_{0},g_{1},g_{\infty}$.  Nous en déduisons un 
rev\^{e}tement galoisien de $\PP^{1}-\{0,1,\infty\}$ de groupe 
d'automorphismes $H$.  Nous complétons ce rev\^{e}tement en un 
rev\^{e}tement ramifié, noté $\pi:\Sigma\rightarrow \PP^{1}$. La 
variété $\Sigma$ est une surface de Riemann compacte.  Considérons 
l'application $\ev:X_{(g_{0},g_{1},g_{\infty})}\rightarrow X$ qui à 
$(x,g_{0},g_{1},g_{\infty})$ associe $x$. L'application $\ev$ est une 
bonne application. Nous renvoyons le lecteur à l'article de Chen et 
Ruan pour la démonstration dans le cas général (cf. le lemme $4.2.3$ 
dans \cite{CRnco}). Cependant, dans le cas des espaces projectifs à 
poids, nous verrons que l'application $\ev$ est une bonne application 
orbifolde (cf. $5$ lignes avant l'équation 
(\ref{eq:rang,fibre,poids})). 
Nous définissons le fibré 
orbifold sur $X_{(g_{0},g_{1},g_{\infty})}$ par 
\begin{align} 
  \label{eq:fibre,defi} 
 E(g_{0},g_{1},g_{\infty})&:=\left(\ev^{\ast} TX \otimes H^{0,1}(\Sigma,\CC)\right)^{H}.  
\end{align} 
Ce fibré est appelé fibré d'obstruction. 
 
D'après la preuve du théorème $4.1.5$ p.$20$ de l'article 
\cite{CRnco}, le rang du fibré $E(g_{0},g_{1},g_{\infty})$ est 
\begin{align} 
  \label{eq:rang,fibre} 
  \dim_{\CC}X_{(g_{0},g_{1},g_{\infty})}-\dim_{\CC}X +\age(g_{0})+\age(g_{1})+\age(g_{\infty}). 
\end{align}

\begin{defi}[cf. définition $4.1.2$ dans \cite{CRnco}]\label{defi:3,tenseur} 
  Soient $\alpha,\beta,\gamma$ dans $H^{\star}_{\orb}(X,\CC)$. On 
  définit 
  \begin{align*} 
  (\alpha,\beta,\gamma):=\sum_{(g_{1},g_{2},g_{3})\in 
    T_{3}}\int^{\orb}_{X_{(g_{1},g_{2},g_{3})}} 
  \ev_{1}^{\ast}\alpha\wedge \ev_{2}^{\ast}\beta \wedge 
  \ev_{3}^{\ast}\gamma\wedge c_{max}(E(g_{1},g_{2},g_{3})) 
  \end{align*} 
  où $c_{max}(E(g_{1},g_{2},g_{3}))\in H^{2\rg(E(\gamma_{0},\gamma_{1},\gamma_{\infty}))}(|X_{(g_{0},g_{1},g_{\infty})}|,\CC)$ est la classe de Chern 
  maximale  du fibré orbifold $E(g_{1},g_{2},g_{3})$ et 
  $ev_{i}:X_{(g_{1},g_{2},g_{3})}\rightarrow \Sigma{X}$ est 
  l'application qui à $(x,(g_{1},g_{2},g_{3}))$ associe $(x,g_{i})$ 
  pour $i\in\{1,2,3\}$. 
\end{defi} 
 
Puis, nous définissons le cup produit orbifold. 
 
\begin{defi}[cf. définition $4.1.3$ dans \cite{CRnco}] 
  \label{defi:Cup-produit-orbifold} 
  On définit le cup produit orbifold sur $H^{\star}_{\orb}(X,\CC)$ par 
  la formule 
\begin{align*} 
  \langle \alpha\cup\beta,\gamma\rangle=(\alpha,\beta,\gamma). 
\end{align*} 
\end{defi} 
 
\begin{expl} 
  \begin{enumerate} 
  \item  Soit $Y$ une variété complexe vue comme orbifold. Le fibré 
  obstruction est de rang $0$ et le cup 
  produit orbifold est simplement le cup produit ordinaire car  
   l'intégrale orbifolde est l'intégrale ordinaire. 
\item Soit $G$ un groupe commutatif qui agit trivialement sur une variété 
complexe $Y$ de dimension $n$. Le quotient $X:=Y/G$ est une orbifold. 
Le rang du fibré obstruction est $0$.  Le cup produit ordinaire sur $|Y|=|X|$  
est le cup produit orbifold restreint à $X_{(id)}=X$. En effet, comme 
dans la définition \ref{defi:3,tenseur}, nous intégrons avec 
l'intégrale orbifolde, nous avons la m\^{e}me constante qui appara\^{i}t 
de part et d'autre de l'égalité dans la définition  \ref{defi:Cup-produit-orbifold}. 
  \end{enumerate} 
\end{expl} 
 
Nous énonçons les propriétés du cup produit orbifold données dans le 
théorème $4.1.5$ de \cite{CRnco} mais nous renvoyons à la preuve dans 
l'article \cite{CRnco}. 
 
\begin{thm} 
  \label{thm:Cup-produit-orbifold} 
  Soit $X$ une orbifold complexe, commutative, compacte et connexe. 
\begin{enumerate} 
\item Le cup produit orbifold respecte la graduation de 
  $H^{\star}_{\orb}(X,\CC)$, c'est-à-dire 
    \begin{align*} 
      \cup : H^{p}_{\orb}(X,\CC)\times H_{\orb}^{q}(X,\CC) \rightarrow H^{p+q}_{\orb}(X,\CC). 
    \end{align*} 
  \item Le cup produit est associatif et son unité est la classe de 
    $1$ dans $H^{0}(|X|,\CC)$. 
  \item Pour tout $(\alpha,\beta) \in H^{d}_{\orb}(X,\CC)\times 
    H^{2n-d}_{\orb}(X,\CC)$, nous avons 
    \begin{align*} 
      \int^{\orb}_{X}\alpha\cup\beta=\langle \alpha,\beta\rangle. 
    \end{align*} 
  \item Le cup produit orbifold restreint à la cohomologie ordinaire 
    $H^{\star}(|X|,\CC)$ est le cup produit ordinaire sur $|X|$. 
  \end{enumerate} 
\end{thm}

\subsection{Calcul du cup produit pour l'orbifold $\PP(w)$}\label{subsec:Calcul-du-cup} 
 
Avant de calculer le fibré obstruction dans le cas de $\PP(w)$, nous 
démontrons quelques résultats préliminaires.

Nous allons calculer le fibré tangent orbifold de 
$\PP(w)$. Considérons l'injection suivante (cf. notation \ref{not:injection} ): 
\begin{align*} 
\alpha : \widetilde{U}& 
\hookrightarrow  \widetilde{V}  \\ 
(y_{0}, \ldots ,1_{i}, \ldots ,y_{n}) & \mapsto  
(y_{0}/y_{j}^{w_{0}/w_{j}},\ldots,1_{j},\ldots,y_{n}/y_{j}^{w_{n}/w_{j}})
\end{align*} 
 où  $\widetilde{U}\subset \widetilde{U}_{i}$ et  $\widetilde{V} \subset \widetilde{U}_{j}$.
Pour $k\neq j$, nous notons $t_{k}:=y_{k}/y_{j}^{w_{k}/w_{j}}$.  La 
fonction de transition sur $U\cap V$ est donnée par la matrice 
: 
\begin{align}\label{eq:25} 
\g_{\alpha}(y_{0}, \ldots ,1_{i}, \ldots 
,y_{n})&=\left(\frac{\partial t_{k}}{\partial 
    y_{\ell}}\right)_{(k,\ell)\in \{0, \ldots ,n\}\times\{0, \ldots 
  ,n\}-\{(j,i)\}} 
\end{align} 
 
D'après le corollaire \ref{cor:fibre,inverse}, un fibré vectoriel 
orbifold sur $\PP(w)$ se restreint à $\PP(w_{I})$. 
\begin{lem}\label{lem:decomposition} 
  Pour tout sous-ensemble $I$ de $\{0, \ldots ,n\}$, on a la 
  décompo\-si\-tion suivante : 
  \begin{align*} 
  T\PP(w)\mid_{\PP(w_{I})}\simeq \left(\bigoplus_{i\in 
    I^{c}}\mathcal{O}_{\PP(w_{I})}(w_{i})\right)\bigoplus T\PP(w_{I}). 
  \end{align*} 
\end{lem} 
 
\begin{proof} Nous allons démontrer le résultat pour $I=\{0, \ldots ,\delta\}$. 
  Soient $i,j$ dans $\{0, \ldots ,n\}$. Nous reprenons les notations 
  ci-dessus.  En coordonnées, la matrice de transition de $T\PP(w)$ 
  est donnée par 
 
\begin{align}\label{eq:changement} 
\ds{\frac{\partial}{\partial y_{\ell}}}&=\begin{cases} 
 \ds{\sum_{\stackrel{k=0}{k\neq j}}^{n} 
-\frac{w_{k}}{w_{j}}\frac{y_{k}}{y_{j}^{\frac{w_{k}}{w_{j}}-1}}\frac{\partial 
  }{\partial t_{k}}} & \mbox{si }  \ell=j\,; \\ 
 \ds{\frac{1}{y_{j}^{w_{\ell}/w_{j}}}  \frac{\partial }{\partial t_{\ell}} }  
 &  \mbox{si }  \ell \neq j.  
\end{cases}  
\end{align} pour $\ell\neq i$.

Restreindre le fibré $T\PP(w)$ à $\PP(w_{I})$ revient à poser 
$y_{\ell}=0$ pour $\ell>\delta$. Ainsi, les fonctions de transition du 
fibré orbifold $T\PP(w)\mid_{\PP(w_{I})}$ sont données par les 
relations (\ref{eq:changement}) où on ne fait la somme que jusqu'à 
$\delta$.  Finalement, la matrice de transition du fibré 
$T\PP(w)\mid_{\PP(w_{I})}$ qui correspond à l'injection $\alpha$ est 
du type suivant 
\begin{align*} 
\left( 
\begin{array}{cc} A_{ij} & 0 \\ 
0 & B_{ij}    \end{array} 
\right) 
\end{align*} 
où $A_{ij}$ est une matrice (de taille $\delta\times \delta$) de transition 
du fibré $T\PP(w_{I})$ et $B_{ij}$ est la matrice diagonale 
\begin{align*}\diag\left( \frac{1}{y_{j}^{w_{\delta+1}/w_{j}}},\ldots, 
  \frac{1}{y_{j}^{w_{n}/w_{j}}}\right).\end{align*} 
 
D'après la remarque \ref{rem:fibre}(\ref{rem:fibre,O(m)}), l'ensemble 
des matrices $\{B_{ij}\mid i,j \in\{0, \ldots ,\delta\}\}$ forme les 
fonctions de transition du fibré 
$\mathcal{O}_{\PP(w_{I})}(w_{\delta+1})\oplus \ldots \oplus 
\mathcal{O}_{\PP(w_{I})}(w_{n})$. 
\end{proof} 
 
Soient $\gamma_{0},\gamma_{1}$ et $\gamma_{\infty}$ dans $\V$ tels que 
$\gamma_{0}+\gamma_{1}+\gamma_{\infty}\in \NN$.  Nous posons 
\begin{align*} 
  I(\gamma_{0},\gamma_{1},\gamma_{\infty})&:=I(\gamma_{0})\cap 
I(\gamma_{1})\cap I(\gamma_{\infty}). 
\end{align*} 
Soit $H$ le groupe abélien 
engendré par $e^{2i\pi\gamma_{0}},e^{2i\pi\gamma_{1}}$ et 
$e^{2i\pi\gamma_{\infty}}$. De la construction expliquée dans le 
début du paragraphe \ref{sec:cup,produit}, nous déduisons un 
rev\^{e}tement ramifié de groupe d'automorphismes $H$, noté 
$\pi:\Sigma\rightarrow \PP^{1}$ où $\Sigma$ est une surface de Riemann 
compacte. 
 
L'application $\ev$ est simplement une inclusion de 
$\PP(w_{I(\gamma_{0},\gamma_{1},\gamma_{\infty})})$ dans $\PP(w)$. 
C'est une bonne application d'après la proposition 
\ref{prop:bonne,application}.  Nous en déduisons que le fibré 
obstruction, défini par la formule (\ref{eq:fibre,defi}) est 
\begin{align*} 
E(\gamma_{0},\gamma_{1},\gamma_{\infty}):=\left(T\PP(w)\mid_{\PP(w_{I(\gamma_{0},\gamma_{1},\gamma_{\infty})})}\otimes 
  H^{0,1}(\Sigma,\CC)\right)^{H}. 
\end{align*} 
 
D'après la formule (\ref{eq:rang,fibre}), le rang de ce fibré est 
\begin{align}\label{eq:rang,fibre,poids} 
\dim_{\CC}\PP(w_{I(\gamma_{0},\gamma_{1},\gamma_{\infty})})-\dim_{\CC}\PP(w)+a(\gamma_{0})+a(\gamma_{1})+a(\gamma_{\infty}). 
\end{align}

\begin{thm}\label{thm:fibre,obstruction}Soient 
  $\gamma_{0},\gamma_{1}$ et $\gamma_{\infty}$ dans $\V$ tels que 
  $\gamma_{0}+\gamma_{1}+\gamma_{\infty}\in \NN$.  Le fibré orbifold 
  $E(\gamma_{0},\gamma_{1},\gamma_{\infty})$ est isomorphe à 
   
  \begin{align*} 
  \bigoplus_{j\in J_{w}(\gamma_{0},\gamma_{1},\gamma_{\infty})} 
  \mathcal{O}_{\PP(w_{I(\gamma_{0},\gamma_{1},\gamma_{\infty})})}(w_{j}) 
  \end{align*} 
  où $J_{w}(\gamma_{0},\gamma_{1},\gamma_{\infty})\:=\{i\in\{0, 
  \ldots ,n\}\mid 
  \{\gamma_{0}w_{i}\}+\{\gamma_{1}w_{i}\}+\{\gamma_{\infty}w_{i}\}=2\}$. 
\end{thm}

\begin{proof}[Démonstration du théorème \ref{thm:fibre,obstruction}] 
 D'après la décomposition du lemme \ref{lem:decomposition}, le fibré
 obstruction  $E(\gamma_{0},\gamma_{1},\gamma_{\infty})$ est
 \begin{align*}
  \left(\bigoplus_{i\in 
    I^{c}(\gamma_{0},\gamma_{1},\gamma_{\infty})}\mathcal{O}_{\PP(w_{I(\gamma_{0},\gamma_{1},\gamma_{\infty})})}(w_{i})\otimes H^{0,1}(\Sigma,\CC)
\right)^{H}\bigoplus \Bigg(T\PP(w_{I(\gamma_{0},\gamma_{1},\gamma_{\infty})}) \otimes H^{0,1}(\Sigma,\CC)\Bigg)^{H}.
 \end{align*}
  Puisque $H$ agit trivialement sur 
  $T\PP(w_{I(\gamma_{0},\gamma_{1},\gamma_{\infty})})$ et 
  \begin{align*} 
 H^{0,1}(\Sigma,\CC)^{H}=0 \mbox{ car } H^{1}(\PP^{1},\CC)=0,    
  \end{align*} 
  nous obtenons la décomposition suivante :
  \begin{align*} 
  \ds{E(\gamma_{0},\gamma_{1},\gamma_{\infty})=\left(\bigoplus_{i\in 
        I^{c}(\gamma_{0},\gamma_{1},\gamma_{\infty})} 
      \mathcal{O}_{\PP(w_{I(\gamma_{0},\gamma_{1},\gamma_{\infty})})}(w_{i})\otimes 
      H^{0,1}(\Sigma,\CC)\right)^{H}}. 
  \end{align*} 
   
  Pour $i\in\{0, \ldots ,n\}$, notons $\chi_{i}$ le caractère du
  groupe $H$ qui à $e^{2\sqrt{-1}\pi\gamma_{j}}$ associe
  $e^{2\sqrt{-1}\pi\gamma_{j}w_{i}}$ pour $j\in\{0,1,\infty\}$. Le groupe $H$
  agit sur
  $\mathcal{O}_{\PP(w_{I(\gamma_{0},\gamma_{1},\gamma_{\infty})})}(w_{i})$
  par multiplication par $\chi_{i}$. Soit $\omega$ une $(1,0)$-forme
  différentielle fermée telle que
  $\varphi_{g}^{\ast}\omega=\chi_{i}(g)\omega$. Nous avons
  \begin{align*} 
    \overline{\varphi_{g}^{\ast}\omega}&=\overline{\chi_{i}(g)\omega}\\ 
&={\chi_{i}(g)}^{-1} {\overline{\omega}}.  
  \end{align*} 
Nous en déduisons que 
$\overline{H^{1,0}(\Sigma,\CC)}_{\chi_{i}}=H^{0,1}(\Sigma,\CC)_{\chi_{i}^{-1}}$. 
Finalement, nous obtenons 
  \begin{align*} 
  E(\gamma_{0},\gamma_{1},\gamma_{\infty})=\bigoplus_{i\in 
    I^{c}(\gamma_{0},\gamma_{1},\gamma_{\infty}) } 
  \left(\mathcal{O}_{\PP(w_{I(\gamma_{0},\gamma_{1},\gamma_{\infty})})}(w_{i})\otimes 
    \overline{H^{1,0}(\Sigma,\CC)_{\chi_{i}}}\right). 
  \end{align*} 
    Pour finir la démonstration, il suffit d'appliquer le lemme suivant. 
\end{proof} 
 
\begin{lem}\label{lem:cohomologie} Pour tout $i\in\{0, \ldots ,n\}$, on a   
  \begin{align*}H^{1,0}(\Sigma,\CC)_{\chi_{i}}= 
\begin{cases}  
\CC & \mbox {si } 
\{\gamma_{0}w_{i}\}+\{\gamma_{1}w_{i}\}+\{\gamma_{\infty}w_{i}\}=2\,;\\ 
0 & \mbox{sinon. } 
 \end{cases} 
\end{align*} 
\end{lem}

\begin{rem}\label{rem:partie,fractionnaire}  
  Nous avons l'inclusion suivante 
  \begin{align*} 
J_{w}(\gamma_{0},\gamma_{1},\gamma_{\infty}) \subset (I(\gamma_{0})\bigcup I(\gamma_{1})\bigcup I(\gamma_{\infty}))^{c}. 
      \end{align*} 
\end{rem} 
 
Soit $\chi$ un caractère du groupe $H$. 
Avant de démontrer ce lemme, nous  allons d'abord calculer la 
caractéristique d'Euler 
$e(\PP^{1},\left(\pi_{\ast}\underline{\CC}_{\Sigma}\right)_{\chi})$ 
dans les lemmes \ref{lem:caracteristique,euler} et \ref{lem:euler} 
puis le lemme \ref{lem:appendice,caractere} calculera $H^{1}(\Sigma,\CC)_{\chi}$. 
 
Pour tout ouvert $U$ dans $\PP^{1}$, notons $h_{0}(\pi^{-1}(U))$ le 
nombre de composantes connexes de $\pi^{-1}(U)$. Le groupe $H$ agit 
transitivement sur 
$\pi_{\ast}\underline{\CC}_{\Sigma}(U)=\CC^{h_{0}(\pi^{-1}(U))}$ en 
permutant les coordonnées de $\CC^{h_{0}(\pi^{-1}(U))}$. Posons 
\begin{align*} 
  (\pi_{\ast}\underline{\CC}_{\Sigma})_{\chi}(U)=\{x\in 
\CC^{h_{0}(\pi^{-1}(U))}\mid \forall h\in H, h\cdot x=\chi(h)x\}. 
\end{align*} 
Notons $\mathcal{F}$ le faisceau 
$(\pi_{\ast}\underline{\CC}_{\Sigma})_{\chi}$.  Comme 
$\pi:\Sigma\rightarrow\PP^{1}$ est un rev\^etement ramifié aux points 
$0,1$ et $\infty$ de $\PP^{1}$, le faisceau 
$\mathcal{F}\mid_{\PP^{1}-\{0,1,\infty\}}$ est un faisceau localement 
constant de rang $1$. En effet, une fonction constante 
$f:\pi^{-1}(U)\rightarrow \CC$, qui vérifie $f(hx)=\chi(h)f(x)$ pour tout 
$(h,x)\in H\times \pi^{-1}(U)$, est déterminée par sa valeur sur une 
composante connexe de $\pi^{-1}(U)$.

\begin{lem}\label{lem:caracteristique,euler} 
  On a l'égalité suivante : 
  \begin{align*} 
  e(\PP¹,\mathcal{F})=e(\PP¹-\{0,1,\infty\},\mathcal{F}\mid_{\PP¹-\{0,1,\infty\}})+ 
  \dim_{\CC}\mathcal{F}_{0}+\dim_{\CC}\mathcal{F}_{1}+\dim_{\CC}\mathcal{F}_{\infty} 
  \end{align*} 
  où $e$ est la caractéristique d'Euler. 
\end{lem} 
 
\begin{proof}[Démonstration du lemme \ref{lem:caracteristique,euler}] 
  On a la suite exacte de Mayer-Vietoris suivante : 
  \begin{align*} 
  \cdots \rightarrow 
  H^{i}_{\{0,1,\infty\}}(\PP¹,\mathcal{F})\rightarrow 
  H^{i}(\PP¹,\mathcal{F}) \rightarrow 
  H^{i}_{\{0,1,\infty\}}(\PP¹-\{0,1,\infty\},\mathcal{F}\mid_{\PP¹-\{0,1,\infty\}}) 
  \rightarrow \cdots 
  \end{align*} 
  où $H^{i}_{\{0,1,\infty\}}(\PP¹,\mathcal{F})$ désigne la 
  cohomologie à support dans $\{0,1,\infty\}$.  Nous en dédui\-sons 
  l'égalité suivante : 
   \begin{align*}e(H^{\star}_{\{0,1,\infty\}}(\PP¹,\mathcal{F}))-e(\PP¹,\mathcal{F})+ 
  e(\PP¹-\{0,1,\infty\},\mathcal{F}\mid_{\PP¹-\{0,1,\infty\}})=0. 
  \end{align*} 
  Comme nous avons 
  \begin{align*}e(H^{\star}_{\{0,1,\infty\}}(\PP¹,\mathcal{F})) = 
  e(H^{\star}_{\{0\}}(\PP¹,\mathcal{F}))+e(H^{\star}_{\{1\}}(\PP¹,\mathcal{F}))+e(H^{\star}_{\{\infty\}}(\PP¹,\mathcal{F})),\end{align*} 
  il suffit de montrer que $ 
  e(H^{\star}_{\{0\}}(\PP¹,\mathcal{F}))=\dim_{\CC}\mathcal{F}_{0}$. 
  Soit $V$ un disque ouvert centré en $0$ ne contenant ni $1$ ni 
  $\infty$ tel que 
  $\dim_{\CC}\mathcal{F}(V)=\dim_{\CC}\mathcal{F}_{0}$. Le lemme 
  d'excision nous donne la suite exacte suivante : 
  \begin{align*} 
  \xymatrix{ \cdots\ar[r] & 
    H^{i}_{\{0\}}(V,\mathcal{F}\mid_{V})\ar[r] & 
    H^{i}(V,\mathcal{F}\mid_{V}) \ar[r] & 
    H^{i}(V-\{0\},\mathcal{F}\mid_{V-\{0\}})\ar[r] & \cdots} 
  \end{align*} 
  Nous obtenons l'égalité 
  \begin{align*} 
  e(H^{\star}_{\{0\}}(V,\mathcal{F}\mid_{V})) - 
  e(V,\mathcal{F}\mid_{V})+ e( V-\{0\},\mathcal{F}\mid_{V-\{0\}})=0. 
  \end{align*} 
   
  Or comme $\mathcal{F}\mid_{V-\{0\}}$ est un faisceau localement 
  constant de rang $1$, d'après \cite{DMmhf} p.$11$, nous avons $e( 
  V-\{0\},\mathcal{F}\mid_{V-\{0\}})=e(V-\{0\})=0$ car $V-\{0\}$ est 
  homéomorphe à $S¹$. 
   
  Finalement, nous obtenons 
  $e(H^{\star}_{\{0\}}(V,\mathcal{F}\mid_{V})) = 
  e(V,\mathcal{F}\mid_{V})$. Le complexe de faisceaux  
  $(\pi_{\ast}\mathcal{E}^{\bullet}_{\Sigma})_{\chi}$, défini dans la 
  démonstration du théorème \ref{thm:cohomologie,caractere}, est une 
  résolution acyclique pour le foncteur $\Gamma(\PP^{1},\cdot)$ du 
  faisceau $\mathcal{F}$.  Nous en déduisons que le complexe 
  $(\pi_{\ast}\mathcal{E}^{\star}_{\Sigma})_{\chi}\mid_{V}$ est une 
  résolution acyclique pour le foncteur $\Gamma(V,\cdot)$ du faisceau 
  $\mathcal{F}\mid_{V}$.  Comme $V$ est convexe, le lemme de Poincaré 
  et le théorème \ref{thm:cohomologie,caractere} impliquent que 
  $H^{i}(V,\mathcal{F}\mid_{V})=0$ pour $i>0$.  Nous en déduisons que 
  $e(V,\mathcal{F}\mid_{V})=\dim_{\CC}H^{0}(V,\mathcal{F}\mid_{V})$. 
\end{proof}

\begin{lem}\label{lem:euler} 
  Soit $\chi$ un caractère non trivial de $H$. On a 
  \begin{align*} 
  e(\PP¹,(\pi_{\ast}\underline{\CC}_{\Sigma})_{\chi})=\begin{cases} 
      -1 & \mbox{ si } \chi(e^{2i\pi\gamma_{0}})\neq 1, 
      \chi(e^{2i\pi\gamma_{1}})\neq 1 \mbox{ et } 
      \chi(e^{2i\pi\gamma_{\infty}})\neq 1\,; \\ 0 & \mbox{sinon.} 
  \end{cases}\end{align*} 
\end{lem} 
 
\begin{rem}\label{rem:caractere}Rappelons que $H$ est engendré par 
  $e^{2i\pi\gamma_{0}},e^{2i\pi\gamma_{1}}$ et 
  $e^{2i\pi\gamma_{\infty}}$ et que 
  $\gamma_{0}+\gamma_{1}+\gamma_{\infty}\in \NN$ c'est-à-dire que 
  $e^{2i\pi\gamma_{0}}.e^{2i\pi\gamma_{1}}.e^{2i\pi\gamma_{\infty}}=1$. Nous en déduisons que  
 si $\chi$ est un caractère non trivial de $H$ alors le cardinal de 
  l'ensemble $\{j\in\{0,1,\infty\}\mid \chi(e^{2i\pi\gamma_{j}})=1\}$ 
  est au plus $1$. 
\end{rem} 
 
\begin{proof}[Démonstration du lemme \ref{lem:euler}] 
  D'après le lemme \ref{lem:caracteristique,euler}, il suffit de 
  calculer les nombres 
  $e(\PP¹-\{0,1,\infty\},\mathcal{F}\mid_{\PP¹-\{0,1,\infty\}}), 
  \dim_{\CC}\mathcal{F}_{0},\dim_{\CC}\mathcal{F}_{1}$ et 
  $\dim_{\CC}\mathcal{F}_{\infty}$.  D'après la relation $(2.2.1)$ de 
  \cite{DMmhf}, on a 
  \begin{align*}
 e(\PP¹-\{0,1,\infty\},\mathcal{F}\mid_{\PP¹-\{0,1,\infty\}})= 
  e(\PP¹)-e(\{0,1,\infty\})=2-3=-1. 
\end{align*}

  Montrons l'égalité suivante  : 
  \begin{align*} 
  \dim_{\CC}\mathcal{F}_{0}=\begin{cases} 1 &\mbox{ si } 
      \chi(e^{2i\pi\gamma_{0}})=1 \,; \\ 0 & \mbox{ si } 
      \chi(e^{2i\pi\gamma_{0}})\neq 1. 
  \end{cases} 
\end{align*} 
 
Soit $V$ un voisinage ouvert de $0$ tel que 
$\dim_{\CC}\mathcal{F}(V)=\dim_{\CC}\mathcal{F}_{0}$ et que le nombre 
de composantes connexes de $\pi^{-1}(V)$ soit $\# H/\langle 
e^{2i\pi\gamma_{0}} \rangle$, où $\langle e^{2i\pi\gamma_{0}}\rangle$ 
est le sous-groupe engendré par $e^{2i\pi\gamma_{0}}$.  La condition 
$\chi(e^{2i\pi\gamma_{0}})=1$ est équivalente à l'existence d'un 
élément non nul dans $\mathcal{F}(V)$. Ainsi nous avons 
\begin{align*} 
\dim_{\CC}\mathcal{F}(V)=1 \Leftrightarrow 
\chi(e^{2i\pi\gamma_{0}})=1. 
\end{align*}

Le même raisonnement pour les points $1$ et $\infty$ et la remarque 
\ref{rem:caractere} nous donne l'alternative suivante : 
\begin{itemize} \item   soit $\chi(e^{2i\pi\gamma_{0}})\neq 1, \chi(e^{2i\pi\gamma_{1}})\neq 1$ et  
  $\chi(e^{2i\pi\gamma_{\infty}})\neq 1$ et alors $\dim 
  \mathcal{F}_{0}=\dim \mathcal{F}_{1}=\dim \mathcal{F}_{\infty}=0$ ;
\item soit il existe 
  $\gamma\in\{\gamma_{0},\gamma_{1},\gamma_{\infty}\}$ tel que 
  $\chi(e^{2i\pi\gamma})=1$ et alors $\dim \mathcal{F}_{0}+\dim 
  \mathcal{F}_{1}+\dim \mathcal{F}_{\infty}=1$. 
 \end{itemize} 
\end{proof} 
  
Les deux lemmes précédents \ref{lem:caracteristique,euler} et 
\ref{lem:euler}, nous permettent de calculer  $H¹(\Sigma,\CC)_{\chi}$. 
 
\begin{lem}\label{lem:appendice,caractere} Soit $\chi$ un caractère 
  non trivial de $H$.  On a 
  \begin{align*} 
H¹(\Sigma,\CC)_{\chi}=  
\begin{cases}  
\CC & \mbox{ si }     \chi(e^{2i\pi\gamma_{0}})\neq 1,\, \chi(e^{2i\pi\gamma_{1}})\neq 1 \mbox{ et } 
    \chi(e^{2i\pi\gamma_{\infty}})\neq 1 \,; \\ 0 & \mbox{sinon.} 
\end{cases} 
\end{align*} 
\end{lem} 
 
\begin{proof}[Démonstration du lemme \ref{lem:appendice,caractere}]  
  Montrons que  
  \begin{align*} 
    H^{0}(\PP¹,\mathcal{F})=H^{2}(\PP¹,\mathcal{F})=0. 
  \end{align*} 
  Le théorème \ref{thm:cohomologie,caractere}  implique que 
  $H^{j}(\PP¹,\mathcal{F})=H^{j}(\Sigma,\underline{\CC}_{\Sigma})_{\chi}$ 
  pour tout $j$. Or $\Sigma$ est connexe donc 
  $H^{0}(\Sigma,\underline{\CC}_{\Sigma})=\CC$ et 
  $H^{2}(\Sigma,\underline{\CC}_{\Sigma})=\CC$. Puisque $\chi$ n'est 
  pas le caractère trivial, nous avons 
  $H^{0}(\Sigma,\underline{\CC}_{\Sigma})_{\chi}=H^{2}(\Sigma,\underline{\CC}_{\Sigma})_{\chi}=0$. 
   
  Finalement, d'après le théorème \ref{thm:cohomologie,caractere} et 
  le lemme \ref{lem:caracteristique,euler}, nous obtenons 
  \begin{align*}\dim_{\CC}H^{1}(\Sigma,\CC_{\Sigma})_{\chi}=\dim_{\CC} 
  H^{1}(\PP^{1},\mathcal{F})= -e(\PP^{1},\mathcal{F}).\end{align*} 
\end{proof} 
 
Il nous reste à démontrer le lemme \ref{lem:cohomologie}. 
 
\begin{proof}[Démonstration du lemme \ref{lem:cohomologie}]
 Rappelons que $\pi:\Sigma\to \PP^{1}$ est un rev\^{e}tement ramifié
 aux points $0,1$ et $\infty$ de $\PP^{1}$ de groupe d'automorphismes 
 $H$ où $\Sigma$ est une surface de Riemann compacte. Nous avons 
 \begin{align*}
   \#\pi^{-1}(0)=\#H/o(e^{2\sqrt{-1}\pi\gamma_{0}}) \,;\\
\#\pi^{-1}(1)=\#H/o(e^{2\sqrt{-1}\pi\gamma_{1}}) \,;\\
\#\pi^{-1}(\infty)=\#H/o(e^{2\sqrt{-1}\pi\gamma_{\infty}}).
 \end{align*}
 
  D'après le lemme \ref{lem:appendice,caractere}, nous pouvons nous 
  restreindre au cas où $i\in (I(\gamma_{0})\bigcup I(\gamma_{1}) \bigcup 
  I(\gamma_{\infty}))^{c}$. 
\begin{itemize} \item   Montrons d'abord que si 
  $H^{1,0}(\Sigma,\CC)_{\chi_{i}}=\CC$ alors $i\in 
  J_{w}(\gamma_{0},\gamma_{1},\gamma_{\infty})$. 
   
  Notons $t$ la coordonnée sur $\PP^{1}-\{\infty\}$. Soit $f$ une 
  fonction holomorphe sur $\Sigma$. Posons 
  $\omega=f\pi^{\ast}(dt/t(t-1))$. 
   
  Nous cherchons une condition sur $f$ pour que 
\begin{enumerate}  
\item[(i)] $\omega$ soit holomorphe ;
\item [(ii)] on ait 
  $\left(h^{-1}\right)^{\ast}\omega=\chi_{i}(h)\omega$ quelque soit 
  $h$ dans $H$. 
  \end{enumerate} 
  Soit $z_{0}$ la coordonnée dans un voisinage d'un point de 
  $\pi^{-1}(0)$.  Dans ces coordonnées, la $1$-forme holomorphe 
  $\pi^{\ast}(dt/t(t-1))$ s'écrit 
  \begin{align*} 
    o(e^{2\sqrt{-1}\pi\gamma_{0}})dz_{0}/z_{0}(z_{0}^{o(e^{2\sqrt{-1}\pi\gamma_{0}})}-1) 
  \end{align*} 
  où $o(e^{2\sqrt{-1}\pi\gamma_{0}})$ est l'ordre de $e^{2\sqrt{-1}\pi\gamma_{0}}$ 
  dans $H$.  La condition $(\rm{ii})$ avec $h=e^{2\sqrt{-1}\pi\gamma_{0}}$ se 
  traduit par 
  \begin{align*} 
    f(\exp(-2\sqrt{-1}\pi/o(e^{2\sqrt{-1}\pi\gamma_{0}}))z_{0})=\exp(2\sqrt{-1}\pi\gamma_{0}w_{i})f(z_{0}). 
  \end{align*} 
  Posons $f(z_{0})=\sum_{n\geq n_{0}} a_{n}z_{0}^{n}$ avec 
  $a_{n_{0}}\neq 0$. Cette condition sur $f$ impose que si $a_{n}\neq 
  0$ alors $\gamma_{0}w_{i}+n/o(e^{2\sqrt{-1}\pi\gamma_{0}})$ est dans $\ZZ$. 
  Ainsi, nous avons 
  \begin{align*} 
    n_{0}/o(e^{2\sqrt{-1}\pi\gamma_{0}})=-\{\gamma_{0}w_{i}\}+\alpha_{0} 
    \mbox{ avec } \alpha_{0}\in \ZZ. 
   \end{align*} 
  La condition $(\rm{i})$ implique que $n_{0}\geq 1$ c'est-à-dire 
  $\alpha_{0}\geq 1$ car $\alpha_{0}$ est un entier supérieur ou égal 
  à $1/o(e^{2\sqrt{-1}\pi\gamma_{0}})+\{\gamma_{0}w_{i}\}$. 
   
  Nous en déduisons que le nombre de zéros de $\omega$ comptés avec 
  multiplicité\footnote{Ne pas confondre la multiplicité des zéros et
    la multiplicité des points de ramifications. Chaque point de
    ramification de  $\pi^{-1}(0)$ a pour multiplicité
    $o(e^{2\sqrt{-1}\pi\gamma_{0}})$. Dans notre calcul, nous comptons
    les multiplicités des zéros sans les multiplicités des points
    de ramifications.} aux $\# H/o(e^{2\sqrt{-1}\pi\gamma_{0}})$ points de $\pi^{-1}(0)$ est 
\begin{align}\label{eq:0} 
(n_{0}-1)\# H/o(e^{2\sqrt{-1}\pi\gamma_{0}})&=\#  
H(-\{\gamma_{0}w_{i}\}+\alpha_{0}-1/o(e^{2\sqrt{-1}\pi\gamma_{0}})). 
\end{align} 
 
Les conditions $(\rm{i})$ et $(\rm{ii})$ appliquées dans un voisinage d'un point 
de $\pi^{-1}(1)$ impliquent que $\alpha_{1}\geq 1$ et que le nombre de 
zéros de $\omega$ comptés avec multiplicité aux $\# H/o(e^{2\sqrt{-1}\pi\gamma_{1}})$ points de $\pi^{-1}(1)$ est 
\begin{align}\label{eq:1} 
\# H(-\{\gamma_{1}w_{i}\}+\alpha_{1}-1/o(e^{2\sqrt{-1}\pi\gamma_{1}})). 
\end{align} 
 
Soit $z_{\infty}$ la coordonnée dans un voisinage d'un point de 
$\pi^{-1}(\infty)$.  Dans ces coordonnées, nous avons 
\begin{align*} 
  \pi^{\ast}(dt/t(t-1))=-o(e^{2\sqrt{-1}\pi\gamma_{\infty}})dz/(1-z^{o(e^{2\sqrt{-1}\pi\gamma_{\infty}})}). 
\end{align*} 
De m\^{e}me, les conditions $(i)$ et $(ii)$ impliquent 
$\alpha_{\infty}\geq 0$ et le nombre de zéros de $\omega$ comptés avec 
multiplicité aux $\# H/o(e^{2\sqrt{-1}\pi\gamma_{\infty}})$ points de $\pi^{-1}(\infty)$ est 
\begin{align}\label{eq:infty} 
\# H(-\{\gamma_{\infty}w_{i}\}+\alpha_{\infty}+1-1/o(e^{2\sqrt{-1}\pi\gamma_{\infty}})). 
\end{align}  
 
La formule de Riemann-Hurwitz nous donne le genre de $\Sigma$, noté 
$g_{\Sigma}$. 
\begin{align}\label{eq:RH} 
2g_{\Sigma}-2&=\# H\left(1-\frac{1}{o(e^{2\sqrt{-1}\pi\gamma_{0}})}-\frac{1}{o(e^{2\sqrt{-1}\pi\gamma_{1}})}-\frac{1}{o(e^{2\sqrt{-1}\pi\gamma_{\infty}})}\right). 
\end{align} 
 
Comme $\omega$ est une $1$-forme holomorphe sur une surface de Riemann 
compacte $\Sigma$, la somme de ses zéros vaut $2g_{\Sigma}-2$. D'après 
les relations (\ref{eq:0}), (\ref{eq:1}), (\ref{eq:infty}) et 
(\ref{eq:RH}), nous en déduisons 
\begin{align*} 
\{\gamma_{0}w_{i}\}+\{\gamma_{2}w_{i}\}+\{\gamma_{\infty}w_{i}\}=\alpha_{0}+\alpha_{1}+\alpha_{\infty}. 
\end{align*} 
Vu les conditions sur $\alpha_{0},\alpha_{1}$ et $\alpha_{\infty}$, 
l'égalité ci-dessus n'est possible que si $\alpha_{0}=\alpha_{1}=1$ et 
$\alpha_{\infty}=0$.

\item Pour finir la démonstration, il suffit de démontrer que si 
  $H^{1,0}(\Sigma,\CC)_{\chi_{i}}=0$ alors 
  $\{w_{i}\gamma_{0}\}+\{w_{i}\gamma_{1}\}+\{w_{i}\gamma_{\infty}\}\neq 
  2$.  D'après le lemme \ref{lem:appendice,caractere}, nous avons 
  l'implication suivante 
  \begin{align*} 
    i\notin 
  I(\gamma_{0})\bigcup I(\gamma_{1}) \bigcup I(\gamma_{\infty}) \Rightarrow  
  H^{1}(\Sigma,\CC)_{\chi_{i}}=\CC. 
  \end{align*} 
  Or, nous avons la décomposition 
  \begin{align*} 
    H^{1}(\Sigma,\CC)_{\chi_{i}}=H^{1,0}(\Sigma,\CC)_{\chi_{i}}\oplus 
  \overline{H^{1,0}(\Sigma,\CC)_{\chi^{-1}_{i}}}. 
  \end{align*} 
Ainsi si 
  $H^{1,0}(\Sigma,\CC)_{\chi_{i}}=0$ alors 
  $H^{1,0}(\Sigma,\CC)_{\chi^{-1}_{i}}=\CC$. Puis la première partie 
  de la démonstration implique que nous avons 
  \begin{align*} 
    \{-w_{i}\gamma_{0}\}+\{-w_{i}\gamma_{1}\}+\{-w_{i}\gamma_{\infty}\}= 
  2  \end{align*} 
c'est-à-dire 
\begin{align*} 
  \{w_{i}\gamma_{0}\}+\{w_{i}\gamma_{1}\}+\{w_{i}\gamma_{\infty}\}=1\neq 
  2 
\end{align*} 
car pour $i\notin I(\gamma_{0})\bigcup I(\gamma_{1}) \bigcup 
  I(\gamma_{\infty})$, nous avons 
  $\{-w_{i}\gamma_{j}\}=1-\{w_{i}\gamma_{j}\}$ pour tout 
  $j\in\{0,1,\infty\}$. 
 \end{itemize} 
\end{proof} 
 
Les applications $\ev_{j}$ de la définition \ref{defi:3,tenseur} sont 
simplement des inclusions naturelles de  
$\PP(w_{I(\gamma_{0},\gamma_{1},\gamma_{\infty})})$ dans 
$\PP(w)_{e^{2i\pi\gamma_{j}}}$, pour $j$ dans $\{0,1,\infty\}$. 
 
Le corollaire suivant calcule l'expression du $3$-tenseur 
$(\cdot,\cdot,\cdot)$ dans la base $\bs{\eta}$ de 
$H^{\star}_{\orb}(\PP(w),\CC)$. 
  
\begin{cor}\label{cor:tenseur} 
  Soient $\eta^{d_{0}}_{\gamma_{0}},\eta^{d_{1}}_{\gamma_{1}}$ et 
  $\eta^{d_{\infty}}_{\gamma_{\infty}}$ des éléments de la base 
  $\bs{ \eta}$. \begin{enumerate} 
  \item Si $\gamma_{0}+\gamma_{1}+\gamma_{\infty}$ n'est pas un 
    entier, alors 
    $(\eta^{d_{0}}_{\gamma_{0}},\eta^{d_{1}}_{\gamma_{1}},\eta^{d_{\infty}}_{\gamma_{\infty}} 
    )= 0$. 
  \item Si $\gamma_{0}+\gamma_{1}+\gamma_{\infty}$ est un entier, 
    alors on a 
\begin{align*} 
  \left( 
    \eta^{d_{0}}_{\gamma_{0}},\eta^{d_{1}}_{\gamma_{1}},\eta^{d_{\infty}}_{\gamma_{\infty}} 
  \right)= \begin{cases} {\frac{\ds{\prod_{i\in 
            J_{w}(\gamma_{0},\gamma_{1},\gamma_{\infty})} 
          w_{i}}}{\ds{\prod_{i\in
            I(\gamma_{0},\gamma_{1},\gamma_{\infty})} w_{i}}} }&
    \mbox{ si }
    \ds{\sum_{i\in\{0,1,\infty\}}\deg^{\orb}(\eta^{d_{i}}_{\gamma_{i}})=2n } \,; \\ 
    0 & \mbox{ sinon.} 
\end{cases} 
\end{align*} 
  \end{enumerate}  
\end{cor}  
 
\begin{proof}\begin{enumerate} \item   
    La première partie du corollaire découle de la défi\-ni\-tion et de la 
    formule \ref{eq:triplet}. 
     
  \item Notons 
    $w_{0,1,\infty}:=w_{I(\gamma_{0},\gamma_{1},\gamma_{\infty})}$. 
    Soient $\gamma_{0},\gamma_{1}$ et $\gamma_{\infty}$ dans $\V$ tels 
    que $\gamma_{0}+\gamma_{1}+\gamma_{\infty}$ soit un entier. Par 
    définition, 
    $(\eta^{d_{0}}_{\gamma_{0}},\eta^{d_{1}}_{\gamma_{1}},\eta^{d_{\infty}}_{\gamma_{\infty}} 
    )$ est non nul si 
 
    \begin{align*}d_{0}+d_{1}+d_{\infty}+rang(E(\gamma_{0},\gamma_{1},\gamma_{\infty}))=\dim_{\CC}\PP(w_{0,1,\infty}).\end{align*} 
    Sous cette condition, qui est exactement 
    $\sum_{i\in\{0,1,\infty\}}\deg^{\orb}(\eta^{d_{i}}_{\gamma_{i}})=2n$ 
    (cf.la formule (\ref{eq:rang,fibre,poids})), nous avons 
     
\begin{align*}(\eta^{d_{0}}_{\gamma_{0}},\eta^{d_{1}}_{\gamma_{1}},\eta^{d_{\infty}}_{\gamma_{\infty}} 
  ) = 
  \ds{\int^{\orb}_{\PP(w_{0,1,\infty})}\iota_{\gamma_{0}}^{\ast}\eta^{d_{0}}_{\gamma_{0}}\wedge 
    \iota_{\gamma_{1}}^{\ast}\eta^{d_{1}}_{\gamma_{1}}\wedge 
    \iota_{\gamma_{\infty}}^{\ast}\eta^{d_{\infty}}_{\gamma_{\infty}}\wedge 
    c_{\max}(E(\gamma_{0},\gamma_{1},\gamma_{\infty}))}\\ 
  = \ds{\prod_{i\in J_{w}(\gamma_{0},\gamma_{1},\gamma_{\infty})} 
    w_{i} \int^{\orb}_{\PP(w_{0,1,\infty})} \left( \frac 
      {c_{1}(\mathcal{O}_{\PP(w_{0,1,\infty})})(\pgcd(w_{0,1,\infty}))} 
      {\pgcd(w_{0,1,\infty})} \right)^{\dim_{\CC}\PP(w_{0,1,\infty})}}. 
    \end{align*} 
     
    Puis la proposition \ref{prop:integrale} termine la 
    démonstration. 
  \end{enumerate} 
\end{proof} 
 
\begin{expl} 
  Pour l'exemple $w=(1,2,2,3,3,3)$, nous allons regrouper les 
  $3$-tenseurs non nuls selon le triplet 
  $(\gamma_{0},\gamma_{1},\gamma_{\infty})$. 
Pour $(\gamma_{0},\gamma_{1},\gamma_{\infty})=(0,0,0)$ le fibré 
obstruction est de rang $0$ et nous avons 
\begin{align*} 
 (\eta_{0}^{0},\eta_{0}^{0},\eta_{10}^{0}) &  = 2^{-2}3^{-3}, &  (\eta^{0}_{0},\eta_{0}^{1},\eta_{0}^{4}) &  = 2^{-2}3^{-3},\\ 
(\eta_{0}^{0},\eta_{0}^{2},\eta_{0}^{3}) &  =2^{-2}3^{-3}, &  (\eta^{1}_{0},\eta^{1}_{0},\eta^{3}_{0}) &  =2^{-2}3^{-3},\\ 
(\eta^{1}_{0},\eta^{2}_{0},\eta^{2}_{0}) &  =2^{-2}3^{-3}.
\end{align*} 
Ceci va nous donner le cup produit sur $H^{\ast}(|\PP(1,2,2,3,3,3)|)$ 
(cf. le théorème \ref{thm:Cup-produit-orbifold}). 
Pour  $(\gamma_{0},\gamma_{1},\gamma_{\infty})=(0,1/2,1/2)$ le fibré 
obstruction est de rang $0$ et nous avons 
\begin{align*} 
  (\eta_{0}^{0},\eta_{1/2}^{0},\eta_{1/2}^{1}) &  = 2^{-2}, & 
(\eta_{0}^{1},\eta_{1/2}^{0},\eta_{1/2}^{0}) &  = 2^{-2}. \\ 
\end{align*} 
Pour  $(\gamma_{0},\gamma_{1},\gamma_{\infty})=(0,1/3,2/3)$ le fibré 
obstruction est de rang $0$ et nous avons 
\begin{align*} 
(\eta_{0}^{0},\eta_{1/3}^{0},\eta_{2/3}^{2}) & = 3^{-3},&  (\eta_{0}^{0},\eta_{1/3}^{2},\eta_{2/3}^{0}) & = 3^{-3},\\ 
(\eta_{0}^{0},\eta_{1/3}^{1},\eta_{2/3}^{1}) & =3^{-3},& (\eta_{0}^{1},\eta_{1/3}^{1},\eta_{2/3}^{0}) & = 3^{-3},\\ 
 (\eta_{0}^{1},\eta_{1/3}^{0},\eta_{2/3}^{1}) & = 3^{-3},& (\eta_{0}^{2},\eta_{1/3}^{0},\eta_{2/3}^{0}) & = 3^{-3}.
 \end{align*} 
Pour  $(\gamma_{0},\gamma_{1},\gamma_{\infty})=(1/3,1/3,1/3)$, le 
fibré obstruction est $\mathcal{O}(2)\oplus \mathcal{O}(2)$ et nous avons 
  \begin{align*} 
    (\eta_{1/3}^{0},\eta_{1/3}^{0},\eta_{1/3}^{0}) & = 4.3^{-3}. 
  \end{align*} 
 Pour  $(\gamma_{0},\gamma_{1},\gamma_{\infty})=(2/3,2/3,2/3)$, le 
fibré obstruction est $\mathcal{O}(1)$ et nous 
 avons 
 \begin{align*} 
   (\eta_{2/3}^{0},\eta_{2/3}^{0},\eta_{2/3}^{1}) & = 1.3^{-3}. 
 \end{align*} 
 
\end{expl} 
 
D'après la définition $4.1.3$ de l'article \cite{CRnco}, le cup 
produit orbifold est défini par l'égalité 
$(\alpha,\beta,\gamma)=\langle \alpha\cup\beta,\gamma \rangle$.  Comme 
$(\cdot,\cdot,\cdot)$ est symétrique en ses $3$ arguments,  
$(H_{\orb}^{\star}(\PP(w),\CC),\cup,\langle\cdot,\cdot\rangle)$ est 
une algèbre de Frobenius graduée.

\begin{cor}\label{cor:cup} Soient $\eta^{d_{0}}_{\gamma_{0}}$ et 
  $\eta^{d_{1}}_{\gamma_{1}}$ des éléments de la base $\bs{\eta}$. 
  On a 
  \begin{align*} 
  \eta^{d_{0}}_{\gamma_{0}}\cup \eta^{d_{1}}_{\gamma_{1}}= 
  \left(\prod_{i\in K(\gamma_{0},\gamma_{1})} w_{i}\right) 
  \eta^{d}_{\{\gamma_{0}+\gamma_{1}\}} 
  \end{align*} où \begin{align*}K(\gamma_{0},\gamma_{1}):= 
        J_{w}\Bigl(\gamma_{0},\gamma_{1},\{1-\{\gamma_{0}+\gamma_{1}\}\}\Bigr)\bigsqcup 
        I(\{\gamma_{0}+\gamma_{1}\})-I(\gamma_{0})\cap 
        I(\gamma_{1})\end{align*} 
  avec 
  $d:=\frac{\deg^{\orb}(\eta^{d_{0}}_{\gamma_{0}})}{2}+\frac{\deg^{\orb}(\eta^{d_{1}}_{\gamma_{1}})}{2}-a(\{\gamma_{0}+\gamma_{1}\})$. 
\end{cor} 
 
\begin{rem} 
\begin{enumerate} 
\item Si $d\geq\delta(\{\gamma_{0}+\gamma_{1}\})$ alors 
  $\eta^{d}_{\{\gamma_{0}+\gamma_{1}\}}=0$. 
\item Pour tout $\gamma\in S_{w}$ et pour tout $d\in\{0, \ldots 
  ,\delta(\gamma)\}$, nous avons 
  \begin{align*} 
  \eta_{0}^{1}\cup \eta_{\gamma}^{d}=\eta_{\gamma}^{d+1}. 
  \end{align*} 
\end{enumerate} 
\end{rem} 
 
\begin{proof}
   
  D'après la définition du cup produit orbifold, nous avons 
  \begin{align*} 
  \eta^{d_{0}}_{\gamma_{0}}\cup \eta^{d_{1}}_{\gamma_{1}}=\sum_{ 
    \stackrel{\gamma\in \V}{\delta\in\{0, \ldots ,\delta(\gamma)-1\} 
      }} 
  (\eta^{d_{0}}_{\gamma_{0}},\eta^{d_{1}}_{\gamma_{1}},\eta^{\delta}_{\gamma}) 
  {\eta^{\delta}_{\gamma}}^{\ast}. 
  \end{align*} 
  La remarque \ref{rem:dual} nous donne le dual de 
  $\eta^{\delta}_{\gamma}$.  Le corollaire \ref{cor:tenseur} et la 
  formule (\ref{eq:triplet}) impliquent que si 
  $(\eta^{d_{0}}_{\gamma_{0}},\eta^{d_{1}}_{\gamma_{1}},\eta^{\delta}_{\gamma})\neq 
  0$ alors $\gamma=1-\{\gamma_{0}+\gamma_{1}\}$. 
   
  Pour finir la démonstration, il suffit de vérifier que la condition 
  \begin{align*} 
    \deg^{\orb}(\eta^{d_{0}}_{\gamma_{0}})+\deg^{\orb}(\eta^{d_{1}}_{\gamma_{1}})+\deg^{\orb}(\eta^{\delta}_{\gamma})=2n 
  \end{align*} 
  est équivalente à 
  \begin{align*} 
    \delta=a(\{\gamma_{0}+\gamma_{1}\})-1+\delta(\{\gamma_{0}+\gamma_{1}\})-\deg^{\orb}(\eta^{d_{0}}_{\gamma_{0}})/2-\deg^{\orb}(\eta^{d_{1}}_{\gamma_{1}})/2. 
  \end{align*} 
\end{proof} 
 
\begin{expl} 
  Pour les poids $w=(1,2,2,3,3,3)$, nous donnons la table du cup produit 
  orbifold dans la base $\bs{\eta}$. 
  \begin{align*} 
    \begin{array}{|c|c|c|c|c|c|c||c|c|c|c|c|c|c|c|} 
\hline &\eta_{0}^{0}&\eta_{0}^{1}&\eta_{0}^{2}&\eta_{0}^{3}&\eta_{0}^{4}&\eta_{0}^{5}&\eta_{1/3}^{0}&\eta_{1/3}^{1}&\eta_{1/3}^{2}&\eta_{1/2}^{0}&\eta_{1/2}^{1}&\eta_{2/3}^{0}&\eta_{2/3}^{1}&\eta_{2/3}^{2} \\ 
\hline \eta^{0}_{0}&\eta_{0}^{0}&\eta_{0}^{1}&\eta_{0}^{2}&\eta_{0}^{3}&\eta_{0}^{4}&\eta_{0}^{5}&\eta_{1/3}^{0}&\eta_{1/3}^{1}&\eta_{1/3}^{2}&\eta_{1/2}^{0}&\eta_{1/2}^{1}&\eta_{2/3}^{0}&\eta_{2/3}^{1}&\eta_{2/3}^{2} \\ 
\hline\eta_{0}^{1}&&\eta_{0}^{2}&\eta_{0}^{3}&\eta_{0}^{4}&\eta_{0}^{5}&0&\eta_{1/3}^{1}&\eta_{1/3}^{2}&0&\eta_{1/2}^{1}&0&\eta_{2/3}^{1}&\eta_{2/3}^{2}&0 \\ 
\hline\eta_{0}^{2}&&&\eta_{0}^{4}&\eta_{0}^{5}&0&0&\eta_{1/3}^{2}&0&0&0&0&\eta_{2/3}^{2}&0&0 \\ 
\hline\eta_{0}^{3}&&&&0&0&0&0&0&0&0&0&0&0&0 \\ 
\hline\eta_{0}^{4}&&&&&0&0&0&0&0&0&0&0&0&0 \\ 
\hline\eta_{0}^{5}&&&&&&0&0&0&0&0&0&0&0&0 \\ 
\hline\hline\eta_{1/3}^{0}&&&&&&&4.\eta_{2/3}^{2}&0&0&0&0&4\eta_{0}^{3}&4\eta_{0}^{4}&4\eta_{0}^{5} \\ 
\hline\eta_{1/3}^{1}&&&&&&&&0&0&0&0&4\eta_{0}^{4}&4\eta_{0}^{5}&0 \\ 
\hline\eta_{1/3}^{2}&&&&&&&&&0&0&0&4\eta_{0}^{5}&0&0 \\ 
\hline\eta_{1/2}^{0}&&&&&&&&&&3^{3}\eta_{0}^{4}&3^{3}\eta_{0}^{5}&0&0&0 \\ 
\hline\eta_{1/2}^{1}&&&&&&&&&&&0&0&0&0 \\ 
\hline\eta_{2/3}^{0}&&&&&&&&&&&&1.\eta_{1/3}^{1}&1.\eta_{1/3}^{2}&0 \\ 
\hline\eta_{2/3}^{1}&&&&&&&&&&&&&0&0\\ 
\hline\eta_{2/3}^{2}&&&&&&&&&&&&&&0 \\ 
 \hline   \end{array} 
  \end{align*} 
  La partie en haut à gauche est simplement le cup produit sur 
  $H^{\ast}(|\PP(1,2,2,3,3,3)|)$. 
Parmi ces produits, le calcul explicite du fibré obstruction est 
nécessaire pour les produits 
\begin{align*} 
  \eta_{1/3}^{0}\cup \eta_{1/3}^{0}&  =4.\eta_{2/3}^{2} & &\mbox{où le fibré 
    obstruction est } \mathcal{O}(2)\oplus\mathcal{O}(2); \\ 
\eta_{2/3}^{0}\cup\eta_{2/3}^{0}&= 1.\eta_{1/3}^{1}  & &\mbox{où le fibré 
    obstruction est } \mathcal{O}(1); \\ 
\eta_{2/3}^{0}\cup\eta_{2/3}^{1}&= 1.\eta_{1/3}^{2}  & &\mbox{où le fibré 
    obstruction est } \mathcal{O}(1). \\ 
\end{align*} 
\end{expl} 
 
\chapter{Cohomologie quantique orbifolde des espaces projectifs à poids}\label{cha:la-cohom-quant} 

\section{Les invariants de Gromov-Witten orbifolds}\label{sec:IGWO} 
 
Nous suivons le paragraphe $2.3$ de l'article \cite{CRogw}. 
\begin{defi}\label{defi:courbe,nodale} 
  \emph{Une courbe nodale de genre $0$ avec $k$ points marqués} est la donnée 
  suivante : 
\begin{itemize} \item  un ensemble fini $I$ et pour tout $i$ 
  dans $I$, une application continue $\varphi_{i}:C_{i}\rightarrow C$ 
  entre une courbe complexe lisse de genre $0$ et un espace 
  topologique ;
\item $k$ points distincts sur $C$, notés $\bs{z}:=(z_{1}, \ldots 
  ,z_{k})$. 
\end{itemize} 
Ces données vérifient les conditions suivantes :
\begin{enumerate}  
\item l'espace topologique $C$ est la réunion des 
  $\varphi_{i}(C_{i})$ ; 
\item pour tout $p_{i} $ dans $C_{i}$, il existe un voisinage 
  $U_{p_{i}}$ de $p_{i}$ tel que l'application 
  $\varphi_{i}\mid_{U_{p_{i}}}:U_{p_{i}}\rightarrow C$ soit un 
  homéomorphisme sur son image ; 
\item pour tout $p$ dans $C$, nous avons 
  $\sum_{i}\varphi_{i}^{-1}(p)\leq 2$ ; 
\item pour tout $z_{j}\in \bs{z}$, nous avons 
  $\sum_{i}\varphi_{i}^{-1}(p)=1$ ;
\item l'ensemble des points nodaux $\{p\in C \mid 
  \sum_{i}\varphi_{i}^{-1}(p)= 2\}$ est fini. 
\end{enumerate} 
\end{defi} 
  
\begin{rem} 
\begin{enumerate}  
\item Les points nodaux ne sont pas marqués. 
\item Nous dirons qu'un point $p\in C_{i}$ est marqué (resp. nodal) si nous avons 
  $\varphi_{i}(p)\in \bs{z}$ (resp. $\varphi_{i}(p)$ est nodal). 
\end{enumerate} 
\end{rem} 
 Une application $\theta:(C,\bs{z})\rightarrow (C',\bs{z'})$ est un 
  \emph{isomorphisme} si $\theta$ est un homéo\-mor\-phisme qui se relève en un 
  biholomorphisme $\theta_{ij}:C_{i}\rightarrow C'_{j}$ sur chaque 
  composante $C_{i}$ de $C$ et si $\theta(z_{i})=z'_{i}$.

Pour tout $(t,r)\in \RR^{+}\times\left(\RR^{+}-\{0\}\right)$, posons 
$X(t,r)=\{(x,y)\in \CC^{2}\mid |x|,|y|\leq r, xy=t\}$.  Le groupe 
$\bs{\mu}_{n}$ agit sur $X(t,r)$ par la formule $\zeta(x,y)=(\zeta 
x,\zeta^{-1}y)$.  Notons $\varphi_{n}:X(t,r)\rightarrow 
X(t^{n},r^{n})$ l'application qui à $(x,y)$ associe $(x^{n},y^{n})$. 
Nous en déduisons que le triplet $(X(t,r),\bs{\mu}_{n},\varphi_{n})$ 
est une carte de $X(t^{n},r^{n})$. 
 
\begin{defi}\label{defi:orbicourbe} 
  \emph{Une orbicourbe nodale} est une courbe nodale marquée $(C,\bs{z})$ 
  avec une structure orbifolde qui vérifie 
\begin{enumerate}  
\item tout point orbifold $p_{i}$ de $C_{i}$, c'est-à-dire $\# 
  G_{p_{i}}>1$, est  un point marqué ou un point nodal ; 
\item pour tout point $z_{j}\in \bs{z}$, il existe une carte  
  d'un voisinage de $p$ qui est donnée par le rev\^etement ramifié 
  $z\mapsto z^{m_{j}}$ ; 
\item pour tout point nodal, il existe une carte donnée par 
  $(X(0,r_{\ell}),\bs{\mu}_{n_{\ell}},\varphi_{n_{\ell}})$. 
\end{enumerate} 
Notons $(C,\bs{z},\bs{m},\bs{n})$, ou plus simplement $(C,\bs{z})$ 
quand il n'y a pas d'ambigu\"{i}té, une telle orbicourbe nodale. 
\end{defi}

  Un \emph{isomorphisme entre deux orbicourbes nodales}  
  \begin{align*} 
    \widetilde{\theta} 
  :(C,\bs{z},\bs{m},\bs{n})\rightarrow (C,\bs{z},\bs{m},\bs{n}) 
  \end{align*} 
est 
  une collection d'isomorphismes $\theta_{ij}$ entre les orbicourbes 
  $C_{i}$ et $C'_{j}$ telle qu'ils induisent un isomorphisme 
  $\theta:(C,\bs{z})\rightarrow (C',\bs{z'})$.

\begin{defi}\label{defi:application,stable}  
 Soit $X$ une orbifold complexe et commutative.
 Une application orbifolde \emph{stable} est la donnée suivante : 
  \begin{itemize} 
  \item une  orbicourbe nodale $(C,\bs{z},\bs{m},\bs{n})$ ;  
  \item une application $f:C\rightarrow X$ continue ; 
\item et une classe d'isomorphisme de structures compatibles notée 
  $\xi$. 
  \end{itemize} 
Ce triplet $(f,(C,\bs{z},\bs{m},\bs{n}),\xi)$ vérifie 
\begin{enumerate} 
\item pour tout $i\in I$, l'application $f_{i}:=f\circ\varphi_{i}$ est holomorphe de $C_{i}$ dans  $X$;  
\item \label{item:monomorphisme} pour tout point $z_{i}$ marqué ou 
  nodal, le morphisme de groupes induit par $\xi$ de $G_{z_{i}}$ dans 
  $G_{f(z_{i})}$ est injectif ;
\item et si $f_{i}$ est l'application constante sur $C_{i}$ alors $C_{i}$ 
  a au moins trois points singuliers (i.e. nodal ou marqué). 
\end{enumerate} 
\end{defi} 
 
Nous munissons l'ensemble des applications orbifoldes stables d'une 
relation d'équivalence de la manière suivante. Nous dirons que les 
deux applications orbifoldes stables 
$(f,(C,\bs{z},\bs{m},\bs{n}),\xi)$ et 
$(f',(C',\bs{z'},\bs{m'},\bs{n'}),\xi')$ sont équivalentes s'il existe 
un isomorphisme 
\begin{align*} 
  \widetilde{\theta}:(C,\bs{z},\bs{m},\bs{n}) 
\rightarrow(C',\bs{z'},\bs{m'},\bs{n'}) 
\end{align*} 
tel que $f'\circ\widetilde{\theta}=f$ et que l'image inverse de $\xi$ 
par $\widetilde{\theta}$ soit isomorphe à $\xi'$.  Notons 
$[f,(C,\bs{z},\bs{m},\bs{n}),\xi]$ la classe d'équivalence de 
l'application orbifolde stable $(f,(C, \bs{z},\bs{m},\bs{n}),\xi)$. 
 
Soit $(f,(C,\bs{z},\bs{m},\bs{n}),\xi)$ une application stable. Nous 
pouvons associer à cette application stable la classe d'homologie dans 
$H_{2}(|X|,\ZZ)$ définie par 
$f_{\ast}([C])=\sum_{i}(f\circ\varphi_{i})_{\ast}[C_{i}]$. 
Cette classe d'homologie ne dépend que de la classe d'équivalence de l'application stable. 
 Pour chaque 
point marqué $z_{i}$, la classe de structure compatible $\xi$ induit 
un monomorphisme de groupes $\kappa_{i}:G_{z_{i}}\hookrightarrow 
G_{f(z_{i})}$. Ce monomorphisme ne dépend que de la classe d'équivalence de l'application stable. 
 
Rappelons que $\Sigma{X}=\bigsqcup_{x\in X}G_{x}=\bigsqcup_{(g)\in T} X_{(g)}$ où $X_{(g)}$ 
est une composante connexe (cf. paragraphe \ref{subsec:Dfinition-de-la}). 
 Nous avons une application d'évaluation, 
notée $\ev$, qui à chaque classe $[f,(C,\bs{z},\bs{m},\bs{n}),\xi]$ d'applications 
orbifoldes stables  associe 
$\left((f(z_{1}),\kappa_{1}(e^{{2i\pi}/{m_{1}}})), \ldots 
  ,(f(z_{k}),\kappa_{k}(e^{{2i\pi}/{m_{k}}}))\right) \in {\Sigma X}^{k}$. 
 
Une application orbifolde stable $(f,(C,\bs{z},\bs{m},\bs{n}),\xi)$ 
est dite de \emph{type} $(g_{1}, \ldots ,g_{k})$ si 
$\left((f(z_{1}),\kappa_{1}(e^{{2i\pi}/{m_{1}}})), \ldots 
  ,(f(z_{k}),\kappa_{k}(e^{{2i\pi}/{m_{k}}}))\right)$ appartient à 
$X_{(g_{1})}\times\cdots X_{(g_{k})}$. S'il n'y a pas d'ambigu\"ité dans la 
notation, nous écrirons $\underline{g}$ pour le $k$-uplet $(g_{1}, \ldots ,g_{k})$. 
 
\begin{defi}  Soit $X$ une orbifold complexe et commutative.
Soit $A\in H_{2}(|X|,\ZZ)$. Nous définissons  
$\overline{\mathcal{M}}_{k}(A,\underline{g})$ l'espace de modules des classes 
d'équivalence des applications orbifoldes stables avec $k$ points 
marqués, de classe d'homologie $A$ et de type $\underline{g}$, c'est-à-dire 
\begin{align*} 
\overline{ \mathcal{M}}_{k}(A,\underline{g})= 
\left\{  
 \begin{array}{c} 
 \mbox{ } [(f,(C,\bs{z},\bs{m},\bs{n}),\xi)]\mid 
   \# z =k, f_{\ast}[C]=A, 
  \\ ev(f,(C,\bs{z},\bs{m},\bs{n}),\xi)\in X_{(g_{1})}\times\cdots\times  
  X_{(g_{k})} 
\end{array}   
\right\}.  
\end{align*}  
 \end{defi} 
 
D'après les résultats de l'article \cite{CRogw} (cf. proposition 
$2.3.6$), l'espace de modules $\overline{\mathcal{M}}_{k}(A,\underline{g})$ 
est un espace topologique compact et métrisable. 
 Chen et Ruan définissent également une structure de 
Kuranishi pour l'espace de modules 
$\overline{\mathcal{M}}_{k}(A,\underline{g})$ dont la dimension est donnée 
par le théorème suivant.

\begin{thm}[cf. théorème $A$ de \cite{CRogw}]\label{thm:dimension} 
 Soit $X$ une orbifold complexe et commutative. La dimension\footnote{Dans la littérature, on la trouve également
 sous le nom de dimension attendue ou dimension virtuelle.} de la
 structure de Kuranishi considérée par Chen et Ruan de
 $\overline{\mathcal{M}}_{k}(A,\underline{g})$ est
 \begin{align*} 
  2\left(\int_{A}c_{1}(TX)+\dim_{\CC}X-3+k-\sum_{i=1}^{k}age(g_{i})\right). 
 \end{align*} 
\end{thm} 
   
Cette structure de Kuranishi définit (cf. théorème
$6.12$ et le paragraphe $17$ de l'article \cite{FOgwi}), une classe d'homologie,
appelée classe fondamentale de la structure de Kuranishi, 
\begin{align}\label{eq:23} 
   \rm{ev}_{\ast}[\overline{\mathcal{M}}_{k}(A,\underline{g})]&\in 
   H_{2\left(\int_{A}c_{1}(TX)+n-3+k-\sum_{i=1}^{k}age(g_{i})\right)} 
 (X_{(g_{1})}\times\cdots\times X_{(g_{k})},\CC). 
 \end{align} 
 
Pour chaque $i\in\{1, \ldots ,k\}$, soit $\alpha_{i}$ une classe dans 
$H^{2(\star-age(g_{i}))}(X_{(g_{i})},\CC) \subset 
H^{2\star}_{\orb}(\PP(w),\CC)$. La formule $(1.3)$ de \cite{CRogw} 
définit les invariants de Gromov-Witten orbifolds par

 \begin{align}\label{eq:defi,IGW} 
\Psi^{A}_{k,\underline{g}} :  
H^{\star}(X_{(g_{1})},\CC)\otimes \cdots \otimes H^{\star}(X_{(g_{k})},\CC) & 
 \longrightarrow  \CC \\ 
 \alpha_{1}\otimes\cdots\otimes\alpha_{k} &  \longmapsto  \int_{\rm{ev}_{\ast}[\overline{\mathcal{M}}_{k}(A,\underline{g})]}\alpha_{1}\wedge\cdots\wedge\alpha_{k}\nonumber. 
 \end{align} 
  
 
\section{Potentiel de Gromov-Witten pour $\PP(w)$} 
\label{sec:pot,GW} 
 
Dans ce paragraphe, nous allons définir le champ d'Euler et le 
potentiel de Gromov-Witten. Dans la proposition 
\ref{prop:Potentiel-de-Gromov}, nous montrons que le potentiel est 
homogène par rapport au champ d'Euler et nous calculons certaines 
conditions initiales de la structure de Frobenius sur $H^{\ast}(\PP(w),\CC)$.   
 
Nous commençons par un lemme qui va nous permettre de calculer la 
classe de Chern orbifolde du fibré $T\PP(w)$. 
 
\begin{lem}\label{lem:suite,exacte} 
  Nous avons une suite exacte de fibrés  
  \begin{align*} 
\xymatrix{0\ar[r]& \underline{\CC} \ar[r]^-{\Phi} & 
    \mathcal{O}_{\PP(w)}(w_{0}) \oplus \cdots \oplus 
    \mathcal{O}_{\PP(w)}(w_{n}) \ar[r]^-{\varphi}& 
    T\PP(w) \ar[r] &0} 
  \end{align*} où $\underline{\CC}$ est le fibré orbifold trivial  
  de rang $1$ sur $\PP(w)$. 
\end{lem} 
 
Si les poids sont tous égaux à $1$, ce lemme est bien connu : on peut 
le trouver au paragraphe $3$ du chapitre $3$ de 
\cite{GHag}. D'ailleurs la preuve de ce lemme s'inspire de la 
preuve donnée dans le livre de Griffiths et Harris.  
 
\begin{proof}[Démonstration du lemme 
  \ref{lem:suite,exacte}]
Pour démontrer ce lemme, nous allons d'abord construire les morphismes de 
fibrés puis montrer que la suite est exacte. 
 
 Pour tout $p\in\PP(w)$, soit $(\widetilde{U}_{p},G_{p},\pi_{p})$ une carte  
d'un ouvert $U_{p}$, qui contient $p$, dans l'atlas $\mathcal{A}(|\PP(w)|)$ qui trivialise les fibrés 
$ \mathcal{O}_{\PP(w)}(w_{0}) \oplus \cdots \oplus 
\mathcal{O}_{\PP(w)}(w_{n})$ et $T\PP(w)$. Il existe un unique 
$i\in\{0, \ldots ,n\}$ tel que $\widetilde{U}_{p}\subset 
\widetilde{U}_{i}$. Dans la suite, nous allons supprimer l'indice 
$p$ à la carte $(\widetilde{U}_{p},G_{p},\pi_{p})$ pour ne pas alourdir les notations. 
 
Nous allons d'abord construire un morphisme de fibrés, noté $\varphi$,  entre 
$\mathcal{O}_{\PP(w)}(w_{0}) \oplus \cdots \oplus 
\mathcal{O}_{\PP(w)}(w_{n})$ et $T\PP(w)$. 
Notons $y_{0}, \ldots ,y_{n}$ les coordonnées sur $\widetilde{U}$. Une 
base du fibré tangent à $\widetilde{U}$ est formé par les champs de 
vecteurs $\partial_{y_{k}}:={\partial}/{\partial y_{k}}$ pour $k\neq i$. 
Considérons l'application linéaire $(\id,\widetilde{\varphi}_{\widetilde{U}})$ 
\begin{align*} 
  \widetilde{U}\times \CC^{n+1} 
&\rightarrow \widetilde{U}\times \CC^{n} \\ 
(y,v_{0}, \ldots ,v_{n})& \mapsto \left(y, 
  \left(-\frac{w_{0}}{w_{i}}y_{0}v_{i}+v_{0}\right){\partial_{y_{0}}}+ \cdots +\widehat{i}+ \cdots + \left(-\frac{w_{n}}{w_{i}}y_{n}v_{i}+v_{n}\right){\partial_{y_{n}}}\right) 
\end{align*} 
où $\hat{i}$ signifie que le terme en position $i$ n'appara\^{i}t 
pas. Remarquons que cette application est surjective et que $\ker\widetilde{\varphi}_{\widetilde{U}}(y)=\CC.(y_{0}w_{0}, \ldots ,y_{n}w_{n})$. 
 
Pour construire ce morphisme, nous allons montrer que pour toute injection 
$\alpha:\widetilde{U}\subset\widetilde{U}_{i}\hookrightarrow 
\widetilde{V}\subset\widetilde{U}_{j}$, nous avons  
\begin{align}\label{eq:27} 
  \psi_{\alpha}^{T\PP(w)}\circ\widetilde{\varphi}_{\widetilde{U}}&=\widetilde{\varphi}_{\widetilde{V}}\circ\psi_{\alpha}^{\oplus\mathcal{O}(w_{i})}. 
\end{align} 
D'après les notations \ref{not:injection}, nous avons 
\begin{enumerate} 
\item soit $i=j$ et alors $\alpha(y)=\zeta y$ avec $\zeta$ dans 
  $\bs{\mu}_{i}$, 
\item soit $i\neq j$ et alors  
  \begin{align*} 
    \alpha(y)=(y_{0}/y_{j}^{w_{0}/w_{j}}, \ldots ,1_{j}, \ldots ,y_{n}/y_{j}^{w_{0}/w_{n}}). 
  \end{align*} 
\end{enumerate} 
 
Dans le cas où $i=j$, nous avons  
\begin{align*} 
  \psi_{\alpha}^{T\PP(w)}(y)(\underline{v}) & = 
   (\zeta^{w_{0}}v_{0}, \ldots, \zeta^{w_{n}}v_{n})= \zeta \underline{v} \\ 
 \psi_{\alpha}^{\oplus\mathcal{O}(w_{i})}(y)(\underline{v}) & 
 = \zeta \underline{v} 
\end{align*} où $y=(y_{0}, \ldots ,1_{i}, \ldots ,y_{n})$ et $\underline{v}=(v_{0}, \ldots ,v_{n})$. 
   
Nous en déduisons que  
\begin{align*} 
&\widetilde{\varphi}_{\widetilde{V}}\circ\psi_{\alpha}^{\oplus\mathcal{O}(w_{i})}(y)(\underline{v})  
\\ 
&= \widetilde{\varphi}_{\widetilde{V}}(\zeta y)(\zeta\underline{v}) \\ 
&= 
\left(\left(-\frac{w_{0}}{w_{i}}\zeta^{w_{0}}y_{0}v_{i}+v_{0}\zeta^{w_{0}}\right), \ldots,\widehat{i}, \ldots,  
  \left(-\frac{w_{n}}{w_{i}}\zeta^{w_{n}}y_{n}v_{i}+\zeta^{w_{n}}v_{n}\right)\right)\\ 
&= \zeta\widetilde{\varphi}_{\widetilde{U}}(y)(\underline{v})\\ 
&= \psi_{\alpha}^{T\PP(w)}\circ\widetilde{\varphi}_{\widetilde{U}}(y)(\underline{v}). 
\end{align*} 
 
Dans le cas où $i\neq j$, nous avons  
 
\begin{align*} 
  \psi_{\alpha}^{T\PP(w)}(y)(\underline{v}) & = 
   \left(\frac{\partial t_{k}}{\partial {y_{\ell}}}\right)_{k\neq j, \ell\neq i} \underline{v}\,; \\ 
       \psi_{\alpha}^{\oplus\mathcal{O}(w_{i})}(y)(\underline{v}) & 
 = (1/y_{j}^{1/w_{j}})\underline{v}=\left(v_{0}/y_{j}^{w_{0}/w_{j}}, \ldots ,v_{n}/y_{j}^{w_{n}/w_{j}}\right) 
\end{align*} 
où $t_{k}=y_{k}/y_{j}^{w_{k}/w_{j}}$ pour $k\in\{0, \ldots 
,n\}-\{j\}$  et  
\begin{align*} 
  \left(\frac{\partial t_{k}}{\partial y_{\ell}}\right)_{k\neq j, 
  \ell\neq i} 
\end{align*} 
est une matrice  de taille $n\times n$ (cf. (\ref{eq:25})). 
Nous en déduisons que  
\begin{align*} 
 & 
 \widetilde{\varphi}_{\widetilde{V}}\circ\psi_{\alpha}^{\oplus\mathcal{O}(w_{i})}(y)(\underline{v})\\ 
 & = \widetilde{\varphi}_{\widetilde{V}}(\alpha(y))(1/y_{j}^{1/w_{j}}\underline{v})\\ 
 &= \left(\left(-\frac{w_{0}}{w_{j}}t_{0}\frac{v_{j}}{y_{j}}+\frac{v_{0}}{y_{j}^{w_{0}/w_{j}}}\right)\partial_{t_{0}}+ \cdots + \widehat{j}+\cdots +\left(-\frac{w_{n}}{w_{j}}t_{n}\frac{v_{j}}{y_{j}}+\frac{v_{n}}{y_{j}^{w_{n}/w_{j}}}\right)\partial_{t_{n}}\right). 
\end{align*} 
D'un autre c\^{o}té, nous avons 
\begin{align*} 
&  \psi_{\alpha}^{T\PP(w)}\circ\widetilde{\varphi}_{\widetilde{U}}(y)(\underline{v}) \\&= \psi_{\alpha}^{T\PP(w)}(y)\left(\left(-\frac{w_{0}}{w_{i}}y_{0}v_{i}+v_{0}\right){\partial_{y_{0}}}+ \cdots +\widehat{i}+ \cdots + \left(-\frac{w_{n}}{w_{i}}y_{n}v_{i}+v_{n}\right){\partial_{y_{n}}}\right). 
\end{align*} 
 
D'après l'égalité (\ref{eq:changement}), pour $\ell\neq i$ nous avons  
\begin{align*} 
\ds{{\partial_{y_{\ell}}}}=\begin{cases} 
 \ds{\sum_{\stackrel{k=0}{k\neq j}}^{n} 
-\frac{w_{k}}{w_{j}}\frac{y_{k}}{y_{j}^{\frac{w_{k}}{w_{j}}-1}} 
 \partial_{ t_{k}}} & \mbox{si }  \ell=j\,;\\ 
 \ds{\frac{1}{y_{j}^{w_{\ell}/w_{j}}}}  \partial_{t_{\ell}}   
 &  \mbox{si }  \ell \neq j.  
\end{cases}  
\end{align*} 
 
Le terme devant $\partial_{t_{0}}$ dans 
\begin{align*} 
\psi_{\alpha}^{T\PP(w)}(y)\left(\left(-\frac{w_{0}}{w_{i}}y_{0}v_{i}+v_{0}\right){\partial 
    y_{0}}+ \cdots +\widehat{i}+ \cdots + 
  \left(-\frac{w_{n}}{w_{i}}y_{n}v_{i}+v_{n}\right){\partial 
    y_{n}}\right) 
\end{align*} 
est  
\begin{align} 
&\left( 
  \left(\frac{-w_{0}}{w_{i}}y_{0}v_{i}+v_{0}\right)1/y_{j}^{w_{0}/w_{j}}\right)  
+\left(-\frac{w_{0}}{w_{j}}\frac{y_{0}}{y_{j}^{w_{0}/w_{j}}}\frac{1}{y_{j}}\left(-\frac{w_{j}}{w_{i}}y_{j}v_{i}+v_{j}\right)\right). 
\label{eq:26} 
\\ 
=& \left(-\frac{w_{0}}{w_{j}}t_{0}\frac{v_{j}}{y_{j}}+\frac{v_{0}}{y_{j}^{w_{0}/w_{j}}}\right)\partial_{t_{0}}.\nonumber 
\end{align} 
Le premier terme de (\ref{eq:26}) vient de l'égalité 
$\partial_{y_{0}}=1/y_{j}^{w_{0}/w_{j}}\partial_{t_{0}}$ et le 
deuxième terme de (\ref{eq:26}) vient du terme $k=0$ dans l'égalité 
\begin{align*} 
 \partial_{y_{j}}&=\sum_{\stackrel{k=0}{k\neq j}}^{n} 
-\frac{w_{k}}{w_{j}}\frac{y_{k}}{y_{j}^{\frac{w_{k}}{w_{j}}-1}} 
 \partial_{ t_{k}}.  
\end{align*} 
Nous faisons le m\^{e}me raisonnement pour les termes devant les 
$\partial_{t_{k}}$ pour $k\neq j$ et nous obtenons l'égalité 
(\ref{eq:27}). L'ensemble des applications 
$\widetilde{\varphi}_{\widetilde{U}}$ définit par passage au quotient 
une application $\varphi :\oplus \mathcal{O}(w_{i})\rightarrow T\PP(w)$. 
Comme nous avons l'égalité (\ref{eq:27}) qui est satisfaite pour toute  
injection $\alpha$, nous en déduisons que $\varphi$ est un morphisme surjectif de fibrés 
entre $\oplus \mathcal{O}(w_{i})$ et $T\PP(w)$. 
 
Maintenant, nous allons construire un morphisme injectif de
fibrés, 
noté $\Phi$, entre le 
fibré trivial de rang $1$, noté $\underline{\CC}$, et le fibré 
$\oplus \mathcal{O}(w_{i})$.   
Considérons l'application $(\id,\widetilde{\Phi}_{\widetilde{U}})$ 
\begin{align*} 
  \widetilde{U}\times \CC 
&\longrightarrow \widetilde{U}\times \CC^{n+1} \\ 
(y,v)& \longmapsto \left(y,vy_{0}w_{0}, \ldots ,vy_{n}w_{n}  \right) 
\end{align*} 
Remarquons que $\im\widetilde{\Phi}_{\widetilde{U}}(y)=\CC.(w_{0}v_{0}, \ldots ,w_{n}v_{n})=\ker\widetilde{\varphi}_{\widetilde{U}}(y)$. 
De m\^{e}me que précédemment, pour construire ce morphisme, nous allons montrer que pour toute injection 
$\alpha:\widetilde{U}\subset\widetilde{U}_{i}\hookrightarrow 
\widetilde{V}\subset\widetilde{U}_{j}$, nous avons  
\begin{align}\label{eq:28} 
  \psi_{\alpha}^{\oplus\mathcal{O}(w_{i})}\circ\widetilde{\Phi}_{\widetilde{U}}&=\widetilde{\Phi}_{\widetilde{V}}\circ\psi_{\alpha}^{\underline{\CC}}. 
\end{align} 
Nous séparons les cas $j=i$ et $j\neq i$. 
Dans le cas $j=i$, nous avons 
\begin{align*} 
  \psi_{\alpha}^{\underline{\CC}}(y)({v}) & = v \,; 
 \\ 
 \psi_{\alpha}^{\oplus\mathcal{O}(w_{i})}(y)({v}) & 
 = \zeta {v} &\mbox{ où } v\in \CC. 
\end{align*} 
Un calcul direct montre l'égalité (\ref{eq:28}). 
 
Dans le cas $j\neq i$, nous avons 
\begin{align*} 
  \psi_{\alpha}^{\underline{\CC}}(y)({v}) & =v\,; \\ 
       \psi_{\alpha}^{\oplus\mathcal{O}(w_{i})}(y)({v}) & 
 = (1/y_{j}^{1/w_{j}}){v}=\left(v_{0}/y_{j}^{w_{0}/w_{j}}, \ldots ,v_{n}/y_{j}^{w_{n}/w_{j}}\right). 
\end{align*} 
Nous en déduisons que  
\begin{align*} 
  \widetilde{\Phi}_{\widetilde{V}}\circ\psi_{\alpha}^{\underline{\CC}}(y)(v)&= 
  \widetilde{\Phi}_{\widetilde{V}}(\alpha(y))(v)\\ 
&=\left(\frac{y_{0}}{y_{j}^{w_{0}/w_{j}}}vw_{0}, \ldots ,\frac{y_{n}}{y_{j}^{w_{n}/w_{j}}}vw_{n}\right). 
\end{align*} 
D'un autre c\^{o}té, nous avons 
\begin{align*} 
   \psi_{\alpha}^{\oplus\mathcal{O}(w_{i})}\circ\widetilde{\Phi}_{\widetilde{U}}(y)(v)&= \psi_{\alpha}^{\oplus\mathcal{O}(w_{i})}(y)(vy_{0}w_{0}, \ldots ,vy_{n}w_{n})\\ 
&= \left(\frac{y_{0}}{y_{j}^{w_{0}/w_{j}}}vw_{0}, \ldots ,\frac{y_{n}}{y_{j}^{w_{n}/w_{j}}}vw_{n}\right). 
\end{align*} 
Nous en déduisons l'égalité (\ref{eq:28}). La même raisonnement que 
précédemment nous montre que nous avons un 
morphisme de fibrés $\Phi:\underline{\CC}\rightarrow\oplus 
\mathcal{O}(w_{i})$. 
Comme l'application $\widetilde{\Phi}_{\widetilde{U}}$ est injective, le morphisme  
de fibrés est injectif. 
 
Comme les applications $\widetilde{\varphi}_{\widetilde{U}}$ et 
$\widetilde{\Phi}_{\widetilde{U}}$ sont de rang constant 
($\rg\widetilde{\varphi}_{\widetilde{U}}=n$ et 
$\rg\widetilde{\Phi}_{\widetilde{U}}=1$), la proposition 
\ref{prop:noyau,image,morphisme,fibre} implique que  les fibrés 
$\ker\varphi$ et $\im\Phi$ sont bien définis. 
Comme $\ker\widetilde{\varphi}_{\widetilde{U}}(y)=\im\widetilde{\Phi}_{\widetilde{U}}(y)$ pour  
tout $y\in \widetilde{U}$, nous avons bien une suite exacte de fibrés 
  \begin{align*} 
\xymatrix{0\ar[r]& \underline{\CC} \ar[r]^-{\Phi}& 
    \mathcal{O}_{\PP(w)}(w_{0}) \oplus \cdots \oplus 
    \mathcal{O}_{\PP(w)}(w_{n}) \ar[r]^-{\varphi}& 
    T\PP(w) \ar[r]& 0} 
  \end{align*} 
\end{proof}

Dans le suite de ce chapitre, 
pour tout $i\in\{0, \ldots ,\mu-1\}$, nous posons 
\begin{align}\label{eq:24} 
\eta_{i}&=\eta_{\Inv{s(i)}}^{i-k^{\min}(s(i))}.  
\end{align}

Le lemme  \ref{lem:suite,exacte} et la proposition \ref{prop:suite,exacte,fibre,chern} implique que  
\begin{align*} 
   c(T\PP(w)) &= c(\mathcal{O}_{\PP(w)}(w_{0}) \oplus \cdots \oplus 
  \mathcal{O}_{\PP(w)}(w_{n}) )  \\ & =
  \prod_{i=0}^{n}\left(1+c_{1}(\mathcal{O}_{\PP(w)}(w_{i}))\right) \\
&= \prod_{i=0}^{n}\left(1+w_{i}\eta_{1}\right). 
\end{align*} 
  Nous en déduisons que  
  \begin{align} 
    \label{eq:C1,tangent} 
   c_{1}(T\PP(w))&=\mu\eta_{1}. 
   \end{align}   
 
Soient $\gamma_{1}, \ldots ,\gamma_{k}$ dans $S_{w}$. Quand il n'y a 
pas d'ambiguïté, nous noterons $\underline{\gamma}=(\gamma_{1}, \ldots ,\gamma_{k})$. 
D'après le théorème \ref{thm:dimension}, nous avons 
\begin{align*} 
  \deg \ev_{\ast}\left[\mathcal{M}_{k}(A,\underline{\gamma})\right]=2\left(
    \mu\int_{A}\eta_{1} +n -3 - \sum_{i=1}^{n}a(\gamma_{i}) \right).
\end{align*}

Soient $t_{0}, \ldots ,t_{\mu-1}$ les coordonnées de 
$H^{2\star}_{\orb}(\PP(w),\CC)$ dans la base $\bs{\eta}$.  Posons 
$T:=\sum_{i=0}^{\mu-1}t_{i}\eta_{i}$.  Le potentiel de Gromov-Witten 
orbifold de l'espace projectif à poids $\PP(w)$, noté $F^{GW}$, est 
défini par 
\begin{align*} 
F^{GW}:=\sum_{k\geq 0} \sum_{\stackrel{A\in 
    H_{2}(\PP(w),\ZZ),}{\bs{\gamma}\in S_{w}^{k}}} 
\frac{\Psi^{A}_{k,\bs{\gamma}}(T, \ldots ,T)}{n!}. 
\end{align*}  
 
\begin{lem}\label{lem:forme,potentiel}  
  Le potentiel de Gromov-Witten de $\PP(w)$ est de la forme suivante 
  \begin{align*} 
  F^{GW}=\sum_{\alpha} 
  \Psi^{A}_{|\bs{\alpha}|,\bs{\gamma}}(\eta_{0}^{\otimes\alpha_{0}}, 
  \ldots 
  ,\eta_{\mu-1}^{\otimes\alpha_{\mu-1}})\bs{\frac{t^{\alpha}}{\alpha!}}. 
  \end{align*} 
\end{lem} 
 
 \begin{proof}  
   Comme les invariants de Gromov-Witten sont linéaires en chaque 
   variable, nous obtenons 
   \begin{align*} 
   F^{GW}=\sum_{\alpha} \sum_{\stackrel{A\in 
       H_{2}(\PP(w),\ZZ),}{\bs{\gamma}\in S_{w}^{|\bs{\alpha}|}}} 
   \Psi^{A}_{|\bs{\alpha}|,\bs{\gamma}}(\eta_{0}^{\otimes\alpha_{0}}, 
   \ldots 
   ,\eta_{\mu-1}^{\otimes\alpha_{\mu-1}})\bs{\frac{t^{\alpha}}{\alpha!}}. 
   \end{align*} 
   D'après la formule (\ref{eq:defi,IGW}), l'invariant 
   $\Psi^{A}_{|\bs{\alpha}|,\bs{\gamma}}(\eta_{0}^{\otimes\alpha_{0}}, 
   \ldots ,\eta_{\mu-1}^{\otimes\alpha_{\mu-1}})$ n'a de sens que si 
   $\bs{\gamma}=(\underbrace{s(0), \ldots ,s(0)}_{\alpha_{0} \mbox{ 
       fois }}, \ldots ,\underbrace{s(\mu-1), \ldots 
     ,s(\mu-1)}_{\alpha_{\mu-1} \mbox{ fois }})$.  Ainsi, il est 
   inutile de sommer sur les $\bs{\gamma}$ dans $S_{w}^{k}$ dans la 
   formule du potentiel.  D'après le théorème \ref{thm:dimension} et 
   la formule (\ref{eq:defi,IGW}), si 
   $\Psi^{A}_{|\bs{\alpha}|,\bs{\gamma}}(\eta_{0}^{\otimes\alpha_{0}}, 
   \ldots ,\eta_{\mu-1}^{\otimes\alpha_{\mu-1}})$ est non nul alors 
   nous avons l'égalité\footnote{Il suffit de v\'erifier que $\deg^{\orb}\eta_{i}=2\sigma(i)$.} 
\begin{align}\label{eq:gw,non,nul}\mu\int_{A}\eta_{1}+n-3&=\sum_{i=0}^{\mu-1}\alpha_{i}(\sigma(i)-1). 
 \end{align} 
 Ainsi, si l'on fixe $\alpha$, la classe $A$ est déterminée par 
 l'équation ci-dessus.  Ce qui termine la démonstration. 
 \end{proof} 
  
 D'après la démonstration ci-dessus, nous pouvons omettre l'indice 
 $\bs{\gamma}$ dans la notation 
 $\Psi^{A}_{|\bs{\alpha}|,\bs{\gamma}}(\eta_{0}^{\otimes\alpha_{0}},\ldots,\eta_{\mu-1}^{\otimes\alpha_{\mu-1}})$.

 \begin{rem}  
   Soit $V$ une variété complexe projective. Posons $H^{\star}(V,\CC)=\oplus_{i} 
   \CC \Delta_{i}$ où $\Delta_{i}$ appartient à $H^{\deg(\Delta_{i})}(V,\CC)$. Soit $(t_{i})$ les coordonnées plates, par rapport  
   à la dualité de Poincaré, sur $H^{\star}(V,\CC)$. 
D'après le paragraphe $I.4.4$ du livre de Manin \cite{Mfm}, le champ d'Euler est défini par  
\begin{align*} 
  \mathfrak{E}= 
  \sum_{i}\left(1-\frac{\deg(\Delta_{i})}{2}\right)t_{i}\partial_{t_{i}}+\sum_{b\mid \deg(\Delta_{b})=2} 
  r^{b} \partial_{t_{b}} 
\end{align*} 
où $r_{b}$ est défini par $c_{1}(TV)=\sum_{b\mid \deg{\Delta}_{b}=2} r^{b}\Delta_{b}$. 
 \end{rem} 
  
Dans notre cas, nous avons $\deg^{\orb}(\eta_{i})=2\sigma(i)$ 
(cf. propositions \ref{prop:base} et \ref{prop:spectre}  ) et 
$c_{1}(T\PP(w))=\mu\eta_{1}$ (cf. égalité  (\ref{eq:C1,tangent})).  Nous définissons le champ d'Euler par 
\begin{align}\label{eq:champ,euler,A} 
  \mathfrak{E}&=\sum_{k=0}^{\mu-1}(1-\sigma(k))t_{k}\partial 
_{t_{k}}+\mu\partial_{t_{1}}. 
\end{align} 
 
\begin{prop}\label{prop:Potentiel-de-Gromov} 
  \begin{enumerate} 
  \item Le potentiel est homogène de degré $3-n$ par rapport au champ
    d'Euler modulo l'addition d'un polyn\^{o}me de degr\'e inf\'erieur ou \'egal
    \`a deux. 
  \item  
La matrice $A_{\infty}:=\id-\nabla \mathfrak{E}$ dans la base $\bs{\eta}$ 
est 
\begin{align*} 
  \diag (\sigma(0), \ldots ,\sigma(\mu-1)). 
\end{align*} 
Cette matrice vérifie $A_{\infty}+A_{\infty}^{\ast}=n\id$. 
  \end{enumerate} 
\end{prop} 
 
\begin{proof} 
  \begin{enumerate} 
  \item Le coefficient devant $\bs{t^{\alpha}/\alpha!}$ du potentiel 
    $F^{GW}$ est 
    \begin{align*} 
      \Psi^{A}_{|\bs{\alpha}|}(\eta_{0}^{\otimes\alpha_{0}},\ldots,\eta_{\mu-1}^{\otimes\alpha_{\mu-1}}). 
    \end{align*} 
Le coefficient devant $\bs{t^{\alpha}/\alpha!}$ de $\mathfrak{E}\cdot F^{GW}$ 
est 
\begin{align*} 
  \Psi^{A}_{|\bs{\alpha}|}(\eta_{0}^{\otimes\alpha_{0}},\ldots,\eta_{\mu-1}^{\otimes\alpha_{\mu-1}})\left(\sum_{k=0}^{\mu-1}\alpha_{k}(1-\sigma(k)) 
\right) + \mu \Psi^{A}_{|\bs{\alpha}|+1}(\eta_{0}^{\otimes\alpha_{0}},\eta_{1}^{\otimes\alpha_{1}+1}\ldots,\eta_{\mu-1}^{\otimes\alpha_{\mu-1}}). 
\end{align*} 
 D'après l'axiome du diviseur (cf. théorème $3.4.2$ de l'article 
 \cite{CRogw}), nous avons  
 \begin{align*} 
   \Psi^{A}_{|\bs{\alpha}|+1}(\eta_{0}^{\otimes\alpha_{0}},\eta_{1}^{\otimes\alpha_{1}+1}\ldots,\eta_{\mu-1}^{\otimes\alpha_{\mu-1}})&= \int_{A}\eta_{1} \Psi^{A}_{|\bs{\alpha}|}(\eta_{0}^{\otimes\alpha_{0}},\eta_{1}^{\otimes\alpha_{1}},\ldots,\eta_{\mu-1}^{\otimes\alpha_{\mu-1}})  
 \end{align*} 
Nous en déduisons que le coefficient devant $\bs{t^{\alpha}/\alpha!}$ 
de $\mathfrak{E}F^{GW}$ est  
\begin{align*} 
  \Psi^{A}_{|\bs{\alpha}|}(\eta_{0}^{\otimes\alpha_{0}},\ldots,\eta_{\mu-1}^{\otimes\alpha_{\mu-1}})\left(\sum_{k=0}^{\mu-1}\alpha_{k}(1-\sigma(k)) 
-\mu\int_{A}\eta_{1}\right). 
\end{align*} 
Puis l'égalité (\ref{eq:gw,non,nul}) finit la démonstration. 
\item Les coordonnées $t_{0}, \ldots ,t_{\mu-1}$ sont plates 
  pour la forme bilinéaire non dégénérée $\langle\cdot,\cdot\rangle$. Comme  $\bs{\nabla}$ 
  est la connexion de Levi-Civita associée à 
  $\langle\cdot,\cdot\rangle$, nous avons  
  \begin{align*} 
    \bs{\nabla}_{\partial_{t_{j}}} \partial_{t_{k}}=0. 
  \end{align*} 
  Ainsi, nous obtenons que 
  $\bs{\nabla}_{\partial_{t_{j}}}\mathfrak{E}=(1-\sigma(j))\partial_{t_{j}}$ c'est-à-dire  
    \begin{align*} 
      A_{\infty}=\diag (\sigma(0), \ldots ,\sigma(\mu-1)). 
    \end{align*} 
 La matrice 
 duale $A^{\ast}_{\infty}$ par rapport à la forme bilinéaire non dégénérée 
 $\langle\cdot,\cdot\rangle$ est  
 \begin{align*} 
      A_{\infty}^{\ast}=\diag (\sigma(n), \ldots ,\sigma(0),\sigma(\mu-1), \ldots ,\sigma(n+1)). 
 \end{align*} 
 Pour finir la démonstration il suffit de remarquer que si
 $\overline{j+k}=n$ alors $\sigma(j)+\sigma(k)=n$ o\'u
 $\overline{j+k}$ est la somme  modulo $\mu$. 
  \end{enumerate} 
\end{proof} 
 
 \section{Calcul de certains invariants de Gromov-Witten orbifolds}\label{sec:calcul-de-certains} 
 
Nous rappelons que 
pour chaque $i\in\{0, \ldots ,\mu-1\}$, nous posons 
\begin{align}\label{eq:3} 
\eta_{i}&=\eta_{\Inv{s(i)}}^{i-k^{\min}(s(i))}.  
\end{align} 
  
 Pour calculer les conditions initiales de la variété de Frobenius, il 
 nous reste à calculer la matrice 
 $A_{0}^{\circ}:=\mathfrak{E}\star\mid_{\bs{t}=0}$. Or le champ 
 d'Euler en $\bs{t}=0$ est simplement $\mu\partial t_{1}$. Il faut 
 calculer les invariants de Gromov-Witten 
 $\Psi_{3}^{A}(\eta_{1},\eta_{j},\eta_{k})$ pour tous $j,k$ et pour 
 tout $A \in H_{2}(\PP(w),\ZZ)$.

 \begin{thm}\label{thm:cup} Soient $j,k$ dans $\{0, \ldots 
   ,\mu-1\}$ tels que $\overline{1+j+k}= n$ et 
   $\sigma(1)+\sigma(j)+\sigma(k)= n$ o\'u
 $\overline{1+j+k}$ est la somme modulo $\mu$. Soit $A(j,k)$ la classe dans 
   $H_{2}(|\PP(w)|,\ZZ)$ définie par l'égalité 
   $\int_{A(j,k)}\eta_{1}=(1+j+k-n)/\mu-s(j)-s(k)$.  Nous avons 
   \begin{align*} 
   \Psi^{A}_{3}(\eta_{1},\eta_{j},\eta_{k})=  
\begin{cases} 
\left(\prod_{i\in I(s(j))}w_{i}\right)^{-1} & 
\mbox{ si }  A=A(j,k) \,;\\ 
0 & \mbox{sinon.} 
 \end{cases} 
\end{align*} 
 \end{thm} 
 
 \begin{proof} 
   Les hypothèses sur $j$ et $k$ impliquent que $A(j,k)=0$. Ainsi, 
   nous obtenons 
   \begin{align*}\Psi^A_{3}(\eta_{1},\eta_{j},\eta_{k})= 
\begin{cases} 
(\eta_{1},\eta_{j},\eta_{k}) & \mbox{ si } A=0 \,;\\ 
 0 & \mbox{ sinon.} 
\end{cases}\end{align*} 
Puis le corollaire \ref{cor:tenseur} termine la démonstration. 
 \end{proof} 
 
 \begin{thm}\label{thm:invariant,0}  
   Supposons que $\mu$ et $p_{w}:=\ppcm (w_{0}, \ldots ,w_{n})$ soient 
   premiers entre eux.  Soient $j,k$ dans $\{0, \ldots ,\mu-1\}$ tels 
   que $\overline{1+j+k}\neq n$.  Pour toute classe $A$ dans 
   $H_{2}(|\PP(w)|,\CC)$, nous avons 
   $\Psi^{A}_{3}(\eta_{1},\eta_{j},\eta_{k})=0$. 
 \end{thm} 
 
Avant de démontrer ce théorème, nous allons énoncer un résultat sur   
 les géné\-ra\-teurs de $H_{2}(|\PP(w)|,\ZZ)$. 
 
 Pour tous $i,j$ dans $\{0, 
\ldots ,\mu-1\}$, posons $p_{ij}:=\ppcm(w_{i},w_{j})$ et $r_{ij}:=p_{w}/p_{ij}$. 
Sur la dernière page de l'article \cite{Kcwps}, nous avons le 
diagramme commutatif suivant 
\begin{align*} 
\xymatrix{H_{2}(\PP^1,\ZZ) \ar[rr] \ar[d]^-{{f_{ij}}_{\ast}} && 
  H_{2}(\PP^n,\ZZ) \ar[d]^-{{f_{w}}_{\ast}\ \ \ \ \  \cong}  && 
  \ZZ \ar[rr]^-{\id} \ar[d]^-{\cdot p_{ij}} && \ZZ \ar[d]^-{\cdot p_{w}} 
 \\ H_{2}(|\PP(w_{i},w_{j})|,\ZZ)  
 \ar[rr]^-{{\ \ \ \iota_{ij}}_{\ast}} & &H_{2}(|\PP(w)|,\ZZ)&& \ZZ \ar[rr]^-{\cdot r_{ij}} && \ZZ} 
\end{align*}

\begin{prop}\label{prop:generateur,H2} 
Il existe un générateur, noté $D_{w}$, de $H_{2}(|\PP(w)|,\ZZ)$ tel
que nous ayons  
\begin{align*} 
\int_{D_{w}}\eta_{1}=\frac{1}{p_{w}}. 
\end{align*}  
\end{prop}

\begin{proof} 
Comme nous nous intéressons à des notions topologiques, nous pouvons 
supposer que les poids sont premiers entre eux. 
Nous en déduisons que les nombres $r_{ij}$ sont premiers entre eux 
c'est-à-dire qu'il existe des entiers $\alpha_{ij}$ tels que 
$\sum_{i,j}\alpha_{ij}r_{ij}=1$. 
D'après l'article \cite{Kcwps}, nous avons  
\begin{align*} 
H_{2}(|\PP(w)|,\ZZ)=\sum_{i,j} 
\mbox{im}\left( H_{2}(|\PP(w_{i},w_{j})|,\ZZ) \rightarrow 
  H_{2}(|\PP(w)|,\ZZ) \right). 
\end{align*} 
Nous posons 
$D_{w}:=\sum_{i,j}\alpha_{ij}[|\PP(w_{i},w_{j})|]$. 
Finalement nous obtenons 
\begin{align*} 
\int_{D_{w}}\eta_{1}=\sum_{i,j}\alpha_{ij}\int_{|\PP(w_{i},w_{j})|}\iota_{ij}^{\ast}\eta_{1}=\sum_{i,j}\alpha_{ij}/p_{ij}=\sum_{i,j}\alpha_{ij}r_{ij}/p_{w}=1/p_{w}. 
\end{align*} 
\end{proof}

\begin{proof}[Démonstration du théorème \ref{thm:invariant,0}] 
Nous allons démontrer la contra\-po\-sé. Soient $A\in H_{2}(|\PP(w)|,\ZZ)$ 
et $j,k$ dans $\{0, \ldots ,\mu-1\}$ tels que 
$\Psi_{3}^{A}(\eta_{1},\eta_{j},\eta_{k})$ soit non nul. D'après 
  l'égalité (\ref{eq:gw,non,nul}), nous obtenons 
\begin{align}\label{eq:gw,dimension} 
\int_{A}\eta_{1}&=\frac{1+j+k-n}{\mu} -s(j)-s(k). 
\end{align} 
Comme la classe $A$ est un multiple entier de $D_{w}$, nous en 
déduisons que 
\begin{align*}\frac{p_{w}(1+j+k-n)}{\mu}-p_{w}(s(j)+s(k)) \in \NN.\end{align*} 
Comme $p_{w}$ 
et $\mu$ sont premiers entre eux, nous en concluons que 
$\overline{1+j+k}=n$. 
\end{proof} 
 
  \begin{conj}\label{conj:invariant,dur} 
   Soient $j,k$ dans $\{0, \ldots ,\mu-1\}$ tels que 
   $\overline{1+j+k}= n$ et $\sigma(1)+\sigma(j)+\sigma(k)\neq n$. 
   Soit $A(j,k)$ la classe dans $H_{2}(|\PP(w)|,\ZZ)$ définie par 
   l'égalité $\int_{A(j,k)}\eta_{1}=(1+j+k-n)/\mu-s(j)-s(k)$.  Nous 
   avons 
  \begin{align*} 
   \Psi^{A}_{3}(\eta_{1},\eta_{j},\eta_{k})=  
     \begin{cases} \left(\ds{\prod_{i\in I(s(j))\bigsqcup 
             I(s(k))}w_{i}}\right)^{-1} & \mbox{ si } A=A(j,k)\,; 
       \\ 
       0 & \mbox{sinon.} 
 \end{cases} 
\end{align*} 
 \end{conj} 
 
Posons 
\begin{align*} 
(\!(\eta_{1},\eta_{j},\eta_{k})\!):=\frac{\partial^{3}F^{GW}}{\partial t_{1}{\partial t_{j}}{\partial t_{k}}}\mid_{\bs{t=0}} 
\end{align*}

Les théorèmes \ref{thm:cup} et \ref{thm:invariant,0} et la conjecture 
\ref{conj:invariant,dur} impliquent le corollaire suivant. 
 
\begin{cor}\label{cor:bilan} Supposons que $\mu$ et le plus petit commum multiple des 
  poids soient premiers entre eux.  Soient $j,k$ dans $\{0, \ldots ,\mu-1\} $. 
\begin{enumerate} \item Si $\overline{1+j+k} \neq n$  alors 
  $(\!(\eta_{1},\eta_{j},\eta_{k})\!)=0$. 
\item Si $\overline{1+j+k} = n$ alors nous avons 
\begin{align*}(\!(\eta_{1},\eta_{j},\eta_{k})\!)= 
 \begin{cases} 
   \left({{\prod_{i\in I(j,k)}w_{i}}}\right)^{-1} & \mbox{ si } 
   \sigma(1)+\sigma(j)+\sigma(k)\neq n\,;\\ 
   \left({{\prod_{i\in I(s(j))}w_{i}}}\right)^{-1} & \mbox{ si } 
   \sigma(1)+\sigma(j)+\sigma(k)=n 
 \end{cases} 
\end{align*} 
où $I(j,k):=I(s(j))\bigsqcup  I(s(k))$. 
\end{enumerate} 
\end{cor} 
  
Pour déterminer la matrice $A_{0}^{\circ}=\mathfrak{E}\star 
\mid_{\bs{t}=0}$, nous utilisons la formule suivante 
\begin{align*} 
 \frac{\partial^3 F^{GW}(t_{0},\ldots,t_{\mu-1})}{\partial 
 t_{i}\partial t_{j} \partial t_{k}}&=\langle \partial t_{i}\star 
 \partial t_{j}, \partial t_{k}\rangle.  
\end{align*} 
Cette formule montre que la donnée des $(\!(\eta_{1},\eta_{j},\eta_{k})\!)$ pour 
tout $j,k\in\{0, \ldots ,\mu-1\}$ et la dualité de Poincaré orbifolde 
$\langle \cdot,\cdot\rangle$ nous permettent de calculer la matrice 
$A_{0}^{\circ}$.

\section{Remarques sur la conjecture  \refhautdepage}
Dans ce paragraphe, nous nous plaçons dans les hypothèses de la 
conjecture \ref{conj:invariant,dur} c'est-à-dire que   $j,k$ sont 
dans $\{ 0, \ldots ,\mu-1 \}$ et sont tels que 
$\overline{1+j+k}= n$ et $\sigma(1)+\sigma(j)+\sigma(k)\neq n$. 
Notons $A(j,k)$ la classe dans $H_{2}(|\PP(w)|,\ZZ)$ définie par 
l'égalité $\int_{A(j,k)}\eta_{1}=(1+j+k-n)/\mu-s(j)-s(k)$. 
Ces hypothèses impliquent que  
\begin{align}\label{eq:20} 
  s(k)&=\{1-s(j+1)\}>0, \hspace{0.5cm} A(j,k)\neq 0,\\ 
1+j+k&=n+\mu. \nonumber 
\end{align} 
 
D'après l'axiome du diviseur (cf. théorème $3.4.2.$ de \cite{CRogw}), nous avons 
\begin{align*} 
\Psi^{A(j,k)}_{3}(\eta_{1},\eta_{j},\eta_{k})=\left(\int_{A(j,k)} 
\eta_{1}\right)\Psi_{2}^{A(j,k)}(\eta_{j},\eta_{k}) 
\end{align*} 
 
\begin{lem}\label{lem:dimension,virtuelle,conj} 
  Le degré de la classe 
  $$\ev_{\ast}\left[\overline{\mathcal{M}}_{2}(\PP(w),A(j,k),(\{1-s(j)\},\{1-s(k)\}))\right]$$
  est
\begin{align*} 
  2 \dim_{\CC} \left(\PP(w)_{I(s(j))}\times \PP(w)_{I(s(k))}\right). 
\end{align*} 
\end{lem} 
 
\begin{rem} 
D'après (\ref{eq:23}) et le lemme précédent nous avons    
\begin{align*} 
   \rm{ev}_{\ast}[\overline{\mathcal{M}}_{2}(\PP(w),A(j,k),(\{1-s(j)\},\{1-s(k)\}))]=\cst [\PP(w)_{I(s(j))}\times \PP(w)_{I(s(k))}]. 
\end{align*} 
Pour montrer la conjecture \ref{conj:invariant,dur}, il suffit de 
montrer que cette constante est  
\begin{align*} 
  \left(\int_{A(j,k)} 
\eta_{1}\right)^{-1}=(s(j+1)-s(j))^{-1} \mbox{ d'après les égalités (\ref{eq:20})}. 
\end{align*} 
\end{rem} 
 
\begin{proof}[Démonstration du lemme \ref{lem:dimension,virtuelle,conj}] 
D'après le théorème \ref{thm:dimension}, ce degré est  
\begin{align*} 
  \mu\int_{A(j,k)}\eta_{1}+n-1-a(\{1-s(j)\})-a(\{1-s(k)\}). 
\end{align*} 
D'après les égalités (\ref{eq:20}), nous avons   
\begin{align*} 
  \mu\int_{A(j,k)}\eta_{1}=\mu(s(j+1)-s(j))>0. 
\end{align*} 
Nous en déduisons que  
\begin{align}\label{eq:21} 
  j&=k^{\max}(s(j)), \hspace{0.5cm} j+1=k^{\min}(s(j+1)).  
\end{align} 
D'après l'égalité (\ref{eq:fondamentale}), le degré cherché est 
\begin{align*} 
 2(\mu(s(j+1)-s(j)) -2+\delta(s(j))+a(s(j))-a(s(j+1))). 
\end{align*} 
Or, nous avons  
\begin{align}\label{eq:22} 
\mu(s(j+1)-s(j))&= 1-\sigma(j+1)+\sigma(j). 
 \end{align} 
 Nous utilisons (\ref{eq:21}) dans l'égalité ci-dessus. Finalement, nous déduisons de la 
 proposition \ref{prop:spectre} que la moitié de ce degré   
 est  $$\delta(s(j))-1+ \delta(s(j+1))-1=\dim_{\CC} \PP(w)_{I(s(j))}\times \PP(w)_{I(s(k))}.$$ 
\end{proof} 
 
\begin{rem}Soient 
  $j,k$  dans $\{ 0, \ldots ,\mu-1 \}$ tels que $\overline{1+j+k}= 
  n$ et $\sigma(1)+\sigma(j)+\sigma(k)\neq n$. 
D'après les égalités (\ref{eq:21}), nous en déduisons que 
\begin{align*} 
  j&=k^{\max}(s(j))\, \mbox{ et } \,k = k^{\max}(s(k)). 
\end{align*} 
Si nous revenons aux notations du chapitre précédent 
(cf. (\ref{eq:24})), nous avons 
\begin{align*} 
  \eta_{j}&=\eta_{1-s(j)}^{\delta(s(j))-1} \in 
  H^{2(\delta(s(j))-1)}(|\PP(w)_{I(s(j))}|)\,;\\  \eta_{k}&=\eta_{1-s(k)}^{\delta(s(k))-1} \in 
  H^{2(\delta(s(k))-1)}(|\PP(w)_{I(s(k))}|). 
\end{align*} 
Ainsi, l'invariant de Gromov-Witten $\psi_{2}^{A(j,k)}(\eta_{j},\eta_{k})$ \og compte\fg le 
nombre de courbes de degré $A(j,k)$ qui passent par un point 
 général dans $\PP(w)_{I(s(j))}$ et un point général dans $\PP(w)_{I(s(k))}$. 
\end{rem} 
 
Dans la suite, nous supposons que les poids sont premiers 
entre eux deux à deux.  Dans ce cas, l'espace des lieux singuliers de 
$\PP(w)$ est réduit à $n+1$ points distincts.  Ainsi, toutes les 
applications orbifoldes non constantes sont bonnes de façon unique 
d'après la proposition \ref{prop:application,regulier}. Nous n'avons 
donc plus de problème avec les classes de systèmes compatibles 
définis au paragraphe \ref{sec:IGWO}.

Notons 
\begin{align}
  \label{not:espace,module}
  \overline{\mathcal{M}}_{j,k} 
:=\overline{\mathcal{M}}_{2}(\PP(w),A(j,k),(\Inv{s(j)},\Inv{s(k)})).
\end{align}
 
D'après le lemme \ref{lem:dimension,virtuelle,conj}, nous avons 
\begin{align*} 
\deg \ev_{\ast}\left[\overline{\mathcal{M}}_{j,k}\right]=  
\begin{cases} 
 0 &  \mbox{ si } s(j)\neq 0\,; \\  
2n & \mbox{ si } s(j)=0. 
   \end{cases} 
 \end{align*} 
  
 Le cas où $s(j)=0$ correspond au couple $(j,k)=(n,\mu-1)$. 
 Comme nous avons 
 \begin{align*} 
   \int_{A(n,\mu-1)}\eta_{1}=1/w_{n} 
 \end{align*} 
où $w_{n}$ est le plus 
 grand poids, nous en déduisons que 
\begin{align*} 
  \Psi_{3}^{A(n,\mu-1)}(\eta_{1},\eta_{n},\eta_{\mu-1})&=\int_{A(n,\mu-1)}\eta_{1}\cst\int_{\PP(w)\times\PP(w_{n})}^{\orb}\eta_{n}\wedge\eta_{\mu-1}\\ 
&= \frac{1}{w_{n}}\frac{\cst}{w_{n}\prod_{i=0}^{n}w_{i}}  \mbox{ \  d'après la 
  proposition \ref{prop:integrale}}.  
\end{align*} 
 
Finalement, pour montrer la conjecture \ref{conj:invariant,dur} dans 
le cas $(j,k)=(n,\mu-1)$, il faut montrer que cette constante vaut $w_{n}$.

\begin{prop}\label{prop:ev,surjective} 
Il existe une unique application $f:\PP(1,w_{n})\rightarrow \PP(w)$ holomorphe 
  telle que  
  \begin{itemize} 
  \item l'application $f$ se relève en une application holomorphe de $\CC^{2}-\{0\}$ dans 
  $\CC^{n+1}-\{0\}$ dont les applications composantes sont des 
  polyn\^{o}mes ; 
\item l'application $f$ envoie les points $[1:0]$ et $[0:1]$ sur 
respectivement $[a_{0}:\ldots:a_{n}]\neq [0:\ldots:0:1]$ et
$[0:\ldots:0:1]$ ;  
\item  la classe $f_{\ast}[\PP(1,w_{n})]$ soit $A(n,\mu-1)$. 
 \end{itemize} 
\end{prop} 
 
Avant de démontrer cette proposition, nous allons montrer le lemme suivant. 
       
\begin{lem}\label{lem:pullback,fibre} 
   Soit l'application 
 \begin{align*} 
\begin{array}{rrcl} 
f: &\PP(1,w_{n}) &  \longrightarrow & \PP(w) \\ 
&  [x:y] & \longmapsto & 
[a_{0}x^{w_{0}}:\ldots:a_{n-1}x^{w_{n-1}}:b_{n}y+a_{n} x^{w_{n}}]     
\end{array} 
\end{align*} 
Le fibré $\mathcal{O}_{\PP(1,w_{n})}(1)$ est isomorphe au fibré 
$f^{\ast}\mathcal{O}_{\PP(w)}(1)$. 
\end{lem} 
 
\begin{proof} 
  Comme les poids sont premiers entre eux deux à deux, l'ensemble 
  $f^{-1}(\widehat{\PP(w)}_{\reg})$ est un ouvert dense et connexe. 
  Nous en déduisons d'après les propositions \ref{prop:pullback,fibre} 
  et \ref{prop:application,regulier} que l'image inverse du fibré 
  orbifold $\mathcal{O}_{\PP(w)}(1)$ existe.   D'après 
les hypothèses sur les nombres complexes $a_{0}, \ldots ,a_{n-1}$, il 
existe un indice $i$ tel que $a_{i}\neq 0$. Pour simplifier les 
notations, nous supposerons que $i=0$ c'est-à-dire que $a_{0}\neq 0$. 
Fixons $a_{0}^{1/w_{0}}$ une racine $w_{0}$-ième de $a_{0}$. Remarquons que 
  l'application 
\begin{align*} 
  \PP(1,w_{n})&\longrightarrow \PP(1,w_{n}) \\ 
[x:y] & \longmapsto [x/a_{0}^{1/w_{0}}:b_{n}y+a_{n}x^{w_{n}}] 
\end{align*} 
est un automorphisme de $\PP(1,w_{n})$.  Quitte à composer par  
l'isomorphisme ci-dessus, nous pouvons donc supposer que 
$a_{n}=0$, $a_{0}=1$ et $b_{n}=1$ dans la formule de l'application $f$.  
Notons $\widetilde{U}_{0}$ et $\widetilde{U}_{1}$ les cartes affines 
de $\PP(1,w_{n})$ et $\widetilde{U}_{0}^{w}, \ldots 
,\widetilde{U}_{n}^{w}$ les cartes affines de $\PP(w)$.  De manière 
générale, nous rajouterons un exposant $w$ pour les objets définis sur 
$\PP(w)$.  Les applications suivantes relèvent $f\mid_{U_{0}}$ et 
$f\mid_{U_{1}}$ : 
\begin{align*} 
  \widetilde{f}_{U_{0}U_{0}^{w}} : \widetilde{U}_{0} &\longrightarrow 
  \widetilde{U}_{0}^{w} \\ 
 (1,y_{1})& \longmapsto (1,a_{1}, \ldots 
 ,a_{n-1},y_{1})\\ 
\widetilde{f}_{U_{1}U_{n}^{w}} : \widetilde{U}_{1} &\longrightarrow 
  \widetilde{U}_{n}^{w} \\ 
 (y_{0},1)& \longmapsto (y_{0}^{w_{0}},\ldots,a_{n-1}y_{0}^{w_{n-1}},1)\\ 
\end{align*} 
Comme $f$ est une bonne application orbifolde, nous avons une 
correspondance bijective, notée $\mathfrak{F}$, entre les ouverts d'un  
recouvrement compatible $\mathcal{U}$ de $|\PP(1,w_{n})|$ et les ouverts 
d'un recouvrement compatible $\mathcal{U}^{w}$ de $|\PP(w)|$. 
Soient $U$ et $V$ deux ouverts du recouvrement $\mathcal{U}$ tels que 
$U\subset V$. 
Soit $\alpha:\widetilde{U}\subset\widetilde{U}_{i}\hookrightarrow 
\widetilde{V}\subset\widetilde{U}_{j}$ une injection où 
$i,j\in\{0,1\}$ et $\widetilde{U},\widetilde{V}$ sont deux cartes de 
respectivement $U$ et $V$. Posons  
\begin{align*} 
  k^{\max}(0) &=0 \mbox{ et } k^{\max}(1)=n. 
\end{align*} 
 Il existe une injection 
$\mathcal{F}(\alpha):\widetilde{\mathfrak{F}(U)}^{w}\hookrightarrow 
\widetilde{\mathfrak{F}(V)}^{w}$ où $\widetilde{\mathfrak{F}(U)}^{w}\subset 
\widetilde{U}_{k^{\max}(i)}^{w}$ et 
$\widetilde{\mathfrak{F}(V)}^{w}\subset\widetilde{U}_{k^{\max}(j)}^{w}$ sont 
deux cartes de respectivement $\mathfrak{F}(U)$ et $\mathfrak{F}(V)$ qui 
font commuter le diagramme suivant 
  \begin{align*} 
  \xymatrix{  \widetilde{U}\subset \widetilde{U}_{i} \ar@{^{(}->}[rr]^-{\alpha} 
    \ar[d]_-{\widetilde{f}_{U_{i}U_{k^{\max}(i)}^{w}}\mid_{\widetilde{U}}} && 
    \widetilde{V}\subset \widetilde{U}_{j} \ar[d]^-{\widetilde{f}_{U_{j}U_{k^{\max}(j)}^{w}}\mid_{\widetilde{V}}} \\ 
    \widetilde{\mathfrak{F}(U)}^{w}\subset\widetilde{U}_{k^{\max}(i)}^{w} \ar@{^{(}->}[rr]^-{\mathcal{F}(\alpha)}&& 
    \widetilde{\mathfrak{F}(V)}^{w}\subset\widetilde{U}_{k^{\max}(j)}^{w}} 
  \end{align*} 
D'après la proposition \ref{prop:pullback,fibre}, les fonctions de 
transition de $f^{\ast}\mathcal{O}_{\PP(w)}(1)$ sont défi\-nies par 
\begin{align}\label{eq:33} 
  \psi_{\alpha}^{f^{\ast}\mathcal{O}_{\PP(w)}(1)}(y)&:=\psi_{\mathfrak{F}(\alpha)}^{\mathcal{O}_{\PP(w)}(1)}(\widetilde{f}_{U_{i}U_{k^{\max}(i)}^{w}}\mid_{\widetilde{U}}(y)) & &\forall y\in \widetilde{U}. 
\end{align} 
Nous avons quatre cas à calculer : $(i,j)=(0,0),(1,1),(0,1),(1,0)$. 
Dans le cas où $(i,j)=(0,0)$, l'injection $\alpha$  est simplement 
l'inclusion. Nous en déduisons que $\mathcal{F}(\alpha)$ est aussi 
l'inclusion. Nous obtenons que  
\begin{align*} 
  \psi_{\alpha}^{f^{\ast}\mathcal{O}_{\PP(w)}(1)}(1,y_{1})(v)=v 
\end{align*} 
pour tout 
$(1,y_{1})\in\widetilde{U}$ et pour tout $v\in\CC$. 
 
Dans le cas où $(i,j)=(1,1)$, l'injection $\alpha$ est  
l'action par un élément $\zeta$ de $\bs{\mu}_{w_{n}}$ (cf. les notations 
\ref{not:injection}). Nous en déduisons que $\mathcal{F}(\alpha)$ est  
aussi l'action par cet élément $\zeta$. Nous obtenons que  
\begin{align*} 
  \psi_{\alpha}^{f^{\ast}\mathcal{O}_{\PP(w)}(1)}(y_{0},1)(v)=\zeta v 
\end{align*} 
pour tout $(y_{0},1)\in \widetilde{U}$ et pour tout $v\in\CC$.

Dans le cas où $(i,j)=(0,1)$, l'injection $\alpha$ envoie $(1,y_{1})$ 
sur $(1/y_{1}^{1/w_{n}},1)$ (cf. les notations 
\ref{not:injection}). Nous en déduisons que l'injection 
$\mathcal{F}(\alpha)$ est  
\begin{align*} 
  (1,x_{1}, \ldots ,x_{n})\longmapsto 1/x_{n}^{1/w_{n}}(1,x_{1}, 
\ldots ,x_{n}). 
\end{align*} 
 D'après l'égalité (\ref{eq:33}), nous obtenons que 
\begin{align*} 
 \psi_{\alpha}^{f^{\ast}\mathcal{O}_{\PP(w)}(1)}(1,y_{1})(v)=v/y_{1}^{1/w_{n}}  
\end{align*} 
 pour tout $(1,y_{1})$ dans $\widetilde{U}$ et pour tout $v\in\CC$. 
 
Dans le cas où $(i,j)=(1,0)$, l'injection $\alpha$ envoie $(y_{0},1)$ 
sur $(1,1/y_{0}^{w_{n}})$ (cf. les notations 
\ref{not:injection}). Nous en déduisons que l'injection $\mathcal{F}(\alpha)$ 
 est 
 \begin{align*} 
   (x_{0}, \ldots ,x_{n-1},1)\longmapsto 1/x_{0}^{1/w_{0}}(x_{0}, 
\ldots ,x_{n-1},1). 
 \end{align*} 
 D'après l'égalité (\ref{eq:33}), nous obtenons 
que 
\begin{align*} 
  \psi_{\alpha}^{f^{\ast}\mathcal{O}_{\PP(w)}(1)}(y_{0},1)(v)=v/y_{0} 
\end{align*} 
pour tout $(1,y_{1})$ dans $\widetilde{U}$ et pour tout $v\in \CC$. 
 
Nous retrouvons bien les fonctions de transition du fibré 
$\mathcal{O}_{\PP(1,w_{n})}(1)$ (cf. la démonstration de la proposition \ref{prop:fibre,O(1)}). 
\end{proof}

\begin{proof}[Démonstration de la proposition \ref{prop:ev,surjective}] 
\textbf{Existence:} 
 L'application 
  \begin{align*} 
\begin{array}{rrcl} 
f: &\PP(1,w_{n}) &  \longrightarrow & \PP(w) \\ 
&  [x:y] & \longmapsto & 
[a_{0}x^{w_{0}}:\ldots:a_{n-1}x^{w_{n-1}}:b_{n}y+a_{n} x^{w_{n}}]     
\end{array} 
\end{align*} 
envoie bien les deux points marqués $[1:0]$ et $[0:1]$ sur 
respectivement $[a_{0}:\ldots:a_{n}]$ et $[0:\ldots:0:1]$. 
Comme $[a_{0}:\ldots:a_{n}]$ et $[0:\ldots:0:1]$ sont deux points 
différents, cette application n'est pas constante.  D'après le lemme 
\ref{lem:pullback,fibre}, nous avons  
\begin{align*} 
  \int_{f_{\ast}[\PP(1,w_{n})]}\eta_{1}=\int_{\PP(1,w_{n})}c_{1}(\mathcal{O}_{\PP(1,w_{n})}(1))=1/w_{n}. 
\end{align*} 
 
\textbf{Unicité:} Soit $h:\PP(1,w_{n}) \longrightarrow \PP(w)$ une 
application qui vérifie les trois conditions du lemme. 
Notons une application 
$\widetilde{h}:\CC^{2}-\{0\}\rightarrow \PP(w)$ qui relève $h$. Nous 
avons 
$\widetilde{h}(\lambda\cdot(x,y))=\lambda\cdot\widetilde{h}(x,y)$ 
c'est-à-dire que pour tout $i$ dans $\{0, \ldots ,n\}$, nous avons 
$\widetilde{h}_{i}(\lambda\cdot(x,y))=\lambda^{w_{i}}\widetilde{h}_{i}(x,y)$. 
Supposons que $\widetilde{h}_{i}(x,y)=c_{i}x^{u}y^{v}$ avec $c_{i}$ 
dans $\CC$. Nous obtenons l'égalité 
\begin{align*} 
u+vw_{n}=w_{i}. 
\end{align*} 
Comme $w_{n}>w_{i}$ pour $i\neq n$, nous en déduisons que 
\begin{itemize}  
\item si $i\neq n$ alors nous avons $v=0$ et $u=w_{i}$\,; 
\item si $i=n$ alors soit $v=1$ et $u=0$ soit $u=w_{n}$ et $v=0$. 
\end{itemize} 
Comme l'application $h$ envoie $[1:0]$ et $[0:1]$ sur respectivement 
$[a_{0}:\ldots:a_{n}]$ et $[0:\ldots:0:1]$, nous en déduisons que 
$h=f$. 
\end{proof} 
 
\begin{rem}\label{rem:espace,module} 
  \begin{enumerate} 
  \item \label{item:24} 
Si $[a_{0}:\ldots:a_{n}]\neq[0:\ldots:0:1]$, 
 l'application  
\begin{align*} 
\begin{array}{rrcl} 
f: &\PP(1,w_{n}) &  \longrightarrow & \PP(w) \\ 
&  [x:y] & \longmapsto & 
[a_{0}x^{w_{0}}:\ldots:a_{n-1}x^{w_{n-1}}:b_{n}y+a_{n} x^{w_{n}}]     
\end{array} 
\end{align*} est dans l'espace de modules $\mathcal{M}_{n,\mu-1}$
(cf. notation \eqref{not:espace,module}).

\item Nous allons décrire une autre famille d'applications dans 
$\mathcal{M}_{n,\mu-1}$. Soit 
$C:=\PP(1,w_{n})\cup\PP^{1}_{w_{n},w_{n}}$ la courbe nodale dont le 
point nodal de $\PP(1,w_{n})$ (resp. $\PP^{1}_{w_{n},w_{n}}$ 
cf. exemple \ref{expl:orbifold}.(\ref{item:17})) et 
$[0:1]$ (resp. $[0:1]$) (cf. figure \ref{fig:courbe2} ). Les deux points marqués sur $C$  
sont les points $z_{1}:=[1:0]$ est $z_{2}:=[a:b]\neq[1:0],[0:1]$ sur $\PP^{1}_{w_{n},w_{n}}$. 
Nous définissons une application $g: C\rightarrow \PP(w)$ sur chacune 
de ses composantes. 
La composante  $\PP^{1}_{w_{n},w_{n}}$ est envoyée sur le point 
$[0:\ldots:0:1]$ c'est-à-dire que $g$ est constante sur cette 
composante. 
\begin{figure}[tbhp] 
\begin{center} 
 \psfrag{x}{$z_{1}$}  
\psfrag{y}{$z_{2}$} 
\psfrag{z}{$[1:0]$} 
\psfrag{A}{$\PP(1,w_{n})$} 
\psfrag{B}{$\PP^{1}_{w_{n},w_{n}}$} 
\psfrag{a}{$\bs{\mu}_{w_{n}}$}  
\psfrag{b}{$\id$} 
\psfrag{c}{$\bs{\mu}_{w_{n}}$} 
  \includegraphics[width=0.4\linewidth,height=0.4\textheight,keepaspectratio]{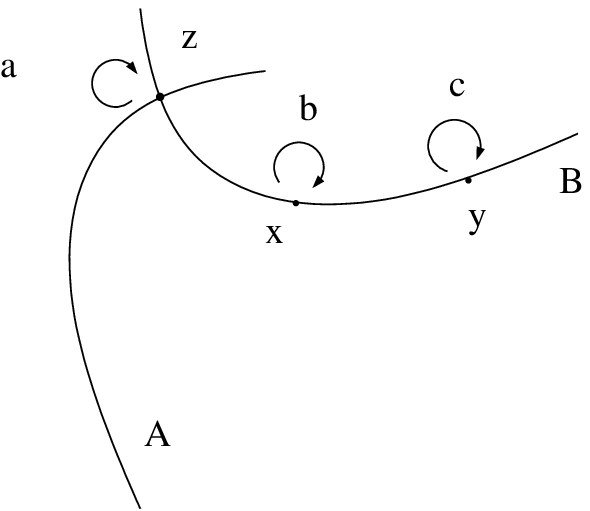} 
\end{center} 
\caption{}\label{fig:courbe2} 
\end{figure}    
 
Sur la composante $\PP(1,w_{n})$, l'application $g$ est définie par  
\begin{align*} 
\begin{array}{rrcl} 
 &\PP(1,w_{n}) &  \longrightarrow & \PP(w) \\ 
&  [x:y] & \longmapsto & 
[a_{0}x^{w_{0}}:\ldots:a_{n-1}x^{w_{n-1}}:b_{n}y+a_{n} x^{w_{n}}]     
\end{array} 
\end{align*} où $[a_{0}:\ldots:a_{n}]\neq[0:\ldots:0:1]$. 
L'application stable $(C,g)$ est dans l'espace de modules $\mathcal{M}_{n,\mu-1}$. 
\end{enumerate} 
\end{rem} 
 
Pour les applications dans $\mathcal{M}_{n,\mu-1}$ qui n'ont qu'une 
seule composante c'est-à-dire celles qui sont comme dans la remarque 
\ref{rem:espace,module}.(\ref{item:24}), nous avons le résultat de 
convexité suivant. 
 
\begin{lem}\label{lem:f,convexite} Soit $(a_{0}, \ldots ,a_{n})\in 
  \CC^{n+1}-\{0\}$ tel que $[a_{0}: \ldots :a_{n}]\neq 
  [0:\ldots:0:1]$ dans $\PP(w)$. Soit $b_{n}$ un nombre complexe non nul. 
  Soit l'application \begin{align*} 
\begin{array}{rrcl} 
f: &\PP(1,w_{n}) &  \longrightarrow & \PP(w) \\ 
&  [x:y] & \longmapsto & 
[a_{0}x^{w_{0}}:\ldots:a_{n-1}x^{w_{n-1}}:b_{n}y+a_{n} x^{w_{n}}]     
\end{array} 
\end{align*} 
 Alors nous avons $H^{1}(\PP(1,w_{n}),f^{\ast}\Theta_{|\PP(w)|})=0$ où  
 $\Theta_{|\PP(w)|}$ est le faisceau des sections du fibré tangent $T\PP(w)$. 
\end{lem}

\begin{proof}[Démonstration du lemme \ref{lem:f,convexite}] 
D'après le lemme \ref{lem:suite,exacte}, nous avons la suite exacte de  
faisceaux \begin{align*} 
    0\rightarrow \mathcal{O}_{|\PP(w)|}\rightarrow 
    \mathcal{O}_{\PP(w)}(w_{0}) \oplus \cdots \oplus \mathcal{O}_{\PP(w)}(w_{n}) \rightarrow 
    \Theta_{|\PP(w)|}\rightarrow 0. 
  \end{align*} 
Comme l'application orbifolde $f:\PP(1,w_{n})\rightarrow \PP(w)$ est 
continue entre les espaces topologiques sous-jacents, nous en déduisons une suite exacte   
\begin{align*} 
    0\rightarrow f^{\ast}\mathcal{O}_{|\PP(w)|}\rightarrow 
    f^{\ast}\mathcal{O}_{\PP(w)}(w_{0}) \oplus \cdots \oplus f^{\ast}\mathcal{O}_{\PP(w)}(w_{n}) \rightarrow 
    f^{\ast}\Theta_{|\PP(w)|}\rightarrow 0. 
  \end{align*} 
Nous avons les égalités entre les faisceaux 
\begin{align*} 
  f^{\ast}\mathcal{O}_{|\PP(w)|}&=\mathcal{O}_{|\PP(1,w_{n})|}\,;\\ 
f^{\ast}\mathcal{O}_{\PP(w)}(w_{k})&=\mathcal{O}_{\PP(1,w_{n})}(w_{k})  
 \mbox{\ \ \ pour  } k \in\{0, \ldots ,n \} \ (\mbox{cf. lemme \ref{lem:pullback,fibre}}). 
\end{align*} 
Comme $ H^{2}\left(|\PP(1,w_{n})|,\mathcal{O}_{|\PP(1,w_{n})|}\right)= 
H^{2}\left(\PP^{1},\mathcal{O}_{\PP^{1}}\right)=0$ (cf. l'exemple 
\ref{expl:espace,annele} et la remarque 
\ref{rem:groupe,isotropie}.(\ref{rem:item:epp,P,twist})), nous en 
déduisons une longue suite exacte en cohomologie 
\begin{align*} 
&  0  
 \rightarrow 
  H^{0}\left(\PP(1,w_{n}),\mathcal{O}_{|\PP(1,w_{n})|}\right)  
\rightarrow 
   H^{0}\left(\PP(1,w_{n}),\mathcal{O}_{\PP(1,w_{n})}(w_{0})\oplus \cdots 
   \oplus \mathcal{O}_{\PP(1,w_{n})}(w_{n})\right) 
\\& \rightarrow  
 H^{0}(\PP(1,w_{n}),f^{\ast}T\PP(w))  
\rightarrow  
 H^{1}\left(\PP(1,w_{n}),\mathcal{O}_{|\PP(1,w_{n})|}\right)\\& 
\rightarrow 
   H^{1}\left(\PP(1,w_{n}),\mathcal{O}_{\PP(1,w_{n})}(w_{0})\oplus \cdots 
   \oplus \mathcal{O}_{\PP(1,w_{n})}(w_{n})\right) 
\rightarrow  
 H^{1}\left(\PP(1,w_{n}),f^{\ast}T\PP(w)\right)  
\rightarrow  
0. 
\end{align*} 
Pour démontrer le lemme, le suffit de montrer que 
$H^{1}\left(\PP(1,w_{n}),\mathcal{O}_{\PP(1,w_{n})}(1)\right)=0$. 
Nous allons utiliser la cohomologie de \v{C}ech.  
Nous recouvrons $|\PP(1,w_{n})|$ avec les ouverts 
$U_{0}:=\{[y_{0},y_{1}]\in|\PP(1,w_{n})|\mid y_{0}\neq 0\}$ et 
$U_{1}:=\{[y_{0},y_{1}]\in|\PP(1,w_{n})|\mid y_{1}\neq 0\}$. 
 
Remarquons que  
\begin{itemize} 
\item le faisceau  
$\mathcal{O}_{\PP(1,w_{n})}(1)\mid_{U_{0}}$ est isomorphe au 
faisceau des fonctions holomorphes sur $\CC$ ; 
\item le faisceau  
$\mathcal{O}_{\PP(1,w_{n})}(1)\mid_{U_{0}\cap U_{1}}$ est isomorphe au 
faisceau des fonctions holomorphes sur $\CC^{\ast}$ ; 
\item le faisceau  
$\mathcal{O}_{\PP(1,w_{n})}(1)\mid_{U_{1}}$ est isomorphe au 
faisceau $\left(\pi_{\ast}\mathcal{O}_{\CC}\right)^{\bs{\mu}_{w_{n}}}$  
où $\pi:\CC\rightarrow \CC$ est l'application qui à $z$ associe 
$z^{w_{n}}$. Le corollaire \ref{cor:coho,faisceau} montre que  
\begin{align*} 
  H^{i}\left(\CC/\bs{\mu}_{w_{n}},\left(\pi_{\ast}\mathcal{O}_{\CC}\right)^{\bs{\mu}_{w_{n}}}\right)= 
H^{i}\left(\CC,\mathcal{O}_{\CC}\right)^{\bs{\mu}_{w_{n}}}=0 \mbox{ 
  pour } i>0. 
\end{align*} 
\end{itemize} 
Le recouvrement $U_{0},U_{1}$ est bien acyclique. 
Soit  $s$ une section de $\mathcal{O}_{\PP(1,w_{n})}(1)(U_{0}\cap U_{1})$. 
La carte $(\widetilde{U}_{0},\id,\pi_{0})$ de $U_{0}$ induit une carte  
$(\widetilde{U_{0}\cap U_{1}}^{0},\id,\pi_{0})$ de $U_{0}\cap U_{1}$. 
De m\^{e}me, la carte $(\widetilde{U}_{1},\bs{\mu}_{w_{n}},\pi_{1})$ de $U_{1}$ induit une carte  
$(\widetilde{U_{0}\cap U_{1}}^{1},\bs{\mu}_{w_{n}},\pi_{1})$ de $U_{0}\cap U_{1}$. 
Notons $y_{0}$ (resp. $y_{1}$) la coordonnée sur $\widetilde{U_{0}\cap  
  U_{1}}^{1}$ 
(resp. $\widetilde{U_{0}\cap U_{1}}^{0}$). Sur $U_{0}\cap U_{1}$, nous avons  
\begin{align*} 
  [1:y_{1}]=[1/y_{1}^{w_{n}}:1]. 
\end{align*} 
C'est-à-dire $y_{0}^{w_{n}}y_{1}=1$. 
La section $s\in\mathcal{O}_{\PP(1,w_{n})}(1)(U_{0}\cap U_{1})$ se relève en une application 
\begin{align*} 
  \widetilde{s}:\widetilde{U_{0}\cap U_{1}}^{1} & \longrightarrow \CC \\ 
y_{0 }& \longmapsto  y_{0}\sum_{p\in \ZZ} a_{p} y_{0}^{w_{n}p} 
\end{align*} 
Posons $\widetilde{s}_{1}(y_{0}):=y_{0}\sum_{p\geq 0} a_{p} 
y_{0}^{w_{n}p}$. L'application $\widetilde{s}_{1}$ est un relevé d'un 
élément, noté $s_{1}$, de $\mathcal{O}_{\PP(1,w_{n})}(1)(U_{1})$. 
Posons $\widetilde{s}_{0}(y_{1}):= \sum_{p> 0} a_{-p} y_{1}^{p}$. 
L'application $\widetilde{s}_{0}$ est un relevé d'un élément, noté 
$s_{0}$, de $\mathcal{O}_{\PP(1,w_{n})}(1)(U_{0})$.  Montrons que 
$s=s_{0}\mid_{U_{0}\cap U_{1}}+s_{1}\mid_{U_{0}\cap U_{1}}$. 
  
Soit $\widetilde{h}: \widetilde{U_{0}\cap U_{1}}^{1}\rightarrow 
\widetilde{U_{0}\cap U_{1}}^{0}$ l'application qui à $y_{0}$ associe $1/y_{0}^{w_{n}}$. 
Nous avons le diagramme commutatif 
\begin{align*} 
\xymatrix{ \widetilde{U_{0}\cap U_{1}}^{1} \times \CC \ar[d] 
  \ar[r]^-{(\widetilde{h},\psi_{\widetilde{h}})} & 
  \widetilde{U_{0}\cap U_{1}}^{0} \times 
  \CC \ar[d] \\ 
  \widetilde{U_{0}\cap U_{1}}^{1} \ar[r]^-{\widetilde{h}} \ar[d]_-{\pi_{1}}& \widetilde{U_{0}\cap U_{1}}^{0}\ar[dl]^-{\pi_{0}} \\ 
  U_{0}\cap U_{1}} 
\end{align*} 
où $\psi_{\widetilde{h}}:\widetilde{U_{0}\cap U_{1}}^{1}\rightarrow GL(1,\CC)$ 
est l'application qui à $y_{0}$ associe la multiplication par 
$1/y_{0}$.  Le diagramme ci-dessus montre que si nous avons un relevé 
$\widetilde{s}'_{1}:\widetilde{U_{0}\cap U_{1}}^{1}\rightarrow \CC$ d'une 
section $s'_{1}$ de fibré $\mathcal{O}_{\PP(1,w_{n})}(1)\mid_{U_{0}\cap 
  U_{1}}$, nous en déduisons un autre relevé 
$\widetilde{s}'_{0}:\widetilde{U_{0}\cap U_{1}}^{0}\rightarrow \CC$ définie par  
\begin{align}\label{eq:4} 
  y_{1} &\longmapsto \psi_{\widetilde{h}}(\widetilde{s}'_{1}(1/y_{1}^{1/w_{n}})) 
\end{align} 
où $1/y_{1}^{1/w_{n}}$ est dans $\widetilde{h}^{-1}(y_{1})$. 
L'application $\widetilde{s}'_{0}$ est bien définie car la formule 
$(\ref{eq:4})$ ne dépend pas du choix de  l'élément dans 
$\widetilde{h}^{-1}(y_{1})$. 
Le calcul ci-dessous montre que le relevé $\widetilde{s}_{1}$ induit le 
relevé $\widetilde{s}_{0}$ :  
\begin{align*} 
  \widetilde{s}_{0}(\widetilde{h}(y_{0}))=\psi_{\widetilde{h}}(y_{0})\left(y_{0}\sum_{p\geq 0} a_{-p} y_{0}^{-w_{n}p}\right)=\sum_{p\geq 0} a_{-p} y_{0}^{-w_{n}p}. 
\end{align*} 
Nous en déduisons que  
\begin{align*} 
  \widetilde{s}_{0}\circ\widetilde{h}+\widetilde{s}_{1}=\widetilde{s} 
\end{align*} 
sur $\widetilde{U_{0}\cap U_{1}}^{1}$. 
Finalement, nous obtenons $s=s_{0}\mid_{U_{0}\cap 
  U_{1}}+s_{1}\mid_{U_{0}\cap U_{1}}$ c'est-à-dire que nous avons \begin{align*}H^{1}\left(\PP(1,w_{n}),\mathcal{O}_{\PP(1,w_{n})}(1)\right)=0.\end{align*} 
\end{proof}

   
\chapter{Structure de Frobenius associée au polyn\^{o}me de Laurent $f$}\label{cha:les-singularites-du} 
 
Nous gardons les notations du chapitre \ref{cha:prel-comb}.  Soit 
$U :=\{ (u_{0}, \ldots ,u_{n})\in \CC^{n+1}\mid u_{0}^{w_{0}}\cdots 
u_{n}^{w_{n}}=1\}$. Soit $f:U\rightarrow \CC$ la fonction définie par 
$f(u_{0},\ldots,u_{n})=u_{0}+\cdots+u_{n}$.  Le polyn\^{o}me $f$ n'est 
pas exactement celui considéré dans \cite{DSgm2} mais nous pouvons 
appliquer les m\^{e}mes méthodes. 
 
Un calcul élémentaire montre que les points critiques de $f$ sont les 
points 
\begin{align*} 
\zeta\left({\prod_{i=0}^{n} 
    w_{i}^{w_{i}}}\right)^{-1/\mu}(w_{0}, \ldots ,w_{n})\in U 
\end{align*} 
où $\zeta$ est une racine $\mu$-ième de l'unité.  Les valeurs 
critiques de $f$ sont les nombres complexes 
\begin{align*} 
{\mu\zeta}\left({\prod_{i=0}^{n}w_{i}^{w_{i}}}\right)^{-1/\mu}. 
\end{align*}

Nous verrons, en utilisant l'article \cite{DSgm1}, qu'il existe 
une structure de Frobenius canonique sur la base de son déploiement 
universel. 
 
Soit $A_{0}^{\circ}$ la matrice $\mu\times\mu$ telle que les seuls éléments 
non nuls soient en position $(\overline{j+1},j)$ (rappelons que 
$\overline{j+1}$ est la somme  modulo $\mu$) et 
\begin{align*} 
(A_{0})_{\overline{j+1},j}= 
\begin{cases}  
\mu & \mbox{ si }  s(\overline{j+1})=s(j)\,; \\ 
{\mu}\left({\prod_{i\in I(s(j))} w_{i}}\right)^{-1} & \mbox{sinon}. 
 \end{cases} 
\end{align*} 
Les valeurs propres de $A_{0}^{\circ}$ sont exactement les valeurs 
critiques de $f$. Ainsi, $A_{0}^{\circ}$ est une matrice semi-simple 
régulière. 
 
Soit $A_{\infty}$ la matrice $\mu\times\mu$ définie par 
\begin{align*} 
A_{\infty}=\diag (\sigma(0), \ldots ,\sigma(\mu-1)). 
\end{align*}

Dans la base canonique $(e_{0}, \ldots ,e_{\mu-1})$ de $\CC^{\mu}$, 
nous définissons la forme bilinéaire non dégénérée $g$ par 
 
  \begin{align*} 
  g(e_{k},e_{\ell}) & =  \begin{cases} \left(\prod_{i\in I(s(k))} 
      w_{i}\right)^{-1} & \mbox{ si } \overline{k+\ell}=n\,; \\ 0 & 
    \mbox{sinon.}  \end{cases} 
    \end{align*} 
     
    La matrice $A_{\infty}$ vérifie $A_{\infty}+ {\ }^{t}\!A_{\infty} 
    =n\cdot{\id}$ où la transposée est définie par rapport à la forme 
    bilinéaire $g$. 
     
    Les données $(A_{0}^{\circ},A_{\infty},g,e_{0})$ définissent 
    (cf. théorème \ref{thm:iso,frobenius,dubrovin}) un unique germe de 
    structure de Frobenius semi-simple au point 
 
    \begin{align*}\left({\mu}\left({\prod_{i=0}^{n}w_{i}^{w_{i}}}\right)^{-1},{\mu\zeta}\left({\prod_{i=0}^{n}w_{i}^{w_{i}}}\right)^{-1}, 
      \ldots 
      ,{\mu\zeta^{\mu-1}}\left({\prod_{i=0}^{n}w_{i}^{w_{i}}}\right)^{-1}\right)
\end{align*} 
    de $\CC^{\mu}$. 
 
\begin{thm}\label{thm:structure,frobenius} 
  La structure de Frobenius canonique de tout germe de 
   déploiement universel du polyn\^{o}me de Laurent $f(u_{0}, \ldots 
  ,u_{n})=u_{0}+\cdots+u_{n}$ sur $U$ est isomorphe au germe de la 
  structure de Frobenius semi-simple définie par les conditions 
  initiales $(A_{0}^{\circ},A_{\infty},e_{0},g)$ au point 
  \begin{align*}\left({\mu}\left({\prod_{i=0}^{n}w_{i}^{w_{i}}}\right)^{-1}, 
    {\mu\zeta}\left({\prod_{i=0}^{n}w_{i}^{w_{i}}}\right)^{-1}, \ldots 
    ,{\mu\zeta^{\mu-1}}\left({\prod_{i=0}^{n}w_{i}^{w_{i}}}\right)^{-1} 
  \right)
\end{align*} 
  de $\CC^{\mu}$. 
\end{thm}

Au paragraphe suivant, nous allons résoudre le problème de 
Birkhoff au point $0$ de l'espace des paramètres d'un déploiement 
universel de $f$. Cette solution sera canonique d'après les résultats 
du paragraphe $5$ de \cite{DSgm2}.  Puis, nous allons calculer les 
conditions initiales de la variété de Frobenius. 

\section[Le système de Gauss-Manin et le réseau de Brieskorn de $f$]{Le système de Gauss-Manin et le réseau de Brieskorn associés au polyn\^{o}me $f$} \label{sec:le-systeme-de} 
Soit l'ensemble 
 $U :=\{ (u_{0}, \ldots ,u_{n})\in \CC^{n+1}\mid 
u_{0}^{w_{0}}\cdots u_{n}^{w_{n}}=1\}$. Soit $f:U\rightarrow \CC$ la 
fonction définie par 
$f(u_{0},\ldots,u_{n})=u_{0}+\cdots+u_{n}$.  
Dans l'article de Douai et Sabbah \cite{DSgm2}, le polyn\^{o}me n'est 
pas exactement le m\^{e}me et les poids sont premiers entre eux dans 
leur ensemble.  
Nous aurons deux types de modifications : 
l'une pour le polyn\^{o}me et l'autre pour tenir compte du pgcd des poids.  
 
Notons $d$ le plus grand 
diviseur commun des entiers $w_{0}, \ldots ,w_{n}$. 
Pour $\alpha=0, \ldots ,d-1$, nous posons 
\begin{align*} 
  U_{\alpha} :=\{ (u_{0,\alpha}, \ldots ,u_{n,\alpha})\in \CC^{n+1}\mid 
u_{0,\alpha}^{w_{0}/d}\cdots u_{n,\alpha}^{w_{n}/d}=\zeta^{\alpha}\} 
\end{align*} 
où $\zeta:=\exp(2i\pi/d)$. 
Nous avons 
\begin{align*} 
  U:= \bigsqcup_{\alpha=0}^{d-1}U_{\alpha}. 
\end{align*} 
Chaque $U_{\alpha}$ est isomorphe à un tore complexe
$(\CC^{\star})^{n}$ et la restriction, noté $f_{\alpha}=u_{0,\alpha}+ \ldots +u_{n,\alpha}$, de $f$ à
$U_{\alpha}$ est un polyn\^{o}me de Laurent. Dans des coordonnées
convenable, son polyèdre de Newton, enveloppe convexe dans $\QQ^{n}$
des exposants de ses mon\^{o}mes, est un simplexe contenant l'origine
dans son intérieur. On peut montrer comme dans \cite{DSgm2} (cf. lemme
$1.2$) qu'il est commode et non dégénéré\footnote{Les définitions ci-dessous viennent de
  l'article de  Kouchnirenko \cite{Knp}. Soit
  $g\in\CC[v,v^{-1}]=\CC[v_{1},v_{1}^{-1}, \ldots
  ,v_{n},v_{n}^{-1}]$. Le polyn\^{o}me de Laurent $g$ s'écrit
  $g=\sum_{i\in\ZZ^{n}}a_{i}v^{i}$. Notons $\rm{Supp}(g):=\{i\in\ZZ^{n}\mid
  a_{i}\neq 0\}$. Notons $\Gamma(g)$ l'enveloppe convexe de
  $\rm{Supp}(g)-\{0\}$. Le polyn\^{o}me de  Laurent $g$ est
  \emph{commode} si le point $0$ n'appartient à aucune des faces
  de dimension $i$ pour $1\leq i\leq n-1$ de $\Gamma(g)$. Le polyn\^{o}me
  de  Laurent $g$ est \emph{non dégénéré} par rapport à  $\Gamma(g)$
  si  pour chaque face $\Delta$ de $\Gamma(g)$ le polynome de Laurent
  $\sum_{i\in\ZZ^{n}\cap\Delta}a_{i}v^{i}$ n'a pas de point critique
  sur $(\CC^{\star})^{n}$.}.
 
Nous rappelons les notations et les principaux résultats des articles 
de Douai et Sabbah \cite{DSgm1} et \cite{DSgm2} adaptés à notre 
polyn\^{o}me de Laurent $f$. 
 
Le système de Gauss-Manin de $f_{\alpha}$ est défini par  : 
\begin{align*} 
G(f_{\alpha}):= \Omega^{n}(U_{\alpha})[\theta,\theta^{-1}]/(\theta d - d 
f_{\alpha}\wedge)\Omega^{n-1}(U_{\alpha})[\theta,\theta^{-1}]. 
\end{align*} 
 
Le système de Gauss-Manin de $f$ est défini par  : 
\begin{align*} 
G:= \Omega^{n}(U)[\theta,\theta^{-1}]/(\theta d - d 
f\wedge)\Omega^{n-1}(U)[\theta,\theta^{-1}]. 
\end{align*} 
  
Comme $U$ a $d$ composantes connexes, nous avons 
\begin{align}\label{eq:9} 
  G&= \bigoplus_{\alpha=0}^{d-1} G(f_{\alpha}).  
\end{align} 
 
Le réseau de Brieskorn de $f_{\alpha}$, 
défini par $G_{0}(f_{\alpha}) := Im(\Omega^{n}(U_{\alpha})[\theta]\rightarrow G(f_{\alpha}))$, est un 
$\CC[\theta]$-module libre de rang $\mu/d$. 
 
Le réseau de Brieskorn de $f$, 
défini par $G_{0} := Im(\Omega^{n}(U)[\theta]\rightarrow G)$, est un 
$\CC[\theta]$-module libre de rang $\mu$. Nous avons 
\begin{align}\label{eq:10} 
   G_{0}&= \bigoplus_{\alpha=0}^{d-1} G_{0}(f_{\alpha}). 
\end{align}

On note $V_\cdot G(f_{\alpha})$ la filtration de Malgrange Kashiwara de 
$G(f_{\alpha})$ (cf.  le paragraphe $2$.e de \cite{DSgm1}). Cette 
filtration est croissante et exhaustive.  La filtration 
$V_{\bullet}G(f_{\alpha})$ induit une filtration sur le réseau de Brieskorn 
définie par $V_{\beta}G_{0}(f_{\alpha}) :=V_{\beta}G(f_{\alpha})\cap 
G_{0}(f_{\alpha})$.  D'après les égalités (\ref{eq:9}) et (\ref{eq:10}), la  
somme directe des filtrations $V_{\bullet}G(f_{\alpha})$ est une filtration croissante et 
exhaustive sur $G$. Nous la notons $V_\cdot G$. Nous en déduisons une filtration notée 
$V_{\bullet}G_{0}$ sur le réseau de Brieskorn $G_{0}$. 
 
Soit $\omega_{0,\alpha}$ la $n$-forme  
\begin{align*} 
\frac{\frac{d{u_{0}}}{u_{0}}\wedge \cdots \wedge \frac{d 
    u_{n}}{u_{n}}}{d\left(\prod_{i}u_{i}^{w_{i}}\right)}\left|_{\prod_{i}u_{i}^{w_{i}/d}=\zeta^{\alpha}} 
\right.  
\end{align*} sur $U_{\alpha}$.
Nous posons  
\begin{align*} 
  \omega_{0}:=(\omega_{0,0},\ldots,\omega_{0,d-1}). 
\end{align*} 
La $n$-forme $\omega_{0}$ est définie sur $U$.

On définit par récurrence la suite $(\underline{a}_{w}(k),i_{w}(k))\in \NN^{n+1}\times 
\{0, \ldots ,n\}$ par 
\begin{align*} 
\underline{a}_{w}(0)&=(0, \ldots ,0),  & i_{w}(0)&=0, \\ 
\underline{a}_{w}(k+1)&=\underline{a}_{w}(k)+1_{i_{w}(k)},  & i_{w}(k+1)&=\min\{j\mid \underline{a}_{w}(k+1)_{j}/w_{j}\}.
\end{align*} 
 
Notons $w/d:=(w_{0}/d, \ldots ,w_{n}/d)$. Comme nous avons  
\begin{align}\label{eq:12} 
  i_{w/d}(\cdot)&=i_{w}(\cdot), & \underline{a}_{w/d}(\cdot)=\underline{a}_{w}(\cdot). 
\end{align} 
nous supprimons l'indice dans les notations $\underline{a}_{w}$ et $i_{w}$. 
 
Pour tout $k$, on a 
$|\underline{a}(k)| :=\sum_{i}\underline{a}(k)_{i}=k$.  Pour tout 
$(q,r)\in \{0,\ldots,d-1\}\times\{0, \ldots ,(\mu/d)-1\}$, nous avons 
les égalités suivantes : 
\begin{align}\label{eq:a} 
i\left(\frac{q\mu}{d}+r\right)&=i(r), &  
\underline{a}\left(\frac{q\mu}{d}+r\right)& 
=\underline{a}\left(\frac{q\mu}{d}\right)+\underline{a}(r), \\ 
\underline{a}\left(\frac{q\mu}{d}\right)&=\left(\frac{qw_{0}}{d}, \ldots 
  ,\frac{qw_{n}}{d}\right). \nonumber 
\end{align} 
 
\begin{notation}Nous avons besoin de préciser les notations du chapitre 
\ref{cha:prel-comb}. 
Nous rajoutons un indice $w$ ou $w/d$ aux notations $s,\sigma$ du chapitre 
\ref{cha:prel-comb} pour préciser que les nombres rationnels
$s_{w}(\cdot)$ et $\sigma_{w}(\cdot)$ (resp. $s_{w/d}(\cdot)$ et $\sigma_{w/d}(\cdot)$)
sont calculés avec les poids $w_{0}, \ldots ,w_{n}$ (resp. avec les poids $w_{0}/d, \ldots ,w_{n}/d$).
\end{notation} 
 
\begin{lem}\label{lem:nombres,spectraux} 
  Pour tout $(q,r)\in \{0,\ldots,d-1\}\times\{0, \ldots 
,(\mu/d)-1\}$, nous avons les égalités   
\begin{align*} 
  s_{w}\left(\frac{q\mu}{d}+r\right)&=\frac{q}{d}+s_{w}(r) \,;& 
 \sigma_{w}\left(\frac{q\mu}{d}+r\right)&=\sigma_{w}(r)=\sigma_{w/d}(r). 
\end{align*} 
\end{lem} 
 
\begin{proof} 
Pour $k\in\{0, \ldots ,\mu-1\}$, nous avons  
$s_{w}(k)=\underline{a}(k)_{i(k)}/w_{i(k)}$. 
Pour $k\in\{0, \ldots ,(\mu/d)-1\}$, nous en déduisons que  
$s_{w/d}(k)=ds_{w}(k)$. 
Comme $\sigma_{w}(k)=k-\mu s_{w}(k)$,  nous obtenons   
\begin{align}\label{eq:14} 
  \sigma_{w/d}(k)&=\sigma_{w}(k) 
\end{align} 
pour $k\in\{0, \ldots ,(\mu/d)-1\}$. 
Les égalités (\ref{eq:a}) et (\ref{eq:14}) terminent la démons\-tra\-tion. 
\end{proof}

Pour tout $\alpha\in\{0,\ldots,d-1\}$, et pour 
$k\in\{0,\ldots,\mu-1\}$, nous définissons les éléments 
$\omega_{k,\alpha} :=u_{\alpha}^{\underline{a}(k)}\omega_{0,\alpha}$ 
où 
$u_{\alpha}^{\underline{a}(k)} :=u_{0,\alpha}^{\underline{a}(k)_{0}}\cdots 
u_{n,\alpha}^{\underline{a}(k)_{n}}$.  
D'après les égalités (\ref{eq:a}),  pour tout $(q,r)\in \{0,\ldots,d-1\}\times\{0, \ldots 
,(\mu/d)-1\}$ nous avons  
\begin{align} 
\label{eq:15} 
\omega_{\frac{q\mu}{d}+r,\alpha} & =\zeta^{\alpha q}\omega_{r,\alpha}. 
\end{align} 
 
  Pour tout $k$, la classe de 
 $\omega_{k,\alpha}$ appartient à $G_{0}(f_{\alpha})$.   
Pour tout $k\in\{0,\ldots,\mu-1\}$, nous posons  
\begin{align*} 
  \omega_{k}:=(\omega_{k,0},\ldots,\omega_{k,d-1})\in G_{0}. 
\end{align*} 
Pour tout $(q,r)\in \{0,\ldots,d-1\}\times\{0, \ldots 
,(\mu/d)-1\}$, nous avons  
\begin{align}\label{eq:11} 
  \omega_{\frac{q\mu}{d}+r}&=(1,\zeta^q,\ldots,\zeta^{q(d-1)})\omega_{r}. 
\end{align} 
 
La proposition suivante est l'analogue de la proposition de $3.2.$ de 
\cite{DSgm2}. 
\begin{prop}\label{prop:douai,sabbah} 
  Les classes des éléments $\omega_{0}, \ldots ,\omega_{\mu-1}$ forment 
  une $\CC[\theta]$-base de $G_{0}$, notée $\bs{\omega}$. De plus, pour $k\in\{0, \ldots 
  ,\mu-1\}$, nous avons les relations 
  \begin{align*} 
  -\frac{1}{\mu}\theta(\sigma_{w}(k)-\theta\partial_{\theta})\omega_{k}=\frac{\omega_{\overline{k+1}}}{w_{i(k)}} 
  \end{align*} 
où $\overline{k+1}$ désigne la réduction modulo $\mu$. 
  De plus, l'ordre de $\omega_{k}$ pour la filtration $V_{\bullet}$ 
  est $\sigma_{w}(k)$ et $\bs{\omega}$ induit une base sur 
  $\oplus_{\alpha}\gr^{V}_{\alpha}(G_{0}/\theta G_{0})$. 
\end{prop} 
 
\begin{proof} 
  Pour démontrer cette proposition, il suffit d'adapter la 
  dé\-mons\-tra\-tion de la proposition de $3.2.$ de \cite{DSgm2}.  
  \begin{itemize} 
  \item  
  Nous allons d'abord démontrer cette proposition pour le polyn\^ome de 
  Laurent $f_{\alpha}$ sur le tore $U_{\alpha}$. 
En particulier, il faut changer la relation $(3.5)$ de \cite{DSgm2} en 
\begin{align}\label{eq:recurrence} 
 \theta\left(\frac{1}{w_{i}/d}u_{i,\alpha}\partial_{u_{i,\alpha}}-\frac{1}{w_{0}/d}u_{0,\alpha}\partial_{u_{0,\alpha}}\right)\varphi\omega_{0,\alpha}&=\left(\frac{u_{i,\alpha}}{w_{i}/d}-\frac{u_{0,\alpha}}{w_{0}/d}\right)\varphi\omega_{0,\alpha} 
\end{align} 
pour tout $\varphi\in 
\CC[u_{\alpha},u_{\alpha}^{-1},\theta,\theta^{-1}]$ et tout $i\in\{0, 
\ldots ,n\}$.   Le reste de la dé\-mons\-tra\-tion est identique 
à celle de \cite{DSgm2}.  Nous en déduisons la proposition pour le 
polyn\^ome de Laurent $f_{\alpha}$ sur le tore $U_{\alpha}$. En particulier, 
pour tout $k \in \{0, \ldots ,(\mu/d)-1\}$, nous avons les relations 
\begin{align}\label{eq:16} 
   -\frac{1}{\mu/d}\theta(\sigma_{w/d}(k)-\theta\partial_{\theta})\omega_{k}&=\frac{\omega_{k+1}}{w_{i(k)}/d} & \mbox{ où } \omega_{\mu/d,\alpha}=\zeta^{\alpha}\omega_{0,\alpha}. 
\end{align}

\item Montrons que $\omega_{0},\ldots,\omega_{\mu-1}$ forment une 
  $\CC[\theta]$-base de $G_{0}$. Nous savons que 
  $\omega_{0,\alpha},\ldots,\omega_{(\mu/d)-1,\alpha}$ forment une 
  $\CC[\theta]$-base de $G_{0}(f_{\alpha})$ pour tout 
  $\alpha\in\{0,\ldots,d\}$.  D'après l'égalité (\ref{eq:10}), une 
  base $\CC[\theta]$-base de $G_{0}$ est formée par les vecteurs suivants 
  \begin{align*} 
    \begin{array}{cccc} 
(\omega_{0,1},0,\ldots,0) & (0,\omega_{0,2},0,\ldots,0) & \cdots & 
(0,\ldots,0,\omega_{0,d})\\ 
(\omega_{1,1},0,\ldots,0) & (0,\omega_{1,2},0,\ldots,0) & \cdots & 
(0,\ldots,0,\omega_{1,d})\\ 
 \vdots& \vdots&\vdots &\vdots \\ 
(\omega_{(\mu/d)-1,1},0,\ldots,0) & (0,\omega_{(\mu/d)-1,2},0,\ldots,0) & \cdots & 
(0,\ldots,0,\omega_{(\mu/d)-1,d})\\ 
 \end{array} 
  \end{align*} 
Pour $(i,j)\in\{1,\ldots,d\}\times\{0,\ldots,(\mu/d)-1\}$, notons 
\begin{align*}
e_{id+j-1}&:=(0,\ldots,0,\omega_{i,j},0,\ldots,0).
\end{align*}
  Pour tout $(q,r)\in \{0,\ldots,d-1\}\times\{0, \ldots 
,(\mu/d)-1\}$, nous 
  avons  
\begin{align*} 
\omega_{\frac{q\mu}{d}+r}=e_{dr}+\zeta^q e_{dr+1}+\cdots+\zeta^{q(d-1)}e_{dr+d-1}. 
\end{align*} 
Nous en déduisons que la matrice de passage entre les $\omega_{\frac{q\mu}{d}+r}$ 
et les $e_{i}$ est une matrice diagonale par bloc où les blocs sont 
tous égaux à la matrice 
\begin{align*} 
\left(  \begin{array}{cccc} 
1 & 1& \cdots & 1 \\ 
1& \zeta & \cdots &\zeta^{d-1} \\ 
\vdots& \vdots & & \vdots \\ 
 1& \zeta^d & \cdots & \zeta^{d(d-1)} 
  \end{array}\right) 
\end{align*} 
Cette matrice est inversible. Nous en déduisons que 
$\omega_{0},\ldots,\omega_{\mu-1}$ est une $\CC[\theta]$-base de 
$G_{0}$.  
 Le reste de la proposition découle de la formule 
 (\ref{eq:recurrence}) et du lemme \ref{lem:nombres,spectraux}. 
\end{itemize} 
\end{proof}

Pour tout $\alpha\in\{0,\ldots,d-1\}$, nous avons un produit sur le quotient  
$G_{0}(f_{\alpha})/\theta G_{0}(f_{\alpha})$ qui provient d'un 
isomorphisme entre le quotient jacobien de $f_{\alpha}$ et
$G_{0}(f_{\alpha})/\theta G_{0}(f_{\alpha})$. Pour tout
$\alpha\in\{0,\ldots,d-1\}$, notons  $[g]$ la classe d'un élément $g \in
G_{0}(f_{\alpha})$ dans le quotient $G_{0}(f_{\alpha})/\theta
G_{0}(f_{\alpha})$.
Ce produit est donné par la formule 
\begin{align}\label{eq:18} 
  [h_{1}\omega_{0,\alpha}] \star_{\alpha}[h_{2}\omega_{0,\alpha}]&=[h_{1}h_{2}\omega_{0,\alpha}] 
\end{align} où $h_{1},h_{2}$ sont dans 
$\CC[u_{\alpha},u^{-1}_{\alpha}]$.  
Ceci  nous permet de définir un produit, noté $\star$, sur le quotient 
  $G_{0}/\theta G_{0}$. D'après la proposition
  \ref{prop:douai,sabbah}, la base $\bs{\omega}$ de $G_{0}$ induit une 
  base $[\bs{\omega}]:=([\omega_{0}], \ldots ,[\omega_{\mu-1}])$ de  $G_{0}/\theta G_{0}$.

\begin{lem}\label{lem:produit,star} 
  Dans la base $[\bs{\omega}]$ de $G_{0}/\theta G_{0}$, le produit est 
  donné par 
  \begin{align*} 
  [\omega_{i}]\star[\omega_{j}]=w^{\underline{a}(i)+\underline{a}(j)-\underline{a}(i+j)}[\omega_{\overline{i+j}}] 
  \end{align*} 
  où $ \overline{i+j}$ désigne la somme modulo $\mu$. 
\end{lem}

\begin{proof}\begin{itemize} \item    
    D'abord, nous allons montrer que  
    \begin{align*} 
      [\omega_{k}] 
    =[w^{\underline{a}(k)}(u_{0}/w_{0})^{k}\omega_{0}]. 
    \end{align*} 
Les relations 
    (\ref{eq:recurrence}) dans le quotient $G_{0}/ \theta G_{0}$ 
    deviennent 
\begin{align}\label{equation:quotient} 
\left[\frac{u_{i}}{w_{i}}\omega_{0}\varphi\right]&=\left[\frac{u_{0}}{w_{0}}\omega_{0}\varphi\right], 
\end{align} 
pour tout $\varphi\in \CC[u,u^{-1},\theta,\theta^{-1}]$ et tout 
$i\in\{0, \ldots ,n\}$.  Nous en déduisons que 
$[\omega_{k}]=[w^{\underline{a}(k)}(u_{0}/w_{0})^{k}\omega_{0}]$.

\item Par définition du produit dans $G_{0}/\theta G_{0}$, nous avons 
  \begin{align*}[\omega_{i}]\star 
  [\omega_{j}]=\left[w^{\underline{a}(i)+\underline{a}(j)}\left(\frac{u_{0}}{w_{0}}\right)^{i+j}\omega_{0}\right]\end{align*} 
  et 
  \begin{align*} 
  \left[\left(\frac{u_{0}}{w_{0}}\right)^{i+j}\omega_{0}\right]=\left[w^{-\underline{a}(i+j)}u^{\underline{a}(i+j)}\omega_{0}\right]. 
  \end{align*} 
   
  Or la relation $\prod u_{i}^{w_{i}}=1$ sur $G$ induit 
  $[u^{\underline{a}(i+j)}\omega_{0}]=[u^{\underline{a}(\overline{i+j})}\omega_{0}]= 
  [\omega_{\overline{i+j}}]$. 
\end{itemize} 
\end{proof} 
 
Pour tout $\gamma\in S_{w}$, nous notons 
 
\begin{align*} 
  k^{\min}(\gamma):=\min\{j\in \{0, \ldots ,(\mu/d)-1\}\mid s_{w}(j)=\gamma\}. 
\end{align*} 
D'après (\ref{eq:min,egalite}), nous avons 
$k^{\min}(\gamma)=k^{\max}(\gamma)-\delta(\gamma)+1$. Pour tout $i\in\{0, 
\ldots ,\mu-1\}$, posons 
\begin{align}\label{eq:19} 
\widetilde{\omega}_{i}&:=\frac{w^{\underline{a}(k^{\min}(s_{w}(i)))}}{w^{\underline{a}(i)}}\omega_{i}. 
\end{align} 
 
\begin{prop} 
  Dans la base $\bs{\widetilde{\omega}}$ de $G_{0}$, nous avons 
  \begin{align*} 
  \theta^{2}\partial_{\theta}\bs{\widetilde{\omega}}=\bs{\widetilde{\omega}}A_{0}^{\circ}+ 
  \theta \bs{\widetilde{\omega}}A_{\infty}. 
  \end{align*} 
\end{prop} 
 
\begin{proof} 
  D'après la proposition \ref{prop:douai,sabbah}, nous avons 
 \begin{align*} 
  -\frac{1}{\mu}\theta(\sigma_{w}(k)-\theta\partial_{\theta})\widetilde{\omega}_{k}&=\widetilde{\omega}_{\overline{k+1}} 
  \frac{w^{\underline{a}(k^{\min}(s_{w}(k)))}}{w^{\underline{a}(k^{\min}(s_{w}(k+1)))}}\\ 
  \theta^{2}\partial_{\theta}\widetilde{\omega}_{k}&= 
  \mu\frac{w^{\underline{a}(k^{\min}(s_{w}(k)))}}{w^{\underline{a}(k^{\min}(s_{w}(k+1)))}}\widetilde{\omega}_{\overline{k+1}}+\theta{\sigma_{w}(k)} \widetilde{\omega}_{k}. 
  \end{align*} 
\end{proof} 
 
Maintenant, nous allons calculer la forme bilinéaire non dégénérée au point $0$ de l'espace  
des paramètres d'un déploiement universel de $f$. 
 
Soit $G$ un $\CC[\theta,\theta^{-1}]$-module. Nous notons 
$\overline{G}$ le $\CC$-espace vectoriel $G$ équipé de la structure de 
module $p(\theta)\cdot g=p(-\theta)g$ où $p(\theta)\in 
\CC[\theta,\theta^{-1}]$ . Nous notons par $\overline{g}$ les éléments 
de $\overline{G}$. Si $G$ est équipé d'un opérateur $\partial_{\theta}$ compatible à la 
multiplication par $\theta$, alors sur $\overline{G}$, nous avons un opérateur 
$\partial_{\theta}\overline{g}:=\overline{-\partial_{\theta}g}$. 
Remarquons que  
$\theta\partial_{\theta}\overline{g}=\overline{\theta\partial_{\theta}g}$. 
 
D'après les résultats du paragraphe $4$ de l'article 
\cite{DSgm2}, pour tout $\alpha\in\{0, \ldots ,d-1\}$, il existe une forme 
$\CC[\theta,\theta^{-1}]$-bilinéaire non dégénérée 
\begin{align*}S_{\alpha}: G(f_{\alpha}) \otimes_{\CC[\theta,\theta^{-1}]} 
\overline{G(f_{\alpha})} & \longrightarrow \CC[\theta,\theta^{-1}] 
\end{align*} 
qui vérifie les propriétés suivantes : 
\begin{enumerate} 
\item\label{item:20} 
  $\theta\partial_{\theta} S_{\alpha}(p_{1},\overline{p_{2}})=S_{\alpha} (\theta\partial_{\theta}p_{1},\overline{p_{2}})+S_{\alpha} (p_{1},\overline{\theta\partial_{\theta}p_{2}})$ ; 
\item \label{item:21} $S_{\alpha} $ envoie $V_{0}G(f_{\alpha})\otimes \overline{V_{<1}G(f_{\alpha})}$ dans 
  $\CC[\theta^{-1}]$ ; 
\item\label{item:22} $S_{\alpha}$ envoie 
  $G_{0}(f_{\alpha})\otimes\overline{G_{0}(f_{\alpha})}$ dans 
  $\theta^{n}\CC[\theta]$ ; 
\item \label{item:23}$S_{\alpha}(p_{1},\overline{p_{2}})=(-1)^{n}\overline{S_{\alpha}(p_{1},\overline{p_{2}})}$. 
\end{enumerate}

\begin{lem}\label{lem:forme,S,alpha} Pour tout $\alpha\in\{0, \ldots 
  ,d-1\}$, il existe 
  une unique, à une constante près, forme bilinéaire non dégénérée $S_{\alpha}$ 
  telle que  les conditions (\ref{item:20}), (\ref{item:21}), 
  (\ref{item:22}) et \eqref{item:23} soient vérifiées. 
Dans la base $\omega_{0,\alpha}, \ldots ,\omega_{(\mu/d)-1,\alpha}$ de  
$G(f_{\alpha})$, nous avons  
\begin{align*} 
  S_{\alpha}(\omega_{r_{j},\alpha},{\omega}_{r_{k},\alpha})= 
\begin{cases} 
C\cdot S_{\alpha}(\omega_{0,\alpha},\omega_{n,\alpha}) & \mbox{ si } 
{r_{j}+r_{k}}= n\,; \\ 
\zeta^{\alpha}C\cdot S_{\alpha}(\omega_{0,\alpha},\omega_{n,\alpha})& \mbox{ si } 
{r_{j}+r_{k}}= n +(\mu/d)\,; \\ 
0 & \mbox{ sinon.} 
\end{cases} 
\end{align*} 
où  $C=w_{n} 
 w^{\underline{a}(r_{j})+\underline{a}(r_{k})-\underline{a}(n+1)-\underline{a}(r_{j}+r_{k}-n)}$. 
\end{lem} 
 
\begin{rem} 
  Comme $|\underline{a}(k)|=k$, nous avons l'égalité suivante 
  \begin{align*} 
    w_{n} 
 w^{\underline{a}(r_{j})+\underline{a}(r_{k})-\underline{a}(n+1)-\underline{a}(r_{j}+r_{k}-n)} &=  dw_{n}({w}/{d})^{\underline{a}(r_{j})+\underline{a}(r_{k})-\underline{a}(n+1)-\underline{a}(r_{j}+r_{k}-n)}. 
  \end{align*} 
\end{rem}

\begin{proof}[Démonstration du lemme \ref{lem:forme,S,alpha}] 
  Nous suivons la preuve du lemme $4.1$ de l'article \cite{DSgm2}. 
  Pour tout $r_{j},r_{k}\in\{0, \ldots ,(\mu/d)-1\}$, d'après (\ref{item:22}), nous avons  
 $S_{\alpha}(\omega_{r_{j},\alpha},\omega_{r_{k},\alpha})\in 
 \theta^{n}\CC[\theta]$ et  $S_{\alpha}(\omega_{0,\alpha},\omega_{r_{k},\alpha})\in 
 \theta^{[\sigma_{w/d}(r_{k})]}\CC[\theta]$ d'après  (\ref{item:21}) où 
 $[\cdot]$ désigne  la partie entière. Si $\sigma_{w/d}(r_{k})=n$ 
 alors nous avons 
 $S_{\alpha}(\omega_{0,\alpha},\omega_{r_{k},\alpha})\neq 0$. 
Or $\sigma_{w/d}(r_{k})=n$ implique que $r_{k}=n$. 
 
D'après les égalités (\ref{item:20}) et (\ref{eq:16}), 
nous obtenons 
\begin{align*} 
  &   -\frac{1}{\mu}\left(-\theta\partial_{\theta}+n\right)S_{\alpha}(\omega_{r_{j},\alpha},\overline{\omega}_{r_{k},\alpha}) \\ 
  &  =  
  \frac{\sigma_{w/d}(r_{j})+\sigma_{w/d}(r_{k})-n}{\mu/d}S(\omega_{r_{j},\alpha},\overline{\omega}_{r_{k},\alpha})\\  
  & + \theta^{-1}\left( 
    S_{\alpha}\left(\omega_{r_{j},\alpha},\frac{\overline{\omega}_{r_{k}+1,\alpha}}{w_{i(r_{k})}/d}\right)-S_{\alpha}\left(\frac{\omega_{r_{j}+1,\alpha}}{w_{i(r_{j})}/d},\overline{\omega}_{r_{k},\alpha}\right)\right) 
  \nonumber 
\end{align*} 
où $\omega_{\mu/d,\alpha}=\zeta^{\alpha}\omega_{0,\alpha}$ (cf. (\ref{eq:15})). 
Comme $S_{\alpha}(\omega_{r_{j},\alpha},\overline{\omega}_{r_{k},\alpha})\in \theta^{n}\CC[\theta]$, 
le membre de gauche de l'égalité ci-dessus est nul. 
Puis, un raisonnement par récurrence montre que si 
$S_{\alpha}(\omega_{r_{j},\alpha},\omega_{r_{k},\alpha})\neq 0$ alors 
$r_{j}+r_{k}\equiv n \mod [\mu/d]$. 
Nous en déduisons  
\begin{itemize} 
\item Si $r_{j}+r_{k} = n $ nous avons  
\begin{align*} 
   S_{\alpha}\left(\omega_{r_{j},\alpha},\overline{\omega}_{r_{k}+1,\alpha}\right)=\frac{w_{i(r_{k})}}{w_{i(r_{j})}}S_{\alpha}\left(\omega_{r_{j}+1,\alpha},\overline{\omega}_{r_{k},\alpha}\right).\end{align*} 
 \item Si $r_{j}+r_{k} = n + (\mu/d) $ nous avons  
\begin{align*} 
   S_{\alpha}\left(\omega_{r_{j},\alpha},\overline{\omega}_{r_{k}+1,\alpha}\right)=\zeta^{\alpha}\frac{w_{i(r_{k})}}{w_{i(r_{j})}}S_{\alpha}\left(\omega_{r_{j}+1,\alpha},\overline{\omega}_{r_{k},\alpha}\right).\end{align*} 
\end{itemize} 
\end{proof}

Pour tout $\alpha\in\{0, \ldots ,d-1\}$, nous posons 
\begin{align*} 
  S_{\alpha}(\omega_{0,\alpha},\omega_{n,\alpha})=\theta^{n}/(dw_{n}). 
\end{align*} 
 
Pour tout $(p_{1},p_{2})\in 
G_{0}(f_{\alpha})\otimes\overline{G_{0}}(f_{\alpha})$, le  
coefficient devant  $\theta^{n}$ de $S_{\alpha}(p_{1},p_{2})$ ne dépend que de la classe de 
$p_{1},p_{2}$ dans $G_{0}(f_{\alpha})/\theta G_{0}(f_{\alpha})$. 
Nous notons $[g_{\alpha}]([p_{1}],[p_{2}])$ ce coefficient. 
Nous en déduisons une forme bilinéaire 
 symétrique et non dégénérée, notée $[g_{\alpha}]$, sur  $G_{0}(f_{\alpha})/\theta G_{0}(f_{\alpha})$. 
Nous posons $[g]:=\sum_{\alpha=0}^{d-1}[g_{\alpha}]$ et  
nous obtenons une forme bilinéaire 
 symétrique et non dégénérée sur  $G_{0}/\theta G_{0}$.

\begin{lem}\label{lem:metrique} 
  Dans la base $[\bs{\widetilde{\omega}}]$ de $G_{0}/\theta G_{0}$, la
  forme bilinéaire non dégénérée est donnée par
  \begin{align*} 
  [g]([\widetilde{\omega}_{j}],[\widetilde{\omega}_{k}])= 
\begin{cases} 
\left(\prod_{i\in I(s_{w}(j))}w_{i}\right)^{-1} & \mbox{ si } \overline{j+k}=n\,;\\ 
0 & \mbox{ sinon} 
\end{cases} 
\end{align*} 
où $\overline{j+k}$ désigne la réduction modulo $\mu$. 
\end{lem} 
 
\begin{proof} 
Pour tout $j,k\in\{0, \ldots ,\mu-1\}$, il existe des uniques couples 
$(q_{j},r_{j})$ et $(q_{k},r_{k})$ dans $\{0, \ldots ,d-1\}\times\{0, 
\ldots ,(\mu/d)-1\}$ tels que 
\begin{align*} 
  j&=q_{j}\frac{\mu}{d}+r_{j}\, \mbox{ et } \, k=q_{k}\frac{\mu}{d}+r_{k}. 
\end{align*} 
Nous en déduisons 
\begin{align*} 
  [g]([\omega_{j}],[\omega_{k}]) & = 
  \sum_{\alpha=0}^{d-1}[g_{\alpha}]([\omega_{j,\alpha}],[\omega_{k,\alpha}]) \\ 
& = 
\sum_{\alpha=0}^{d-1}\zeta^{\alpha(q_{j}+q_{k})}[g_{\alpha}]([\omega_{r_{j},\alpha}],[\omega_{r_{k},\alpha}])  
  & \mbox{ d'après (\ref{eq:15})}.  
\end{align*} 
D'après le lemme \ref{lem:forme,S,alpha} et un calcul direct,  nous obtenons 
\begin{align*} 
  [g]([\omega_{j}],[\omega_{k}])= 
  \begin{cases} 
    w^{\alpha} & 
      \mbox{ si }\begin{cases} \left(r_{j}+r_{k}=n\right) \\ \mbox{ et }\left( q_{j}+q_{k}=0\mbox{ ou }d\right)\,;\end{cases}\\ 
  w^{\alpha} & \mbox{ si } \begin{cases}\left(r_{j}+r_{k}=n +(\mu/d)\right)\\ \mbox{ et } \left(q_{j}+q_{k}+1=d\right)\,;\end{cases}\\ 
   0 & \mbox{ sinon } 
  \end{cases} 
\end{align*} 
où  
\begin{align*} 
  \alpha=\underline{a}(r_{j})+\underline{a}(r_{k})-\underline{a}(n+1)-\underline{a}(r_{j}+r_{k}-n). 
\end{align*} 
Comme $\overline{j+k}=n$ est équivalent à $s_{w}(j)=\{1-s_{w}(k)\}$, 
nous avons les équivalences suivantes 
\begin{align} 
  j+k=n &\Leftrightarrow r_{j}+r_{k}=n \mbox{ et } q_{j}+q_{k}=0 \label{eq:13}\,;\\ 
 j+k=n+\mu &\Leftrightarrow  
 \begin{cases} 
   r_{j}+r_{k}=n \mbox{ et } q_{j}+q_{k}=d\\ 
 \mbox{ ou } \\ 
r_{j}+r_{k}=n +(\mu/d)\mbox{ et }q_{j}+q_{k}+1=d 
 \end{cases}\label{eq:17} 
\end{align} 
Dans la base $[\bs{\omega}]$ la forme bilinéaire $[g]$ est donnée 
  par : 
\begin{align*} 
  [g]([\omega_{j}],[\omega_{k}])= 
\begin{cases}   
{w^{\underline{a}(r_{j})+\underline{a}(r_{k})-\underline{a}(n+1)-\underline{a}(r_{j}+r_{k}-n)}} & \mbox{ si } \overline{j+k}=n\,;\\ 
0 & \mbox{ sinon.} 
\end{cases} 
\end{align*} 
 
Puis dans la base $[\bs{\widetilde{\omega}}]$ (cf. égalité (\ref{eq:19})), nous avons 
\begin{align*} 
[g]([\widetilde{\omega}_{j}],[\widetilde{\omega}_{k}])= 
\begin{cases}   
\ds{\frac{w^{\alpha}}{w^{\beta}}} & \mbox{ si } \overline{j+k}=n\,;\\ 
0 & \mbox{ sinon.} 
\end{cases} 
\end{align*} 
où 
  \begin{align*} 
\alpha&=\underline{a}(k^{\min}(s_{w}(j)))+\underline{a}(k^{\min}(s_{w}(k)))\,;\\ 
\beta&=\underline{a}(n+1)+\underline{a}(r_{j}+r_{k}-n)+\underline{a}(j)-\underline{a}(r_{j})+\underline{a}(k)-\underline{a}(r_{k}) 
  \end{align*} 
D'après les égalités (\ref{eq:a}), nous avons 
\begin{align*} 
  \underline{a}(j)-\underline{a}(r_{j})+\underline{a}(k)-\underline{a}(r_{k}) & = a\left((q_{j}+q_{k})\frac{\mu}{d}\right). 
\end{align*} 
L'hypothèse $\overline{j+k}=n$, vue à travers les équivalences 
(\ref{eq:13}) et (\ref{eq:17}), implique que $\beta=\underline{a}(n+1)+\underline{a}(j+k-n)$. 
 
Si $j+k=n$ alors $s_{w}(j)=s_{w}(k)=0$ et nous avons 
$[g]([\widetilde{\omega}_{j}],[\widetilde{\omega}_{k}])=\left(\prod_{i=0}^{n}w_{i}\right)^{-1}$. 
 
Supposons que $j+k=n+\mu$. D'abord nous allons montrer que pour tout 
$\gamma>0$ dans $S_{w}$, nous avons 
\begin{align}\label{eq:a,gamma} 
\underline{a}(k^{\min}(\gamma))+\underline{a}(k^{\max}(\Inv{\gamma})+1)&=\underline{a}(\mu)+\underline{a}(n+1). 
\end{align} 
Pour tout $\alpha \in\{0, \ldots ,\mu-1\}$, nous avons 
\begin{align*} 
\underline{a}(k^{\min}(\gamma))_{\alpha}= 
\begin{cases} 
 [ w_{\alpha}\gamma ] & \mbox{ si } \alpha \in I(\gamma)\,;  \\    
 [ w_{\alpha}\gamma ] +1 & \mbox{ sinon}. 
\end{cases} 
\end{align*} 
Comme $\underline{a}(k^{\max}(\gamma)+1)_{\alpha}=[w_{\alpha}\gamma]+1$, nous obtenons 
\begin{align*} 
\underline{a}(k^{\min}(\gamma))_{\alpha}+\underline{a}(k^{\max}(\Inv{\gamma})+1)_{\alpha}= 
\begin{cases} 
w_{\alpha} + [\gamma w_{\alpha}] + [-\gamma w_{\alpha}] +1  & \mbox{ si } \alpha \in I(\gamma)\,;  \\    
 w_{\alpha} + [\gamma w_{\alpha}] + [-\gamma w_{\alpha}] +2 & \mbox{ sinon}. 
\end{cases} 
\end{align*} 
L'égalité 
\begin{align*} 
[\gamma w_{\alpha}]+[-\gamma w_{\alpha}]=\begin{cases} 0 
    & \mbox{si } \alpha\in 
    I(\gamma)\,;\\ 
    -1 & \mbox{sinon.}  \end{cases} 
\end{align*} 
montre que 
\begin{align*}\underline{a}(k^{\min}(\gamma))_{\alpha}+\underline{a}(k^{\max}(\Inv{\gamma})+1)_{\alpha}=w_{\alpha}+1.\end{align*} 
Puis, les égalités (\ref{eq:a}) impliquent la formule 
(\ref{eq:a,gamma}). 
 
Comme $j+k=n+\mu$, nous avons $s_{w}(j)=\Inv{s_{w}(k)}$. Nous 
appliquons la formule (\ref{eq:a,gamma}) à $\gamma=s_{w}(k)$ et nous en 
déduisons que 
\begin{align*} 
[g]([\widetilde{\omega}_{j}],[\widetilde{\omega}_{k}])=w^{\underline{a}(k^{\min}(s_{w}(j)))-\underline{a}(k^{\max}(s_{w}(j))+1)}=\left(\prod_{i\in 
  I(s_{w}(j))}w_{i}\right)^{-1}. 
\end{align*} 
\end{proof}  
 
\begin{notation} 
  Dorénavant, nous supprimons les indices $w$ dans $s_{w}$ et 
  $\sigma_{w}$ car nous travaillerons toujours avec les poids $w_{0}, \ldots ,w_{n}$. 
\end{notation}

Soit $F:U\times X\rightarrow\CC$ un déploiement universel de $f$.  La 
forme $[\omega_{0}]\in G_{0}/\theta G_{0}$ est une section primitive
homogène et canonique
c'est-à-dire que $[\omega_{0}]$ vérifie les propriétés suivantes :
\begin{itemize}
\item (\emph{primitive})  $[\omega_{0}]$ induit un isomorphisme entre
  $G_{0}/\theta G_{0}$
  est le quotient jacobien de $f$ ;
\item (\emph{homogène})  $[\omega_{0}]$ est un vecteur propre de
  l'endomorphisme dont la matrice est $-A_{\infty}$ dans la base
  $[\bs{\omega}]$ ;
\item (\emph{canonique}, cf. paragraphe $3$.c de \cite{DSgm1})
  $[\omega_{0}]$ est un vecteur propre de la matrice $-A_{\infty}$
  pour la valeur propre maximale (pour le polyn\^{o}me de Laurent $f$
  cette valeur propre est $0$).
\end{itemize}
 
Notons $\varphi_{[\omega_{0}]}^{\circ}:T_{0}X\rightarrow G_{0}/\theta 
G_{0}$ l'application de période infinitésimale en 
$\underline{x}=\underline{0}$ défini par K.\kern2pt Saito
\cite{KSpmpf} (voir aussi \cite{Sdivf} p.$244$). Nos conditions initiales de la variété 
de Frobenius sont donc 
$(A_{0}^{\circ},A_{\infty},\varphi_{[\omega_{0}]}^{-1}[g],\varphi_{[\omega_{0}]}^{-1}[\omega_{0}])$. 
 
La matrice $A_{0}^{\circ}$ est la multiplication par le champ d'Euler 
$\mathfrak{E}$ en $\underline{x}=\underline{0}$. Pour faciliter la 
correspondance entre le c\^{o}té $A$ et le c\^{o}té $B$, nous allons définir 
un objet qui va coder cette multiplication par le champ d'Euler  
en  $\underline{x}=\underline{0}$. 
Nous définissons $(\!([a],[b],[c])\!):=[g]([a]\star [b],[c])$ pour 
tout $[a],[b],[c]$ dans $G_{0}/\theta G_{0}$.

\begin{prop}\label{prop:3tenseur} 
  Soient $j,k$ dans $\{0, \ldots ,\mu-1\} $. 
\begin{enumerate} \item Si $\overline{1+j+k} \neq n$  alors 
  $(\!([\widetilde{\omega}_{1}],[\widetilde{\omega}_{j}],[\widetilde{\omega}_{k}])\!)=0$. 
\item Si $\overline{1+j+k} = n$ alors nous avons 
\begin{align*}(\!([\widetilde{\omega}_{1}],[\widetilde{\omega}_{j}],[\widetilde{\omega}_{k}])\!)= 
 \begin{cases} 
   \left({{\prod_{i\in I(j,k)}w_{i}}}\right)^{-1} & \mbox{ si } 
   \sigma(1)+\sigma(j)+\sigma(k)\neq n\,;\\ 
   \left({{\prod_{i\in I(s(j))}w_{i}}}\right)^{-1} & \mbox{ si } 
   \sigma(1)+\sigma(j)+\sigma(k)=n
 \end{cases} 
\end{align*} 
où $I(j,k):=I(s(j))\bigsqcup  I(s(k))$. 
\end{enumerate} 
\end{prop}

\begin{proof} Quitte à changer $j$ et $k$, on peut supposer que $j\leq  
  k$. 
\begin{itemize} 
\item Dans un premier temps, nous allons démontrer la formule suivante 
\begin{align}\label{eq:3,tenseur} 
(\!([\widetilde{\omega}_{1}],[\widetilde{\omega}_{j}],[\widetilde{\omega}_{k}])\!)&=  
\begin{cases} 
\ds{\frac{w^{\underline{a}(k^{\min}(s(j)))}}{w^{\underline{a}(k^{\max}(\Inv{s(k)})+1)}}}  
  & \mbox{ si } \overline{1+j+k}=n\,; \\ 
0 & \mbox{ sinon.} 
\end{cases} 
\end{align} 
D'après le lemme \ref{lem:produit,star} et le lemme 
\ref{lem:metrique}, nous obtenons 
\begin{align*}(\!([\widetilde{\omega}_{1}],[\widetilde{\omega}_{j}],[\widetilde{\omega}_{k}])\!)= 
\begin{cases} 
\ds{\frac{w^{\alpha}}{w^{\beta}} } 
& \mbox{ si } \overline{1+j+k}=n \,;\\ 0 & \mbox{ sinon} 
\end{cases} 
\end{align*} 
où 
  \begin{align*} 
\alpha&=\underline{a}(k^{\min}(s(j)))+\underline{a}(k^{\min}(s(k))) \,;\\ 
\beta&=\underline{a}(n+1)+\underline{a}(1+j+k-n). 
  \end{align*} 
Si $1+j+k=n$ alors $s(j)=s(k)=0$. Nous en déduisons la formule 
(\ref{eq:3,tenseur}). 
 
Supposons que $1+j+k=n+\mu$. 
 
Si $(w_{0}, \ldots ,w_{n})=(1, \ldots ,1)$, la seule possibilité est 
$j=k=n$ et la formule (\ref{eq:3,tenseur}) est vraie. 
 
Supposons que $(w_{0}, \ldots ,w_{n})\neq(1, \ldots ,1)$. Nous en 
déduisons que les inégalités suivantes  
\begin{align*} 
  \mu&\geq n+2,& n+1& \leq j+1\leq 
\mu-1\\ 
 \mbox{et }  s(j+1)&>0. 
\end{align*} 
 La condition $1+j+k=n+\mu$ implique que 
$s(k)=\Inv{s(\overline{j+1})}>0$.  Puis, nous appliquons la 
formule (\ref{eq:a,gamma}) à $\gamma=s(k)$ et nous en déduisons la 
formule (\ref{eq:3,tenseur}) 
 
\item Soient $j,k$ tels que $1+j+k=\varepsilon+n$ où $\varepsilon=0$ 
  ou $1$. Nous avons $s(k)=\Inv{s(\overline{j+1})}$.  Pour 
  finir la démonstration, il suffit de vérifier que 
\begin{align}\label{eq:inverse,s(k)} 
\Inv{s(k)}= 
\begin{cases} s(j) & \mbox{si } 
    \sigma(1)+\sigma(j)+\sigma(k)=n \,;\\ 
s(\overline{j+1})\neq s(j) & \mbox{sinon.} 
\end{cases} 
\end{align} 
Or, nous avons 
\begin{align*} 
\frac{1}{\mu}( 
\sigma(1)+\sigma(j)+\sigma(k)-n)=\varepsilon-s(j)-s(k) \in ] 
-1,1[. 
\end{align*} 
Finalement, les équivalences suivantes démontrent l'égalité 
(\ref{eq:inverse,s(k)}) 
\begin{align*} 
\sigma(1)+\sigma(j)+\sigma(k)=n \Leftrightarrow s(j)+s(k) \in \NN 
\Leftrightarrow s(j)=\Inv{s(k)}. 
\end{align*} 
\end{itemize} 
\end{proof}

\section{Les conditions initiales du potentiel}\label{sec:les-cond-init} 
Dans ce paragraphe, nous allons montrer que les nombres 
$(\!([\widetilde{\omega}_{1}],[\widetilde{\omega}_{j}],[\widetilde{\omega}_{k}])\!)$ 
engendrent le potentiel de la structure de Frobenius du polyn\^{o}me 
de Laurent $f$. 
 
Soit $X$ l'espace de base d'un déploiement universel de $f$. Soit 
$t_{0}, \ldots ,t_{\mu-1}$ des coordonnées plates au voisinage de $0$ 
dans $X$. 
 
Le champ d'Euler est défini par 
\begin{align}\label{eq:champ,euler} 
\mathfrak{E}&=\sum_{k=0}^{\mu-1}(1-\sigma(k))t_{k}\partial 
_{t_{k}}+\mu\partial_{t_{1}}. 
\end{align} 
Nous développons le potentiel de la structure de Frobenius en série 
entière et nous le notons 
\begin{align*} 
F^{sing}(\bs{t})&=\sum_{\alpha_{0}, \ldots ,\alpha_{\mu-1}\geq 0} 
A(\bs{\alpha}) \frac{t^{\bs{\alpha}}}{\bs{\alpha}!} 
\end{align*} 
où $\bs{\alpha}:=(\alpha_{0}, \ldots ,\alpha_{\mu-1})$ et 
$\frac{t^{\bs{\alpha}}}{\bs{\alpha}!}:=\frac{t_{0}^{\alpha_{0}}}{\alpha_{0}!}\cdots\frac{t_{\mu-1}^{\alpha_{\mu-1}}}{\alpha_{\mu-1}!}$. 
Nous appelons $|\alpha|:=\alpha_{0}+\cdots+\alpha_{\mu-1}$, la 
longueur de $A(\bs{\alpha})$.

Notons $(g^{ab})$ la matrice inverse de la forme bilinéaire non dégénérée dans les 
coordonnées $\bs{t}$.  Pour tout $a\in\{0, \ldots ,\mu-1\}$, notons 
par $a^{\star}$ l'unique élément de $\{0, \ldots ,\mu-1\}$ tel que 
$g^{aa^{\star}}\neq 0$. Pour tout $i,j,k,\ell \in \{0, \ldots 
,\mu-1\}$, le potentiel vérifie les équations WDVV

\begin{align}\label{eq:WDVV} 
(i,j,k,\ell):\ \ \  \sum_{a=0}^{\mu-1}F^{sing}_{ija} g^{aa^{\star}}F^{sing}_{a^{\star}k\ell}&=\sum_{a=0}^{\mu-1}F^{sing}_{jka} g^{aa^{\star}}F^{sing}_{a^{\star}i\ell}, 
\end{align} 
la condition d'homogénéité par rapport au champ d'Euler 
\begin{align}\label{eq:homo,potentiel} 
\mathfrak{E}\cdot F^{sing}&= (3-n)F^{sing} 
\end{align} 
et les conditions 
\begin{align}\label{eq:condition,initiale,potentiel} 
F_{ijk}^{sing}(\bs{0})&=g\mid_{\bs{t}=0}(\partial_{t_{i}}\star \partial_{t_{j}},\partial_{t_{k}}) \mbox{  
  où } F^{sing}_{ijk}:=\frac{\partial^{3} F^{sing}}{\partial 
  t_{i}\partial t_{j} \partial t_{k}}. 
\end{align} 
 
Remarquons que $g\mid_{\bs{t}=0}(\partial_{t_{i}}\star 
\partial_{t_{j}},\partial_{t_{k}})=(\!(\widetilde{\omega}_{i},\widetilde{\omega}_{j},\widetilde{\omega}_{k})\!)$. 
Notons $A_{ijk}(\bs{\alpha})$ le nombre $A(\alpha_{0}, \ldots 
,\alpha_{i}+1, \ldots ,\alpha_{j}+1, \ldots ,\alpha_{k}+1, \ldots 
,\alpha_{\mu-1})$. 
 
\begin{thm}\label{thm:basique,potentiel} 
  Le potentiel $F^{sing}$ est déterminé par les nombres 
  $A_{1jk}(\bs{0})$ avec $j,k\in\{0, \ldots ,\mu-1\}$ tels que 
  $\overline{1+j+k}=n$. 
\end{thm} 
 
\begin{rem}\label{rem:basique,potentiel} 
\begin{enumerate}  
\item Nous avons 
  $A_{1jk}(\bs{0})=(\!(\widetilde{\omega}_{1},\widetilde{\omega}_{j},\widetilde{\omega}_{k})\!)$. 
\item Si $\overline{1+j+k}\neq n$, alors $A_{1jk}(\bs{0})=0$ d'après 
  l'égalité (\ref{eq:condition,initiale,potentiel}). 
  \end{enumerate} 

\end{rem}

\begin{lem} \label{lem:condition,initiale} 
  Le potentiel $F^{sing}$ est déterminé par les nombres 
  $A_{ijk}(\bs{0})$ avec $i,j,k\in \{0, \ldots ,\mu-1\}.$ 
\end{lem} 
 
 \begin{proof} 
   Nous allons démontrer le lemme par récurrence sur la longueur des 
   nombres $A(\bs{\alpha})$.  Pour tout $i,j,k,\ell\in \{0, \ldots 
   ,\mu-1\}$, le terme de 
   $F^{sing}_{ija}g^{aa^{\star}}F^{sing}_{a^{\star}k\ell}$ devant 
   $\bs{\frac{\bs{t}^{{\alpha}}}{\alpha!}}$ est 
    
   \begin{align*} 
   g^{aa^{\star}}\sum_{\bs{\beta}+\bs{\gamma}=\bs{\alpha}} {\beta_{0} \choose 
     \alpha_{0}}\cdots {\beta_{\mu-1} \choose \alpha_{\mu-1}} 
   A_{ija}(\bs{\beta})A_{a^{\star}k\ell}(\bs{\gamma}). 
   \end{align*} 
   Ainsi, les termes de plus grande longueur, c'est-à-dire de 
   longueur $|\bs{\alpha}|+3$, dans la somme ci-dessus sont $ 
   g^{aa^{\star}}A_{ija}(\bs{\alpha})A_{a^{\star}k\ell}(\bs{0})$ et 
   $g^{aa^{\star}}A_{ija}(\bs{0})A_{a^{\star}k\ell}(\bs{\alpha})$.  Comme le 
   potentiel vérifie les conditions 
   (\ref{eq:condition,initiale,potentiel}), nous en déduisons que 
   $A_{ija}(\bs{0})\neq 0$ si et seulement si 
   $a=\overline{i+j}^{\star}$.

   Dans l'équation WDVV $(1,j,k,\ell)$, les termes de longueur 
   $|\bs{\alpha}|+3$ devant $\frac{t^{\bs{\alpha}}}{\bs{\alpha}!}$ 
   sont 
 \begin{itemize} 
 \item 
   $g^{\overline{1+j}, \overline{1+j}^{\star}}A_{1j\overline{1+j}^{\star}}(\bs{0})A_{\overline{1+j}k\ell}(\bs{\alpha})$ ; 
 \item 
   $g^{\overline{k+\ell}, \overline{k+\ell}^{\star}}A_{1j\overline{k+\ell}}(\bs{\alpha})A_{\overline{k+\ell}^{\star}k\ell}(\bs{0})$ ;
 \item 
   $g^{\overline{j+k}, \overline{j+k}^{\star}}A_{jk\overline{j+k}^{\star}}(\bs{0})A_{\overline{j+k}1\ell}(\bs{\alpha})$  ;
   et 
 \item 
   $g^{\overline{1+\ell}, \overline{1+\ell}^{\star}}A_{jk\overline{1+\ell}}(\bs{\alpha})A_{\overline{1+\ell}^{\star}1\ell}(\bs{0})$. 
 \end{itemize}  
 Les termes du type $A_{???}(\bs{0})$ se calculent par 
 (\ref{eq:condition,initiale,potentiel}) et la condition d'homogénéité 
 (\ref{eq:homo,potentiel}) implique 
 \begin{align*} 
 A(\alpha_{0},\alpha_{1}+1,\alpha_{2},\ldots,\alpha_{\mu-1})=\frac{1}{\mu}A(\bs{\alpha})d(\bs{\alpha}) 
 \mbox{ pour } |\bs{\alpha}|\geq 3 
 \end{align*} 
 où $d(\bs{\alpha})=3-n+\sum_{k=0}^{\mu-1}\alpha_{k}(\sigma(k)-1)$. 
 Nous en déduisons que les nombres $A_{1??}(\bs{\alpha})$ s'expriment 
 en fonction de nombres de longueur strictement plus petite.  Ainsi, 
 l'équation WDVV $(1,j,k,\ell)$ permet d'obtenir une relation entre 
 $A_{\overline{1+j}k\ell}(\bs{\alpha})$ et 
 $A_{jk\overline{1+\ell}}(\bs{\alpha})$. 
 \end{proof}

 \begin{proof}[Démonstration du théorème \ref{thm:basique,potentiel}] 
   D'après le lemme \ref{lem:condition,initiale}, il suffit de montrer 
   que les nombres $A_{ijk}(\bs{0})$ se calculent en fonction des 
   termes du type $A_{1??}(\bs{0})$.  Les nombres de longueur $3$ dans 
   l'équation $(1,j,k,\ell)$ sont non nuls si et seulement si 
   $\overline{1+j+k+\ell}=n$.  Sous cette condition, nous avons 
   $\overline{1+j}=\overline{k+\ell}^{\star}$ et 
   $\overline{j+k}^{\star}=\overline{1+\ell}$.  Ainsi, les termes de 
   longueur $3$ dans l'équation $(1,j,k,\ell)$ sont 
 \begin{itemize}  
 \item 
   $A_{1j\overline{1+j}^{\star}}(\bs{0})A_{\overline{1+j}k\ell}(\bs{0})$ 
   et 
 \item 
   $A_{jk\overline{1+\ell}}(\bs{0})A_{\overline{1+\ell}^{\star}1\ell}(\bs{0})$. 
 \end{itemize} 
 En considérant successivement les équations WDVV 
 $(1,j-1,k,\ell+1)$,$(1,j-2,k,\ell+2)$,..., nous pouvons exprimer 
 $A_{\overline{1+j}k\ell}(\bs{0})$ en fonction des nombres du type 
 $A_{1??}(\bs{0})$. 
\end{proof}

\section[Une structure d'algèbre de Frobenius graduée sur 
  $\gr^{V}_{\star}(G_{0}/\theta G_{0})$]{Définition d'une structure d'algèbre de Frobenius graduée sur 
  $\gr^{V}_{\star}(G_{0}/\theta G_{0})$}\label{sec:defin-dune-algebre} 
Dans ce paragraphe nous allons quotienter le produit $\star$ et la 
forme bilinéaire $[g]$ obtenus sur 
$G_{0}/\theta G_{0}$ par la filtration $V_{\bullet}G_{0}$. 
 
\begin{prop}\label{prop:produit,gradue} 
  Le produit que nous avons défini sur $G_{0}/\theta G_{0}$ est 
  compatible avec la filtration $V_{\bullet}(G_{0}/\theta G_{0})$, c'est-à-dire qu'on a 
  \begin{align*} 
  V_{\beta_{1}}(G_{0}/\theta G_{0}) \star 
  V_{\beta_{2}}(G_{0}/\theta G_{0}) \subset 
  V_{\beta_{1}+\beta_{2}}(G_{0}/\theta G_{0}). 
  \end{align*} 
\end{prop} 
 
Notons $\gr^{V}_{\star}(G_{0}/\theta G_{0})$ le gradué 
$\oplus_{\beta}\gr^{V}_{\beta}(G_{0}/\theta G_{0})$ et notons $\cup$ 
le gradué du produit $\star$ sur $\gr^{V}_{\star}(G_{0}/\theta 
G_{0})$. 
 
\begin{rem}Soit $\alpha\in\{0, \ldots ,d\}$. 
  Nous allons rappeler les notations et les résultats du paragraphe 
  $4$.a de \cite{DSgm1}.  Le polyn\^{o}me de Laurent  
  \begin{align*} 
    f_{\alpha}(u_{0,\alpha}, \ldots 
  ,u_{n,a})=\sum_{i=0}^{n}u_{i,\alpha} 
  \end{align*} 
sur le tore $U_{\alpha}$ est commode et non dégénéré
par rapport à son polyèdre de Newton (cf. bas de page p.\pageref{eq:9}), noté
$\Gamma(f_{\alpha})$.
  Soit $\sigma$ une face de $\partial \Gamma(f_{\alpha})$. Soit $L_{\sigma}$ la 
  forme linéaire telle que $L_{\sigma}=1$ sur $\sigma$.  Pour 
  $g\in\CC[u_{\alpha},u_{\alpha}^{-1}]$, on pose 
  $\phi_{\alpha}(g) :=\max_{\sigma}\max\{L_{\sigma}(a)\mid a \in \rm{Supp}(g)\}$. 
  Dans la figure \ref{fig:1}, on peut voir $\phi_{\alpha}$ comme la jauge du 
  convexe $\Gamma(f_{\alpha})$. Comme la jauge d'un convexe vérifie l'inégalité 
  triangulaire, on en déduit l'inégalité suivante  : 
\begin{align}\label{relation:convexite} 
\phi_{\alpha}(gh)&\leq \phi_{\alpha}(g)+\phi_{\alpha}(h). 
\end{align}   
 
\begin{figure}[tbhp] 
\begin{center} 
  \psfrag{x}{$u_{1}$} \psfrag{y}{$u_{2}$} \psfrag{O}{$0$} 
  \psfrag{G}{$\partial\Gamma(f)$} \psfrag{g}{$g$} \psfrag{h}{$h$} 
  \psfrag{p}{$gh$} 
  \includegraphics[width=0.4\linewidth,height=0.4\textheight,keepaspectratio]{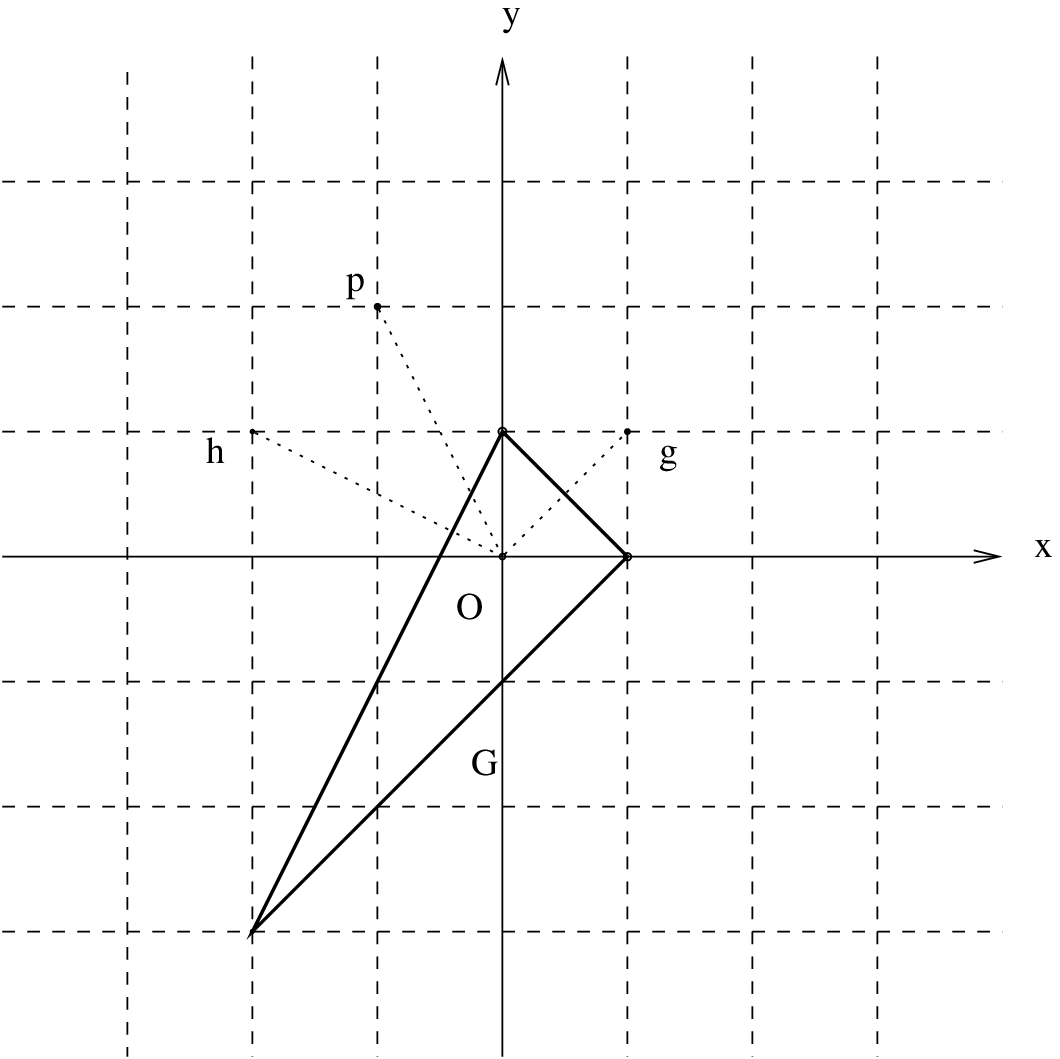} 
\end{center} 
\caption{Exemple pour  $f=u_{1}+u_{2}+u_{1}^{-2}u_{2}^{-3}$, 
  $g=u_{1}u_{2}$ et $h=u_{1}^{-2}u_{2}$.}\label{fig:1} 
\end{figure}    
 
Soit $v_{1,\alpha}, \ldots ,v_{n,\alpha}$ les coordonnées de $(\CC^\star)^n$ qui 
paramètrent le tore $U_{\alpha}$. On pose 
$\frac{dv_{\alpha}}{v_{\alpha}} :=\frac{dv_{1,\alpha}}{v_{1,\alpha}}\wedge \cdots \wedge 
\frac{dv_{n,\alpha}}{v_{n,\alpha}}$.  On définit la filtration de 
$\Omega^{n}(U)[\theta]$ par 
\begin{align*}\mathcal{N}_{\beta}\Omega^{n}(U_{\alpha})[\theta]=\mathcal{N}_{\beta}\Omega^{n}(U_{\alpha})+\theta\mathcal{N}_{\beta-1}\Omega^{n}(U_{\alpha})+  
\theta^{2}\mathcal{N}_{\beta-2}\Omega^{n}(U_{\alpha})\cdots\end{align*} 
où 
$\mathcal{N}_{\beta}\Omega^{n}(U_{\alpha}) :=\{g \frac{dv_{\alpha}}{v_{\alpha}}\in 
\Omega^{n}(U_{\alpha})\mid \phi_{\alpha}(g)\leq \beta\}$. 
 
Cela induit une filtration, notée $\mathcal{N}_{\bullet}G_{0}(f_{\alpha})$, du 
réseau de Brieskorn qui est défini par 
\begin{align*} 
\mathcal{N}_{\beta}G_{0}(f_{\alpha}) :=\mathcal{N}_{\beta}\Omega^{n}(U_{\alpha})[\theta]/(\theta 
d- df_{\alpha}\wedge)\Omega^{n-1}(U_{\alpha})[\theta]\cap 
\mathcal{N}_{\beta}\Omega^{n}(U_{\alpha})[\theta]. 
\end{align*}   
Pour plus de précision sur la filtration de Newton, on renvoie au 
paragraphe $4$ de \cite{DSgm1}.  D'après le théorème $4.5$ de 
l'article \cite{DSgm1}, on a 
$\mathcal{N}_{\beta}G_{0}(f_{\alpha})=V_{\beta}G_{0}(f_{\alpha})$. Par 
passage aux quotients, on obtient $\mathcal{N}_{\beta}(G_{0}(f_{\alpha})/\theta 
G_{0}(f_{\alpha}))=V_{\beta}(G_{0}(f_{\alpha})/\theta G_{0}(f_{\alpha}))$. 
\end{rem} 
 
\begin{proof}[Démonstration de la proposition \ref{prop:produit,gradue}] 
Pour démontrer cette proposition, il suffit de vérifier que 
 \begin{align*} 
  V_{\beta_{1}}(G_{0}(f_{\alpha})/\theta G_{0}(f_{\alpha})) \star_{\alpha} 
  V_{\beta_{2}}(G_{0}(f_{\alpha})/\theta G_{0}(f_{\alpha})) \subset 
  V_{\beta_{1}+\beta_{2}}(G_{0}(f_{\alpha})/\theta G_{0}(f_{\alpha})). 
  \end{align*} pour tout $\alpha\in\{0, \ldots ,d-1\}$. 
 
 Pour $i\in\{1,2\}$, soit $[h_{i}]$ dans $V_{\beta_{i}}(G_{0}(f_{\alpha})/\theta G_{0}(f_{\alpha}))$. 
Nous avons 
\begin{align*} 
  [h_{1}]\star_{\alpha}[h_{2}]& = 
  \left[g_{1}\frac{dv_{\alpha}}{v_{\alpha}}\right]\star_{\alpha} 
  \left[g_{2}\frac{dv_{\alpha}}{v_{\alpha}}\right] && \mbox{ où } 
  g_{1},g_{2}\in \CC[v_{\alpha},v_{\alpha}^{-1}] \mbox{ tels que 
    }\phi_{\alpha}(g_{i})\leq \beta_{i}\\ 
 & = \left[g_{1}g_{2}\frac{dv_{\alpha}}{v_{\alpha}}\right] && \mbox{ 
   d'après (\ref{eq:18})} 
\end{align*} 
L'inégalité (\ref{relation:convexite}) montre que $\phi_{\alpha}(g_{1} 
  g_{2})\leq \beta_{1}+\beta_{2}$  c'est-à-dire que 
  \begin{align*} 
    [h_{1}]\star[h_{2}]\in V_{\beta_{1}+\beta_{2}}(G_{0}(f_{\alpha})/\theta 
  G_{0}(f_{\alpha})). 
  \end{align*} 
\end{proof} 
 
Notons $[\![g]\!]$ la classe d'un élément $g \in G_{0}$ dans $\gr^{V}_{\star}(G_{0}/\theta G_{0})$.
D'après le lemme \ref{lem:produit,star}, la base $\bs{\omega}$
de $G_{0}$ induit une base  $[\![\bs{\omega}]\!]:=([\![{\omega}_{0}]\!], \ldots ,[\![{\omega}_{\mu-1}]\!])$ de  $\gr^{V}_{\star}(G_{0}/\theta G_{0})$.
 
\begin{prop}\label{prop:cup} Dans la base $[\![\bs{\omega}]\!]$ de 
  $\gr^{V}_{\star}(G_{0}/\theta G_{0})$, le produit $\cup$ est donné 
  par la formule 
  \begin{align*} 
  [\! [\omega_{i}]\!]\cup [\![\omega_{j}]\!] =  
\begin{cases} 
 \ds{\frac{w^{\underline{a}(i)+\underline{a}(j)}}{w^{\underline{a}(i+j)}}} [\![\omega_{\overline{i+j}} ]\!] & \mbox{ si } 
 \sigma(\overline{i+j})=\sigma(i)+\sigma(j) \,;
 \\ 
0  & \mbox{ si } 
\sigma(\overline{i+j})<\sigma(i)+\sigma(j). 
\end{cases} 
\end{align*} 
\end{prop}

La forme bilinéaire non dégénérée $[g]$, définie sur $G_{0}/\theta G_{0}$, est compatible à 
la filtration $V_{\bullet}(G_{0}/\theta G_{0})$. Notons 
$[\![g]\!](\cdot,\cdot)$ la forme bilinéaire non dégénérée induite sur 
$\gr^{V}_{\star}(G_{0}/\theta G_{0})$. Le lemme \ref{lem:metrique}
implique la proposition suivante.

\begin{prop}\label{prop:metrique}  Soient $\gamma$ et $\gamma'$ dans $\V$. Soient $d$ et $d'$ dans
  respectivement $\{0, \ldots ,\delta(\gamma)-1\}$ et $\{0, \ldots 
  ,\delta(\gamma')-1\}$. 
\begin{enumerate}  
\item Si $\gamma'\neq\{1-\gamma\}$ alors on a 
  $[\![g]\!]([\![\widetilde{\omega}_{k^{\max}(\gamma)-d}]\!],[\![\widetilde{\omega}_{k^{\max}(\gamma')-d'}]\!])= 
  0$. 
   
\item Si $\gamma'=\{1-\gamma\}$ alors on a 
  \begin{align*} 
  [\![g]\!]([\![\widetilde{\omega}_{k^{\max}(\gamma)-d}]\!],[\![\widetilde{\omega}_{k^{\max}(\{1-\gamma\})-d'}]\!])\end{align*} 
  \begin{align*} 
\begin{cases} 
\left(\prod_{k\in I(\gamma)} w_{k}\right)^{-1} & \mbox{ si } \sigma(k^{\max}(\gamma)-d)+\sigma(k^{\max}(\{1-\gamma\})-d')=n\,;\\ 
0 & \mbox{ sinon.} 
\end{cases} 
\end{align*} 
\end{enumerate} 
\end{prop}

\begin{proof} 
  Soient $d_{i}:=k^{\max}(s(i))-i$ et $d_{j}:=k^{\max}(s(j))-j$. Nous avons 
  l'équivalence suivante : 
  \begin{align*} 
    \overline{i+j}=n \Leftrightarrow \begin{cases} 
      s(j)=\{1-s(i)\} \\ 
      \mbox{et } d_{i}+d_{j}=\delta(s(i))-1.  \end{cases} 
  \end{align*} 
  D'apr\`es le lemme \ref{lem:metrique}, l'expression de
  $[\![g]\!](\cdot,\cdot)$ dans la base $[\![\omega]\!]$ est 
  : 
  \begin{align*} 
  [\![g]\!]([\![\omega_{k^{\max}(\gamma)-d}]\!],[\![\omega_{k^{\max}(\gamma')-d'}]\!])= 
  \end{align*} 
  \begin{align*} 
\begin{cases} 
  \ds{\frac{w^{\widetilde{\alpha}}}{w^{\widetilde{\beta}}}} & \mbox{si 
    }\begin{cases} \gamma'=\{1-\gamma\}\\ \mbox{et } 
    d+d'=\delta(\gamma)-1 \,;
    \end{cases}\\ 
0 & \mbox{ sinon} 
\end{cases} 
\end{align*} 
où 
  \begin{align*} 
\widetilde{\alpha} & =\underline{a}(k^{\max}(\gamma)-d)+\underline{a}(k^{\max}(\gamma')-d')  \\ \widetilde{\beta} &=\underline{a}(n+1)+\underline{a}(k^{\max}(\gamma)-d+k^{\max}(\gamma')-d'-n)  
  \end{align*}

Si nous avons $\gamma'=\{1-\gamma\}$, la remarque 
\ref{rem:symetrie} montre que nous avons l'équi\-va\-len\-ce entre 
$d+d'=\delta(\gamma)-1$ et 
$\sigma(k^{\max}(\gamma)-d)+\sigma(k^{\max}(\{1-\gamma\}-d')=n$.  Dans la base 
$[\![\bs{\widetilde{\omega}}]\!]$, nous obtenons 
\begin{align*} 
[\![g]\!]([\![\widetilde{\omega}_{k^{\min}(\gamma)+d}]\!],[\![\widetilde{\omega}_{k^{\max}(\{1-\gamma\})-d}]\!])= 
\frac{w^{\alpha}}{w^{\beta}} 
\end{align*} 
où 
  \begin{align*} 
\alpha & = 
\underline{a}(k^{\min}(\gamma))+\underline{a}(k^{\min}(\{1-\gamma\}) 
\\ \beta & =\underline{a}(n+1)+\underline{a}(k^{\max}(\gamma)+k^{\max}(\{1-\gamma\})-\delta(\gamma)+1-n)  
  \end{align*} 
 
Pour tout $j$ dans $\{0, \ldots ,n\}$, nous avons 
$\underline{a}(k^{\max}(\gamma)-\delta(\gamma)+1)_{j}=\#\{\ell\in\{0, \ldots 
,w_{j}-1\}\mid \ell<\gamma w_{j}\}$. Nous en déduisons 
\begin{align*} 
\underline{a}(k^{\max}(\gamma)-\delta(\gamma)+1)_{j}+\underline{a}(k^{\max}(\{1-\gamma\})-\delta(\{1-\gamma\})+1)_{j}=\end{align*} 
\begin{align*} 
\begin{cases} 
0 & \mbox{ si } \gamma=0\,;\\ 
w_{j} & \mbox{ si }\gamma>0  \mbox{ et } j\in I(\gamma)\,;\\  
w_{j}+1 & \mbox{ si }  \gamma>0\mbox{ et }j\in I^{c}(\gamma). 
\end{cases} 
\end{align*} 
Pour finir la démonstration, il suffit d'utiliser le corollaire 
\ref{cor:position,inverse} et l'égalité $\underline{a}(n+1)=(1, \ldots ,1)$. 
\end{proof} 
 
Nous définissons une forme $3$-linéaire, notée 
$(\!(\cdot,\cdot,\cdot)\!)$, sur $\gr_{\star}(G_{0}/ \theta G_{0})$ 
par la formule $(\!(a,b,c)\!):=[\![g]\!](a\cup b , c)$. 
 
\begin{prop}\label{prop:3,tenseur} 
  Soient $\gamma_{0},\gamma_{1},\gamma_{\infty}$ dans $\V$. Soient 
  $d_{0},d_{1},d_{\infty}$ dans respectivement $\{0, \ldots 
  ,\delta(\gamma_{0})-1\},\{0, \ldots ,\delta(\gamma_{1})-1\}$ et 
  $\{0, \ldots ,\delta(\gamma_{\infty})-1\}$. 
\begin{enumerate} \item Si  $\gamma_{0}+\gamma_{1}+\gamma_{\infty}$ 
  n'est pas un entier alors on a 
  \begin{align*}(\!([\![\widetilde{\omega}_{k^{\max}(\gamma_{0})-d_{0}}]\!],[\![\widetilde{\omega}_{k^{\max}(\gamma_{1})-d_{1}}]\!],[\![\widetilde{\omega}_{k^{\max}(\gamma_{\infty})-d_{\infty}}]\!])\!)= 
  0.\end{align*} 
\item Si $\gamma_{0}+\gamma_{1}+\gamma_{\infty}$ est un entier alors 
  nous avons : 
\begin{align*} (\!( 
   & [\![\widetilde{\omega}_{k^{\max}(\gamma_{0})-d_{0}}]\!],[\![\widetilde{\omega}_{k^{\max}(\gamma_{1})-d_{1}}]\!],[\![\widetilde{\omega}_{k^{\max}(\gamma_{\infty})-d_{\infty}}]\!])\!) 
    \\ 
&    = 
\begin{cases}{ 
\frac{\ds{\prod_{i\in J_{w}(\{1-\gamma_{0}\},\{1-\gamma_{1}\},\{1-\gamma_{\infty}\}) } w_{i}}}{\ds{\prod_{i\in I(\gamma_{0},\gamma_{1},\gamma_{\infty})} w_{i}}}}&  
      \mbox{ si }
      \sum_{i\in\{0,1,\infty\}}\sigma(k(\gamma_{i})-d_{i})=n \,;\\ 
 0 & \mbox{ sinon} 
\end{cases} 
\end{align*} où 
$J_{w}(\{1-\gamma_{0}\},\{1-\gamma_{1}\},\{1-\gamma_{\infty}\}):=\{i\mid 
\{\{1-\gamma_{0}\}w_{i}\}+\{\{1-\gamma_{1}\}w_{i}\}+\{\{1-\gamma_{\infty}\}w_{i}\}=2\}$. 
\end{enumerate} 
\end{prop} 
 
\begin{rem} Comme $(\!(\cdot,\cdot,\cdot)\!)$ est symétrique en ses 
  trois arguments, l'algèbre $(\gr_{\star}^{V}(G_{0}/ \theta 
  G_{0}),\cup,[\![g]\!](\cdot,\cdot))$ est une algèbre de Frobenius 
  graduée. 
\end{rem} 
  
\begin{proof}[Démonstration de la proposition \ref{prop:3,tenseur}] 
  Les propositions \ref{prop:cup} et \newline \ref{prop:metrique}  impliquent 
  l'égalité suivante 
  \begin{align*} 
  (\!([\![\omega_{i}]\!],[\![\omega_{j}]\!],[\![\omega_{k}]\!])\!)= 
    \begin{cases} 
      \ds{\frac{w^{\underline{a}(i)+\underline{a}(j)+\underline{a}(k)}}{w^{\underline{a}(n+1)+\underline{a}(i+j+k-n)}}} 
      & 
\begin{cases} \mbox{ si } \overline{i+j+k}=n\\ 
   \mbox{ et si } \sigma(i)+\sigma(j)+\sigma(k)=n\,;  
\end{cases}\\ 
0 & \mbox{ sinon.}  
 \end{cases} 
\end{align*}

Par définition des nombres spectraux, nous avons 
\begin{align*} 
  \sum_{i\in\{0,1,\infty\}}\sigma(k^{\max}(\gamma_{i})-d_{i})&=\sum_{i\in 
  \{0,1,\infty\}} k^{\max}(\gamma_{i})-d_{i} - 
\mu(\gamma_{0}+\gamma_{1}+\gamma_{\infty}). 
\end{align*} 
 Nous en déduisons 
l'équivalence 
\begin{align*} 
\begin{cases} 
    \overline{\sum_{i\in\{0,1,\infty\}}k^{\max}(\gamma_{i})-d_{i}}=n  \\ 
    \mbox{et } 
    \sum_{i\in\{0,1,\infty\}}\sigma(k^{\max}(\gamma_{i})-d_{i})=n 
\end{cases}  
\Leftrightarrow  
\begin{cases} 
\gamma_{0}+\gamma_{1}+\gamma_{\infty} \in\NN \\ 
\mbox{et } \sum_{i\in\{0,1,\infty\}}\sigma(k^{\max}(\gamma_{i})-d_{i})=n 
\end{cases} 
\end{align*}

Nous avons l'équivalence 
\begin{align*}\begin{array}{c} 
  (\!([\![\widetilde{\omega}_{k^{\max}(\gamma_{0})-d_{0}}]\!],[\![\widetilde{\omega}_{k^{\max}(\gamma_{1})-d_{1}}]\!],[\![\widetilde{\omega}_{k^{\max}(\gamma_{\infty})-d_{\infty}}]\!])\!) 
  \neq 0 \\ \Leftrightarrow  
    \begin{cases} 
\gamma_{0}+\gamma_{1}+\gamma_{\infty} \in \NN  \\ 
\mbox{ et } \sum_{i\in\{0,1,\infty\}}\sigma(k^{\max}(\gamma_{i})-d_{i})=n. 
    \end{cases} 
    \end{array}\end{align*} 
     
    Sous les conditions $\gamma_{0}+\gamma_{1}+\gamma_{\infty} \in 
    \NN$ et 
    $\sum_{i\in\{0,1,\infty\}}\sigma(k^{\max}(\gamma_{i})-d_{i})=n$, nous 
    avons 
    \begin{align*}(\!([\![\widetilde{\omega}_{k^{\max}(\gamma_{0})-d_{0}}]\!],[\![\widetilde{\omega}_{k^{\max}(\gamma_{1})-d_{1}}]\!],[\![\widetilde{\omega}_{k^{\max}(\gamma_{\infty})-d_{\infty}}]\!])\!)=\frac{w^{\alpha}}{w_{\beta}}\end{align*} 
où 
  \begin{align*} 
\alpha &=\underline{a}(k^{\min}(\gamma_{0}))+\underline{a}(k^{\min}(\gamma_{1}))+\underline{a}(k^{\min}(\gamma_{\infty})) \,;
\\ 
\beta &=\underline{a}(n+1)+\underline{a}((\gamma_{0}+\gamma_{1}+\gamma_{\infty})\mu). 
  \end{align*}

    Posons 
    $\alpha_{j}:=\underline{a}(k^{\min}(\gamma_{0}))_{j}+\underline{a}(k^{\min}(\gamma_{1}))_{j}+\underline{a}(k^{\min}(\gamma_{\infty}))_{j}-1$. 
     
    Pour tout $j\in \{0, \ldots ,n\}$, nous avons 
\begin{align}\label{eq:multi} 
\underline{a}(k^{\max}(\gamma)-\delta(\gamma)+1)_{j}&= 
\begin{cases}  
\ [\gamma w_{j}] & \mbox{ si } \gamma\in I(\gamma)\,; \\ 
\ [\gamma w_{j}]+1 & \mbox{ si } \gamma \in I^{c}(\gamma). \\ 
\end{cases} 
\end{align} 
 
Comme $\gamma_{0}+\gamma_{1}+\gamma_{\infty}\in \NN$, nous avons 
\begin{align*}\{0,\ldots ,n\}=I(\gamma_{0},\gamma_{1},\gamma_{\infty})\bigsqcup 
J_{w}(\gamma_{0},\gamma_{1},\gamma_{\infty})\bigsqcup 
J_{w}(\Inv{\gamma_{0}},\Inv{\gamma_{1}},\Inv{\gamma_{\infty}})\end{align*} 
\begin{align*}\bigsqcup \{i\mid \exists ! k\in\{0,1,\infty\} \mbox{ tel que }i\in 
I(\gamma_{k})\}.\end{align*} 
La formule (\ref{eq:multi}) permet de donner la 
valeur de $\alpha_{j}$ dans les quatre cas suivants. 
\begin{itemize} \item Si $j\in 
  I(\gamma_{0},\gamma_{1},\gamma_{\infty})$ alors $\alpha_{j}= 
  w_{j}(\gamma_{0}+\gamma_{1}+\gamma_{\infty})-1$. 
\item Si $j\in J_{w}(\gamma_{0},\gamma_{1},\gamma_{\infty})$ alors 
  $\alpha_{j}= w_{j}(\gamma_{0}+\gamma_{1}+\gamma_{\infty})$. 
   
\item Si $j\in 
  J_{w}(\Inv{\gamma_{0}},\Inv{\gamma_{1}},\Inv{\gamma_{\infty}})$ 
  alors $\alpha_{j}= w_{j}(\gamma_{0}+\gamma_{1}+\gamma_{\infty})+1$. 
  
\item Si $j\in \{i\mid \exists ! k\in\{0,1,\infty\} \mbox{ tel que 
    }i\in I(\gamma_{k})\}$ alors $\alpha_{j}= 
  w_{j}(\gamma_{0}+\gamma_{1}+\gamma_{\infty})$. 
\end{itemize} 
Pour terminer la démonstration, il suffit de remarquer que 
$\underline{a}((\gamma_{0}+\gamma_{1}+\gamma_{\infty})\mu)=((\gamma_{0}+\gamma_{1}+\gamma_{\infty})w_{0}, 
\ldots ,(\gamma_{0}+\gamma_{1}+\gamma_{\infty})w_{n})$. 
\end{proof}

\chapter{Correspondances}\label{cha:correspondance} 

\section{Démonstration de la correspondance classique}\label{section:correspondance} 
Soit $\Xi$ l'application $\CC$-linéaire définie par 
 \begin{align*} 
    \Xi  : H^{2\star}_{\orb}(\PP(w),\CC) & \longrightarrow 
    \gr^{\mathcal{N}}_{\star}\left(G_{0}/\theta G_{0}\right) \\ 
    \eta^{d}_{\gamma} & \longmapsto 
    [\![\widetilde{\omega}_{k^{\min}(\{1-\gamma\})+d}]\!] 
    \end{align*} 
     
 \begin{thm}\label{thm:corres,classique} L'application $\Xi$ est un isomorphisme gradué  
    entre les algèbres de Frobenius graduées entre
    \begin{align*} 
      (H^{2\star}_{\orb}(\PP(w),\CC),\cup,\langle\cdot,\cdot\rangle) 
    \end{align*} 
    et 
    \begin{align*} 
      (\gr^{\mathcal{N}}_{\star}\left(G_{0}/\theta G_{0}\right),\cup,[\![g]\!](\cdot,\cdot)). 
    \end{align*} 
 \end{thm} 
 
\begin{proof}  
\begin{itemize}   
\item D'après la proposition \ref{prop:base}, nous avons  
  \begin{align*} 
    \deg^{\orb}(\eta^{d}_{\gamma})=2(d+a(\gamma)). 
  \end{align*} 
La proposition 
  \ref{prop:spectre} et l'égalité (\ref{eq:min,egalite}) impliquent 
  que $\Xi(\eta^{d}_{\gamma})$ est dans le gradué 
  $\gr^{\mathcal{N}}_{d+a(\gamma)}\left(G_{0}/\theta G_{0}\right)$. 
  Nous en concluons que $\Xi$ est gradué. 
\item En comparant les propositions \ref{prop:dualite} et 
  \ref{prop:metrique}, nous en déduisons que \begin{align*}\langle 
  \eta^{d}_{\gamma},\eta^{d'}_{\gamma'}\rangle=[\![g]\!](\Xi(\eta^{d}_{\gamma}),\Xi(\eta^{d'}_{\gamma'})).\end{align*} 
  
\item De m\^{e}me, le corollaire \ref{cor:tenseur} et la proposition 
  \ref{prop:3,tenseur} impliquent l'égalité 
  \begin{align*}(\eta^{d_{0}}_{\gamma_{0}},\eta^{d_{1}}_{\gamma_{1}},\eta^{d_{\infty}}_{\gamma_{\infty}})=(\!(\Xi(\eta^{d_{0}}_{\gamma_{0}}),\Xi(\eta^{d_{1}}_{\gamma_{1}}),\Xi(\eta^{d_{\infty}}_{\gamma_{\infty}}))\!).\end{align*} 
 \end{itemize}   
\end{proof}  
 
\section{Démonstration de la correspondance quantique} \label{sec:corres-quantique} 
Si nous admettons la conjecture \ref{conj:invariant,dur}, nous 
obtenons le corollaire suivant. 
 
\begin{cor}\label{cor:correspondance,quantique} Soient $w_{0}, \ldots 
  ,w_{n}$ tels $w_{0}+\cdots+w_{n}$ et $\ppcm(w_{0}, \ldots 
  ,w_{n})$ soient premiers entre eux. Les variétés de Frobenius associées au polyn\^{o}me de Laurent $f$ et à $\PP(w)$ sont isomorphes. 
 \end{cor}  
La condition sur les poids provient du corollaire \ref{cor:bilan} 
que nous utilisons dans la démonstration ci-dessous. 
 
\begin{proof}[Démonstration du corolaire 
  \ref{cor:correspondance,quantique}] Nous allons utiliser le théorème 
  \ref{thm:iso,frobenius,dubrovin}. Montrons que les 
  variétés de Frobenius ont  
  les m\^{e}mes conditions initiales. Puis, d'après 
  le début du chapitre \ref{cha:les-singularites-du}, ces conditions 
  initiales vérifient les hypothèses du théorème 
  \ref{thm:iso,frobenius,dubrovin} ce qui montrera que les variétés de  
  Frobenius sont isomorphes.  
 
Le théorème 
  \ref{thm:structure,frobenius} nous donne les conditions initiales 
  $(A_{0}^{\circ},A_{\infty},e_{0},g)$ de la variété de Frobenius pour 
  le polyn\^{o}me de Laurent $f$.  D'après les propositions 
  \ref{prop:Potentiel-de-Gromov} et \ref{prop:dualite} nous avons la 
  m\^{e}me matrice $A_{\infty}$, le m\^{e}me vecteur propre $e_{0}$ 
  pour la valeur propre $q=0$ et la m\^{e}me forme bilinéaire non 
  dégénérée.  Il reste à comparer les matrices $A_{0}^{\circ}$ qui 
  correspondent aux multiplications par les champs d'Euler à 
  l'origine. 
Les formules (\ref{eq:champ,euler}) et (\ref{eq:champ,euler,A}) montrent  
que les champs d'Euler sont égaux. A l'origine, ils sont simplement 
$\mu\partial_{t_{1}}$ où $t_{0}, \ldots ,t_{\mu-1}$ sont les 
coordonnées plates de la variété de Frobenius. 
 Le corollaire \ref{cor:bilan},via la conjecture \ref{conj:invariant,dur}, et  
 la proposition \ref{prop:3tenseur} impliquent l'égalité  
  \begin{align*} 
   (\!(\eta_{1},\eta_{j},\eta_{k})\!)=(\!(\widetilde{\omega}_{1},\widetilde{\omega}_{j},\widetilde{\omega}_{k})\!)\end{align*} 
  pour tous $j,k$. Comme les formes bilinéaires non 
  dégénérées sont les m\^{e}mes, la formule ci-dessus montre que les 
  multiplications par les champs d'Euler à l'origine sont identiques. 
  Ceci montre que les deux variétés de Frobenius ont les m\^{e}mes 
  conditions initiales $(A_{0}^{\circ},A_{\infty},e_{0},g)$.   
\end{proof}  
\chapter{Annexe}\label{cha:annexe} 
\numberwithin{equation}{chapter}

Soit $Y$ une variété $C^{\infty}$ de dimension $n$. Soit $H$ un groupe 
fini qui agit sur $Y$. Soit $\chi$ un caractère de $H$. Notons 
$\pi:Y\rightarrow Y/H$. 
 
Le complexe de de Rham  
\begin{align*} 
\xymatrix{\mathcal{E}^{\bullet}_{Y}  :& 
  \mathcal{E}^{0}_{Y} \ar[r]^-{d}& \mathcal{E}^{1}_{Y} \ar[r]^-{d} & 
  \cdots \ar[r]^-{d} & \mathcal{E}^{n}_{Y} \ar[r]^-{d}& 0} 
\end{align*} 
est une résolution du faisceau $\underline{\CC}_{Y}$.  On a un 
isomorphisme d'espaces vectoriels 
\begin{align}\label{eq:isomorphisme} 
H^{i}(\mathcal{E}^{\bullet}_{Y}(Y)) &\stackrel{\sim}{\longrightarrow} H^{i}(Y,\underline{\CC}_{Y}). 
\end{align}

Comme $H$ agit sur $Y$, nous avons une action sur 
$H^{i}(\mathcal{E}^{\bullet}_{Y}(Y))$ définie par 
$h\cdot[\omega] :=[{h^{-1}}^{\ast}\omega]$. Cette action ne dépend 
pas du représentant choisi. Soit 
$Z^{i}_{\chi}(\mathcal{E}^{\bullet}_{Y}(Y))$ l'ensemble des $\omega$ 
dans $\mathcal{E}^{i}_{Y}(Y)$ tels que nous ayons 
\begin{itemize} \item   $d\omega=0$, 
\item pour tout $h$ dans $H$, il existe $\eta(h)$ dans 
  $\mathcal{E}^{i-1}_{Y}(Y)$ tel que 
  ${h^{-1}}^{\ast}\omega=\chi(h)\omega+d\eta(h)$. 
 \end{itemize} 
 Soit $B^{i}(\mathcal{E}^{\bullet}_{Y}(Y))$ les formes différentielles 
 exactes c'est-à-dire  
 \begin{align*} 
   B^{i}(\mathcal{E}^{\bullet}_{Y}(Y))=\{\omega\in \mathcal{E}^{i}_{Y}(Y)\mid \exists 
 \eta\in \mathcal{E}^{i-1}_{Y}(Y), d\eta=\omega\}. 
 \end{align*} 
  Posons 
\begin{align*} 
  H^{i}(\mathcal{E}^{\bullet}_{Y}(Y))_{\chi} :=Z^{i}_{\chi}(\mathcal{E}^{\bullet}_{Y}(Y))/Z^{i}_{\chi}(\mathcal{E}^{\bullet}_{Y}(Y))\bigcap 
 B^{i}(\mathcal{E}^{\bullet}_{Y}(Y)). 
\end{align*} 
Notons 
 $H^{i}(Y,\underline{\CC}_{Y})_{\chi}\subset H^{i}(Y,\underline{\CC}_{Y})$ l'image de 
 $H^{i}(\mathcal{E}^{\bullet}_{Y}(Y))_{\chi}$  par l'isomorphisme 
 (\ref{eq:isomorphisme}). 
 
  Pour tout ouvert $U$ de $Y/H$, posons 
  \begin{align*} 
    (\pi_{\ast}\underline{\CC}_{Y})_{\chi}(U)=\{x\in\underline{\CC}(\pi^{-1}(U))\mid 
 \forall h\in H, h\cdot x=\chi(h)x\}. 
  \end{align*} 
Notons 
 $(\pi_{\ast}\underline{\CC}_{Y})_{\chi}$ le faisceau ainsi défini sur 
 $Y/H$.

\begin{thm}\label{thm:cohomologie,caractere} Soit $Y$ une variété paracompacte $C^{\infty}$ de dimension $n$. Soit 
  $H$ un groupe fini qui agit sur $Y$. Soit $\chi$ un caractère du 
  groupe $H$.  Pour $i\in\{0, \ldots ,n\}$, nous avons un isomorphisme 
  \begin{align*} 
  H^{i}(Y,\underline{\CC}_{Y})_{\chi}\simeq 
  H^{i}(Y/H,(\pi_{\ast}\underline{\CC}_{{Y}})_{\chi}). 
  \end{align*} 
\end{thm}

\begin{proof} 
 \begin{itemize} \item Nous allons d'abord exprimer la cohomologie \newline 
   $H^{\star}(Y/H,(\pi_{\ast}\underline{\CC}_{{Y}})_{\chi})$ comme la 
   cohomologie d'un complexe. 
    
   Le complexe de faisceaux 
   \begin{align*} 
   \xymatrix{\pi_{\ast}\mathcal{E}^{\bullet}_{Y}  :& \pi_{\ast} 
     \mathcal{E}^{0}_{Y} \ar[r]^-{d}& \pi_{\ast}\mathcal{E}^{1}_{Y} 
     \ar[r]^-{d} & \cdots \ar[r]^-{d} & \pi_{\ast}\mathcal{E}^{n}_{Y} 
     \ar[r]^-{d}& 0} 
   \end{align*} 
   est une résolution de $\pi_{\ast}\underline{\CC}_{Y}$ (cf. la 
   démonstration de la proposition \ref{prop:cohomologie} 
   (\ref{prop:cohomologie,acyclique})). Considérons le foncteur qui à un
   faisceau $\mathcal{F}$ sur lequel $H$ agit associe le faisceau 
   $\mathcal{F}_{\chi}$ défini par $\mathcal{F}_{\chi}(U)=\{s\in 
   \mathcal{F}(U)\mid h\cdot s =\chi(h) s\}$.  Nous appliquons ce 
   foncteur au complexe $\pi_{\ast}\mathcal{E}^{\bullet}_{Y}$ et nous 
   obtenons le complexe 
    
   \begin{align*} 
   \xymatrix{(\pi_{\ast}\mathcal{E}^{\bullet}_{Y})_{\chi}  :& 
     (\pi_{\ast} \mathcal{E}^{0}_{Y})_{\chi} \ar[r]^-{\delta}& 
     (\pi_{\ast}\mathcal{E}^{1}_{Y})_{\chi} \ar[r]^-{\delta} & \cdots 
     \ar[r]^-{\delta} & (\pi_{\ast}\mathcal{E}^{n}_{Y})_{\chi} 
     \ar[r]^-{\delta}& 0} 
   \end{align*} 
    
   Montrons que $(\pi_{\ast} \mathcal{E}^{0}_{Y})_{\chi}$ est fin. 
   Comme $Y$ est paracompacte, il existe des partitions de l'unité sur 
   $Y/H$. Soit $(f_{i})_{i\in I}$ une partition de l'unité de $Y/H$. 
   Posons $\widetilde{f}_{i}:=\frac{1}{\sum_{H}\chi(h)}\sum_{h\in 
     H}\chi(h)h^{\ast}f_{i}$. Ainsi, $(\widetilde{f}_{i})_{i\in I}$ 
   est une partition de l'unité de $Y/H$ et $\widetilde{f}_{i}\in 
   (\pi_{\ast} \mathcal{E}^{0}_{Y})_{\chi}$. 
    
   Comme $(\pi_{\ast} \mathcal{E}^{i}_{Y})_{\chi}$ est un $(\pi_{\ast} 
   \mathcal{E}^{0}_{Y})_{\chi}$-module, le faisceau $(\pi_{\ast} 
   \mathcal{E}^{i}_{Y})_{\chi}$ est fin.  Nous en déduisons que la 
   cohomologie du complexe 
   $(\pi_{\ast}\mathcal{E}^{\bullet}_{Y})_{\chi}(Y/H)$ est la cohomologie 
   $H^{\star}(Y/H,(\pi_{\ast}\underline{\CC}_{{Y}})_{\chi})$. 
 
 \item Posons $Z^{i}((\pi_{\ast}\mathcal{E}^{\bullet}_{Y})_{\chi}(Y/H))$ 
   l'ensemble des $\omega$ dans 
   $(\pi_{\ast}\mathcal{E}^{i}_{Y})_{\chi}(Y/H)$ tels qu'on ait 
 \begin{itemize}  
 \item $d\omega=0$, 
 \item pour tout $h$ dans $H$, on ait 
   ${h^{-1}}^{\ast}\omega=\chi(h)\omega$. 
 \end{itemize} 
 Définissons l'application suivante  : 
  
\begin{align*} 
f^{i} :  Z^{i}_{\chi}(\mathcal{E}^{i}_{Y}(Y))& \longrightarrow  
Z^{i}((\pi_{\ast}\mathcal{E}^{\bullet}_{Y})_{\chi}(Y/H)) \\ 
 \omega & \longmapsto  \frac{1}{\sum_{H}\chi(h)}\sum_{h\in 
   H}\chi(h)h^{\ast}\omega   
\end{align*} 
 
Cette application est bien définie car nous avons 
\begin{itemize} \item   $d(f^{i}(\omega))=0$ et 
\item pour tout $h'\in H$, nous avons 
\begin{align*} {h'^{-1}}^{\ast}\left(\frac{1}{\sum_{H}\chi(h)} \sum_{h\in 
      H}\chi(h)h^{\ast}\omega\right) & =\frac{1}{\sum_{H}\chi(h)} \sum_{h\in 
      H}\chi(h)(hh'^{-1})^{\ast}\omega \\ 
& =   \frac{1}{\sum_{H}\chi(h)} \sum_{u\in H}\chi(h')\chi(u)u^{\ast}\omega  \\ 
&=\chi(h')f^{i}(\omega)  
 \end{align*} 
 \end{itemize} 
  
 Comme l'application $f^{i}$ est surjective, l'application 
\begin{align*} 
\widetilde{f}^{i} :  Z_{\chi}^{i}(\mathcal{E}^{i}_{Y}(Y))& \longrightarrow  
 H^{i}(Y/H,(\pi_{\ast}\underline{\CC}_{{Y}})_{\chi}) \\ 
 \omega & \longmapsto  \left[ \frac{1}{\sum_{H}\chi(h)}\sum_{h\in 
    H}\chi(h)h^{\ast}\omega  \right] 
\end{align*} 
est aussi surjective.  Pour finir la démonstration, il suffit de 
montrer que 
\begin{align*} 
  \left(\widetilde{f}^{i}\right)^{-1}([0])=Z_{\chi}^{i}(\mathcal{E}^{\bullet}_{Y}(Y))\bigcap 
B^{i}(\mathcal{E}^{i}_{Y}(Y)). 
\end{align*} 
Comme $\widetilde{f}^{i}(d\eta)=[0]$, 
nous avons  
\begin{align*} 
  Z_{\chi}^{i}(\mathcal{E}^{\bullet}_{Y}(Y))\bigcap 
B^{i}(\mathcal{E}^{i}_{Y}(Y))\subset 
\left(\widetilde{f}^{i}\right)^{-1}([0]). 
\end{align*} 
 
Inversement, soit $\omega\in Z_{\chi}^{i}(\mathcal{E}^{\bullet}_{Y}(Y))$ 
tel que $\widetilde{f}^{i}(\omega)=[0]$ c'est-à-dire qu'il existe 
$\eta\in \mathcal{E}_{Y}^{i-1}(Y)$ tel que  
\begin{align*} 
  d\eta= 
\frac{1}{\sum_{H}\chi(h)}\sum_{h\in H}\chi(h)h^{\ast}\omega. 
\end{align*} 
Or 
$h^{\ast}\omega=\chi(h^{-1})\omega+d\eta(h^{-1})$ avec 
$\eta(h^{-1})\in \mathcal{E}^{i-1}_{Y}(Y)$, donc 
\begin{align*}d\eta= \frac{1}{\sum_{H}\chi(h)}\sum_{h\in 
  H}\chi(h)\chi(h^{-1})\omega+ \frac{1}{\sum_{H}\chi(h)}\sum_{h\in 
  H}\chi(h)d\eta(h^{-1}).\end{align*} 
Finalement, nous en déduisons que 
$\omega\in B^{i}(\mathcal{E}^{\bullet}_{Y}(Y))$. 
  \end{itemize} 
\end{proof} 
 
Nous pouvons appliquer la m\^{e}me technique pour montrer un résultat 
sur une variété holomorphe $Y$. Soit $\Omega^{p}_{Y}$ le faisceau des 
$p$-formes différentielles holomorphes sur $Y$. Soit 
$\mathcal{A}^{p,q}$ le faisceau des formes différentielles 
$C^{\infty}$ de type $(p,q)$ sur $Y$.  Le complexe 
\begin{align*}\xymatrix{ \mathcal{A}^{p,\bullet}  : \mathcal{A}^{p,0} \ar[r]^-{d''} 
  & \mathcal{A}^{p,1}\ar[r]^-{d''}& \cdots \ar[r]^-{d''} & 
  \mathcal{A}^{p,n} \ar[r]^-{d''}& 0 }\end{align*} 
est une résolution de 
$\Omega^{p}_{Y}$.  On peut définir les m\^{e}mes objets que 
précédemment et nous obtenons le résultat suivant.

\begin{cor}\label{cor:coho,faisceau}Soit $Y$ une variété complexe paracompacte de dimension $n$. Soit 
  $H$ un groupe fini qui agit sur $Y$. Soit $\chi$ un caractère du 
  groupe $H$.  Pour $i\in\{0, \ldots ,n\}$, on a un isomorphisme 
  \begin{align*} 
  H^{i}(Y,\Omega^{p}_{Y})_{\chi}\simeq 
  H^{i}(Y/H,(\pi_{\ast}\Omega^{p}_{Y})_{\chi}). 
  \end{align*} 
\end{cor} 
      
 
\bibliographystyle{alpha} 
\bibliography{biblio}

\end{document}